# TWO COMPLETE INVARIANTS OF ORIENTED KNOTS

DIMITRIOS KODOKOSTAS

ABSTRACT. We construct two complete invariants of oriented classical knots in space. The value of each invariant on any knot is a set, infinite for the first invariant and finite for the second. The finite set is computed algorithmically from any knot diagram in a finite number of steps. The algorithm reveals the knotting number of the knot as well.

## Contents

INTRODUCTION

We are going to construct two complete invariants $S$ and $R$ of oriented knots. They are complete in the sense that each one distinguishes any two distinct oriented knots. Both of them are computed algorithmically by any given knot diagram of a knot. Each one is a set of objects of a specific kind called based symbols. We call $S$ as the symbol invariant and its value $S(\varkappa)$ for any knot $\varkappa$ contains infinitely many based symbols. $R$ is a refinement of $S$ which we call as the reduced symbol invariant, and $R(\kappa)$ contains only finitely many based symbols. The algorithm computing the elements of $R(\kappa)$ terminates after a finite number of steps, so we say the invariant is effectively computable. As a byproduct, the algorithm computes the knotting number of any knot.

We work with smooth and pl knots in the 3-euclidean space $\mathbb{E}^3$ which we identify with $\mathbb{R}^3$. As far as knot classification concerns, we equivalently work in the 3-sphere $S^3$. Wild topological knots do not succumb to this invariant. We do not consider links of more than one component, although a great deal of the work can be generalized to them.

The main idea is that we get all the information needed to construct a representative loop in space of some knot $\varkappa$, by performing a trip around any diagram $K$ of $\varkappa$ on a plane, and recording the crossings we meet along with their over and under information. This recording has to be a kind of unique signature for $\varkappa$ distinguishing it from any other knot, and also it should be something we can retrieve from any diagram of $\varkappa$.

We materialize this idea choosing some base point $p$ other than a crossing to start the trip around $K$ along its orientation, and recording the crossings as we meet them successively by the labels $1, 2, \ldots, 2n$, where say $n$ is the number of all crossings of $K$. This way each crossing point gets two distinct labels. We also record the over and under information for each crossing by a signeed number $+1$ or $-1$ taking into account the orientation of $K$ in the usual way (§1.2). We put the collected set of information in the convenient form $_pK = (a_1, b_1)^{\pi_1}(a_2, b_2)^{\pi_2} \cdots (a_n, b_n)^{\pi_n}$ where the order of the pairs is not important. Each pair of labels is ordered though. The first label in the pair is the one the crossing gets as we pass through this point along an overpass arc. We define $_pK = \varnothing$ whenever $K$ has no crossings, that is, whenever $K$ is a circle. We call the symbol $_pK$ as a based symbol with the exact definition given in §2. So the main idea is that $_pK$ carries all the information one needs to construct a knot representative in space $\mathbb{R}^3$ or $S^3$ of the knot $\kappa$ whose $K$ is a diagram.

$_pK$ depends on both $K$ and $p$. One needs a kind of universal version of such a symbol which remains unaltered under the change of the position for the base point $p$ and also unchanged under any finite sequence of changes of $K$ via the Alexander-Briggs-Reidemeister moves $\Omega_1, \Omega_2, \Omega_3$ moves, since these produce all possible knot diagrams of $\varkappa$.

To achieve this we transform everything to algebraic language. Especially, we define algebraic moves between based symbols (§3), mirroring the usual topological ones between knot diagrams. This mirroring is exact both ways in the following sense: first, any topological move $K \xrightarrow{\omega} K'$ implies the existence of a sequence of algebraic moves from $_pK$ to $_qK'$ for some choices of base points $p, q$ (Lemmata 16, 17, 18, 19, 20). Also, any algebraic move $_pK \xrightarrow{al_\omega} {_qK'}$ implies the existence of a sequence of topological moves from $K$ to $K'$ (Lemmata 22, 23, 25, 26, 27).

We also show that the algebraic moves generate an equivalence relation $\underset{al}{\sim}$ in the set $\Sigma_D$ of all based symbols (Lemma 28). The value $S(\varkappa)$ for a knot $\kappa$ is essentially defined as the equivalence class under $\underset{al}{\sim}$ of the based symbols $_pK$ of the knot diagrams $K$ of $\varkappa$ (§5). As a set, this class contains infinitely many members because we can always perform algebraic moves increasing the number of crossings of a based symbol. On the other hand, $R(\varkappa)$ defined in §6 is a subset of $S(\varkappa)$ containing only the reduced based symbols. Reduced are those based symbols with minimum number of crossings among all based symbols in an equivalence class of $\underset{al}{\sim}$. We prove that $R(\varkappa)$ is finite, and that we can produce all its elements via a finite algorithmic procedure starting from any knot diagram of $\varkappa$ (Corollary 5). In the reasoning we present, we make use of the locally reduced based symbols which are those based symbols on which no negative algebraic move can be performed, that is, no immediate move can reduce their order. Reduced and locally reduced based symbols strongly resemble the total and local minima of graphs, but unlike them, it turns out that in the case of based symbols they coincide (Lemma 35). This fact forces all desirable properties for the invariant $R$ to hold.

In §1 we recall a few facts about oriented knots and knot diagrams and we build on some of the notation to be followed throughout the paper.

In order to ease up the passing to algebraic language, we pay attention to the orientations of the diagrams, to the way a diagram is partitioned by its crossing points and to the way the diagram partitions the plane on which it lies (§2.1). Thus we define a plethora of notions like that of a basic arc, main arc, room, height of a label, and others. We subsequently analyze the topological situation before and after a topological move between diagrams and we use the findings for the gradual but complete passing to algebraic language in the next section §3 (ch. especially Definition 20 and the definitions of the algebraic moves). The definitions of the algebraic moves do not come for free though. Namely, it is not immediate that the result of an algebraic move on a based symbol is always another based symbol. For this reason we define a larger set of symbols which we call numbered symbols (Definition 20) and which for sure includes the result of our moves. The fact that the result is always a based symbol is eventually certified at the Lemmata 16, 17, 18, 19, 20.

The mirroring of the topological moves by algebraic ones is more or less immediate by the definitions. Nevertheless, the converse mirroring of the algebraic moves by topological ones is not similarly immediate. The reasons are topological: although a based symbol can inform us for a loop of some diagram made up by one, two or three arcs, it cannot inform us equally well about the whereabouts of the rest of the arcs of this diagram. In particular it cannot inform us about if the interior of this loop on the plane is void of other points of the diagram or not. So we cannot be sure we can perform a usual $\Omega_1, \Omega_2$ or $\Omega_3$ move using this 1-gon, 2-gon or 3-gon of arcs to mirror the algebraic move. To overcome the difficulty we enrich the set of topological moves defining a bunch of new ones. These generalize the $\Omega_i$'s but they are not necessarily local in character, and in general they can best be characterized as long-range moves (§4.4). Section §4 is devoted to the proof of the above mentioned converse mirroring with the help of the newly defined moves. An additional topological obstacle to the mirroring occurs whenever one tries to decide about if two diagrams with exactly the same based symbol for some choices of their base points, are located on the plane as expected in a first though, namely so as to differ by some isotopy of the plane. Unfortunately the answer is no, but nevertheless the way the diagrams are located can be described (§4.2) with the help of some other moves which we define in §4.1. As a byproduct, this effort leads to the notion of the decomposition of a diagram which is a kind of coded information on a circle capable to recreate the diagram. (§4.1).

The departure from the usual $\Omega$ moves frame, is also coupled with the occasional use of the 2-sphere $S^2 = \rho \cup \{\infty\}$ in place of a plane $\rho$ as the projection surface for the knot diagrams. We justify this change of viewpoint with the proof of a suitable version of the Alexander-Briggs-Reidemeister Theorem (Theorem 3). We work with the sphere only occasionally and although this use of the sphere is convenient it can be avoided altogether.

Some other choices had to be made as well, for example regarding the category of knots we work in and the kind of moves we work with. We opted for those choices that seem to homogenize the study and make it a bit clearer. So early on (§1.4) we consider an isotopy move as a companion to the usual $\Omega$ moves, homogenizing the set of moves for the smooth and the pl cases we work with. In the algebraic setting, we still keep an algebraic version of this isotopy move (§3.2), although this version is actually a special case of another algebraic move which we call as change of base move (§3.1).

The whole work is an interplay between Topology and Algebra aiming to transform the first to the second. A crucial point is Lemma 10 which essentially transcripts to algebraic language the deeply topological property of the transverse intersection at the crossing points of a regular diagram.

The fundamental definition of the algebraic moves is flexible enough to allow a smaller or bigger number of them. The important thing though is to allow enough of them to be able to mirror all topological moves. The usual topological endpoint in the rigorousness of a proof is the Jordan ([3, 8]) or the Schonflies Theorem ([7]) for curves on a plane or the 2-sphere. Some proofs are quite demanding, with a notable example the proof of Lemma 22 for the shake of which a whole preparatory paragraph (§4.1) was devoted.

A considerable number of figures is included, accompanied in many cases by extensive captions. We choose to work with a version of $\mathbb{Z}_{2n}$ containing as elements the integers $1, 2, \ldots, 2n$, suitable for cyclic labelings and relabelings of the crossings.

# 1. A reminder on oriented knots, knot diagrams and Alexander-Briggs-Reidemeister moves

## 1.1. Oriented smooth and pl knots.

We are interested in smooth and pl oriented knots in the three dimensional euclidean space $\mathbb{E}^3$ or in $S^3$. Each such knot is a set of equivalent smooth or pl oriented knotted loops in the relevant ambient space under the equivalence relation generated by smooth or pl isotopies of the ambient space which start with the identity map. These are the orientation preserving isotopies of the space. Below we present some more details.

We identify $\mathbb{E}^3$ with $\mathbb{R}^3$ by equipping $\mathbb{E}^3$ with an orthonormal Cartesian system $Oxyz$ and we consider $S^3$ as the one point compactification $S^3 = \mathbb{R}^2 \cup \{\infty\}$ of $\mathbb{R}^3$ with some new point $\infty$. With respect to $Oxyz$, the positive and negative orientation of $\mathbb{R}^3$ has the usual meaning. Similarly for the two orientations of $S^3$.

Let $S^1$ be the usual euclidean circle of radius 1, centered at $O$ on the $Oxy$ plane, equipped with the positive orientation of $Oxy$.

**Definition 1.** *A topological oriented knotted curve or loop $\kappa$ in $X = \mathbb{R}^3$ or $S^3$ is any homeomorphic image of $S^1$ under some homeomorphism $f : S^1 \to f(S^1) = \kappa \subset X$, equipped with the induced orientation of $S^1$ by the map $f$.*

As a set of points the loop $\kappa$ is described via many homeomorphisms $f$, but each one of them induces just one of the two possible orientations on the curve $\kappa$. Strictly speaking, as an oriented loop, $\kappa$ is always equipped by a homeomorphism $f$, but we adopt the usual convention to omit $f$ or any other kind of notational reference to the orientation of $\kappa$. In the auxiliary figures we denote the orientations as usual with curved arrows.

Whenever $f$ is smooth, we call $\kappa$ as a *smooth oriented loop*. Whenever $S^1$ is partitioned into a finite number of arcs mapped under $f$ to line segments in space, we call $\kappa$ as a polygonal or *pl (piecewise linear) loop*. Of course every smooth or pl loop is also a topological loop.

For an isotopy $F : X \times [0, 1] \to X \times [0, 1]$ of the ambient space $X = \mathbb{R}^3$ or $S^3$ where our loops reside, we will be denoting as $F_t : X \times \{t\} \to X \times \{t\}$ the instances of $F$ for the instances of time $t \in [0, 1]$. $F_t$ are homeomorphisms. If we call $pr : X \times [0, 1] \to X$ the projection of $X \times [0, 1]$ to the first factor, and assuming that no confusion arises, we also denote by $F_t$ the map $X \to X$, $x \mapsto pr(F(x, t))$. Especially for $t = 1$ the map $F_1$ will be called as *the last instance* of $F$. It is clear that at any time $t \in [0, 1]$ an isotopy $F$ sends any knotted loop $\kappa \subset X$ to another knotted loop $F_t(\kappa)$, since $F_t(\kappa) = (F_t \circ f)(S^1)$ and $F_t \circ f$ is a homeomorphism. The loop $F_t(\kappa)$ is considered oriented via the induced homeomorphism by the map $F_t \circ f$.

An isotopy $F$ of $X$ will have for us the usual meaning of an orientation preserving isotopy starting with the identity on $X$, i.e. $F_0 : X \to X$, $F_0(x, 0) = x, \forall x \in X$. We say that any oriented knotted loop $\kappa$ moves inside $X$ to the oriented knotted loop $F_1(\kappa)$ by $F$. We also say that $\kappa$ and $F_1(\kappa)$ are *ambient isotopic* in $X$. This relation of oriented loops via the isotopies is an equivalence relation in the set of oriented loops, hence the usual definition of a knot is as follows:

**Definition 2.** *We denote the equivalence relation in the set of the oriented loops in $X = \mathbb{R}^3$ or $S^3$ via isotopies in $X$ (in the usual sense, starting from the identity) as $\underset{X}{\sim}$. The oriented knots in $X$ are the classes of this equivalence relation.*

For an oriented loop $\kappa$ we will typically denote its knot class as $\varkappa$.

The above setting refers to the general topological case. Whenever for an oriented knot there exists some smooth or pl representative, we call the knot as a *smooth* or *pl oriented knot* respectively. It is not true that all topological knots admit smooth or pl representatives and those which do not are called wild knots. Here we do not consider wild knots. On the other hand, any smooth knot admits pl representatives and any pl knot admits smooth representatives. In order to distinguish between the smooth and the pl categories, we also ask that any two smooth representatives $\kappa, \kappa'$ of a smooth knot are connected via a smooth isotopy $F$ ($\kappa' = F_1(\kappa)$), in which case all instances of the isotopy are smooth homeomorphisms of the ambient space $X$ and all loops through which $\kappa$ moves to arrive eventually at $\kappa'$ are smooth as well.

Similarly in the pl case we ask for $F$ to be a pl map. Such maps can destroy the apparent pl structure of a pl loop $\kappa$. For example they can create more vertices for $\kappa'$ than those of $\kappa$. What we actually ask from $F$ is to respect some finite subdivision of $\kappa$ carrying it to a finite subdivision of $\kappa'$.

If confusion does not arise, for the rest of the paper we systematically omit references to the category we work in and to the kind of isotopy we deal with, unless an explicit mention to them adds to the comprehension of the situation. This convention comes as an addition to the set of standard notational conventions regarding the orientations which we mentioned earlier.

To study knots in $S^3$ one pushes the knotted loops of $S^3$ inside its subset $\mathbb{E}^3$. This is legitimate: almost trivially, any two oriented smooth or pl loops in $S^3 = \mathbb{R}^3 \cup \{\infty\}$ are ambient isotopic in $S^3$ if and only if any two (smooth or pl) pushings of them to $\mathbb{R}^3$ are ambient isotopic in $\mathbb{R}^3$. If for a knot $\kappa$ in $S^3$ we forget those representative loops that go through $\infty$ we get all representative loops for a knot $\kappa'$ in $\mathbb{R}^3$ and the correspondence $\kappa \mapsto \kappa'$ between the sets of knots in the two manifolds is one to one and onto. This makes knot theory in $S^3$ coincide with knot theory in $\mathbb{R}^3$ and we frequently merge the two cases to one, referring to both manifolds as "the space".

Our point of view will be given in a while in §1.5. It is the usual point of view where one considers knots in $\mathbb{R}^3$ and works their equivalence via projections on a plane and moves among the projections. We pay some extra attention though to isotopies of the plane. This attention is not needed if one wants only to establish an equivalence between the study of knots in space and the study of knot projections on a plane.

1.2. **Knot diagrams.** An extremely convenient feature of $\mathbb{R}^3$ in the study of its knotted loops is not its metric structure, but rather its well known projection map onto its planes which preserves the smooth and pl structures of the loops.

So, let $\rho$ be a fixed plane of $\mathbb{R}^3$ and $pr : \mathbb{R}^3 \to \rho$ be the perpendicular projection on $\rho$.

For an oriented smooth (pl) loop $\kappa = f(S^1)$, its projection $K = pr(\kappa) = (pr \circ f)(S^1)$ on $\rho$ is a smooth (pl) closed curve, oriented by the orientation of $S^1$ induced by the map $pr \circ f$. Whenever $K$ has only a finite number of self intersections, and locally around them only two simple arcs of the curve intersect transversally, we say that the projection is in general position. We call each double point $\Delta$ of $K$ as a crossing point or just a crossing of $\kappa'$. In some small 2-disk neighborhood $U_\Delta$ of $\Delta$ in $\rho$, the part of $K$ inside $U_\Delta$ looks like in Figure 1 (a) where one should note the transversal nature of the self intersection of $K$ at $\Delta$.

Neighborhoods like $U_\Delta$ are not necessarily round disks, but in general we consider topological such 2-disks. We call them as canonical for $K$ at $\Delta$ and if necessary we consider them to be smooth or pl in the sense that their boundaries are smooth or pl simple closed curves. Although not of utmost importance, we consider these disks as closed, that is, including their boundaries. We call the four regions at which $U$ is partitioned by $K$ as the angles of $U$. We call the two arcs $\ell_1, \ell_2$ of $K$ in $U_\Delta$ as its legs in this neighborhood. In Figure 1 it is $\ell_1 = s_1 \cup s_2$ and $\ell_2 = s_3 \cup s_4$. We call the arcs $s_1, s_2, s_3, s_4$ as the half-legs of $\ell_1$ and $\ell_2$. There exist exactly two points $\delta_1, \delta_2$ of $\kappa$ which project to $\Delta$ and exactly two arcs $\lambda_1, \lambda_2$ of $\kappa$ through $\delta_1$ and $\delta_2$ respectively which project onto $\ell_1$ and $\ell_2$. Our point of view of the plane $\rho$ and the projection $K$ on it, is from a point far away in one of the half-spaces defined by $\rho$. One of the points $\delta_1, \delta_2$ has a higher attitude from $\rho$ as this is measured on an axis perpendicular to the plane with algebraic values (attitudes from $\rho$) increasing towards us. We say that the corresponding arc $\lambda_1$ or $\lambda_2$ is higher than the other and we pass this terminology to the corresponding projections $\ell_1, \ell_2$ as well. We call the higher of them as the over-crossing arc and the other as the under-crossing arc. We call this information as *over/under information* for the diagram $K$ at its crossing $\Delta$. We denote this information as in Figure 1(b),(c). We say that $\kappa$ projects to $K$ with the accompanied over/under information at each crossing point.

We must stress that the over and under information changes whenever our position with respect to $\rho$ changes from the one to the other half-space. So from now on we adopt the usual tacit convention that we get the over/under information of the diagrams viewing our drawings from an arbitrary chosen but nevertheless always fixed half-space of $\mathbb{R}^3$ with respect to $\rho$. For definiteness, we agree from now that for our position in space, a positive frame $[v_1, v_2]$ of two vectors on $\rho$ is any frame that completed to a triad $[v_1, v_2, v_3]$ with a vector $v_3$ perpendicular to $\rho$ with origin on the plane and pointing to us, gives a positive frame of $\mathbb{R}^3$. It will be no loss of generality to assume $\rho = Oxy$ and our position on the positive $z$ half-axis, in which case a frame with $v_1$ on the positive $x$ axis and $v_2$ on the positive $y$ axis is positive.

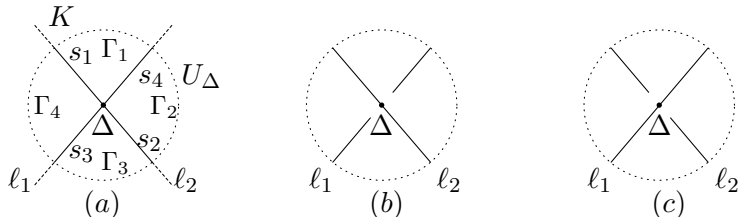

FIGURE 1. (a) A canonical neighborhood $U_\Delta$ of a crossing point $\Delta$ of a diagram $K$. $\ell_1, \ell_2$ are the legs and $s_1, s_2, s_3, s_4$ are the half-legs of $K$ in the neighborhood, $\Gamma_i$'s are the angles of the neighborhood. (b), (c) The possible over/under information at a crossing $\Delta$. In (b) $\ell_1$ is the over-crossing and $\ell_2$ the under-crossing, whereas in (c) they are the other way around.

**Definition 3.** *We call as a knot diagram or just diagram, any smooth or pl oriented closed curve $K$ on a plane $\rho$ with a finite number of self intersections (maybe none), all of them double points where the curve crosses itself (the self intersection is transverse), along with additional over/under information at each crossing point. If $\kappa$ is a smooth or pl oriented loop in space which projects to $K$ with all its orientation and over/under information, we say that $K$ is a regular projection, a knot diagram or just diagram of $\kappa$. We also say that $K$ is a diagram of the knot $\varkappa$ for which $\kappa$ is a representative. And we say that $K$ represents the knot $\varkappa$.*

The projections of an oriented loop $\kappa$ on the various planes $\rho$ are almost always regular in the following sense [2]:

The projection of $\kappa$ on any two parallel planes has the same kind and the same number of singularities. Thus we represent all parallel planes of a given direction by the planes of this direction tangent to a unit sphere $S^2$. Equivalently, we represent them by the couple of tangent points on the unit sphere of the two such tangent planes. Then the set of points on $S^2$ which correspond to non-regular projections on the corresponding tangent planes, is of measure 0 (for the usual measures in $\mathbb{R}^3$).

We usually insist to keep the projection plane $\rho$ fixed. If the projection of $\kappa$ on $\rho$ is not regular, this does not hinder our the study of the knot $\varkappa$ which $\kappa$ represents. The reason is that by the above, there exists a plane $\rho'$ very close to $\rho$ to which $\kappa$ projects regularly. If $\epsilon$ is the common line of the two planes and $F$ the rotation around $\epsilon$ which brings $\rho'$ onto $\rho$, then $F$ is an isotopy in the whole space $\mathbb{R}^3$ which brings $\kappa$ to $F(\kappa) = \kappa'$, a nearby curve with regular projection on $\rho$. Moreover, since $F$ is an isotopy, $\kappa'$ represents the same knot $\varkappa$. Its projection can be as close as we wish to the one of $\kappa$, but it is free from the singularities of the latter.

Given a diagram $K$ on $\rho$ it is clear that if there exists some loop $\kappa$ in space that projects to it, then there exists many loops with the same property: $\kappa + v$ projects to $K$ for any vector $v$ normal to $\rho$. So a diagram cannot inform us about which loop produced it via projection on the plane. Nevertheless, it informs us about which knot produced it. This is a result derived immediately by the Alexander-Briggs-Reidemeister Theorem on the ambient isotopy equivalence of loops in space via the equivalence of their knot diagrams by diagrammatic moves on the plane. We present theses moves in the next section.

As for the existence of oriented loops in space projecting to a given $K$ it is easy to construct some: choose disjoint canonical neighborhoods of the crossings (this is possible since there exist only finitely many crossings) and push vertically with respect to $\rho$ all the over-crossing arcs a bit above the plane $\rho$ fixing their boundaries. In the smooth case take care to produce a smooth result $\kappa$ in space (this can easily be done even at the endpoints of the arcs), whereas in the pl case take care to produce a pl result. Give to $\kappa$ the obvious orientation coming from $K$. By construction $\kappa$ is a smooth or pl oriented loop that projects to $K$ as wanted.

1.3. **Moves on the plane.** Due to the work of Alexander, Briggs and Reidemeister ([1, 6]) we know that two diagrams $K, K'$ on the plane $\rho$ represent the same knot in $\mathbb{R}^3$ iff one of them comes from the other after a finite number of consecutive alterations, each one happening in the interior of a bounded 2-disk:

**Theorem 1.** *(Alexander, Briggs, Reidemeister) Two smooth (or pl) knotted loops in $\mathbb{R}^3$ are representatives of the same knot iff their projections on a plane $\rho$ come one from the other after a finite number of consecutive alterations of types $\Omega_1, \Omega_2, \Omega_3$, each one in the interior of a bounded 2-disk $U$ in the way shown in Figure 2. The alterations include the shown over/under information. In the pl case we need to consider the last two alterations in the bottom row of Figure 2 of types $\Omega_4, \Omega_5$ also as moves (local isotopy moves).*

The $\Omega_4, \Omega_5$ moves in the pl case are indeed necessary ([5]): pl loops in space which differ by a pl isotopy of the space posses subdivisions whose basic elements of vertices and edges correspond in 1-1 fashion under the isotopy. This passes over to the projections on a plane, but the $\Omega_1, \Omega_2, \Omega_3$ moves alone cannot always capture it. For example these moves cannot relate two diagrams with the same underlying set whose apparent pl structure differs, with the second diagram possessing an extra vertex at the interior of an edge of the first.

The moves $\Omega_i$ transform a diagram $K$ to another diagram $K'$ and they generate an equivalence relation on the set of diagrams on $\rho$, for each one of the pl and sooth cases.

**Definition 4.** *Let $\mathbb{D}_s, \mathbb{D}_{pl}$ be respectively the set of smooth and pl diagrams in the plane $\rho$. We denote any one of them as $\mathbb{D}$ and the equivalence relation in it generated by the moves $\Omega_i$ as $\underset{top}{\sim}$.*

We denote an alteration of a diagram $K$ in a 2-disk $U$ by a move $\omega$ to a new diagram $K'$ as $K \xrightarrow{\omega} K'$. We call $U$ as the disk of the move.

According to our symbolism, the last Theorem says that if $\kappa, \kappa'$ are two knotted loops in $\mathbb{R}^3$ and $K, K'$ their projections on $\rho$ then $\kappa \underset{\mathbb{R}^3}{\sim} \kappa' \Leftrightarrow K \underset{top}{\sim} K'$.

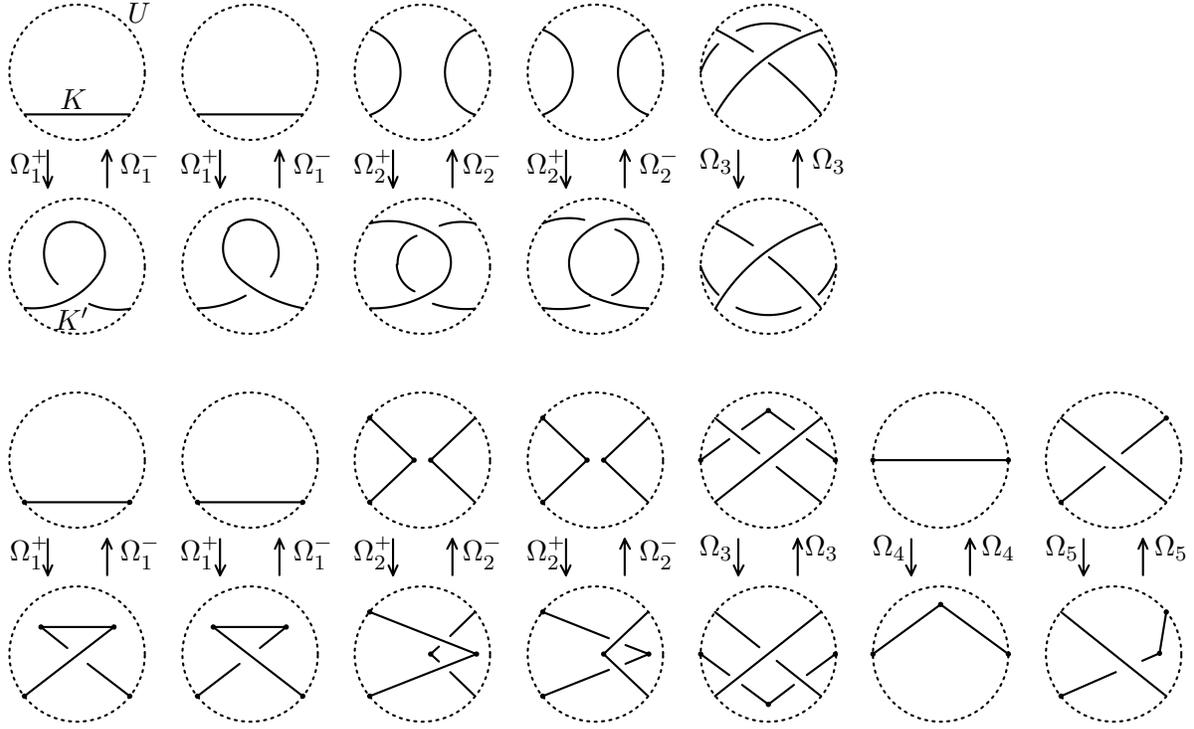

FIGURE 2. Top: smooth $\Omega_1, \Omega_2, \Omega_3$ moves. Bottom: pl $\Omega_1 - \Omega_5$ moves. The dots are vertices. In both smooth and pl cases, each move is actually defined for all possible orientations for the top diagram and with the obvious corresponding orientation for the bottom diagram. The over/under information at the crossing points are shown explicitly. The $\Omega_1, \Omega_2$ moves come into two versions $\Omega_1^+, \Omega_1^-$ and $\Omega_2^+, \Omega_2^-$ respectively and depending on the over/under information each + version comes into two varieties. The "moving" arc in the $\Omega_3$ move can be any one of the three arcs inside the 2-disk, but here only the bottom one is shown moving.

The disks $U$ in Theorem 1 are chosen with smooth or pl boundary $\partial U$ and so that the diagrams intersect $\partial U$ transversely.

Depending on the move, each one of $K \cap U, K' \cap U$ consists of a few oriented arcs of the corresponding diagram $K$ or $K'$ respectively. The arcs here have the usual meaning: we consider the diagrams as parametrized closed curves $f : S^1 \to K$, $g : S^1 \to K'$, and then their arcs are images $f(s), g(s)$ where $s$ are arcs of $S^1$. The orientation of these arcs is the one induced to them as subarcs of the oriented diagrams. By definition $K'$ keeps the orientation of $K$ outside $U$ while its arcs in $U$ are oriented in the obvious way so that it becomes an oriented closed curve with transverse intersections (as opposed for example to becoming a closed curve with double points but not transversing itself at them).

Figure 2 presents the moves without the extra information regarding the orientations. Each move is actually defined for all possible orientations for the top diagram $K$ and the obvious corresponding orientation on the bottom diagram $K'$.

The $\Omega_1$ move comes in two versions: one version creates one more crossing point whereas the other version eliminates one. We denote them as $\Omega_1^+, \Omega_1^-$ respectively. Similarly, the $\Omega_2$ move comes in two versions: one of them creates two crossing points and the other eliminates two crossing points and we denote the version as $\Omega_2^+, \Omega_2^-$ respectively.

The $\Omega_1^+$ move which creates a crossing point actually comes in two varieties which "twist" a small part of the projection like in the first two disks in the second row of Figure 2. The net topological result forgetting the over and under information is the same for the resulting diagram. There is no need right now to assign different names to these two varieties. But we shall work on this matter once again later on (ch. Definition 16).

If an $\Omega_1^s, s = \pm 1$ move alters $K$ inside the 2-disk $U$ to become $K'$, then an $\Omega_1^{-s}$ alters $K'$ inside $U$ to become $K$ once again. with the same disk $U$, we can alter $U \cap K'$ back to $U \cap K$. This way the resulting diagram is $K$ once again. In such a situation we call the two successive moves $\Omega_1^+, \Omega_1^-$ as inverses of each other.

Similar remarks hold for the $\Omega_2$ moves. The ways of placing the two arcs one on top of the other in an $\Omega_2^+$ move is a matter we take up in more detail later on (ch. Definition ??.

The $\Omega_3$ move does not create or eliminate crossing points but it rather redistributes them along the projection. $\Omega_3$ actually comes in three versions according to which one of the three arcs involved is relocated. Nevertheless, no matter if the moving arc is the top, the bottom or the middle one, the gross result of the move is topologically the same for the resulting diagram (two such diagrams as point sets differ by an isotopy of the plane), and remains so even if we include the over/under information. There will be no need for us to assign different names to the three versions. An $\Omega_3$ move that relocates some arc $\ell$ of a diagram $K$ inside a 2-disk $U$, say to some new position $\ell'$, has as an inverse move another $\Omega_3$ move performed in the same 2-disk $U$ and which returns $\ell'$ back to $\ell$ producing once again as a result the diagram $K$.

1.4. **Isotopy or not?** Usually when we perform a move on a diagram we are not very meticulous about the exact location of the arcs inside the move's disk $U$ after the move is over. The reason is that the resulting diagram is always topologically the same: by the definition of the moves, any two locations of the arcs in $U$ give two diagrams which differ by an ambient isotopy of $\rho$ starting from the identity. As long as we keep the same under and over information at corresponding crossing points we can talk about "the same" result since for example we can get from the one diagram to the other by a sequence of two inverse moves.

This lack of attention on isotopies is reinforced by the Alexander-Briggs-Reidemeister Theorem itself: there is no need to include an isotopy move to compare two diagrams. The way to go from a diagram $K$ to an equivalent one $K'$ is by consecutive small local changes inside bounded 2-disks. These changes add up making $K$ move just a bit at a time, ending up to its final position coinciding point-wise with $K'$ (also keeping orientation and over/under information).

Of course, whenever one is given two diagrams $K_1, K_2$ that differ by an isotopy on the plane $\rho$, it is much easier to use this isotopy to make them coincide, than to gradually change one of them using the $\Omega$ moves. So in the literature one can find some versions of the Theorem where isotopies are included next to the classical $\Omega$ moves, sometimes restricted to a single isotopy and usually after the other moves are performed.

Another benefit of the use of isotopies is that they eliminate the technical distinction between the smooth and pl cases: the pl $\Omega_4, \Omega_5$ moves are essentially just isotopies in the sense discussed above of retaining the pl structure of suitable subdivisions of the diagrams. Hence in the presence of an isotopy move of the plane the $\Omega_4, \Omega_5$ moves become obsolete.

In order to produce a homogeneous setting, we give a rigorous definition for the new kind of move for diagrams on the plane:

**Definition 5.** *An $\Omega_{iso}$ move on a diagram $K$ is the replacement of the diagram $K$ by the diagram $K'$ resulting by the movement of $K$ via an ambient isotopy (accordingly smooth or pl) of the plane $\rho$ which starts with the identity. If $F : \rho \times [0,1] \to \rho \times [0,1]$ is this isotopy ($F(x,0) = (x,0), \forall x \in \rho$), then $(K',1) = F_1(K,1)$ where $F_1 : \rho \times 1 \to \rho \times 1$, $F_1(x,1) = F(x,1), \forall x$ is the last moment of $F$.*

*Similarly to the case of isotopies of the ambient space, we denote here for simplicity the composition of $F_1$ with the projection map $\rho \times [0,1] \to \rho$ to the first factor again as $F_1$. We write $K' = F_1(K)$ and we say $F$ carries $K$ to $K'$. If $K$ is parametrized as $K = f(S^1)$ then $K' = (F_1 \circ f)(S^1)$ and $K'$ has the induced orientation of $S^1$ via the map $F_1 \circ f$. At a canonical neighborhood $U'$ of the crossing point $\Delta'$ of $K'$, we define the arc $\ell'_1$ of $K'$ to be over/under the arc $\ell'_2$, exactly when the arc $\ell_1 = F_1^{-1}(\ell'_1)$ of $K = F_1^{-1}(K')$ is over/under the arc $\ell_2 = F_1^{-1}(\ell'_2)$ in the canonical neighborhood $U = F_1^{-1}(U')$ of the crossing $\Delta = F_1^{-1}(\Delta')$ of $K$.*

The orientation we gave to $K'$ is also the induced orientation of $K$ via $F_1$.

If we wish we can say that the whole plane $\rho$ is an open (unbounded) 2-disk for $\Omega_{iso}$. But we can also consider usual bounded 2-disks as disks for this move in the following way: if $K_1$ is taken by the isotopy to $K_2$, then the union of all positions of $K_1$ on the plane during the isotopy is bounded ($K_1$ is compact and the isotopy continuous). Thus for some 2-disk $U$ that contains these positions we can replace the original isotopy with one that fixes everything outside $U$ and behaves as the original one in its interior. So we can consider in our definition, only those isotopies of the plane that fix everything outside a 2-disk which contain the diagram on their interior.

And we can say that an $\Omega_{iso}$ move with $F$ as its associated isotopy, has as inverse another $\Omega_{iso}$ move whose associated isotopy is $G : \rho \times [0,1] \to \rho \times [0,1], G_t = (F_t)^{-1}, \forall t$ for the obvious reason that performing the two of them consecutively we return at our original diagram.

Now let $\kappa_1, \kappa_2$ be two loops in $\mathbb{R}^3 = \rho \times \mathbb{R}$ whose projections on the plane $\rho$ are isotopic via some isotopy $F$. We extend $F$ to an isotopy $H$ of $\mathbb{R}^3$, repeating the movement of $\rho$ to the planes parallel to it: $H((a,z),t) = (F(a,t),z), \forall a \in \rho, z \in \mathbb{R}, t \in [0,1]$. Then $\kappa_1$ moves to a position with the same projection as $\kappa_2$ and the same over/under information. Consequently, a vertical (with respect to $\rho$) isotopy of $\mathbb{R}^3$ brings this position exactly on $\kappa_2$. Thus $\kappa_1, \kappa_2$ are isotopic in space and by Theorem 1 their diagrams are related by the moves $\Omega_i$ for $i = 1,2,3,4,5$. In other words the $\Omega_{iso}$ move can be added to the $\Omega_i$ moves without disturbing the equivalence classes of diagrams in in $\underset{top}{\sim}$. Since it includes the $\Omega_4, \Omega_5$ moves, Theorem 1 becomes:

**Theorem 2.** *(Alexander, Briggs, Reidemeister) (version 2) Two pl (or smooth) knotted loops in $\mathbb{R}^3$ are representatives of the same knot iff their projections on a plane $\rho$ come one from the other after a finite number of consecutive alterations of types $\Omega_1, \Omega_2, \Omega_3$, each one in the interior of a bounded 2-disk $U$ in the way shown in Figure 2, along with some isotopy moves $\Omega_{iso}$ of the whole plane.*

**1.5. Our setting, $\Omega_{iso\sigma}$ moves on the sphere.** Our point of view for the study of knotted loops in space is the usual one with a slight preference to homogeneity:

The ambient space is $\mathbb{R}^3$. The knotted loops are smooth or pl, oriented and with with regular projections (diagrams) on a fixed plane $\rho$. Two loops are equivalent whenever one become the other after an isotopy of $\mathbb{R}^3$. We repeat that isotopies here are considered orientation preserving starting with the identity. Knots in $\mathbb{R}^3$ are the classes of equivalent loops. The Alexander-Briggs-Reidemeister Theorem transforms equivalence of loops in space to equivalence of their diagrams on the plane. This equivalence on the set $\mathbb{D}$ of diagrams is generated as in the second version of the Theorem, i.e. by the moves $\Omega_1, \Omega_2, \Omega_3, \Omega_{iso}$. The first two come into two versions $\Omega_1^+, \Omega_1^-$ and $\Omega_2^+, \Omega_2^-$ respectively.

In what follows, we work almost exclusively with the $\Omega_{iso}$ moves. But at some instances it will pay off if we go a step further and replace them by a set of isotopy moves happening on the 2-sphere instead of the plane. The resulting diagram is always considered on the plane:

We recall the bigger space $S^3 = \mathbb{R}^3 \cup \{\infty\}$. And we consider our diagrams on $\rho$ as also lying on the extended plane $S^2 = \rho \cup \{\infty\}$ which topologically is a 2-sphere. It is not hard to prove that two diagrams on $\rho$ are equivalent via the $\Omega_1, \Omega_2, \Omega_3, \Omega_{iso}$ moves iff they are equivalent via the $\Omega_1, \Omega_2, \Omega_3$ moves on the plane and a version $\Omega_{iso\sigma}$ of $\Omega_{iso}$ in which a diagram can be isotoped to another place on $\rho$ via an isotopy of $S^2$ (which starts from the identity). Indeed:

During an isotopy of $\rho$ bringing $K$ to $K'$, the first diagram moves in a bounded region of the plane, so the isotopy can be replaced by another one which is constant outside this bounded region and still brings $K$ to $K'$. This new isotopy can then be extends to an isotopy of $S^2$ by defining at all times the point $\infty$ as image of $\infty$. So if two diagrams on $\rho$ are equivalent via the $\Omega_1, \Omega_2, \Omega_3, \Omega_{iso}$ moves, they are also equivalent via the $\Omega_1, \Omega_2, \Omega_3$ moves on the plane and the $\Omega_{iso\sigma}$ moves on $S^2$.

Conversely, let an isotopy $F$ of $S^2$ bring the diagram $K$ of the plane to the diagram $K'$ of the plane. And let $\kappa$ be a loop in $\mathbb{R}^3$ whose diagram on $\rho$ is $K$. We write $\mathbb{R}^3 = \rho \times \mathbb{R}$. Then $S^3 = (\rho \times \mathbb{R}) \cup \{\infty\}$. We extend $F : \mathbb{S}^2 \times [0,1] \to \mathbb{S}^2 \times [0,1]$ to a map $G : S^3 \times [0,1] \to S^3 \times [0,1]$ as: $G((a,z),t) = (F(a,t),z)$ for $a \in \rho, z \in \mathbb{R}, t \in [0,1]$ and $G((\infty,t)) = \infty$ for $t \in [0,1]$. This map is a pl and smooth isotopy of $S^3$ starting from the identity. $\kappa$ is isotoped by $G$ to a loop $\kappa'$. The projection of $\kappa'$ on $\rho$ as an oriented closed curve is $K'$. Also, the over/under information at the double points of this projection are clearly those of the diagram $K'$. Thus $K'$ is the diagram of $\kappa'$ on $\rho$. Since the loops $\kappa, \kappa'$ are connected via $G$, they represent the same knot in $S^3$. We already mentioned (close to the end of §1.1) that then they also represent the same knot in $\mathbb{R}^3$. So then by Theorem 2 their diagrams $K, K'$ on $\rho$ are equivalent via $\Omega_1, \Omega_2, \Omega_3, \Omega_{iso\sigma}$ moves. Hence if two diagrams on $\rho$ are equivalent via the $\Omega_1, \Omega_2, \Omega_3$ moves on $\rho$ and the $\Omega_{iso\sigma}$ move in $S^2$, they are also equivalent via the $\Omega_1, \Omega_2, \Omega_3, \Omega_{iso\sigma}$ moves on $\rho$.

As a result we get:

**Theorem 3.** *(Alexander-Briggs-Reidemeister) (version 3) Two pl (or smooth) knotted loops in $\mathbb{R}^3$ are representatives of the same knot iff their projections on a plane $\rho$ come one from the other after a finite number of consecutive alterations of types $\Omega_1, \Omega_2, \Omega_3$, each one in the interior of a bounded 2-disk $U$ of $\rho$ in the way shown in Figure 2, along with some isotopy moves $\Omega_{iso\sigma}$ on the 2-sphere $S^2 = \rho \cup \{\infty\}$; here $S^3 = \mathbb{R}^3 \cup \{\infty\}$ is the 3-sphere considered as the one-point compactification of $\mathbb{R}^3$.*

**Remark 1.** *An isotopy $\Omega_{iso\sigma}$ on the sphere is a kind of long-range move for diagrams on a plane. Since its resulting diagram lies again on the plane, it can be thought as happening inside a 2-disk of the plane; of course then it is a new move on the plane and not an isotopy on the sphere.*

*An alteration of a diagram $K$ outside a 2-disk $U$ as shown in Figure 3 is justified immediately as the result of an $\Omega_{iso\sigma}$ move: fix $U$ point-wise, move arc $\ell$ in the exterior of $U$ far away to the left to go over $\infty$, and then moving again in the exterior of $U$ bring it to the position $\ell'$. Of course one can prove this without the help of the new move: consider a loop $\kappa$ in $\mathbb{R}^3$ whose diagram on $\rho$ is $K$. Say $s$ is the arc of $\kappa$ projecting onto $\ell$. Push $s$ parallel to $\rho$ to the two sides of the rest of $\kappa$ away from it (top and bottom in the Figure), then push it almost vertically with respect to $\rho$ way higher than the rest of the points of $\kappa$ and finally push it parallel to $\rho$ to pass it over the back (left in the Figure) of $\kappa$. This way, the final position of $\kappa$ has $K'$ as its diagram on $\rho$.*

*This argument can be generalized to prove that considering $\Omega_{iso\sigma}$ as a move happening totally in the interior of a 2-disk on the plane, produces a diagram equivalent in $\underset{top}{\sim}$ to the given one without the help of $S^2$ and Theorem 3.*

*In §4.4 we shall construct some more long-range moves $\Omega_{1\gamma}^-, \Omega_{2\gamma}^-, \Omega_{3\gamma}$ which include $\Omega_1^-, \Omega_2^-$ and $\Omega_3$ moves as special cases respectively.*

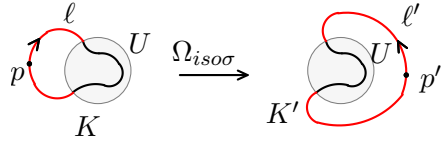

Figure 3. An example of an $\Omega_{iso\sigma}$ move.

Let us notice that unlike the $\Omega_{iso}$ moves, and although the resulting diagram of an $\Omega_{iso\sigma}$ move lies inside some disl on the plane, it is not always true that the isotopy of the move happens to some 2-disk of $\rho$ or even of on $S^2$ fixing throughout all points outside it: the places of a diagram (a compact set) during an isotopy (a continuous map) on $S^2$ (a compact set), can literally occupy the hole $S^2$!

We shall use the $\Omega_{iso\sigma}$ moves to our benefit much farther below in the proof of Lemma 22.

It seems that the power of the $\Omega_{iso}^s$ moves in our study of knots springs from the example in Figure 3. The gist of it is that although there exist exactly two oriented diagrams of a circle on the planecnon-equivalent by $\Omega_{iso}$ moves, these two become equivalent by an $\Omega_{iso}^s$ move on the sphere, thus there exists just a unique oriented diagram of a circle on the sphere.

## 2. Based symbols

We work in the set $\mathbb{D}$ of all regular oriented smooth (pl) diagrams $K$ on a plane $\rho$, which come as projections on $\rho$ of the oriented smooth (pl) loops $\kappa$ in $\mathbb{R}^3$.[1]

As discussed in §1.2 the elements of $\mathbb{D}$ are the smooth (pl) closed curves $K = f(S^1), f : S^1 \to K \subset \rho$ on the plane $\rho$ with only a finite number of self-intersections (possibly none) which we call as crossings, all of them double points in which the curve intersects itself transversely. The orientation of $K$ is the one of $S^1$ induced to it by $f$, or equivalently, the one of the loop $\kappa$ induced to it by the projection $pr : \mathbb{R}^3 \to \rho$. For $\kappa = g(S^1), g : S^1 \to g(S^1) \subset \mathbb{R}^3$ it is $K = (pr \circ g)(S^1), f = pr \circ g$.

Below we give a bunch of definitions aiming at a better description of the $\Omega$ moves. We can think that $\Omega_1^+, \Omega_2^+$ moves actually happen in one of the regions in which $K$ partitions the plane. Below we call these regions as rooms of the diagram. And similarly, the $\Omega_1^-, \Omega_2^-$ can be though as happening in a room which they make disappear, whereas $\Omega_3$ as a move which rearranges the position of a room and the borders of some neighboring to it rooms.

In a nutshell, a great deal of the definitions and results we shall deal with in the sequel originate from the following observations: the boundaries of the rooms are arcs of $K$ and each room together with its boundary is a closed region on the plane. The closed rooms are bounded regions except for one of them which is unbounded. A bounded closed room $\Pi$ can be described as a 2-disk pinched in some of its boundary points so that they become double points at which the boundary $\partial\Pi$ does not cross itself. So these pinched boundaries as closed parametrized curves are not arcs of $K$, although they fully lie on it. As closed curves, they posses two orientations. A fixed orientation on $\partial\Pi$ induces to its boundary arcs an orientation which for some arcs may coincide with the one they have from $K$ but for some other arcs it may not coincide. These remarks are valid even for the unbounded closed room except from the fact that topologically it is a pinched 2-disk without an interior point. The distinction between bounded and unbounded regions disappears when we consider the diagram on the 2-sphere instead on the plane.

---

[1] For a plethora of objects defined in the paper we use the symbols $\alpha, \delta, \phi, \tau, \upsilon, \pi, \Delta, \theta, \rho, \sigma, \epsilon, \gamma, \mu$ coming from the first letters of the greek words $\alpha\rho\iota\sigma\tau\epsilon\rho\acute{o}\varsigma$=left, $\delta\epsilon\xi\iota\acute{o}\varsigma$=right, $\phi o\rho\acute{\alpha}$= direction, geometric sense, $\tau\acute{\alpha}\xi\eta$= order, $\acute{\upsilon}\psi o\varsigma$= height, $\pi\rho\acute{o}\sigma\eta\mu o$= sign, $\delta\iota\alpha\sigma\tau\alpha\acute{\upsilon}\rho\omega\sigma\eta$= crossing, $\theta\acute{\epsilon}\sigma\eta$= place, $\rho o\pi\acute{\eta}$ in $o\mu\acute{o}\rho\rho o\pi o\varsigma$=same (geometric) sense, $\sigma\phi\alpha\acute{\iota}\rho\alpha$=sphere, $\epsilon\acute{\iota}\delta o\varsigma$=kind, $\gamma\epsilon\nu\iota\kappa\acute{o}\varsigma$=general, generic, $\mu\epsilon\gamma\acute{\alpha}\lambda o\varsigma$=big.

## 2.1. Rooms, orientations of arcs, basic arcs, main arcs, cycles.

**Definition 6.** *Let $K = f(S^1), f : S^1 \to K$ be a diagram on the plane $\rho$.*

*We have already discussed (§1.2) the way to judge positive and negative orientations on $\rho$ from our fixed position in one of the two half-spaces with respect to $\rho$: if a frame $(e_1, e_2)$ on $\rho$ is completed to a positive frame $(e_1, e_2, e_3)$ of $\mathbb{R}^3$ by a vector $e_3$ with first point on $\rho$ and final point on the same half-space with us, then we say that $(e_1, e_2)$ is a positive frame or orientation on $\rho$.*

*We denote the positive and the negative orientation of the plane $\rho$ by $\alpha, \delta$ respectively. We denote an arbitrary orientation by $\phi \in \{\alpha, \delta\}$. We write $\delta = -\alpha, \alpha = -\delta$ and as usual we call the two orientations as opposite to each other. We give a sign $-, +$ or rather a number $-1, +1$ to the orientations as: $\pi_\alpha = +1, \pi_\delta = -1$.*

*Let $s = t_1 t_2$ be a non trivial arc ($t_1 \neq t_2$) of $S^1$ with the induced orientation of $S^1$ and let $t_1, t_2$ be the first and second endpoint of $s$ respectively. As usually for parametrized curves like $K$, we say that $e = f(s)$ is an arc of $K$, considered oriented by the induced orientation of $K$, or equivalently oriented by the induced orientation of $s$ via the map $f$. The first and second endpoint of $e$ is $p = f(t_1)$ and $q = f(t_2)$ respectively. We are going to use the notation $e = \overrightarrow{pq}$ to denote $e$ oriented like this.*

*It can be $p = q$, in which case $e$ is actually a loop on the plane, possibly with self intersections parametrized by $f$ and considered as a closed curve with beginning and ending point the point $\Delta = p = q$. This point $\Delta$ is nothing but a crossing point of $K$. We continue calling $e = \overrightarrow{pq}$ as an arc of $K$.*

*We denote the arc $e$ with the opposite orientation as $-e = -\overrightarrow{pq}$.*

*We also denote $-(-e) = e$ and we say that $e$ and $-e$ are opposite to each other.*

*If $K$ has $n \in \mathbb{N}$ crossings we write $\tau(K) = n$ and we call $\tau(K)$ as the order of $K$.*

An alternative appropriate notation for $-e = -\overrightarrow{pq}$ could be $-e = \overleftarrow{pq}$ but we'll not make use of it.

**Definition 7.** *For a diagram $K$ on the plane $\rho$, the set $\rho - K$ is the union of a finite number of disjoint, open, connected regions which we call as the open rooms of $K$. We call the closures $\Pi_1, \Pi_2, \ldots, \Pi_k$ of the open rooms as closed rooms or just rooms of $K$. The open regions are the interiors of the rooms and together with $K$ they form a partition of the plane into disjoint connected sets. We also say that the rooms are the closed subsets into which $K$ divides the plane.*

*For $\tau(K) = 0$, $K$ is homeomorphic to a circle ($K \approx S^1$) and by Jordan's Theorem there exist exactly two rooms, each with $K$ as its boundary.*

*For $\tau(K) > 0$, the boundary $\partial \Pi_i$ of any room $\Pi_i$ is the union of arcs $e$ of $K$, each extending from a crossing point $f(t_1)$ to another crossing point $f(t_2)$ and containing no other crossings in its interior. Whenever $f(t_1) = f(t_2)$ the arc is actually a loop. We call the arcs $e$ with the induced orientation of $K$ as basic arcs of $K$.*

*We call the opposite $-e = -\overrightarrow{pq}$ of a basic arc $e = \overrightarrow{pq}$ as opposite or negative basic arc of $K$.*

*We call basic and negative basic arcs collectively as main arcs.*

*We indicate the difference of the orientations between basic and opposite basic arcs by the sign, or rather the signed number of the arcs defined as: $\pi_e = 1, \pi_{-e} = -1$, for any basic arc $e$.*

Notice that we do not define basic arcs for $K \approx S^1$ that is when $\tau(K) = 0$. According to our terminology, for a basic arc $e$ of a diagram $K = f(S^1)$, it is $e = \overrightarrow{pq} = \overrightarrow{f(t_1)f(t_2)} = f(t_1 t_2)$ where $p = f(t_1), q = f(t_2)$ is respectively the first and second endpoint of $e$ and $s = t_1 t_2$ is an oriented arc of $S^1$ with the induced orientation of $S^1$ ($t_1, t_2$ are the first and second endpoints of $s$ respectively).

The basic arcs $e$ of $K$ are nothing else than the simple pieces of $K$ between its crossings, considered oriented by the induced orientation of $K$.

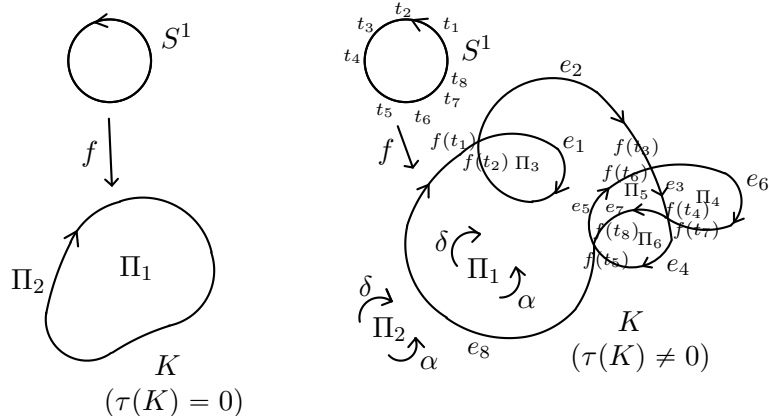

FIGURE 4. Rooms and simple arcs of a diagram $K$. For the diagram on the right it is for example: $(\Pi_1)_{\delta e_8}$, $(\Pi_2)_{\alpha e_8}$ and $(\Pi_1)_{\alpha(-e_8)}$, $(\Pi_2)_{\delta(-e_8)}$.

The boundaries of the rooms are closed curves when considered as point sets bounding the rooms. These curves have self-intersections, but not transversal ones (Figure 5). So considered this way, $\partial\Pi$ are not parametrized curves as subarcs of the parametrized curve $K$.

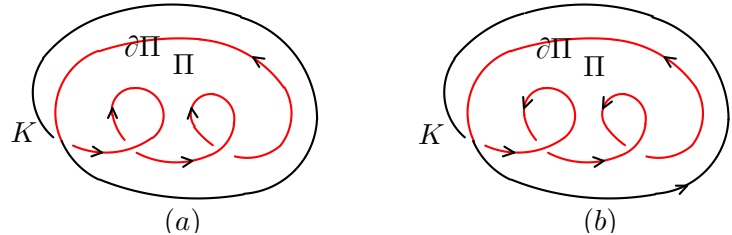

FIGURE 5. $\Pi$ is a room of a diagram $K$. The boundary $\partial\Pi$ of $\Pi$ is the red curve. (a) $\partial\Pi$ as a bounding closed curve, depicted with one of its two possible orientations. (b) $\partial\Pi$ as a closed arc of $K$ depicted with the induced orientation of $K$.

Some more details about the boundaries of the rooms are given in the following Lemma:

**Lemma 1.** *If $\Pi$ is a room of the diagram $K$ then its boundary $\partial\Pi$ as a point set is a closed curve. Removing from $\Pi$ some small disjoint open canonical neighborhoods of all the crossings that might lie on $\partial\Pi$, remains a set $D_\Pi$ which is a 2-disk or the closed complement of a 2-disk on the plane.*

*Proof.* For $\tau(K) = 0$, it is $K \approx S^1$ with no crossings and the result is immediate.

For $\tau(K) > 0$, in some small disjoint canonical neighborhoods $U$ for all crossings of $K$ lying on $\partial\Pi$ we perform a change on $\Pi$ as shown in Figure 6. We end up with a connected region $\Pi'$ in place of $\Pi$, whose boundary $c$ is a simple closed curve on the plane, thus a homeomorphic image of $S^1$. For an arbitrary basic arc $e$, let us denote by $e'$ be its connected part left in $\partial\Pi'$, i.e. its part outside the neighborhoods. We consider such a subarc just as a point set without caring about orientations, and we call it as a truncated basic arc.

Fix an orientation of $\partial\Pi'$ and make a trip along this orientation meeting truncated basic arcs in the order $e'_1, e'_2, \ldots e'_k, e'_1$. Let this trip starts and ends at an interior point $A$ of $e'_1$. Now repeat the trip, but inside each neighborhood $U$ replace the part of $\partial\Pi'$ by the respective part of the basic arcs depicted in Figure 6 ($\sigma_1$ is replaced by $s_1 \cup s_3$ and $\sigma_2$ is replaced by $s_2 \cup s_4$). The new trip is performed on $\partial\Pi$ and in it we meet basic arcs (which we consider just as point sets with no care about their orientations) in the order $e_1, e_2, \ldots e_k, e_1$ finishing at our starting point $A$. So the boundary $\partial \Pi$ of $P$ is indeed a closed curve.

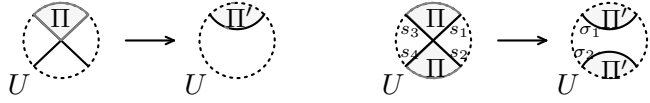

FIGURE 6. Changes on a canonical neighborhood $U$ of a crossing on the boundary of a room $\Pi$. A choice of such a change for each crossing on $\partial \Pi$ transforms the room to a connected region $\Pi_1$ with boundary a simple closed curve on the plane.

By construction $\Pi'$ results by throwing away from $\Pi$ some small neighborhoods of the crossings on $\partial \Pi$ (Figure 7). Since $\partial \Pi'$ is a circle, by Schonflies Theorem $\Pi'$ is either a 2-disk or the closed complement of a 2-disk thus $\Pi'$ is one of the $D_\Pi$'s claimed to exist. $\square$

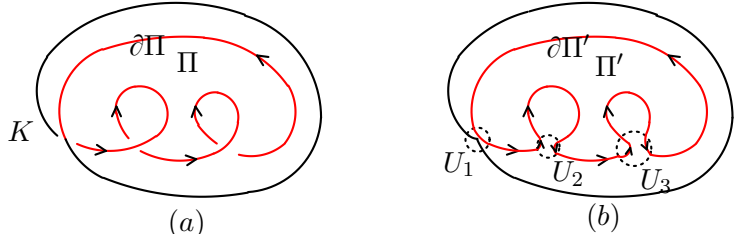

FIGURE 7. (a) A room $\Pi$ of a diagram $K$. The boundary $\partial \Pi$ of $\Pi$ is the red closed curve depicted with one of its two orientations. (b) $U_i$'s are small canonical neighborhoods of the crossings of $K$ that lie on $\partial \Pi$. $\Pi'$ is the result of modifying $\Pi$ inside the neighborhoods $U_i$ as shown in Figure 6. $\partial \Pi'$ is a circle and $\Pi'$ is a 2-disk $D_\Pi$. $\partial \Pi'$ is depicted with the induced orientation of $\Pi$ (the two orientations induce to the common arcs of $\Pi, \Pi'$ the same orientation).

**Definition 8.** *For any room $\Pi$ of a diagram $K$ we call a $D_\Pi$ as in Lemma 1, as a truncated disk of $\Pi$. For $\tau(K) = 0$, there are no crossings on any of its two rooms $\Pi, P$ and we set $D_\Pi = \Pi, D_P = P$ still calling the as the truncated disks of $K$.*

*For any room $\Pi$ of a diagram $K$, we consider the two orientations which the open region $\Pi - \partial \Pi$ of the plane can have, as the two orientations of $\Pi$. We denote $\Pi$ oriented with an orientation $\phi \in \{\alpha, \delta\}$ as $\Pi_\phi$.*

*$\Pi_\phi$ induces to $S = \partial \Pi - crossing$ points on $\partial \Pi$, an orientation in the obvious way. We extend this orientation of the open arcs on $S$ to the closed arcs on $\partial \Pi$ as well, and we get an orientation of the closed by Lemma 1 curve $\partial \Pi_\phi$. We say that this is the boundary orientation of $\partial \Pi$ induced by the orientation of $\Pi_\phi$. We denote $\partial \Pi$ oriented this way as $\partial \Pi_\phi$.*

*For $\tau(K) = 0$, it is $\partial \Pi = K$. The induced orientation by $\Pi_\phi$ may or may not coincide with the one of $K$, and we denote this fact as $\Pi_{\phi K}, \Pi'_{-\phi K}$ respectively. For $\tau(K) > 0$, $\Pi_\phi$ induces on any basic arc $f = e$ or on any opposite basic arc $f = -e$ of its boundary, an orientation which may or may not coincide with its own orientation as an arc of $K$. If the two orientations coincide, we say that $\Pi_\phi$ is a $\phi$-room for $f$ and $\partial \Pi_\phi$ is a $\phi$-cycle for $f$, and we write $\Pi_{\phi f}$ and $\partial \Pi_{\phi f}$. Otherwise we say that that $\Pi_\phi$ is a $-\phi$-room for $f$ and $\partial \Pi_\phi$ is a $-\phi$-cycle for $f$ in $\Pi_\phi$, and we write $\Pi_{-\phi f}$ and $\partial \Pi_{-\phi f}$. We call $\partial \Pi_\phi, \partial \Pi_{-\phi}$ as oriented cycles, and the boundary $\partial \Pi$ without any specific mention to an orientation just as a cycle of $K$.*

*For $\tau(K) = 0$ the cycles of $K$ are $K$ with its orientation as a diagram, and also the point-set $K$ with the opposite orientation, a diagram we denote as $-K$. For $\tau(K) > 0$, we denote a $\phi$-cycle also as $C = (f_1, f_2, \ldots, f_n)_{cycl}$, where the $f_i$'s are the consecutive main arcs (basic or negative basic*

arcs) forming the cycle as an oriented closed curve. In this notation the index cycl indicates we care only about the cyclic order of the arcs, so $C = (f_1, f_2, \ldots, f_n)_{cycl} = (f_2, f_3, \ldots, f_n, f_1)_{cycl} = \cdots$

As we have set $-\alpha = \delta$ and $-\delta = \alpha$, the above imply that for $\phi \in \{\alpha, \delta\}$ the notation $\Pi_{\phi f}$ (or $\Pi_{\phi K}$) always means that $\phi$ is the orientation of $\Pi$ which induces on its boundary arc $f$ the same orientation as naturally $f$ has: that of $K$ for $f$ a basic arc, the opposite of $K$ for $f$ an opposite basic arc. And $\partial \Pi_{\phi f}$ always means that the oriented closed curve $\partial \Pi_\phi$ contains $e$ as an oriented subarc, that is, as an arc whose induced orientation by $\partial \Pi_\phi$ coincides with the natural orientation of $f$.

Now, a truncated disk $D_\Pi$ of a room $\Pi$ is a 2-disk or the closed complement of a 2-disk. So any orientation $\phi$ of $D_\Pi$ induces an orientation to its boundary $\partial D_\Pi$ which is just a circle. This will be $\phi$ whenever $D_\Pi$ is a 2-disk, whereas it will be $-\phi$ whenever $D_\Pi$ is the closed complement of a 2-disk.

Also, $\Pi$ induces anyone of its orientations $\phi$ to its boundary $\partial \Pi$ which is a closed curve.

$D_\Pi$ as a subset of $\Pi$ has its interior $(D_\Pi)^0$ in the interior $(\Pi)^0$ of $\Pi$. So the open set $(D_\Pi)^0$ gets its two orientations on the plane as induced orientations as a subset of the two orientations of the open set $(\Pi)^0$. So the following lemma does not come as a surprise.

**Lemma 2.** *A room $\Pi$ and anyone of its truncated disks $D_\Pi$ equipped with the same orientation $\phi \in \{\alpha, \delta\}$, induce on the common part of their boundaries the same orientation.*

*Proof.* For an arc $f$ of $\partial \Pi \cap \partial D_\Pi$, a small open neighborhood $U_f$ of $f$ in $D_\Pi$ is also an open neighborhood of $f$ in $\Pi$. So the $\phi$ orientation of the interior of $\Pi$ and its induced orientation in the interior of its subset $D_\Pi$, both induce the same orientation on their common open subset $U_f - f$, thus they both induce the same orientation on the boundary arc $f$. □

Returning to the relation of the orientations of rooms and the arcs on their boundaries, it is clear that reversing the orientation of a room reverses the orientation of its boundary. Thus the orientation of any point-set of some arc which lies on the boundary is reversed as well. So any main arc $f$ which is part of the oriented boundary of the room with the first orientation, reverses its orientation becoming $-f$ and is part of the oriented boundary for the second orientation of the room. So:

**Lemma 3.** *For a room $\Pi$ and a main arc $f$ of the diagram $K$ it holds: $\Pi_{\phi f} \Leftrightarrow \Pi_{(-\phi)(-f)}$.*

Let us note that each basic arc $e$ (or an opposite basic arc $-e$) belongs to exactly two rooms, say $\Pi, \Pi'$ and that each one of the rooms equipped with one of its orientations is an $\alpha$-room for $e$ while equipped with the other it is a $\delta$-room for $e$. In order for $\Pi_{\phi_1}, \Pi'_{\phi_2}$ to be the same kind of room for $e$ (or $-e$), they must be oriented in opposite ways, i.e. $\phi_1 = -\phi_2$. The same argument holds for negative basic arcs $-e$ as well. So:

**Lemma 4.** *For the two rooms $\Pi, \Pi'$ containing a main arc $f$ of the diagram $K$ it holds: $(\Pi_{\phi_1 f}, \Pi'_{\phi_2 f})$ $\Leftrightarrow \phi_1 = -\phi_2$.*

Now for two main arcs $e, f$ which as point-sets are parts of the boundary of the same room $\Pi$, we can compare their orientations by comparing each one against an orientation of $\Pi$, or equivalently of $\partial \Pi$:

**Definition 9.** *If for the main arcs $e, f$ and a room $\Pi$ of some diagram it is $\partial \Pi_{\phi e}$ and $\partial \Pi_{\phi' f}$, we say that $e, f$ are oriented in the same sense (or that they have the same sense) in $\partial \Pi$ whenever $\phi = \phi'$, and we say they are oriented in opposite sense (or that they have the opposite sense) in $\partial \Pi$ whenever $\phi = -\phi'$.*

Since every main arc lying as a point-set on the boundary of a room has exactly one of the two induced orientations of this boundary, we can also say that $e, f$ have respectively the same or opposite sense in $\Pi$.

In other words, for $e, f$ which as point-sets lie on the boundary of the same room we say they have the same or opposite sense depending on if they both become parts of the same oriented boundary of the room or not. Notice that we do not define the notion of same or opposite sense for arcs that as point sets do not both lie on the boundary of the same room. Being able to compare the kind of sense for the two main arcs $e, f$ allows us to make similar comparisons for any arc in the set $\{e, -e\}$ with any arc in the set $\{-f, f\}$ since by Lemma 3 we have:

**Lemma 5.** *If for the main arcs $e, f$ the notion of same or opposite sense of orientation is defined, then $-e, -f$ have the same kind of orientation sense as $e, f$, whereas $-e, f$ and $e, -f$ have the other kind of orientation sense than that of $e, f$.*

Taking into account Lemma 5, we can restrict our effort for comparisons only to basic arcs instead of to main arcs in general.

Let us now notice that each one of two given basic arcs $e, f$ as an oriented arc lies on the oriented boundary of exactly two rooms. So we can compare the orientations of $e, f$ at most two times. The two times case happens whenever both arcs as point-sets lie in common on the boundary of two rooms. According to Lemma 6 that follows, the same or opposite way of sense for them does not change for both rooms. So the definition below has a meaning:

**Definition 10.** *We say that two basic arcs $e, f$ of the diagram $K$ have the same or opposite boundary sense in $K$ whenever they are respectively oriented in the same sense or opposite sense in the boundary of some room of $K$. We also say respectively that they are similarly or oppositely boundary oriented in $K$. And we write respectively $e \upuparrows f$ and $e \downuparrows f$.*

*We also define then the number of similarity of boundary sense as $\rho_{e,f} = 1$ whenever they have the same boundary sense in $K$ and $\rho_{e,f} = -1$ whenever they have opposite boundary sense in $K$.*

By the very definition, the order of $e, f$ is irrelevant, i.e. $e, f$ have the same boundary sense iff $f, e$ have the same boundary sense and it holds $\rho_{e,f} = 1 \Leftrightarrow \rho_{f,e} = 1$. Note that for two basic arcs that do not lie on some same room, we do not define a boundary sense number, nor we call them boundary oriented either similarly or oppositely.

The definition and the accompanied notation for the number of sense do not involve any mention to a specific room that contains both $e, f$, since as we said the exact room is irrelevant:

**Lemma 6.** *If both basic arcs $e, f$ of the diagram $K$ lie as point-sets on the boundaries of two rooms $\Pi, P$ of $K$, then they have the same sense on both boundaries of the rooms or else they have opposite sense on both boundaries. In terms of the notation of Definition 8: if $\partial \Pi_{\phi e}, \partial \Pi_{\psi f}, \partial P_{\phi' e}, \partial P_{\psi' f}$ then $\phi = \psi \Leftrightarrow \phi' = \psi'$.*

*Proof.* Let $D_\Pi, D_P$ be truncated disks of $\Pi, P$ produced using the same canonical neighborhoods $U_i$ on the common crossing points of $\Pi, P$ (Figure 8). And let $e_1, f_1$ be the connected subarcs of $e, f$ remaining on the boundaries of both $D_\Pi, D_P$. We consider $e_1, f_1$ oriented by the induced orientations of $e, f$.

By definition, $e, f$ have the same sense on $\partial \Pi$ iff some orientation of $\partial \Pi$ induces to both $e, f$ the orientations they have as basic arcs. This happens iff some orientation of $\partial \Pi$ induces to both $e_1, f_1$ the orientations they have as subarcs of the oriented basic arcs $e, f$. By Lemma 2 this happens whenever for some orientation $\phi$ of $\partial D_\Pi$, the induced orientation on $e_1, f_1$ is the one they have as subarcs of the basic arcs $e, f$. The induced orientation of $\partial D_\Pi$ makes $e_1, f_1$ to have the same sense in the usual way for arcs on a circle (here $\partial D_\Pi$).

And similarly, $e, f$ have the same sense on $\partial D_P$ iff $e_1, f_1$ have the same sense in the usual way for arcs on a circle (here $\partial D_P$).

So we wish to show that the arcs $e_1, f_1$ with a given orientation (namely oriented as subarcs of the basic arcs $e, f$), have the same sense on $\partial D_\Pi$ iff they have the same sense on $\partial D_P$. This is not hard to show:

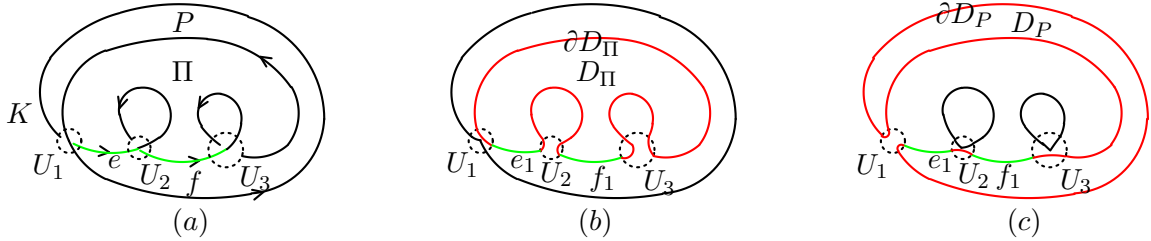

FIGURE 8. The basic arcs $e, f$ of the diagram $K$ lie on the boundary of two rooms $\Pi$ and $P$ of $K$. $D_\Pi, D_P$ are truncated disks of $\Pi, P$ which used the same canonical neighborhoods $U_i$ for all common crossing points on their boundaries. $e_1, f_1$ are the parts of $e, f$ outside the $U_i$'s. $e_1, f_1$ lie on both $\partial D_\Pi, \partial D_P$. In general, $e, f$ have the same or opposite sense on $\partial \Pi$ iff they have the same or opposite sense respectively on $\partial P$. In this figure they have the same sense on both $\partial D_\Pi, \partial D_P$. In general, $e, f$ have the same sense on $\partial \Pi$ iff $e_1, f_1$ have the same sense on the circle $\partial D_\Pi$. Similarly, $e, f$ have the same sense on $\partial P$ iff $e_1, f_1$ have the same sense on the circle $\partial D_P$.

At least one of $D_\Pi, D_P$ is a 2-disk. Say $D_\Pi, D_P$ is a 2-disk. Let us assume that $e_1, f_1$ have an opposite sense in $\partial D_\Pi$. Say $a, b$ be the first and second point of $e_1$ and $c, d$ the first and second point of $f_1$. And let us move on $\partial D_\Pi$ (Figure 9) along its orientation starting at $b$ until we arrive at $f_1$. Say $x$ be the arc we have traveled on. The opposite of sense orientation makes us arrive at $d$. Then exiting at $c$ we need to eventually arrive at $a$. Nevertheless, if we call $z$ the arc on $\partial D_\Pi$ between $b$ and $c$ as we move along the curve's orientation, then $\gamma = x \cup f_1 \cup z$ is a simple closed curve, so by Jordan's Theorem it has an interior $D$ (a 2-disk). Since $x$ is outside $D_\Pi$, the disk $D$ has to be outside $D_\Pi$ as well. If we move on the boundary of $D_P$ after point $c$ along its orientation, say on an arc $\ell$, we need to return to point $a$. But the boundary of $D_P$ except for the arcs $e_1, f_1$ is in the exterior of $D_\Pi$. So $\ell$ near $c$ has to enter in the interior $D$ of $\gamma$. In order to arrive then at the exterior point $a$ we have to cross the boundary circle $\gamma$. Since we cannot cross arcs $e_1, f_1$ again we should either meet arc $z$ or else arc $x$ (by Jordan's Theorem). The first cannot happen because $D_\Pi, D_P$ have no common interior points. The second cannot happen either since the curve $\partial D_P$ we move on is a circle. □

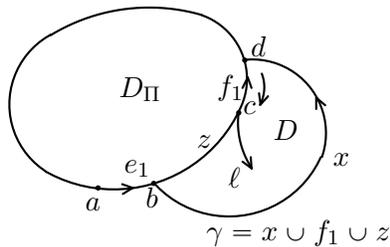

FIGURE 9. Moving on the arc $\ell$ along its orientation and in order to arrive at $a$ avoiding $e_1, f_1$ we need either to enter the 2-disk $D_\Pi$ or else to cross $x$ in some point.

Now, when we perform a trip on the boundary of a room along one of the two possible orientations, every time we arrive at a vertex we have to leave the edge $e$ we are currently walking on, to turn as the orientation dictates, and then to continue walking on a new edge $e'$. Orienting $e$ and $e'$ with the orientation of our trip, each becomes a basic arc or a negative basic arc, not necessarily both of the same kind. We want to say that $e'$ comes after $e$ in such a trip sos we give an official definition which we'll find useful later on (§2.4):

**Definition 11.** *Let $K$ be a diagram with $\tau(K) > 0$ and $\Pi$ be a room of $K$. For $\Pi_\phi$ ($\Pi$ oriented by $\phi \in \{\alpha, \delta\}$), the $\phi$-cycle $\partial\Pi_\phi$ (the oriented boundary of $\Pi_\phi$) is a union of consecutive basic arcs or negative basic arcs $e_i$ (it can contain both): $\partial\Pi_\phi = e_1 \cup e_2 \cup \cdots \cup e_k$ with $\partial\Pi_{\phi e_i}, \forall i$ and final endpoint of $e_i$=first endpoint of $e_{i+1}$ (indices mod $k$). We say that each $e_{i+1}$ is the $\phi$ turn of $e_i$ at their common endpoint and we write $\phi(e_i) = e_{i+1}$.*

Some more information regarding the rooms of a diagram can be found in [4].

### 2.2. Definition of based symbols.

**Definition 12.** *Let $K$ be a diagram with $n \in \mathbb{N}$ crossings.*

*We call any point of $K$ other than a crossing as a base point or origin of $K$.*

*For a chosen base point $p$ we perform a trip around $K$ following its orientation until we return to $p$. We name the crossing points as we meet them successively as $1, 2, \ldots$ This way each crossing $\Delta$ gets two distinct names, say $i$ and $j$ which we also call as its labels. It gets one of them, say $i$ as we cross it in a canonical neighborhood walking on an over-crossing (§1.2), while it gets the other name $j$ when we meet $\Delta$ walking on an under-crossing. Then we denote $\Delta = (i, j)$: no matter if $\Delta$ gets the label $i$ in our first visit or in the second, the number $i$ is always the first entry of the pair.*

*We also give to $\Delta$ a sign $+, -$ or rather a signed number $\pi_\Delta = +1, -1$ in the usual way as in Figure 10: when we cross $\Delta$ walking in a canonical neighborhood $U$ along the over-crossing we have to turn either positively or negatively on $\rho$ in order to find ourselves walking on the under-crossing arc along the orientation of $K$ and then we give to $\Delta$ the sign $\pi = +1$ or $\pi = -1$ respectively.*

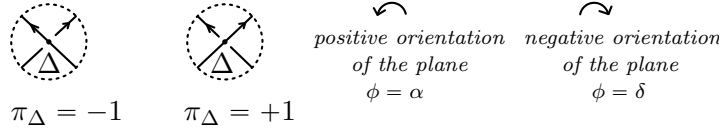

FIGURE 10. Sign $\pi_\Delta$ at a crossing point $\Delta$.

*We denote $\Delta$ with all the above information as $\Delta = (i, j)^\pi$. Whenever the order in the pair is unknown or not important we agree to write $\Delta = i : j = j : i = (i : j)^\pi = (j : i)^\pi = (i : j) = (j : i)$.*

We are ready for the main definition:

**Definition 13.** *For a diagram $K$ with $n > 0$ crossings we call as its based symbol at its base point $p$ the set $\{\Delta_1, \Delta_2, \ldots, \Delta_n\} = \{(i_1, i_2)^{\pi_1}, (i_3, i_4)^{\pi_1}, \ldots, (i_{2n-1}, i_{2n})^{\pi_n}\}$ of its crossings $\Delta_i$ equipped with all of their information. We conveniently denote the based symbol at $p$ as:*

$$_pK = (i_1, i_2)^{\pi_1} \ (i_3, i_4)^{\pi_1} \ldots (i_{2n-1}, i_{2n})^{\pi_n},$$

*where the order of the pairs is not important.*

*For $K$ the unknot, we define as its based symbol at $p$ the empty set $_pK = \varnothing$ for all $p \in K$. We call $\varnothing$ as the empty based symbol.*

*Let $\Sigma_\Delta$ be the set of all based symbols of all diagrams $K$ on the plane.*

*For any based symbol $_pK$ we write $\tau(_pK) = \tau(K) = n$ and we call $n$ as the order of $_pK$.*

*For a crossing $\Delta = (i, j)^\pi$ we define as height of $i$ and $j$ the numbers $v_i = 1, v_j = -1$. We recall that $\pi \in \{-1, +1\}$ is the sign of $\Delta$.*

Height indicates which of the two arcs in a canonical neighborhood of $\Delta$ that give the labels $i, j$ to the crossing is the over-arc and which is the under-arc. In turn, for a loop in space which projects to $K$, this indicates which of the two arcs that project inside the neighborhood is the higher over $\Delta$ and which the lower.

Recall that for a diagram $K$ with $\tau(K) = 0$ we have not defined basic arcs. On the other hand, for a diagram $K$ with $\tau(K) = n > 0$ its basic arcs according to our notation are $e_1 = \overrightarrow{12}, e_2 = \overrightarrow{23}, \ldots, e_{2n-1} = \overrightarrow{(2n-1)(2n)}, e_{2n} = \overrightarrow{(2n)1}$. Recall also that whenever the labels $i, i+1$ are given to the same crossing of $K$ then $\overrightarrow{i(i+1)}$ is actually a loop. We shall also call it as a *basic loop* of $K$.

Finally recall that if $e = \overrightarrow{i(i+1)}$ ($i, i+1$ are (mod $n$) numbers in $\mathbb{Z}_{2n}$) is a basic arc then $-e = -\overrightarrow{i(i+1)}$ is its opposite, i.e. the same arc with the opposite orientation.

Regarding the notation, let us observe that since no crossing of $K$ lies in the interior of any of its basic arcs $e = \overrightarrow{i(i+1)}$, then in case $\tau(K) > 1$ the complementary arc of $e$ contains some crossing and it cannot be basic. Hence we can even use the notation $-e = \overrightarrow{(i+1)i}$ for an opposite basic arc with no ambiguity on its meaning.

2.3. **Modifications of based symbols under topological actions.** Below we describe the way the based symbol ${}_pK$ relates to some based symbol ${}_qK'$ when we perform a change of base or an $\omega \in \{\Omega_{iso}, \Omega_1, \Omega_2, \Omega_3\}$ move $K \xrightarrow{\omega} K'$ on $K$. We need this information in order to define algebraic analogues of these topological modifications later on in §3.

We are going to deal with names of crossings and relabelings of them. So it is a good idea to fix a bit of notation to be followed throughout:

**Definition 14.** $\mathbb{N}_0 = \emptyset$ and $\mathbb{N}_{2n} = \{1, 2, \ldots, 2n\}$ for $n > 0, n \in \mathbb{N}$. For $n > 0$ we consider $\mathbb{N}_{2n}$ as the usual abelian group $(\mathbb{Z}_{2n}, +)$ but we insist on writing its elements as the natural numbers $1, 2, \ldots, 2n$, hence its new notation $\mathbb{N}_{2n}$ instead.

*We consider $\mathbb{N}_{2n}$ cyclically ordered as $1 < 2$, $2 < 3, \ldots, 2n-1 < 2n$, $2n < 1$ and with no other pair of elements related in this order.*

As an additive group it allows expressions like for example $2n + 1 = 1$, and for $a \in \mathbb{N}_{2n}$ the number $a + 1$ is again an element of $\mathbb{N}_{2n}$. For $0 < n < m$ we have $\mathbb{N}_{2n} \subset \mathbb{N}_{2m}$ and for example for $a \in \mathbb{N}_{2n}$ the number $a + 1$ is also considered as an element of $\mathbb{N}_{2m}$.

Below let $K = f(S^1), f : S^1 \to K \subset \rho$ be a diagram on the plane $\rho$, $p$ a base point of $K$ and $\phi$ the orientation of $K$.

2.3.1. *Change of base.* Let q be another base point of $K$.

If $\tau(K) = 0$ we have ${}_pK = \emptyset$ and also ${}_qK = \emptyset$.

Otherwise, let ${}_pK = (i_1, i_2)^{\pi_1}, (i_3, i_4)^{\pi_1}, \ldots, (i_{2n-1}, i_{2n})^{\pi_n}$. Let $p' = f^{-1}(p), q' = f^{-1}(q)$ be the points of $S^1$ that map to $p, q$ under the parametrization $f$ of $K$. And let $s_1 = pq, s_2 = qp$ be the two arcs of $S^1$ between $p'$ and $q'$ which have the induced orientation of $S^1$. Let also $t_1, \ldots t_{a-1}$ be the points of $s_1$ which become crossings of $K$, as we meet them moving along the orientation of $S^1$ starting at point $p'$ of $s_1$. And similarly, let $t_a, \ldots t_{2n}$ be the consecutive points of $s_2$ that become crossings of $K$. $t_i$ induces the label $i$ into the corresponding crossing of $K$.

Changing the base point of $K$ from $p$ to $q$ means that in order to encounter crossings on $K = f(S^1)$ and assign labels to them we first move on the arc $s_2$ and then on $s_1$ along the orientation of $S^1$. Hence $t_{a+}, \ldots t_{2n}$ (in this order) are the first to assign a label to a crossing and then $t_1, \ldots t_{a-1}$ (in this order) will be the next ones. This amounts to the following changes of names to the crossings of $K$:

(2.1)
| base is $p$: | 1 | 2 | $\ldots$ | $a-1$ | $a$ | $a+1$ | $\ldots$ | $2n$ |
|---|---|---|---|---|---|---|---|---|
| base is $q$: | $2n+2-a$ | $2n-a+1$ | $\ldots$ | $2n$ | $t$ | $2$ | $\ldots$ | $2n+1-a$ |

In other words ${}_qK$ comes from ${}_pK$ after the change of labels in the pairs of ${}_pK$ by the correspondence shown in the above table. The signs for each pair remains the same.

So if we call $\mu$ the relabeling function corresponding the numbers in the first row to those of the second, the pairs in ${}_pK$ are replaced as:

$$(x, y)^\pi \mapsto (\mu(x), \mu(y))^\pi.$$

And ${}_pK$ transforms to ${}_qK$ as:

(2.2) $$\quad {}_pK = \cdots (x, y)^\pi \cdots \mapsto {}_qK = \cdots (\mu(x), \mu(y))^\pi \cdots.$$

The table given above is to be understood as a generic one where $a$ can be any number in $\mathbb{N}_{2n} = \{1, 2, \ldots, 2n\}$ and the numbers in the second line are elements of $\mathbb{N}_{2n}$, so the expressions in this line should be understood mod $2n$.

2.3.2. $\Omega_{iso}$ *moves.* Let $K \xrightarrow{\Omega_0} K'$ be an $\Omega_{iso}$ move by an isotopy $F : \rho \times [0, 1] \to \rho \times [0, 1]$ of the plane $\rho$.

Then recall (§1.3) $K' = F_1(K) = (F_1 \circ f)(K)$, $F_1 \circ f : S^1 \to K'$ and $K'$ is considered as a diagram with the induced orientation of $S^1$ by the map $F_1 \circ f$ (which is also the induced orientation of $K$ by the map $F_1$) with crossing points $F(\Delta)$ where $\Delta$ are the crossing points of $K$, and which has the same over/under information at each $F_1(\Delta)$ as those for $K$ at $\Delta$. Recall this means that if $\ell_1, \ell_2$ are the under- and over-crossing arcs in a canonical neighborhood $U$ of $\Delta$ for $K$, then $F_1(\ell_1), F_1(\ell_2)$ are respectively the under- and over-crossing arcs in the canonical neighborhood $F_1(U)$ of $F_1(K) = K'$ of $F_1(\Delta)$.

The points $x_1, x_2$ of $s$ corresponding by $f$ to $\Delta$ are those that correspond by $F_1 \circ f$ to $(F_1 \circ f)(\Delta)$. In two trips along the orientation of $S^1$ from the points $f^{-1}(p)$, $(F_1 \circ f)^{-1}(F_1(p))$, we meet the $x_i$'s in the same order (since $f^{-1}(p) = (F_1 \circ f)^{-1}(F_1(p))$). So in a trip on $K$ and another one on $K'$ along their orientations starting at $p$ and $F_1(p)$ respectively, we label a crossing $F_1(\Delta)$ of $K'$ by $i$ exactly when we label the crossing $\Delta$ of $K$ by $i$. So in the based symbols ${}_pK$, ${}_{F_1(p)}K'$ the pairs $(i : j)$ are common.

Moreover, because of the remark on the over/under information above, $i$ is given as a label moving in both cases along an over-crossing arc, or in both cases moving along an under-crossing arc. So in the based symbols ${}_pK, {}_{F_1(p)}K'$ the ordered pairs $(i, j)$ are common.

Finally, since $F$ starts from the identity, $F_1$ is an orientation preserving homeomorphism of $\rho$, so the over-crossing $\ell_1$ at $\Delta$ turns positively (negatively) at $\Delta$ to become the under-crossing $\ell_2$ iff the over-crossing $F_1(\ell_1)$ at $F_1(\Delta)$ turns positively (negatively) at $F_1(\Delta)$ to become the under-crossing $F_1(\ell_2)$. This means that the sign of $\ell_1, \ell_2$ at $\Delta$ coincides with the sign of $F_1(\ell_1), F_1(\ell_2)$ at $F_1(\Delta)$.

So in the based symbols ${}_pK, {}_{F_1(p)}K'$ the ordered pairs $(i, j)^\pi$ are common, including their signs $\pi$. Thus we proved rigorously the expected relation:

(2.3) $$\quad {}_pK = {}_{F_1(p)}F_1(K).$$

2.3.3. $\Omega_1$ *moves.* Let $K \xrightarrow{\Omega_1} K'$ be an $\Omega_1$ move, $D$ its 2-disk, $\ell$ its arc for $K$ and $\ell'$ its arc for $K'$ (Figure 11). The two arcs share their first and last endpoints (points $t_1, t_2$ in the Figure).

Let us assume without loss of generality that $K \xrightarrow{\Omega_1^+} K'$ and $K' \xrightarrow{\Omega_1^-} K$.

Then $\ell$ does not contain any crossing of $K$ whereas $\ell'$ contains a crossing $\Delta$ of $K'$. Outside $D$ the diagrams $K, K'$ coincide, so they share all crossings of $K$ (if any).

We are going to relate ${}_pK, {}_qK'$ for some convenient choices of base points $p, q$. Then §2.3.1 explains the way to relate them for all possible choices of base points.

First for $\tau(K) = 0$: for any choice of base points $p$ and $q$ of $K$ and $K'$ we have (Figure 11) ${}_pK = \varnothing$ and ${}_qK' = (1, 2)^{\pm 1}$ or ${}_qK' = (2, 1)^{\pm 1}$.

Now for $\tau(K) = n > 0$, that is, when $K$ has some crossing:

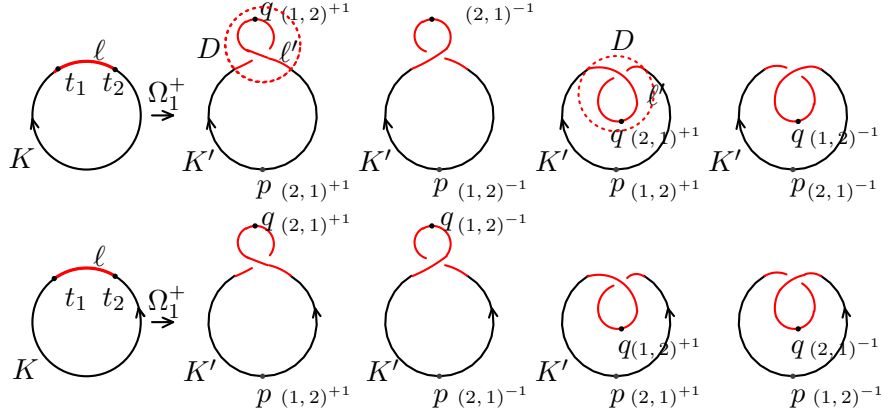

FIGURE 11. All possible topologically distinct ways to perform an $\Omega_1^+$ move on a diagram $K$ with $\tau(K) = 0$. In all cases we calculate the single crossing for all possible positions for the base point of the resulting diagram $K'$ ($q$ and $q'$ are arbitrary base points on the two arcs of the diagram). By §2.3.1 it is enough to know the crossings just for a single position of the base point(say for $q$).

Since $\ell$ contains no crossings, it is contained in the interior of a basic arc $e$. We consider four cases depending on the exact position of $p$ with respect to $e$ and $\ell$ (Figure 12). (A) $p$ is outside $e$, (B) $p$ is between the right endpoints of $e, \ell$, (C) $p$ is on $\ell$, (D) $p$ is between the left endpoints of $e, \ell$.

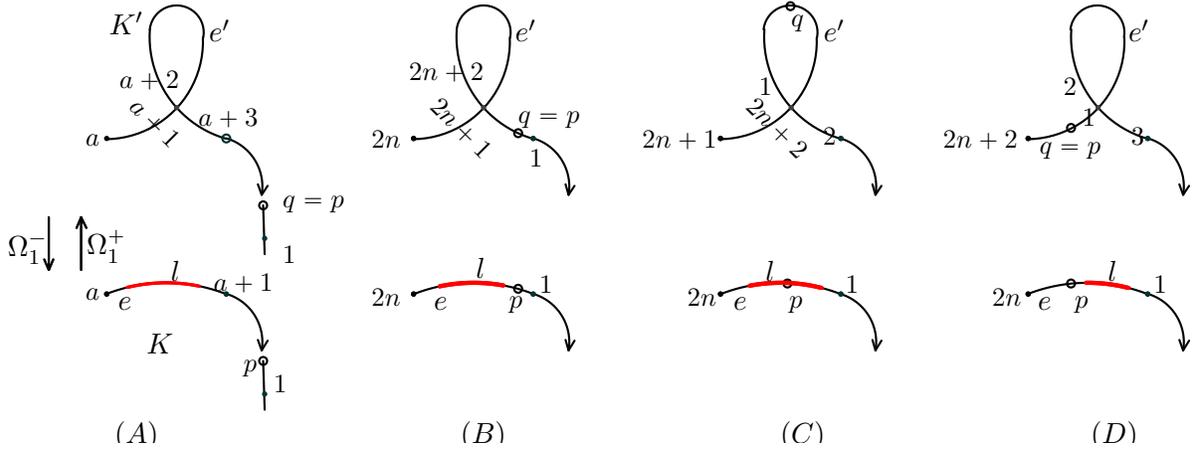

FIGURE 12. The relabelings during $\Omega_1^+, \Omega_1^-$ moves, close to the basic arcs $e = a(a+1)$ and $e'$ of the moves. For the $\Omega_1^+$ move, the position of the base point $p$ of $K$ is described with respect to $e$ and to its subarc $\ell$ on which the move actually takes place as follows: (A) $a \neq 2n$ and $p \notin e$, (B) $a = 2n$ and $p \in e$ but not on $\ell$, lying closer to the right endpoint of $e$, (C) $a = 2n$ and $p \in e$ $p$ on $\ell$, (D) $a = 2n$ and $p \in e$ $p \in e$ but not on $\ell$, lying closer to the left endpoint of $e$. The base point $q$ of $K'$ is chosen as depicted. For the $\Omega_1^-$ move, the position of the base point $q$ of $K'$ is described with respect to the basic arc $e'$ which forms the loop of the crossing that disappears in the move, as follows: (A) $q \notin e'$, (B)-(D) $q \in (2n+2)(1)$. The base point $p$ of $K$ is then chosen as depicted.

In cases (A), (B), (C) we set as base point $q$ of $K'$ again the point $p$, whereas in case $D$ we set as $q$ a point on the final position $\ell'$ of $\ell$ after the move, as depicted in Figure 12.

We work thoroughly case (A) and then we mention briefly the small changes in calculations for the rest of the cases.

In (A):

Let $e = \overrightarrow{a(a+1)}$ for some $a \in \mathbb{N}_{2n} - \{2n\}$.

The labels $1, 2, \ldots, a$ given to crossings of $K$ do not change for the same points considered now as crossings of $K'$. The labels $a+1, a+2$ will be given to the new crossing formed by the move. and the old labels will now become $a+3, \ldots, 2n+2$. The complete correspondence of the labels for the crossings of $_pK$ and $_pK'$ is given in the following table:

(2.4)
$$\begin{array}{llllllllll} K: & 1 & 2 & \ldots & a & & & & a+1 & \ldots & 2n \\ K': & 1 & 2 & \ldots & a & a+1 & a+2 & a+3 & \ldots & 2n+2 \end{array}$$

If $\mu$ is the relabeling function which corresponds the numbers in the first row to those of the second, the pairs in $_pK$ are replaced as:

$$(x_i, y_i)^{\pi_i} \stackrel{\Omega_1^+}{\mapsto} (\mu(x_i), \mu(y_i))^{\pi_i}.$$

And $_pK$ becomes $_pK'$ as follows:

(2.5)
$$_pK = \cdots (x_i, y_i)^{\pi_i} \cdots \stackrel{\Omega_1^+}{\mapsto} {_pK'} = \cdots (\mu(x_i), \mu(y_i))^{\pi_i} \cdots \Delta,$$
$$\Delta = (a+1 : a+2).$$

We can describe completely all relevant information for $\Delta$ in a quite compact form.

For this we give first a definition:

**Definition 15.** *For a basic arc $e$ of a diagram $K$ with $\tau(K) > 0$, we call placements or twists of $e$, the modification of an interior subarc of $e$ as in Figure 13 along with the depicted over/under information. We call the depicted number $\theta_e \in \{-1, 1\}$ as the placement number or twist-number of the corresponding twist. The changes happen inside a 2-disk. We call the placements with $\theta_e = 1$ as over-placements or over-twists and the rest as under-placements or under-twists.*

*Quite similarly we define over-, under-placements and numbers $\theta_e \in \{-1, 1\}$ for a simple arc $e$ of a diagram $K \approx S^1$.*

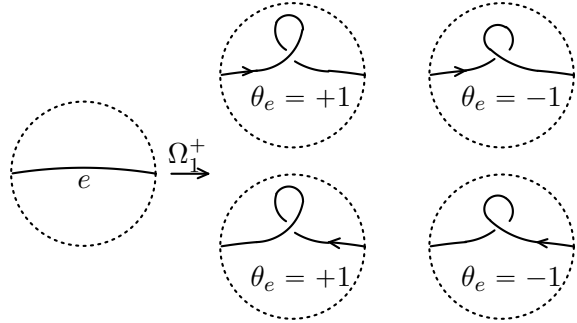

FIGURE 13. The under- and over-twists of a basic arc $e$ (or of a simple arc of a diagram homeomorphic to a circle) and their placement numbers $\theta_e$.

Now, the $\Omega_1^+$ move can be performed in the interior of any one of the two rooms $\Pi$ containing $e = \overrightarrow{a(a+1)}$. Let $\Pi_\phi$ for some $\phi \in \{\alpha, \delta\}$ be the orientation of the room $\Pi$ which induces the orientation of $e$ as a basic arc of $K$. In the notation we have introduced it is $\Pi_{\phi e}$. $\Pi$ can be a bounded or an unbounded room, but this does not troubles us. All possible cases for $\Pi, \phi$ and for the kind of twist are shown in Figure 14.

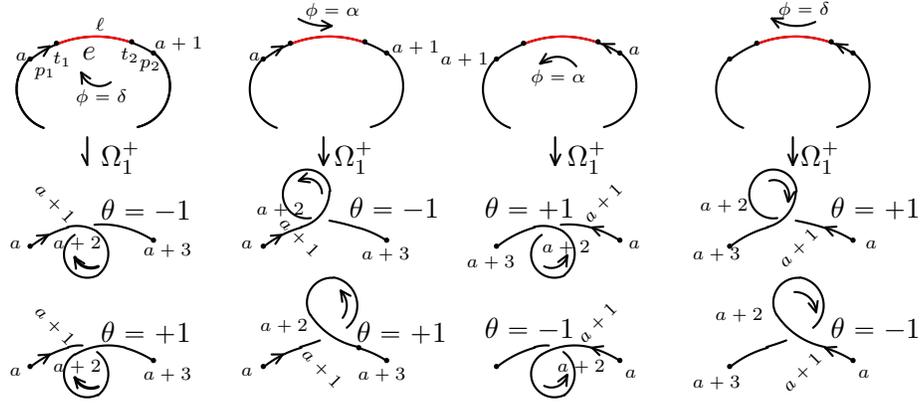

FIGURE 14. All possible topological ways to perform an $\Omega_1^+$ move on a basic arc $e$ on a diagram $K$ with $\tau(K) > 0$. $p_1, p_2$ are the endpoints of $e$, and $a, a+1$ are their labels for a chosen base point $p$ of $K$. The 2-disk of the move contains a subarc $\ell$ of $e$ with endpoints $t_1, t_2$. We perform the move inside a room $\Pi$ which contains $e$ on its boundary. We orient $\Pi$ by the orientation $\phi \in \{\alpha, \delta\}$ which induces on $e$ the orientation it has as a basic arc of $K$. After the move a small loop $e'$ is created. As a simple curve on the plane, $e'$ possesses an interior which by construction contains no other point of the resulting diagram $K'$. So this interior is a room $\Pi'$ of of $K'$. The orientation of $\Pi'$ which induces to its boundary $e'$ the orientation this arc has as a basic arc of $K'$, is always $\phi$ as is readily checked.

Let us also note that if we change the room in which we perform the move but still use the same kind of over- or under-placement, then in the description of $\Delta$ the placement number $\theta_e$ does not change, but the sign of the orientation of the room changes to $\pi_{-\phi} = -\pi_\phi$.

In all cases, it holds:

$$(2.6) \qquad \Delta = \left( \frac{1-\theta_e}{2}(a+1) + \frac{1+\theta_e}{2}(a+2), \frac{1+\theta_e}{2}(a+1) + \frac{1-\theta_e}{2}(a+2) \right)^{\pi_\phi \theta_e}.$$

This also holds in case $\tau(K) = 0$, i.e. whenever $K \approx S^1$, for an $\Omega_1^+$ move performed on a simple arc $e$ of $K$, if we put $a = 0$.

In cases (B),(C),(D) all the above hold with the only changes regarding the exact relabelings. Figure 12 helps us determine each relabeling. Briefly, the following hold:

In (B):

$$(2.7) \quad \begin{array}{l} K: \ 1 \quad 2 \quad \ldots \quad \quad 2n \\ K': \ 1 \quad 2 \quad \ldots \quad \quad 2n \quad 2n+1 \quad 2n+2 \\ {}_pK = \cdots (x,y)^\pi \cdots \mapsto {}_pK' = \cdots (\mu(x_i), \mu(y_i))^{\pi_i} \cdots \Delta, \\ \Delta = \left( \dfrac{1-\theta_e}{2}(2n+1) + \dfrac{1+\theta_e}{2}(2n+2), \dfrac{1+\theta_e}{2}(2n+1) + \dfrac{1-\theta_e}{2}(2n+2) \right)^{\pi_\phi \theta_e}. \end{array}$$

In (C):

(2.8)
$$K: \quad 1 \quad 2 \quad \ldots \quad 2n$$
$$K': \quad 2 \quad 3 \quad \ldots \quad 2n+1 \quad 2n+2 \quad 1$$
$$_pK = \cdots (x,y)^\pi \cdots \mapsto \; _pK' = \cdots (\mu(x_i), \mu(_i))^{\pi_i} \cdots \Delta,$$
$$\Delta = \left( \frac{1-\theta_e}{2}(2n+2) + \frac{1+\theta_e}{2} \cdot 1, \frac{1+\theta_e}{2}(2n+2) + \frac{1-\theta_e}{2} \cdot 1 \right)^{\pi_\phi \theta_e}.$$

In (D):

(2.9)
$$K: \quad 1 \quad 2 \quad \ldots \quad 2n$$
$$K': \quad 3 \quad 4 \quad \ldots \quad 2n+2 \quad 1 \quad 2$$
$$_pK = \cdots (x,y)^\pi \cdots \mapsto \; _pK' = \cdots (\mu(x_i), \mu(y_i))^{\pi_i} \cdots \Delta,$$
$$\Delta = \left( \frac{1-\theta_e}{2} \cdot 1 + \frac{1+\theta_e}{2} \cdot 2, \frac{1+\theta_e}{2} \cdot 1 + \frac{1-\theta_e}{2} \cdot 2 \right)^{\pi_\phi \theta_e}.$$

As noted above, we can dispense cases (B)-(D) altogether but it is of some interest to keep them and to go to some extra detail regarding the relabeling of each case. Figure 15 presents in a more compact way these relabelings and assigns to them respectively for cases (A)-(D) the names:

$$\mu^+_{2n,a,(-1)} \; (a \in \mathbb{N}_{2n} - \{2n\}), \; \mu^+_{2n,2n(-1)}, \; \mu^+_{2n,2n(0)}, \mu^+_{2n,2n(1)}.$$

Let us define for any $m > 0$ and any $x \in \mathbb{N}_{2m}$ the set $\mathbb{N}_{2m,x} = \mathbb{N}_{2m} - \{x, x+1\}$.

Then $\mu^+_{2n,a,(-1)} : \mathbb{N}_{2n} \to \mathbb{N}_{2(n+1),a+1}$ for $a \in \mathbb{N}_{2n} - \{2n\}$ and $\mu^+_{2n,2n,(x)} : \mathbb{N}_{2n} \to \mathbb{N}_{2(n+1),2(n+1)+x}$ for $x \in \{-1, 0, +1\}$. A slightly more unifying description of these relabelings is:

(2.10)
$$\mu^+_{2n,a,(x)} : \mathbb{N}_{2n} \to \mathbb{N}_{2(n+1), a+1+(x+1)}, \; a \in \mathbb{N}_{2n}$$
$$x = -1 \text{ if } a < 2n, x \in \{-1, 0, 1\} \text{ if } a = 2n.$$

The entry $x$ inside the parenthesis is a direct reference to the position of the base point $p \in K$.

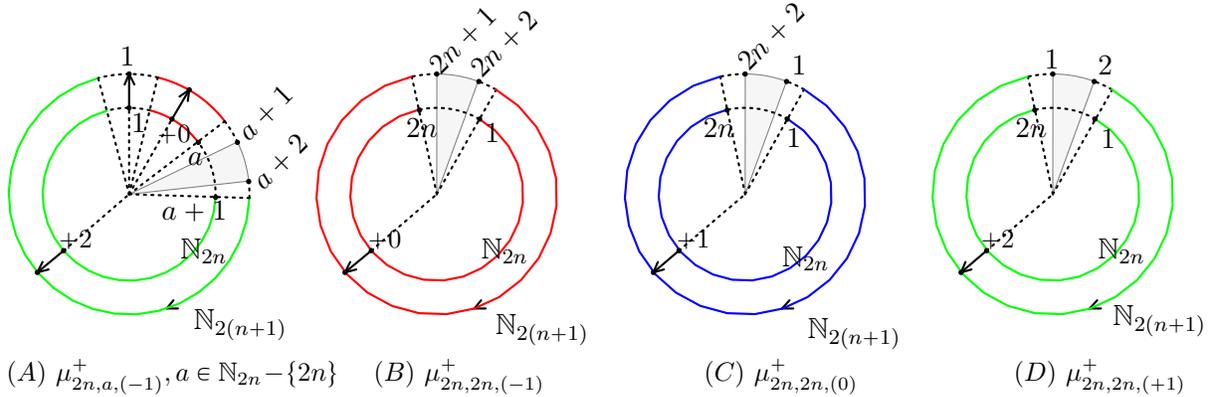

(A) $\mu^+_{2n,a,(-1)}, a \in \mathbb{N}_{2n} - \{2n\}$  (B) $\mu^+_{2n,2n,(-1)}$  (C) $\mu^+_{2n,2n,(0)}$  (D) $\mu^+_{2n,2n,(+1)}$

FIGURE 15. The relabelings of $\Omega^+_1$ moves on diagrams with order $n > 0$. Labels $z$ in a colored arc get the image $z+s$ where $s$ is the depicted signed integer $+0, +1$ or $+2$. The endpoints of the gray arcs (shaded sectors) in the big circle are not used as images.

Now for the inverse move $K' \xrightarrow{\Omega^-_1} K$ the converse of the above happen:

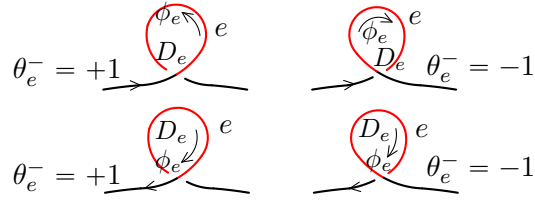

FIGURE 16. A 1-gon $e$ and its attached diagram placement number.

Inside the 2-disk of the move there exists a basic arc $e'$ of $K'$ whose endpoints coincide. The endpoints form the crossing $\Delta$ of $K'$ which disappears during the move. $e'$ forms a loop whose interior is a 2-disk (by Schonflies Theorem) and which contains no other points of $K'$. Actually he 2-disk of the move is just an enlargement of this interior. We say that $e'$ is the 1-loop of $\Delta$ and of the $\Omega_1^-$ move. The notion of 1-loops belongs to the diagrams, so we give the following definition:

**Definition 16.** *We call as a 1-loop or 1-gon of a diagram $K$, any basic arc $e$ whose endpoints coincide and which as a closed curve on the plane contains no other points of $K$ in its interior.*

*We write $D_e$ for the 2-disk enclosed by $e$. This is a room of $K$ which we call as the 2-disk of the 1-gon. We write $\phi_e$ for the orientation of $D_e$ which induces on its boundary $e$ the orientation this arc has as basic arc of $K$.*

*For a 1-gon $e$ of a diagram $K$, we call diagram placement number or diagram twist of $e$ the integer $\theta_e^- \in \{-1, 1\}$ depicted in Figure 16.*

Definitions 15, 16 and a check in Figure 14 imply immediately that:

**Lemma 7.** *Let $e$ be a basic arc of a diagram which produces after an $\Omega_1^+$ move an 1-gon $e'$. Let $\theta_e$ be the chosen placement number of the move, and $\phi$ the orientation of the room in which the move takes place and which induces to $4e4$ its orientation as a basic arc of $K$. Then (a) $\theta_{e'}^- = \theta_e$, (b) $\phi_{e'} = \phi$.*

Now, let us choose the base points $q \in K'$ and $p \in K$ in the way depicted in Figure 12; all possible positions of $q$ in $K'$ are covered in this figure, but this is not really important for our purposes. $e'$ is written as $e' = \overrightarrow{(a+1),(a+2)}$ for some $a \in \mathbb{N}_{2(n+1)}$. The crossings of $K'$ other than $\Delta$, remain in $K$ and these are all the crossings that $K$ has. Each one of them retains its old sign but now it gets new labels as dictated by the inverses $\mu^{-1}$ of the relabeling functions $\mu$ we met above.

Hence the pairs in $_pK'$ other than $\Delta$ are replaced as follows:

$$(x_i, y_i)^{\pi_i} \stackrel{\Omega_1^-}{\mapsto} (\mu^{-1}(x_i), \mu^{-1}(y_i))^{\pi_i}.$$

And $_pK'$ becomes $_pK$ as follows:

(2.11)
$$_pK' = \cdots (x_i, y_i)^{\pi_i} \cdots \Delta \stackrel{\Omega_1^-}{\mapsto} {}_pK = \cdots (\mu^{-1}(x_i), \mu^{-1}(y_i))^{\pi_i} \cdots$$
$$\Delta = \left(\tfrac{1-\theta}{2}(a+1) + \tfrac{1+\theta}{2}(a+2), \tfrac{1+\theta}{2}(a+1) + \tfrac{1-\theta}{2}(a+2)\right)^{\pi\theta},$$

where $\theta = \theta_{e'}$, $\pi = \pi_\phi$ and $\phi = \phi_{e'}$ which, recall, is the orientation of the interior of the loop $e'$ on the plane, which induces to its boundary $e'$ the orientation it has as a basic arc of $K'$.

In this description of $\Delta$ our choice $\theta = \theta_{e'}^-$ is expected since this way by Lemma 7 we just retain the placement number of the basic arc $e$ of $K$ which gives birth to $e'$ via an $\Omega_1^+$ move. Then the use of the orientation $\phi_{e'}$ is indeed the correct one as suggested by Lemma 7: the $\Omega_1^+$ move on the arc $e$ of $K$ takes place in a room $\Pi$ of $K$ (Figure 14) whose orientation, let $\phi_0$, induces on $e$ its orientation as a basic arc of $K$. $\phi_0$ is the orientation we used in Relations 2.6 - 2.9 for the description of $\Delta$ created by that move. But $e'$ lies in $\Pi$ too and $\phi = \phi_0$ as we can easily check in Figure 14 in all cases.

We can give to the relabelings the notation:

$$(\mu^+_{2n,a,(x)})^{-1} = \mu^-_{2(n+1),(a+1)+(x+1)} : \mathbb{N}_{2(n+1),(a+1)+(x+1)} \to \mathbb{N}_{2n},$$

(2.12)     so all of them get the common notation :

$$\mu^-_{2(n+1),z} : \mathbb{N}_{2(n+1),z} \to \mathbb{N}_{2n}, \ z \in \mathbb{N}_{2(n+1)}$$

since for the allowable values of $a$ and $x$, the indices $(a+1)+(x+1)$ get all values of $\mathbb{N}_{2(n+1)}$ exactly once (recall: $x = -1$ if $a < 2n$ and $x \in \{-1, 0, 1\}$ if $a = 2n$).

All these inverse relabelings were indicated in Figure 12 and they are presented again in Figure 17 in a more algebraic form suppressing any reference to diagrams $K$.

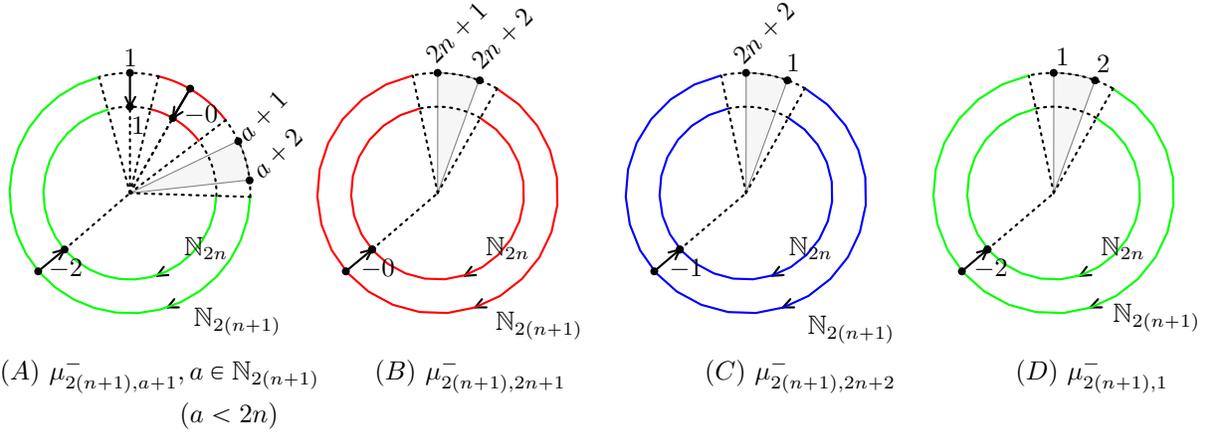

(A) $\mu^-_{2(n+1),a+1}, a \in \mathbb{N}_{2(n+1)}$     (B) $\mu^-_{2(n+1),2n+1}$     (C) $\mu^-_{2(n+1),2n+2}$     (D) $\mu^-_{2(n+1),1}$
($a < 2n$)

FIGURE 17. The relabelings of $\Omega^-_1$ moves on diagrams of order necessarily $\geq 1$. Labels $z$ in a colored arc get the image $z + s$ where $s$ is the depicted signed integer $-0, -1$ or $-2$. The endpoints of the gray arcs (shaded sectors) get no images.

We conveniently merge the relabelings of both $\Omega^+_1$ and $\Omega^-_i$ in a single Figure 18.

2.3.4. $\Omega_2$ *moves.* Let $K \xrightarrow{\Omega_2} K'$ be an $\Omega_2$ move, $D$ its 2-disk, $\ell_1, \ell_2$ the arcs of $K$ in $D$ and $\ell'_1, \ell'_2$ the arcs of $K'$ in $D$ with $\ell_1, \ell'_1$ sharing their endpoints and $\ell_2, \ell'_2$ sharing theirs.

Let us assume without loss of generality that $K \xrightarrow{\Omega^+_2} K'$ and $K' \xrightarrow{\Omega^-_2} K$.

Then $\ell = \ell_1 \cup \ell_2$ does not contain any crossing of $K$, thus its two subarcs $\ell_1, \ell_2$ are disjoint, whereas $\ell' = \ell'_1 \cup \ell'_2$ contains two crossings $\Delta_1, \Delta_2$ of $K'$, which are the intersection points of the two subarcs $\ell'_1, \ell'_2$. Outside $D$ the diagrams $K, K'$ coincide, so they share all crossings of $K$ (if any).

We are going to relate $_pK, _qK'$ for some convenient choices of base points $p, q$. Then §2.3.1 explains the way to relate them for all possible choices of base points.

First for $\tau(K) = 0$: for any choice of base point $p$ of $K$ it is $_pK = \emptyset$. For the choice $q$ in Figure 19 as base point for $K'$, we have $_qK' = (1,4)^{\pm 1}(2,3)^{\mp 1}$ or $_qK' = (4,1)^{\pm 1}(3,2)^{\mp 1}$.

Now for $\tau(K) = n > 0$, that is, when $K$ has some crossing:

Since $\ell$ contains no crossings, the subarcs $\ell_1, \ell_2$ of $\ell$ are contained in the interior of some basic arcs $e_1$ and $e_2$ respectively. It can be $e_1 = e_2$. Due to the cyclic ordering of the labels as elements of $\mathbb{N}_{2n}$, some cases in Figure 20 can merge to a single one.

Ignoring the "a unique choice of $p, q$" convenience mentioned above, in cases cases (I)-(J) of Figure 20 we place the base point $p$ of $K$ in various places. For the position of the base point $q$ of $K'$ we choose the point $p$ again, except when $p$ lies on $\ell_1$ in which case we set as $q$ a point on the final position $\ell'$ of $\ell$ after the move (cases (E), (F), (J)). In cases (A), (B) the base point $p \in K$

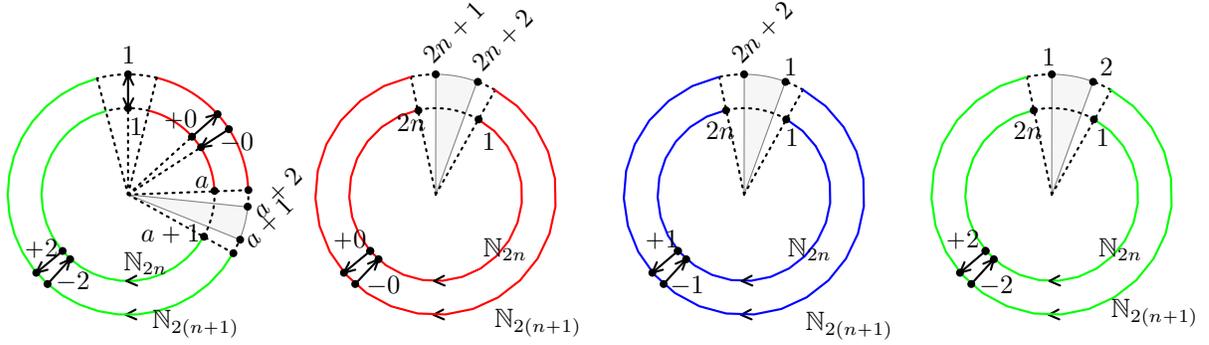

$\Omega_1^+$ : (A) $\mu_{2n,a,(-1)}^+, a \in \mathbb{N}_{2n} - \{2n\}$ (B) $\mu_{2n,2n,(-1)}^+$ (C) $\mu_{2n,2n,(0)}^+$ (D) $\mu_{2n,2n,(+1)}^+$
$\Omega_1^-$ : (A) $\mu_{2(n+1),a+1}^-, a \in \mathbb{N}_{2(n+1)}$ (B) $\mu_{2(n+1),2n+1}^-$ (C) $\mu_{2(n+1),2n+2}^-$ (D) $\mu_{2(n+1),1}^-$
   $(a < 2n)$

FIGURE 18. The relabelings for both versions of $\Omega_1$ moves are presented here in common. It is $n > 0$ for relabelings of $\Omega_i^+$ moves and $n \geqslant 0$ for relabelings of of $\Omega_i^-$ moves. For $n = 0$ we get in each case the trivial inverse relabeling by forgetting all information in the smaller circle. Labels $z$ in a colored arc get the image $z + s$ where $s$ is the depicted signed number $+0, +1, +2, -0, -1, -2$. The endpoints of the gray arcs (shaded sectors) in the big circles are not used as images of the relabelings and get no images in the relabelings of $\Omega_i^-$ moves.

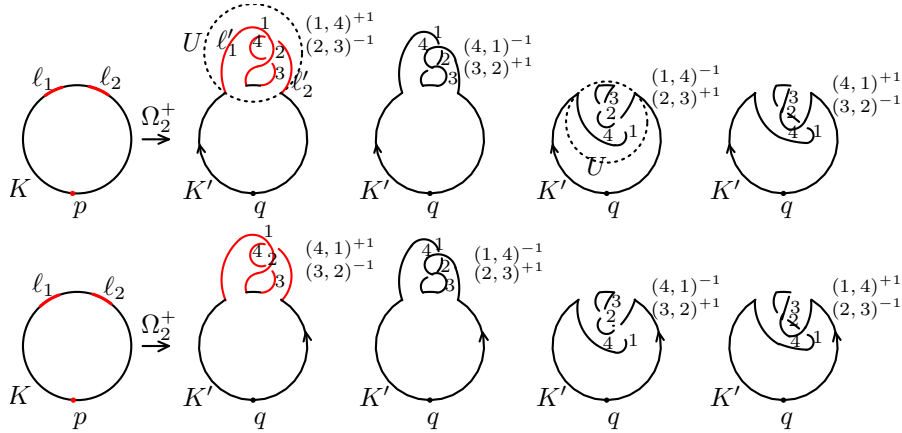

FIGURE 19. All possible topological ways to perform an $\Omega_2^+$ move on a diagram $K$ with $\tau(K) = 0$. We are allowed to extend the red arcs $\ell_1, \ell_2$ outside $K$ to the left and right respectively and clasp them on the bottom so as to engulf point $q$ and at the same time perform an $\Omega_2^+$ move since this can happen inside an elongated version of the 2-disk $U$, but the setting we produce differs to the pictures given here only by an isotopy of the plane. In all cases we calculate the two crossings for $q$ as the base point of $K'$. Changing the orientation of $K$ (and then of $K'$) the calculations of crossings would change by reversing the entries and replacing the signs to the opposite number. Thus the possible descriptions of the crossings are again the shown ones.

does not lie on the basic arcs $e_1, e_2$ and these two arcs are distinct. In cases (C)-(H) point $p$ lies in one of $e_1, e_2$ which continue to be distinct. In cases (I)-(L) $e_1, e_2$ coincide to a single arc $e$, and in case (I) the base point does not lie on it whereas in the rest of the cases it does.

We work thoroughly cases (A), (B) and then we mention briefly the small changes in calculations for the rest of the cases.

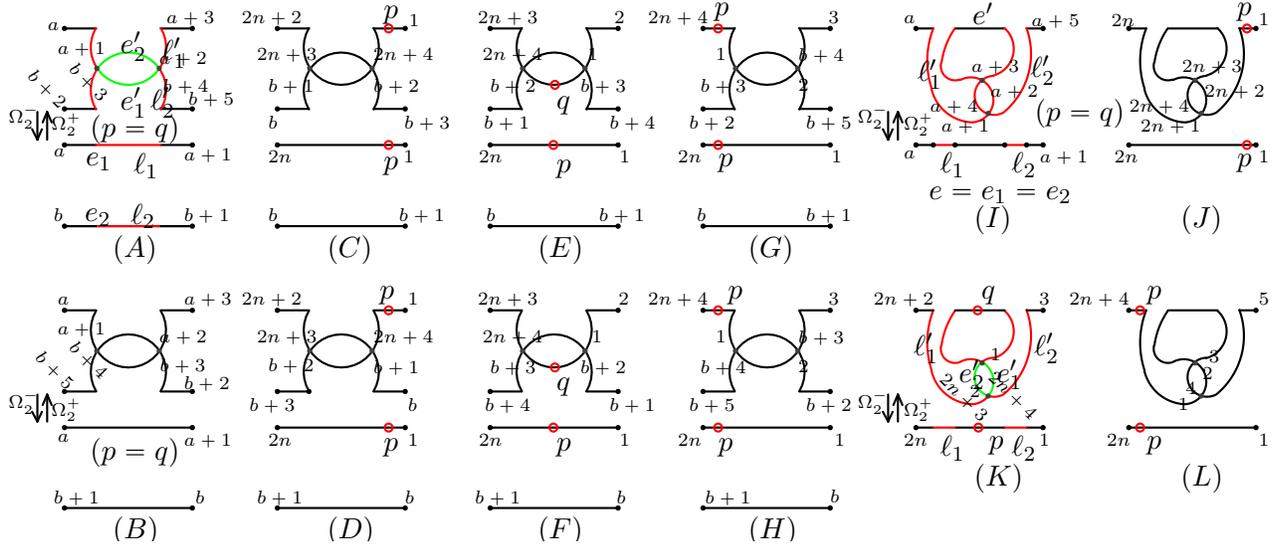

FIGURE 20. The relabelings during $\Omega_2^+, \Omega_2^-$ moves, close to the basic arcs $e_1, e_2$ or $e_1', e_2'$ respectively on which they are performed. For the $\Omega_2^+$ move, the position of the base point $p$ of $K$ with respect to $e_1, e_2$ and to their subarcs $\ell_1, \ell_2$ which the move changes is described as follows: In (A), (B) and (I) $p \notin e_1, e_2$, whereas $p \in e_1$ in the other cases. It is closer to the right endpoint in (C), (D), (K), it lies exactly on $\ell_1$ in cases (E), (F), (J), and it is closer to the left endpoint in (G), (H) and (L). In cases (A)-(H) the arcs $e_1, e_2$ are distinct, whereas in (I)-(L) they coincide. The base point $q$ of $K'$ is chosen as depicted. $q$ and $p$ coincide whenever $q$ is not mentioned. For the $\Omega_2^-$ move, the position of the base point $q$ of $K$ is described with respect to the basic arcs $e_1', e_2'$ which form a 2-gon (ch. Definition 18 below) for the two crossings that disappears in the move, as follows: In (A), (B), (I) $q \notin e_1', e_2'$, whereas in the other cases $q \in e_1' = (2n+4)(1)$. The base point $p$ of $K$ is chosen as depicted. $q$ and $p$ coincide whenever $q$ is not mentioned. We do not provide names and colors in all cases for $e_1, e_2, \ell_1, \ell_2, \ell_1', \ell_2'$ ($\Omega_2^+$ move) or for $e_1', e_2'$ ($\Omega_2^-$ move) since their positions are inferred from the shown ones.

In (A), (B):

It is $e_1 = \overrightarrow{a(a+1)}$, $e_2 = \overrightarrow{b(b+1)}$ for some $a, b \in \mathbb{N}_{2n} - 2n$, $a \neq b$. Without any loss of generality let $a < b$. This condition on $a, b$ creates an artificial asymmetry between $e_1, e_2$ which we shall keep in the sequel. Thus from now on $e_1$ is considered as the first and $e_2$ the second arc among the two.

The correspondence of the labels for the crossings of $K$ and $K'$ is given in the following table:

(2.13) 
| $K:$ | 1 | 2 | ... | $a$ | | | $a+1$ | $a+2$ | ... | $b$ | | | | $b+1$ | ... | $2n$ |
|---|---|---|---|---|---|---|---|---|---|---|---|---|---|---|---|---|
| $K':$ | 1 | 2 | ... | $a$ | $a+1$ | $a+2$ | $a+3$ | $a+4$ | ... | $b+2$ | $b+3$ | $b+4$ | $b+5$ | ... | | $2n+4$ |

If $\mu$ is the relabeling function corresponding the numbers in the first row of the tables to those of the second, then the pairs in ${}_pK$ are replaced as:

$$(x,y)^\pi \overset{\Omega_2^+}{\mapsto} (\mu(x_i), \mu(y_i))^{\pi_i}.$$

And ${}_pK$ becomes ${}_pK'$ as follows:

$${}_pK = \cdots (x_i, y_i)^{\pi_i} \cdots \overset{\Omega_2^+}{\mapsto} {}_pK' = \cdots (\mu(x_i), \mu(y_i))^{\pi_i} \cdots \Delta_1, \Delta_2,$$

where $\Delta_1 = (a+1:b+4), \Delta_2 = (a+2:b+3)$ or $\Delta_1 = (a+1:b+3), \Delta_2 = (a+2:b+4)$.

In Relation (2.15) below, we describe completely all relevant information for $\Delta_1, \Delta_2$. First we give a definition similar to the one of twist placements for the $\Omega_1^+$ moves.

**Definition 17.** *For two distinct basic arcs $e_1, e_2$ of a diagram $K$ with $\tau(K) > 0$, we call over- and under-placements of $e_1$ with respect to $e_2$ the modifications of some interior subarcs of them as in Figure 21 (a) and (b) along with the depicted over/under information. The changes happen inside a 2-disk. We also call them as under- and over- placements of $e_2$ with respect to $e_1$. We denote $\theta_{e_1,e_2} = +1$ and $\theta_{e_1,e_2} = -1$ respectively. We also denote $\theta_{e_2,e_1} = -1$ and $\theta_{e_2,e_1} = +1$ respectively. We call $\theta_{e_1,e_2}, \theta_{e_1,e_2}$ as the placement number of the pair $(e_1, e_2)$ or $(e_2, e_1)$ respectively.*

*For a single basic arc $e$ of some diagram we call the self-intersected arcs in Figure 21 (c) and (d) along with their over/under information respectively as an over- or under-placement of $e$ with respect to itself. We define as placement number of $e$ the depicted number $\theta_{e,e} = +1$ or $\theta_{e,e} = -1$ respectively.*

*For $K \approx S^1$ we define an over- and an under-placement of a simple arc $e$ of $K$ with respect to itself as in Figure 21 (e) and (f). The changes happen inside a 2-disk. We define then the placement number of $e$ with respect to itself as $\theta_{e,e} = +1, \theta_{e,e} = -1$ respectively.*

According to the definition it is always $\theta_{e_2,e_1} = -\theta_{e_1,e_2}$.

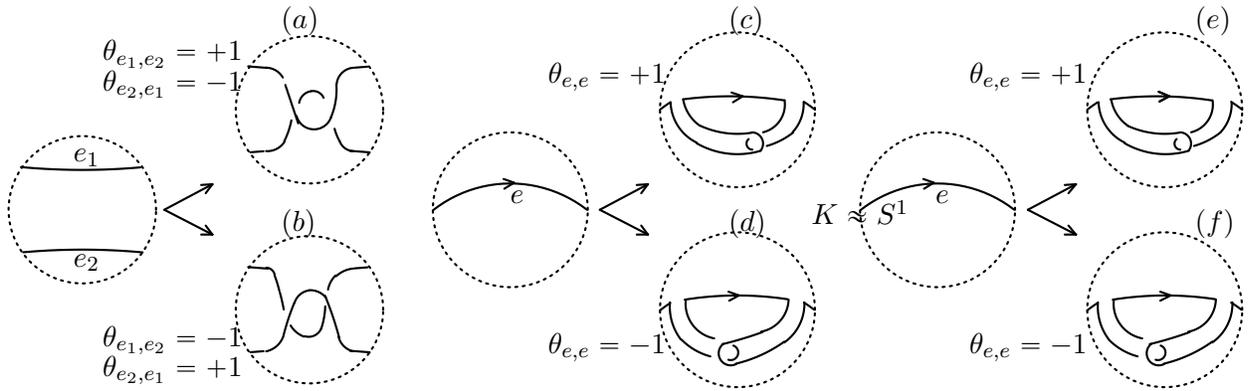

FIGURE 21. The under- and over-placements of a basic arc with respect to another basic arc (in (a),(b)) or to itself (in (c),(d)). Also shown are the under-, over-placements of a simple arc of a circle with respect to itself (in (e),(f)). In each case the corresponding placement number is also denoted.

**Lemma 8.** *If an $\Omega_2^+$ move can be performed on the basic arcs $e_1 = \overrightarrow{a(a+1)}, e_2 = \overrightarrow{b(b+1)}$ of $K$, then $e_1, e_2$ lie on some (unoriented) cycle of $K$ and there exist at most two cycles of $K$ on which $e_1, e_2$ lie simultaneously (the latter happens whenever $e_1, e_2$ share the two rooms in which each one of them lies). Let $D$ be the 2-disk of the $\Omega_2^+$ move and let $\partial\Pi$ be a common (unoriented) cycle of $e_1, e_2$. If the move happens in $\Pi$ then the topological possibilities for the placement of $e_1, e_2$ as chords of $D$ and the information of their orientations are the ones in Figure 23, where the notation $e_1 \uparrow\uparrow e_2, e_1 \downarrow\uparrow e_2$ introduced in Definition 10 indicates whether $e_1, e_2$ have the same sense on $\partial\Pi$ or not, and where $\phi$ is the orientation of $\Pi$ that induces on $e_1$ its orientation as a basic arc of $K$ (in the notation of Definition 8: $\partial\Pi_{\phi e_1}$). The topological possibilities for the orientations of the involved arcs before and after the move and for the crossings created are the ones in Figure 24.*

*Proof.* $e_1, e_2$ are disjoint topological chords of $D$ dividing it into three parts $D_1, D_2, D_3$. Say $D_2$ is the middle part with both $e_1, e_2$ on its boundary. Since no other part of $K$ lies in $D_2$, it must be that $D_2$ is part of a room of $K$, revealing that $e_1, e_2$ are basic arcs on the boundary of the same room, so they belong to the same cycle of $K$. Since each one of $e_1, e_2$ lies on the boundary of two rooms, they can belong simultaneously in at most two cycles (this happens whenever they share the two rooms in which they lie), as wanted.

To decide about the topological possibilities, we follow the technique in the proof of Lemma 1 and delete small disjoint canonical neighborhoods of the crossings of $\Pi$ to transform it to $\Pi'$ which is a 2-disk or the closed complement of a 2-disk. We choose neighborhoods that avoid $D$. The modifications transform $e_1, e_2$ to slightly smaller arcs $e'_1 = \overrightarrow{a'(a+1)'}, e'_2 = \overrightarrow{b'(b+1)'}$ of $K$ (not basic arcs). $D_2$ lies on $\Pi'$ and the common part of $\partial D_2$ with $\partial \Pi'$ is contained in the arcs $e'_1, e'_2$.

Now all possible topological ways for a 2-disk like $D_2$ to sit in another 2-disk (or the complement of a 2-disk) like $\Pi'$ in the above way, are as depicted in Figure 22. So $e_1, e_2$ sit in $D$ as described in Figure 23. And this gives that whenever the move happens in $\Pi$ the orientation information for the arcs before and after the move is as described in Figure 24, that is, as wanted. $\square$

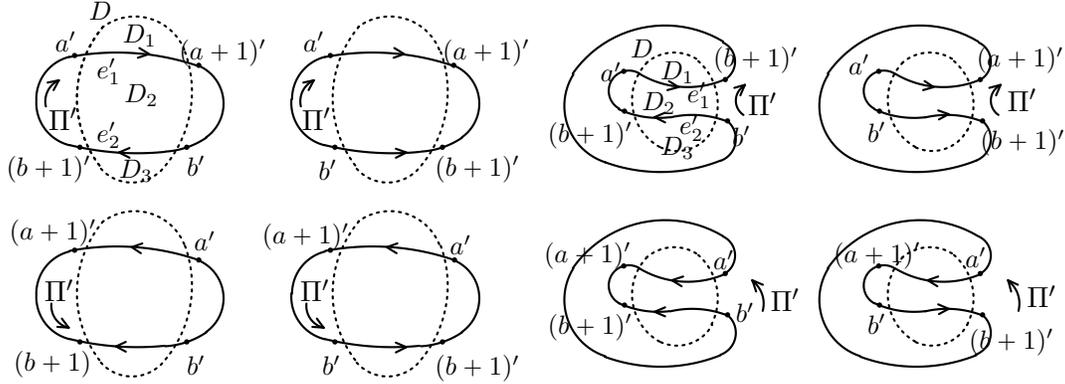

FIGURE 22. $D$ is the 2-disk of an $\Omega_2^+$ move on the basic arcs $e_1 = a(a+1), e_2 = b(b+1)$ $e_1, e_2$ of some diagram $K$. The arcs $e_1, e_2$ split $D$ in three parts $D_1, D_2, D_3$ and $D_2$ is the middle one among them. The clasp of arcs that will be created by the move will lie in $D_2$. $D_2$ lies in a room $\Pi$ of $K$ and $e_1, e_2$ lie on the boundary of $\Pi$. We modify $\Pi$ as in Figure 6 away from $D$, and we get a 2-disk $\Pi'$ or the closed complement of a 2-disk on the plane. The modification transforms $e_1, e_2$ of $\partial \Pi$ to the arcs $e'_1 = a'(a+1)', e'_2 = b'(b+1)'$ of $\Pi'$. $e_1, e'_1$ share all points except from small subarcs at their endpoints which were exchanged during the modification; they certainly share the arcs on the boundary of $D_2$ separating it from $D_1$. Similarly for $e_2, e'_2$. The figure shows the possible topological placements of the 2-disk $D_2$ with respect to the 2-disk $\Pi'$, or with respect to the complement $\Pi'$ on the plane of a 2-disk. In all cases the orientation shown for $\Pi'$ is the one that induces on the arc $e'_1 = a'(a+1')$ the orientation it has as a basic arc of $K$.

All the $\Omega_2^+$ moves in cases (A) - (B) we work on, which can be performed in a room $\Pi$ on the two basic arcs $e_1, e_2$ are depicted with all their details at the crossings and the orientation information of the involved arcs in Figure 24. Setting $\theta = \theta_{e_1, e_2}, \rho = \rho_{e_1, e_2}, \pi = \pi_\phi$ it holds:

(2.14)
$$\Delta_1 = \left(\tfrac{1+\theta}{2}(a+1) + \tfrac{1-\theta}{2}\left(b + \tfrac{7+\rho}{2}\right), \tfrac{1-\theta}{2}(a+1) + \tfrac{1+\theta}{2}\left(b + \tfrac{7+\rho}{2}\right)\right)^{\pi\theta\rho}$$
$$\Delta_2 = \left(\tfrac{1+\theta}{2}(a+2) + \tfrac{1-\theta}{2}\left(b + \tfrac{7-\rho}{2}\right), \tfrac{1-\theta}{2}(a+2) + \tfrac{1+\theta}{2}\left(b + \tfrac{7-\rho}{2}\right)\right)^{-\pi\theta\rho}$$

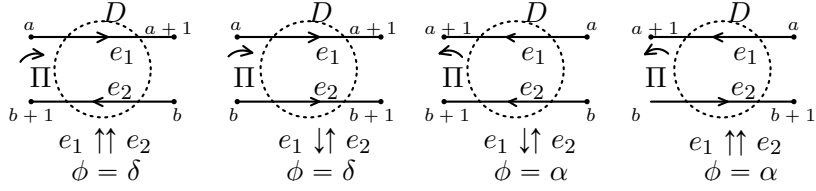

FIGURE 23. All topologically distinct relative positions of $e_1, e_2, D, \Pi$. They coincide with the relative position of $e'_1, e'_2, D, \Pi'$ in Figure 22. In all cases the orientation $\phi$ shown for $\Pi$ is the one that induces on the arc $e_1$ the orientation it has as a basic arc of $K$, and this coincides with the orientation of $\Pi'$ that induces on the arc $e'_1$ the orientation it has as an arc of $K$.

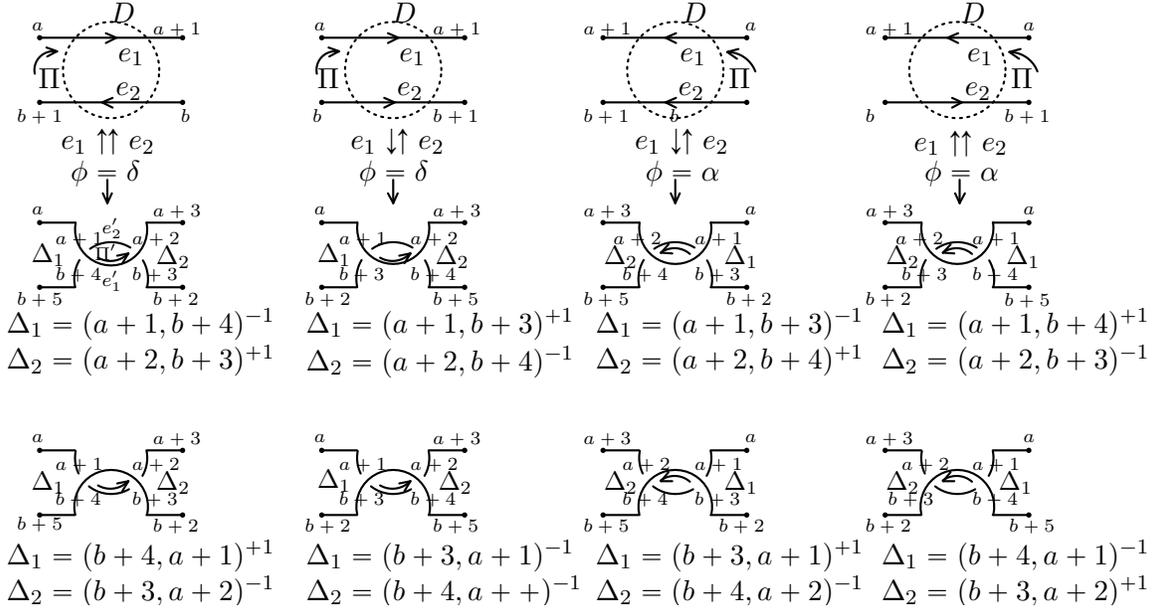

FIGURE 24. All possible topological ways to perform an $\Omega_2^+$ move on two basic arcs $e_1 = a(a+1), e_2 = b(b+1)$ ($a, b \in \mathbb{N}_{2n}$, $a < b$) on a diagram $K$ with $\tau(K) > 0$. Regarding $\phi$, the arcs $e_1, e_2$ appear asymmetric with respect to their order, with $e_1$ considered as first among them: the move takes place inside a room $\Pi$ and $\phi$ is the orientation of $\Pi$ which induces to $e_1$ the orientation it has as a basic arc of $K$ (it holds $\partial \Pi_{\phi e_1}$). After the move two new crossings $\Delta_1, \Delta_2$ are created and the two new basic arcs $e'_1, e'_2$ among them form a small loop $e'_1 \cup e'_2$. As a simple curve on the plane $e'_1 \cup e'_2$ possesses an interior which by construction contains no other point of the resulting diagram $K'$. So this interior is a room $\Pi'$ of $K'$. The orientation of $\Pi'$ which induces to its boundary arc $e_1$ the orientation it has as a basic arc of $K'$ is always $\phi$ as is immediately checked. In each case the crossings $\Delta_1, \Delta_2$ are calculated.

Let us also note that if we change the room in which we perform the move but still use the same kind of over- or under-placement, then in the description of $\Delta_1, \Delta_2$ the placement number $\theta_{e_1, e_2}$ does not change, but the sign of the rooms orientation changes to $\pi_{-\phi} = -\pi_\phi$.

This completes the analysis of cases (A), (B).

In the rest of the cases, the above hold essentially unaltered:

The expressions in relation (2.14) for the new crossings $\Delta_1, \Delta_2$ continue to hold with the values of $a$ and $b$ at hand (which now can take explicitly the value 1 or $2n$), a claim easily checked in Figure 24. The only important changes have to do with the exact relabelings when we pass from $K$ to $K'$. All these relabelings are shown in the more convenient cyclic form of Figure 25 rather than in the line form of table (2.13).

In cases (C)-(H) the last two remarks suffice.

In case (I) we should mention the following particulars:

The two arcs $e_1, e_2$ are now the same arc, call it $e$, and the base point $p \in K$ does not lie on it. Hence $e = \overrightarrow{a(a+1)}$ for some $a \in \mathbb{N}_{2n}$, $a \neq 2n$. Relation (2.14) continues to hold but now it is always $a = b$ and $\rho = 1$ (as one would expect), and also $\theta = \theta_{e_1,e_2}$ can get any value of $\{-1, 1\}$ depending on the choice of placement of $e$. It is:

$$
\begin{aligned}
\Delta_1 &= \left(\tfrac{1+\theta}{2}(a+1) + \tfrac{1-\theta}{2}\left(a + \tfrac{7+1}{2}\right), \tfrac{1-\theta}{2}(a+1) + \tfrac{1+\theta}{2}\left(a + \tfrac{7+1}{2}\right)\right)^{\pi\theta} \\
\Delta_2 &= \left(\tfrac{1+\theta}{2}(a+2) + \tfrac{1-\theta}{2}\left(a + \tfrac{7-1}{2}\right), \tfrac{1-\theta}{2}(a+2) + \tfrac{1+\theta}{2}\left(a + \tfrac{7-1}{2}\right)\right)^{-\pi\theta}
\end{aligned}
\quad (2.15)
$$

The complete relabeling function is the following:

(2.16)
$$
\begin{array}{llllllllllll}
K: & 1 & 2 & \ldots & a & & & & & b+1 & \ldots & 2n \\
K': & 1 & 2 & \ldots & a & a+1 & a+2 & a+3 & a+4 & a+5 & \ldots & 2n+4
\end{array}
$$

As is trivially checked, (2.15) also holds for the diagrams $K$ with $\tau(K) = 0$, i.e. whenever $K \approx S^1$, for $\Omega_2^+$ moves happening on a simple arc $e$ of $K$, if we put $a = b = 0$ and $\theta = \theta_{e,e}$. This is the same as to say that relation (2.14) continues to hold where now $a = b = 0$, $\rho = 1$ and $\theta = \theta_{e,e}$.

In cases (J)-(L) relation (2.15) holds as in (I) with the expected values of the labels at the crossings where now $a = 2n$.

Although of importance to the sequel, it is interesting to give a few comments regarding a more compact notation of the relabeling functions:

For any $m > 0$ and any $x, y \in \mathbb{N}_{2m}$ with $|x - y| \neq 1$ we define the set $\mathbb{N}_{2m,x,y} = \mathbb{N}_{2m} - \{x, x+1, y, y+1\}\}$ for $|x - y| \geq 2$ and $\mathbb{N}_{2m,x,x} = \mathbb{N}_{2m} - \{x, x+1, x+2, x+3\}\}$ for $|x - y| = 0$ i.e. for $x = y$.

Then in the cases (A)-(B) the relabelings are: $\mu^+_{2n,a,b} : \mathbb{N}_{2n} \to \mathbb{N}_{2(n+1),a+1,b+3}$ for $a, b \in \mathbb{N}_{2n}$ so that $a < b < 2n$. And in cases (C)-(H) they are $\mu^+_{2n,a,2n,(x)} : \mathbb{N}_{2n} \to \mathbb{N}_{2(n+2),a+1+(1+x),2n+3+(1+x)}$ for $a \in \mathbb{N}_{2n}$, $a < b = 2n$ and $x \in \{-1, 0, +1\}$.

A more unifying notation for all relabelings is: $\mu^+_{2n,a,b,(x)} : \mathbb{N}_{2n} \to \mathbb{N}_{2(n+2),a+1+(1+x),b+3+(1+x)}$ for $a, b \in \mathbb{N}_{2n}$ such that $a < b$, where $x = -1$ if $b < 2n$ and $x \in \{-1, 0, 1\}$ if $b = 2n$. The entry inside the parenthesis is a direct reference to the position of the base point $p \in K$.

Also, in cases (A)-(I) the relabelings are: $\mu^+_{2n,a,a} : \mathbb{N}_{2n} \to \mathbb{N}_{2(n+1),a+1,a+3}$ for $a(=b) \in \mathbb{N}_{2n}$, $a < 2n$. And in cases (J)-(L) $\mu^+_{2n,2n,2n,(x)} : \mathbb{N}_{2n} \to \mathbb{N}_{2(n+2),2n+1+2(1+x),2n+3+2(1+x)}$ for $a(=b) = 2n$ and $x \in \{-1, 0, +1\}$.

A more unifying notation for all cases (A)-(L) is: $\mu^+_{2n,a,a,(x)} : \mathbb{N}_{2n} \to \mathbb{N}_{2(n+2),a+1+2(1+x),b+3+2(1+x)}$ for $a(=b) < 2n$ and $x = -1$, or for $a(=b) = 2n$ and $x \in \{-1, 0, 1\}$. The entry inside the parenthesis is again a direct reference to the position of the base point $p \in K$.

In compact form, the relabelings in all cases (A)-(L) are:

(2.17)
$$\mu^+_{2n,a,b,(x),(y)} : \mathbb{N}_{2n} \to \mathbb{N}_{2(n+2),a+1+(1+y)(1+x),b+3+(1+y)(1+x)}, \quad a, b \in \mathbb{N}_{2n} \text{ such that } a \leq b,$$

$x = -1$ if $b < 2n$ but $x \in \{-1, 0, +1\}$ if $b = 2n$, $y = +1$ if $a = b = 2n$ but $y = 0$ otherwise.

We have at last completed the analysis of the $K \xrightarrow{\Omega_2^+} K'$ move.

For the inverse move $K' \xrightarrow{\Omega_2^-} K$ the converse of the above happen: two distinct basic arcs $e'_1, e'_2$ of $K'$ share their two distinct endpoints, which are a couple of distinct crossing $\Delta_1, \Delta_2$ of $K$. Thus

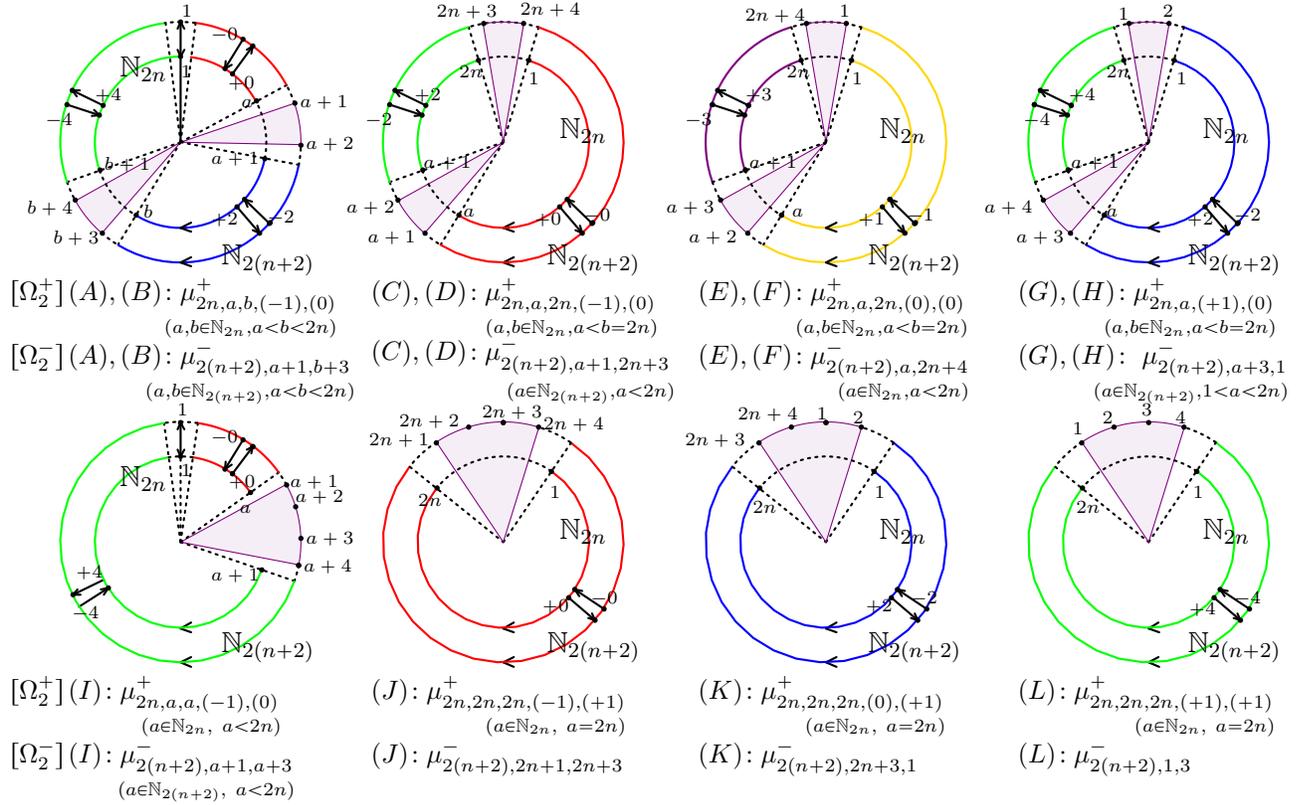

FIGURE 25. The relabelings of the $\Omega_2^+$ moves on diagrams (of order $n \geqslant 0$), shown in cyclic representation together with the relabelings of the $\Omega_2^-$ moves on diagrams (of order necessarily $\geqslant 2$). The $\Omega_2^+$ relabelings are the arrows from the small circle to the big one, whereas the $\Omega_2^-$ relabelings are the arrows from the big circle to the small one. In all cases, labels $z$ in a colored arc get the image $z + s$ where $s$ is the depicted signed integer $+0, +1, +2$ for $\Omega_2^+$ relabelings, or $-0, -1, -2$ for $\Omega_2^-$ relabelings. The endpoints of the grey arcs (shaded sectors) are not used as images and they get no images as well. For $\Omega_2^+$ moves on diagrams with $n = 0$ and $\Omega_1^-$ moves on diagrams with $n = 2$, the inner circle has no labels on it.

the union $e'_1 \cup e'_2$ forms a loop, and moreover this loop bounds by Schonflies Theorem a 2-disk. One of the arcs is over the other in both endpoints in the sense that some small extensions of the arcs at the crossings make one of them the over-arc at both crossings. The move eliminates the crossings. Its 2-disk is am enlargement of the 2-disk bounded by $e_1 \cup e_2$. We say that $e'_1, e'_2$ form a *2-gon* for $\Delta_1, \Delta_2$ and for the $\Omega_2^-$ move. The notion of 2-gons belongs to the diagrams, so we give the following definition:

**Definition 18.** *We call as a 2-gon $e_1, e_2$ of a diagram $K$ with $e_1$ over $e_2$, any two distinct basic arcs $e_1, e_2$ of $K$, none a 1-gon, which share their two distinct endpoints so that $e_1$ is over $e_2$ on both endpoints, and for which the simple closed curve $e_1 \cup e_2$ contains no other points of $K$ in its interior on the plane (Figure 26).*

We write $D_{e_1,e_2}$ for the 2-disk enclosed by $e_1 \cup e_2$ which we call as the *2-disk of the 2-gon*. This is a room of $K$. We write $\phi_{e_1,e_2}, \phi_{e_2,e_1}$ for the orientation of $D_{e_1,e_2}$ which induces respectively in its boundary component $e_1$ or $e_2$, the opposite of the orientation it has as a basic arc of $K$.

For a 2-gon $e_1, e_2$ of $K$, we call *diagram placement number* of the ordered pairs $(e_1, e_2)$, $(e_2, e_1)$ the integers $\theta^-_{e_1,e_2}, \theta^-_{e_2,e_1} \in \{-1, 1\}$ depicted in Figure 26. It is always $\theta^-_{e_2,e_1} = -\theta^-_{e_1,e_2}$.

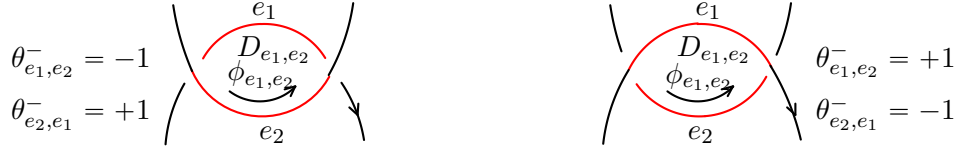

FIGURE 26. A 2-gon $e_1, e_2$, the room $D_{e_1,e_2} = D_{e_2,e_1}$ it creates, and the attached diagram placement numbers $\theta^-_{e_1,e_2}, \theta^-_{e_2,e_1}$.

Of course any one of $\phi_{e_1,e_2} = \phi_{e_2,e_1}$, $\phi_{e_1,e_2} = -\phi_{e_2,e_1}$ might happen.

Definitions 17, 18 and a check in Figure 24 imply immediately that:

**Lemma 9.** *Let $e_1, e_2$ be basic arcs of a diagram which produce after an $\Omega_2^+$ move an 2-gon $e'_1, e'_2$. Let $\theta_{e_1,e_2}$ be the chosen placement number of the move, and $\phi$ the orientation of the room in which the move takes place and which induces to $e_1$ its orientation as a basic arc of $K$. Then (a) $\theta^-_{e'_1,e'_2} = \theta_{e_1,e_2}$, (b) $\phi_{e'_1,e'_2} = \phi$.*

Now, we choose the base point $q \in K'$ as in Figure 20 and then we choose the base point $p \in K$ as in the same Figure. These positions of $q$ do not cover all possible places the base point can have, but this is not important for our purposes. $e'_1, e'_2$ get the form $e'_1 = \overrightarrow{z(z+1)}$, $e'_2 = \overrightarrow{w(w+1)}$ for some $z, w \in \mathbb{N}_{2(n+2)}$ with $z \neq w$. The crossings of $K$ are those of $K'$ other than $\Delta_1, \Delta_2$. Each one of them retains its old sign but now it gets new labels as dictated by the inverses $\mu^{-1}$ of the relabeling functions $\mu$ we met above. All these relabelings are shown in the cyclic form of Figure 25 alongside with the relabelings of the $\Omega_2^+$ moves.

Hence the pairs in $_pK'$ other than $\Delta_1, \Delta_2$ are replaced as follows:

$$(x_i, y_i)^{\pi_i} \stackrel{\Omega_2^-}{\mapsto} (\mu^{-1}(x_i), \mu^{-1}(y_i))^{\pi_i}.$$

And $_pK'$ becomes $_pK$ as follows:

(2.18)
$$pK' = \cdots (x_i, y_i)^{\pi} \cdots \Delta_1 \Delta_2 \stackrel{\Omega_2^-}{\mapsto} {}_pK = \cdots (\mu^{-1}(x_i), \mu^{-1}(y_i))^{\pi_i} \cdots.$$
$$\Delta_1 = \left(\tfrac{1-\theta}{2}(a+1) + \tfrac{1+\theta}{2}\left(b + \tfrac{7+\rho}{2}\right), \tfrac{1+\theta}{2}(a+1) + \tfrac{1-\theta}{2}\left(b + \tfrac{7+\rho}{2}\right)\right)^{\pi\theta\rho}$$
$$\Delta_2 = \left(\tfrac{1-\theta}{2}(a+2) + \tfrac{1+\theta}{2}\left(b + \tfrac{7-\rho}{2}\right), \tfrac{1+\theta}{2}(a+2) + \tfrac{1-\theta}{2}\left(b + \tfrac{7-\rho}{2}\right)\right)^{-\pi\theta\rho},$$

where $\rho = \rho_{e'_1,e'_2}$, $\theta = \theta^-_{e'_1,e'_2}$, $\pi = \pi_\phi$ and $\phi = -\phi_{e'_1,e'_2}$, which, recall, is the orientation of the interior of the simple closed curve $e'_1 \cup e'_2$ on the plane, which induces to its boundary arc $e'_1$ the opposite orientation it has as a basic arc of $K'$.

In this description of $\Delta_1, \Delta_2$, the choice $\theta = \theta^-_{e'_1,e'_2}$ is expected since this way we just retain the placement number of the ordered pair $(e_1, e_2)$ of the basic arcs $e_1, e_2$ of $K$ which gives birth to $e'_1, e'_2$ via an $\Omega_2^+$ move with $e_1$ as first arc among the two. Then the use of the orientation $\phi = -\phi_{e'_1,e'_2}$ as described is indeed the correct one as suggested by Lemma 9: the $\Omega_2^+$ move on the pair $e_1, e_2$ of $K$ takes place in a room $\Pi$ of $K$ (Figure 24) whose orientation, let $\phi_0$, induces on $e_1$ its orientation as a basic arc of $K$. $\phi_0$ is the orientation used in the description of $\Delta_1, \Delta_2$ created by that move. But $e'_1, e'_2$ lie in $\Pi$ too and $\phi = -\phi_0$ as we can easily check in Figure 24 in all cases.

Finally, and although not imperative, let us give to the relabelings of the move the notation:

$$(\mu^+_{2n,a,b,(x),(y)})^{-1} = \mu^-_{2(n+2),a+1+(y+1)(x+1),b+3+(y+1)(x+1)} : \mathbb{N}_{2(n+2),a+1+(y+1)(x+1),b+3+(y+1)(x+1)} \to \mathbb{N}_{2n}$$

$$a, b \in \mathbb{N}_{2n}, a \leqslant b,$$

$x = -1$ if $b < 2n$, but $x \in \{-1, 0, 1\}$ if $b = 2n$; also $y = 1$ if $a = b = 2n$, but $y = 0$ otherwise.

So all relabelings get the common notation:

$$\mu^-_{2(n+2),z,w} : \mathbb{N}_{2(n+2),z,w} \to \mathbb{N}_{2n}, \ z, w \in \mathbb{N}_{2(n+2)}, \ |w - z| \neq 1,$$

since for the allowable values of $a, b, x$ and $y$, the indices $z = a + 1 + (y+1)(x+1)$, $w = b + 3 + (y+1)(x+1)$ get all pairs of values of $\mathbb{N}_{2(n+1)}$ exactly once, except for those with absolute difference 1.

2.3.5. $\Omega_3$ *moves.* Let $K \xrightarrow{\Omega_3} K'$ be an $\Omega_3$ move and $D$ be its 2-disk. It is $n = \tau(K) = \tau(K') \geqslant 3$.

We are going to relate $_pK, _qK'$ for some convenient choices of base points $p, q$. Then as we also mentioned for $\Omega_1, \Omega_2$ moves, §2.3.1 explains the way to relate them for all possible choices of base points.

The arcs $f_1, f_2, f_3$ of $K$ inside $D$ are not basic arcs, but each contains a basic arc whose endpoints are the two common points the overlying $f_i$ has with the other two of the $f_i$'s. Let $e_1, e_2, e_3$ be the basic arcs lying on $f_1, f_2, f_3$ respectively. Let $f'_i, e'_i$ be the positions of $f_i, e_i$ after the move ($i = 1, 2, 3$) and let $\Delta_1 = f_1$. The $f_i$'s are the arcs of $K'$ inside $D$ and the $e'_i$s are the basic arcs of $K'$ lying completely inside $D$. Let $\Delta_1 = f_1 \cap f_2, \Delta_2 = f_1 \cap f_3, \Delta_3 = f_2 \cap f_3$ and $\Delta'_1 = f'_1 \cap f'_2, \Delta'_2 = f'_1 \cap f'_3, \Delta'_3 = f'_2 \cap f'_3$.

Similarly to the cases of the $\Omega_1^-, \Omega_2^-$ moves, by the definition of the $\Omega_3$ move, the unions $e_1 \cup e_2 \cup e_3$ and $e'_1 \cup e'_2 \cup e'_3$ are closed curves on the plane that moreover bound rooms of $K, K'$ respectively, which by Schonflies Theorem are 2-disks $D_0, D'_0$. Unlike the $\Omega_1^-, \Omega_2^-$ cases, these 2-disks do not disappear but rather exchange positions, making $K \xrightarrow{\Omega_3} K'$ to have as an inverse move another one of the same type, namely the move $K' \xrightarrow{\Omega_3} K$ which inside $D$ brings $e'_1, e'_2, e'_3$ back to $e_1, e_2, e_3$. As in the $\Omega_1^-, \Omega_2^-$ cases we give a definition:

**Definition 19.** *We call as a 3-gon $e, f, g$ of a diagram $K$ with $e$ as the over-arc, $f$ as the middle arc and $g$ as the under-arc, any three distinct basic arcs $e, f, g$ of $K$, none a 1-gon and no two forming a 2-gon, each sharing one endpoint with another arc and the second endpoint with the third arc, so that $e$ is over $f, g$ at its endpoints, whereas $g$ is under $e, f$ at its endpoints, and also so that the closed curve $e_1 \cup e_2 \cup e_3$ of their union contains no other points of $K$ in its interior on the plane.*

Now, we choose a base point $p$ of $K$ not on $f_1, f_2, f_3$ and we set $q = p$ as the base point of $K'$. Alternatively, we choose $p \in e_i$ and then we put $q$ on $e'_i$.

Then after the move, no relabeling happens at the crossings outside $D$, and in all cases the old and new crossing are labeled as follows.

$$\Delta_1 = (\alpha : \gamma)^{\pi_1}, \Delta_2 = (\beta : \epsilon)^{\pi_2}, \Delta_3 = (\delta : \zeta)^{\pi_3},$$
$$\Delta'_1 = (\beta : \delta)^{\pi_1}, \ \Delta'_2 = (\alpha : \zeta)^{\pi_2}, \ \Delta'_3 = (\gamma : \epsilon)^{\pi_3}$$

where $\alpha, \gamma, \epsilon, \beta, \delta, \epsilon$ distinct integers in $\mathbb{N}_{2n}$ and $|\alpha - \beta| = |\gamma - \delta| = |\epsilon - \zeta| = 1$ modulo $2n$.

Let us notice, that all possible choices of the signs $\pi_1, \pi_2, \pi_3 \in \{-1, 1\}$ can happen.

Moreover, the height (first or second place) of the labels $\alpha, \gamma, \epsilon, \beta, \delta, \epsilon$ in the ordered pairs of the $\Delta'_i$'s is always the one that their respective companion $\beta, \delta, \zeta, \alpha, \gamma, \epsilon$ had in the ordered pairs of the $\Delta$'s. So if without any loss of generality we assume (as in Figure 27) $f_1$ over $f_2, f_3$ and $f_3$ under

$f_1, f_2$ then:

(2.19) $$\Delta_1 = (\alpha, \gamma)^{\pi_1}, \Delta_2 = (\beta, \epsilon)^{\pi_2}, \Delta_3 = (\delta, \zeta)^{\pi_3},$$
$$\Delta'_1 = (\beta, \delta)^{\pi_1}, \Delta'_2 = (\alpha, \zeta)^{\pi_2}, \Delta'_3 = (\gamma, \epsilon)^{\pi_3}.$$

And $_pK$ becomes $_pK'$ as follows:

(2.20) $$_pK = \cdots (x,y)^\pi \cdots \Delta_1, \Delta_2, \Delta_3 \overset{\Omega_3}{\mapsto} {_qK'} = \cdots (x,y)^\pi \cdots \Delta'_1, \Delta'_2, \Delta'_3.$$

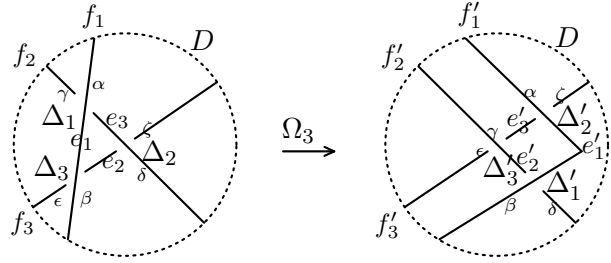

FIGURE 27. Labels of crossings for an $\Omega_3$ move, before and after the move. $K, K'$ have a common base point $p = q$ outside the $e_i$'s, or they have $p \in e_i$ and $q \in e'_1$ for some $i = 1, 2, 3$.

2.4. **Based symbols determine the cycles of the diagrams.** In the current section we intend to convince ourselves for the importance of the based symbols in describing a diagram $K$: we are going to show that knowledge of $_pK$ implies an algorithmic description of all oriented rooms of $K$, or equivalently of all oriented boundaries (oriented cycles) of the rooms (Corollary 1). For $\tau(K) > 0$, the algorithm retrieves information directly from $_pK$ and explains for any given main (basic or negative basic) arc $f$ which are its two cycles $C_1 = (f_1 = f, f_2, \ldots, f_k)_{cycl}$, $C_1 = (g_1 = f, g_2, \ldots, g_m)_{cycl}$.

Knowing $_pK$ we can recover almost all information about the various notions we have defined for $K$ so far, but we cannot recover the exact topological placement of $K$ on the plane as an oriented closed curve. For example, for $_pK = \emptyset$ we know that $K$ is a circle on the plane, but we do not know its exact orientation. So even if we know that there exist exactly two rooms $\Pi_{\alpha K}, \Pi'_{\delta K}$ each with $K$ as their oriented boundary, we cannot tell which of the rooms is the bounded and which is the unbounded one (Figure 28).

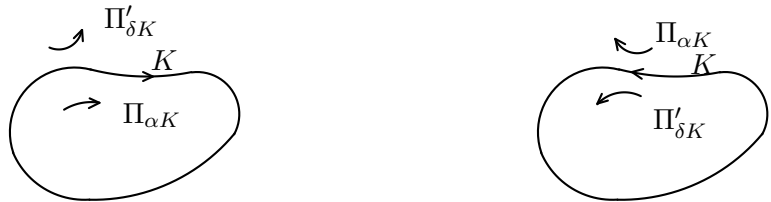

FIGURE 28. The distinct topological placements on the plane of an oriented diagram $K$ for which $_pK = \emptyset$ for some base point $p$. The two rooms of $K$ are also oriented so that they induce on their boundary $K$ the orientation it has as a diagram. The base point is irrelevant and not shown.

Similarly, for $_pK = (1,2)^{+1}$ we know that $K$ has a unique crossing point, and we know the way it crosses itself at that, but the exact topological position of $K$ on the plane could be any one of those placements shown in Figure 29. There depicted are all topologically non-equivalent diagrams with a based symbol $(1,2)^{+1}$. This means that no isotopy of the plane carries one diagram onto the other with all its over/under information.

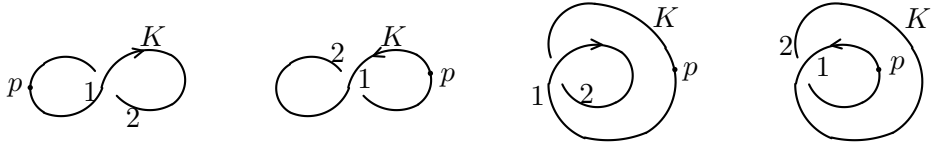

FIGURE 29. The distinct topological placements on the plane of an oriented diagram $K$ and a base point $p$ of it, for which $_pK = (1,2)^+$. They are distinct in the sense that there exists no isotopy of the plane bringing the pair $(K,p)$ from one of the depicted places onto another.

For $\tau(K) > 0$ in general, knowledge of $_pK$ informs us about the number of crossing points, the names of the basic arcs, the signs at the crossing points, the under- and over-arcs for each crossing point, and as we'll soon see, it informs us about the cycles of $K$ as well. So, although we cannot be sure beforehand about the exact placement of $K$ on the plane, we can be sure for example about the rooms which can accept an $\Omega_1^+$ move on a chosen basic arc $e$, or about the existence and the exact room or rooms that can accept an $\Omega_2^+$ move on two chosen basic arcs $e_1, e_2$.

By the end of the paper we will realize that although the exact placement of $K$ on the plane cannot be retrieved from $_pK$, nevertheless it contains enough information to prohibit a non-equivalent diagram to posses the same symbol. That is, it contains enough information to relate any two space loops projecting to diagrams with the same symbol (for some base choices) by the existence of an isotopy between the loops!

In the sequel we consider a diagram $K$ with $\tau(K) > 0$ and some fixed base point $p$ of it. Let $f = e$ be a basic arc or a negative basic arc $f = -e$ of $K$. And let us perform a trip on $K$ from the starting endpoint of $f$ until its finishing endpoint, turning there negatively or positively at will. As we have already mentioned (just before Definition 11) we will then find ourselves are on the boundary of the two rooms $\Pi, \Pi'$ which contain $f$ and which when oriented respectively negatively or positively on the plane induce to $f$ its orientation. In the notation of Definition 8 for the two rooms we have $\partial \Pi_{\alpha f}, \partial \Pi'_{\delta f}$ and according to the notation of Definition 11 the next main arc which we traverse will be $\alpha(f), \delta(f)$ respectively. The lemma that follows explains that knowledge of the based symbol $_pK$ reveals the exact basic arc or the exact negative basic arc into which we walk on after our turn.

**Lemma 10.** *For a diagram $K$ with $\tau(K) > 0$, a base point $p$ of $K$, a room $\Pi$ of $K$ and a basic arc $f = e = \overrightarrow{i(i+1)}$ or a negative basic arc $f = -e = -\overrightarrow{i(i+1)}$ which is part of the $\phi$-oriented boundary of $\Pi$, $\phi \in \{\alpha, \delta\}$ (that is: it holds $\partial \Pi_{\phi f}$), the $\phi$-turn $\phi(f)$ (the next edge of $e$ on the boundary of $\Pi_\phi$, ch. Definition 11) is determined as follows:*

*If $f = e = \overrightarrow{i(i+1)}$, let $\Delta = (i+1 : y) \in \, _pK$. If $f = -e = -\overrightarrow{i(i+1)}$, let $\Delta = (i : y) \in \, _pK$.*

*Then $\phi(f)$ has first endpoint $y$ and last endpoint $y' = y - \pi_\phi \pi_f \pi_\Delta v_y$ with this integer considered modulo $2n$ as an element of $\mathbb{N}_{2n}$. So if $y \prec y'$ it is $\phi(f) = \overrightarrow{yy'}$, whereas if $y' \prec y$ it is $\phi(f) = -\overrightarrow{yy'}$.*

We remind (ch. Definition 14) that $\prec$ denotes the cyclic ordering of $\mathbb{N}_{2n}$.

*Proof.* The proof is immediate by checking the validity of the formula in all cases in Figures 30, 31. □

**Corollary 1.** *Let $_pK$ be a given based symbol and $f$ a given main arc of the diagram $K$ for which $\tau(K) > 0$. Without need of knowledge of the actual placement of the diagram $K$ on the plane but with knowledge of $_pK$, there exists a finite algorithm which describes the cycles (oriented boundaries) $\partial \Pi_{\alpha f}, \partial \Pi'_{\delta f}$ of the two rooms $\Pi, \Pi'$ of $K$ that contain $f$ as an oriented arc on their boundaries. For $\phi \in \{\alpha, \delta\}$ the algorithm is the following: $f_1 = f, f_{i+1} = \phi(f_i), \forall i \in \mathbb{N}$ with $\phi(f_i)$ as described in*

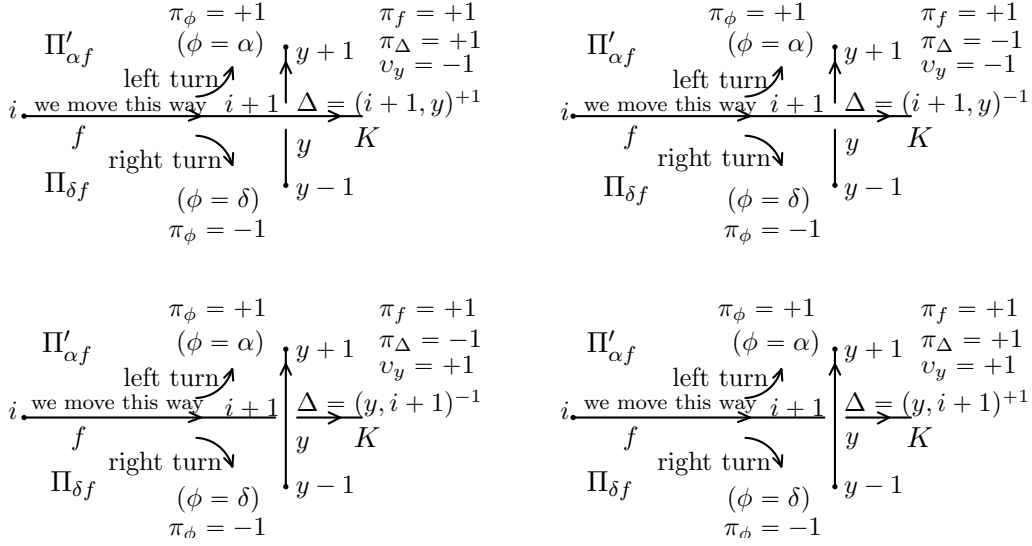

FIGURE 30. $f = i(i+1)$ is a basic arc of a diagram $K$. We make a trip on $\partial \Pi_{\delta f}, \partial \Pi'_{\alpha f}$ along their orientations starting at $i$. Let $y$ be the second label of the crossing point $\Delta$ at the end of $f$. The first label of $\Delta$ is $i+1$. All possible ways for continuing after $\Delta$ are given here.

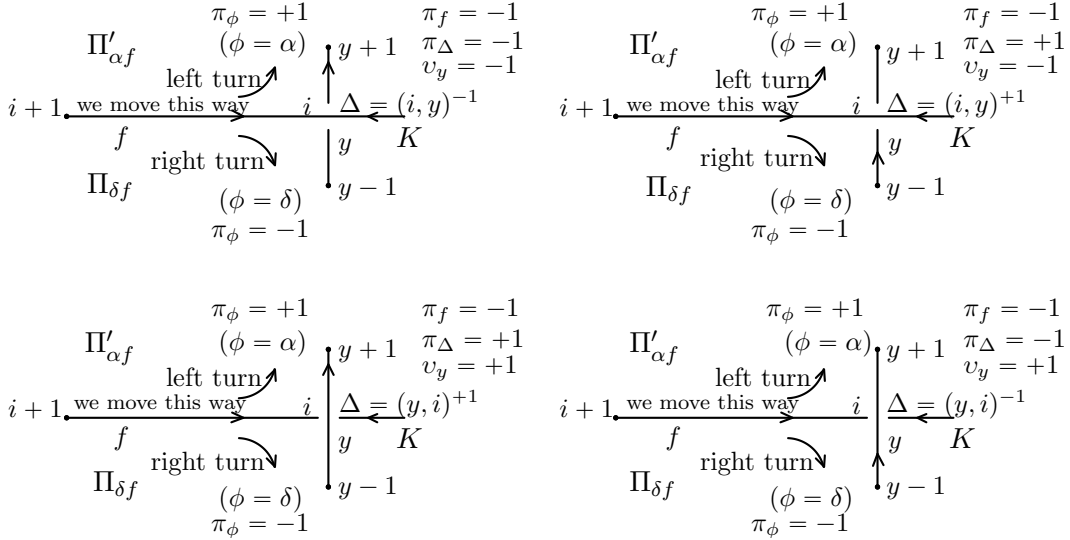

FIGURE 31. $f = -i(i+1)$ is a basic arc of a diagram $K$. We make a trip on $\partial \Pi_{\delta f}, \partial \Pi'_{\alpha f}$ along their orientations starting at $i$. Let $y$ be the second label of the crossing point $\Delta$ at the end of $f$. The first label of $\Delta$ is $i+1$. All possible ways for continuing after $\Delta$ are given here.

**Lemma 10.** $\partial \Pi_{\phi f} = (f_1 = f, f_2, \ldots, f_m)_{cycl}$ where $m+1$ is the smallest integer for which $f_{m+1}$ is an arc already appearing as $f_1$.

*Proof.* Let us fix an orientation, say $\phi = \alpha$ and similarly then for the orientation $\delta$. Let $x, y$ be the first and last endpoints of $f$ which we also call $f_1$. We begin a trip on $K$ with starting point its crossing $\Delta_1 = (x : x')$ and along the orientation of $f$. The next crossing we meet is at the endpoint of $f$, say $\Delta_2 = (y, y')$. There we make an $\alpha$-turn (left turn) traveling along a new main arc $f_2$ at the end of which we turn left again and so on.

This way we get a sequence of main arcs $f_1 = f, f_2, \ldots$ related as $f_{i+1} = \alpha(f_i)$, $\forall i \in \mathbb{N}$. Since there exist only finitely many main arcs, there exist an arc $f_m$ in the sequence that is repeated as $f_{m+k}$ for some $k > 0$. Let $m$ be the smallest index for which this happens. We claim that $m = 1$. To prove it we use the property $\phi(f) = g \Leftrightarrow -f = (-\phi(-g))$ which says that if in a $\phi$-cycle the main arc $g$ follows $f$, then reversing the orientation of the cycle we get a new one in which the main arcs $-f, -g$ lie where $-f$ follows $-g$. Topologically this is immediate, but of course here we need this result purely derived from the knowledge of $_pK$ alone. So we prove it as an independent Lemma at the end of the current Corollary (Lemma 11).

Now if it were $1 < m$, it would be $0 < m - 1$ and the arc $f_{m-1}$ would exist. But for the two arcs $f_{m-1}, f_{m+k-1}$ then we would have $\alpha(f_{m-1}) = f_m$ and $\alpha(f_{m+k-1}) = f_{m+k} = f_m$, thus $\alpha(f_{m-1}) = \alpha(f_{m+k-1})$. Then $-\alpha(f_{m-1}) = -\alpha(f_{m+k-1})$ hence $(-\alpha)(-\alpha(f_{m-1})) = (-\alpha)(-\alpha(f_{m+k-1}))$ which by Lemma 11 gives $-f_{m-1} = -f_{m+k-1}$ and this gives $f_{m-1} = f_{m+k-1}$ implying that $f_{m-1}$ is repeated too, a contradiction to the choice of $m$.

So it is indeed $m = 1$. Then the sequence of the $f_i$'s is actually a finite sequence $f = f_1, f_2, \ldots, f_m$ for some $m > 1$ in which all $f_i$'s are distinct which is repeated thereafter, i.e. $f_{\nu+m} = f_\nu, \forall \nu \in \mathbb{N}$. So $\partial \Pi_{\alpha f} = (f = f_1, f_2, \ldots, f_m)_{cycl}$.

Now by Lemma 10, knowledge of $_pK$ implies that given $f_1$ we algorithmically determine $f_2 = \alpha(f_1)$. The same Lemma implies an algorithmic knowledge of $f_3 = \alpha(f_2)$ and so on, and we have finished. □

Here is the Lemma needed in the last proof.

**Lemma 11.** *Let $_pK$ be a given based symbol and $f, g$ main arcs of $K$. Without need of knowledge of the actual placement of the diagram $K$ on the plane but with knowledge of $_pK$ alone, we can be sure that for any orientation $\phi \in \{\alpha, \delta\}$ it holds $\phi(f) = g \Rightarrow -f = (-\phi)(-g) = (-\phi)(-\phi(f))$.*

*Proof.* We shall prove $\phi(f) = g \Rightarrow -f = (-\phi)(-g)$ and then recalling that $g = \phi(f)$ we get $-f = (-\phi)(-g) = (-\phi)(-\phi(f))$ as wanted. For the relation $\phi(f) = g \Leftrightarrow -f$:

Let $a$ and $b$ be the first and last endpoints of $f$. Let $c$ and $d$ be the first and last endpoints of $g$.
$$a, b, c, d \in \mathbb{N}_{2\tau(K)} \text{ with } |a - b| = |c - d| = 1. \tag{1}$$
$\phi(f) = g \Rightarrow \Delta = (b, c)$ is a crossing of $K$ and by Lemma 10:
$$d = c - \pi_\phi \pi_f \pi_\Delta v_c. \tag{2}$$
Also, $-f$ has $b$ and $a$ as first and last endpoints respectively, whereas $-g$ has $d$ and $c$ as first and last endpoints respectively.

The crossing in $_pK$ with the endpoint $c$ of $-g$ as one label is $\Delta = (b, c)$. By Lemma 10 again we get: $(-\phi)(-g) =$ a main arc with first endpoint $b$ and last endpoint:
$$y = b - \pi_{-\phi} \pi_{-g} \pi_\delta v_b = b - (-\pi_\phi)(-\pi_g)(-v_c) = b + \pi_\phi \pi_g v_c \stackrel{(2)}{=} b + \frac{\pi_g(c-d)}{\pi_f}. \tag{3}$$
It also holds:
$$\pi_f = b - a, \ \pi_g = d - c \tag{4}$$
Then (3),(4) give:
$$y = b + \frac{(d-c)(c-d)}{\pi_\phi} = b + \frac{-(d-c)^2}{\pi_\phi} \stackrel{(1)}{=} b - \frac{1}{\pi_f} \stackrel{\pi_f = \pm 1}{=} b - \pi_f \stackrel{(4)}{=} a.$$
So $(-\phi)(-g)$ has the same first and last endpoints as $-f$ and then it oughts to be $-f$ as wanted. □

We illustrate the content of Lemma 10 and of Corollary 1 with a couple of examples:

**Example 1.** *For the diagram $K$ of Figure 32 with the denoted base point $p$ and for the basic arc $f = \overrightarrow{(15)(16)}$, we can calculate $\partial \Pi_{\alpha f}, \partial \Pi'_{\delta f}$ directly from what we see in the figure. We can also calculate them assuming we were given only the based symbol $_pK$ instead. The based symbol is:*
$$_pK = (1, 8)^{+1}(2, 9)^{-1}(3, 12)^{+1}(13, 4)^{+1}(5, 14)^{+1}(15, 6)^{+1}(7, 19)^{+1}(10, 17)^{-1}(11, 16)^{-1}$$

and we calculate $\partial \Pi'_{\delta f}, \partial \Pi_{\alpha f}$ as follows (after fixing $\phi$, all information is retrieved from $_pK$):

$f = \overrightarrow{(15)(16)}, \phi = \delta,$

$$
\begin{array}{l}
\phantom{(15,)}\phantom{(13,)}\phantom{(11,)}\phantom{(3,)}\;\;15 \\
\phantom{(15,)}\phantom{(13,)}\phantom{(11,)}\phantom{(3,)}\downarrow f_1 \\
\phantom{(15,)}\phantom{(13,)}(11, \phantom{..}16)^{+1} = \Delta_1 \\
\phantom{(15,)}\phantom{(13,)}\phantom{(11,)}\downarrow f_2 \\
\phantom{(15,)}\phantom{(13,)}(3, \phantom{..}12)^{+1} = \Delta_2 \\
\phantom{(15,)}\phantom{(13,)}\phantom{(3,)}\downarrow f_3 \\
\phantom{(15,)}(13, \phantom{..}4)^{+1} = \Delta_3 \\
\phantom{(15,)}\phantom{(13,)}\downarrow f_4 \\
\phantom{(15,)}(5, \phantom{..}14)^{+1} = \Delta_4 \\
\phantom{(15,)}\downarrow f_5 \\
(15, \phantom{..}6)^{+1} = \Delta_5 \\
\downarrow f_6 = f_1 \\
16
\end{array}
$$

$-\pi_\delta \pi_{f_1} \pi_{\Delta_1} v_{11} = -(-1)(+1)(+1)(+1) = +1$

$-\pi_\delta \pi_{f_2} \pi_{\Delta_2} v_3 = -(-1)(+1)(+1)(+1) = +1$

$-\pi_\delta \pi_{f_3} \pi_{\Delta_3} v_{13} = -(-1)(+1)(+1)(+1) = +1$

$-\pi_\delta \pi_{f_4} \pi_{\Delta_4} v_5 = -(-1)(+1)(+1)(+1) = +1$

$-\pi_\delta \pi_{f_1} \pi_{\Delta_5} v_{15} = -(-1)(+1)(+1)(+1) = +1$

So:
$$\partial \Pi'_{\delta f} = \cup f_i = \overrightarrow{(15)(16)} \cup \overrightarrow{(11)(12)} \cup \overrightarrow{(3)(4)} \cup \overrightarrow{(13)(14)} \cup \overrightarrow{(5)(6)}.$$

*Similarly:*

$f = \overrightarrow{(15)(16)}, \phi = \alpha,$

$$
\begin{array}{l}
\phantom{f_1 =}\phantom{(15,)}\phantom{(7,)}\;\;15 \\
\phantom{f_1 =}\phantom{(15,)}\phantom{(7,)}\downarrow f_1 \\
\phantom{f_1 =}\phantom{(15,)}(11, \phantom{..}16)^{+1} = \Delta_1 \\
\phantom{f_1 =}\phantom{(15,)}\downarrow f_2 \\
\phantom{f_1 =}\phantom{(15,)}(10, \phantom{..}17)^{-1} = \Delta_2 \\
\phantom{f_1 =}\phantom{(15,)}\downarrow f_3 \\
\phantom{f_1 =}\phantom{(15,)}(7, \phantom{..}18)^{+1} = \Delta_3 \\
\phantom{f_1 =}\phantom{(15,)}\downarrow f_4 \\
\phantom{f_1 =}(15, \phantom{..}6)^{+1} = \Delta_4 \\
f_1 = \;\downarrow f_5 \\
\phantom{f_1 =}16
\end{array}
$$

$-\pi_\alpha \pi_{f_1} \pi_{\Delta_1} v_{11} = -(+1)(+1)(+1)(+1) = -1$

$-\pi_\alpha \pi_{f_2} \pi_{\Delta_2} v_{17} = -(+1)(-1)(-1)(-1) = +1$

$-\pi_\alpha \pi_{f_3} \pi_{\Delta_3} v_7 = -(+1)(+1)(+1)(+1) = -1$

$-\pi_\alpha \pi_{f_4} \pi_{\Delta_4} v_{15} = -(+1)(-1)(+1)(+1) = +1$

So:
$$\partial \Pi_{\alpha f} = \cup f_i = \overrightarrow{(15)(16)} \cup (-\overrightarrow{(10)(11)}) \cup \overrightarrow{(17)(18)} \cup (-\overrightarrow{(6)(7)}).$$

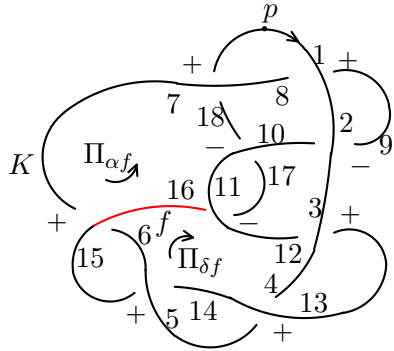

FIGURE 32. A diagram $K$ and a base point $p$ of it. The symbol $_pK$ is calculated according to the definition. Knowledge of $_pK$ alone, suffices to describe all oriented cycles (and oriented rooms) like for example the left and right cycles $\partial \Pi_{\alpha f}, \partial \Pi'_{\delta f}$ of the basic arc $f = (15)(16)$.

**Example 2.** *For the diagram $K$ of Figure 33 with the denoted base point $p$ and for the basic arc $f = \overrightarrow{(2)(3)}$, we can calculate $\partial \Pi_{\delta f}, \partial \Pi'_{\alpha f}$ directly from what we see in the figure. We can also calculate them assuming we were given the based symbol $_pK$ instead. The based symbol is:*

$$_pK = (1,2)^{-1}(3,4)^{-1}(5,6)^{-1}$$

*and we calculate $\partial \Pi'_{\delta f}, \partial \Pi_{\alpha f}$ as follows (after fixing $\phi$, all information is retrieved from $_pK$):*

$f = \overrightarrow{(2)(3)}, \phi = \delta,$

$$\begin{array}{l}
2 \\
\downarrow \\
(3 \quad 4)^{-1} \\
\quad \downarrow \\
\quad (5, \quad 6)^{-1} \\
\qquad \downarrow \\
\qquad (1, \quad 2)^{-1} \\
\qquad \quad \downarrow \\
\qquad \quad 3
\end{array}
\qquad
\begin{array}{l}
-\pi_\delta \pi_{(2)(3)} \pi_{(3,4)^{-1}} v_4 = -(-1)(+1)(-1)(-1) = +1 \\[4pt]
-\pi_\delta \pi_{(4)(5)} \pi_{(5,6)^{-1}} v_6 = -(-1)(+1)(-1)(-1) = +1 \\[4pt]
-\pi_\delta \pi_{(6)(1)} \pi_{(1,2)^{-1}} v_2 = -(-1)(+1)(-1)(-1) = +1
\end{array}$$

*So:*

$$\partial \Pi'_{\delta f} = \cup f_i = \overrightarrow{(2)(3)} \cup \overrightarrow{(4)(5)} \cup \overrightarrow{(6)(1)}.$$

*Similarly:*

$f = \overrightarrow{(2)(3)}, \phi = \alpha,$

$$\begin{array}{l}
2 \\
\downarrow \\
(3, \quad 4)^{-1} \\
\quad \downarrow \\
\quad (3, \quad 4)^{-1} \\
\qquad \downarrow \\
\qquad (5, \quad 6)^{-1} \\
\qquad \quad \downarrow \\
\qquad \quad (5, \quad 6)^{-1} \\
\qquad \qquad \downarrow \\
\qquad \qquad (1, \quad 2)^{-1} \\
\qquad \qquad \quad \downarrow \\
\qquad \qquad \quad (1, \quad 2)^{-1} \\
\qquad \qquad \qquad \downarrow \\
\qquad \qquad \qquad 3
\end{array}
\qquad
\begin{array}{l}
-\pi_\alpha \pi_{(2)(3)} \pi_{(3,4)^{-1}} v_4 = -(+1)(+1)(-1)(-1) = -1 \\[4pt]
-\pi_\alpha \pi_{(4)(3)} \pi_{(3,4)^{-1}} v_4 = -(+1)(-1)(-1)(-1) = +1 \\[4pt]
-\pi_\alpha \pi_{(4)(5)} \pi_{(5,6)^{-1}} v_6 = -(+1)(+1)(-1)(-1) = -1 \\[4pt]
-\pi_\alpha \pi_{(6)(5)} \pi_{(5,6)^{-1}} v_6 = -(+1)(-1)(-1)(-1) = +1 \\[4pt]
-\pi_\alpha \pi_{(6)(1)} \pi_{(1,2)^{-1}} v_2 = -(+1)(+1)(-1)(-1) = -1 \\[4pt]
-\pi_\alpha \pi_{(2)(1)} \pi_{(1,2)^{-1}} v_2 = -(+1)(-1)(-1)(-1) = +1
\end{array}$$

*So:*

$$\partial \Pi_{\alpha f} = \cup f_i = \overrightarrow{(2)(3)} \cup (-\overrightarrow{(3)(4)}) \cup \overrightarrow{(4)(5)} \cup (-\overrightarrow{(6)(5)}) \cup \overrightarrow{(6)(1)} \cup (-\overrightarrow{(1)(2)}).$$

**Remark 2.** *A subtle but crucial point in Lemma 10 may go unnoticed in a first reading. Namely, the labels at the crossings in Figures 30, 31 incorporate the deeply topological property that at the crossing points the intersection of a diagram with itself is transverse. So for example the labels in Figure 34 (a) are permitted whereas those in Figure 34(b) are not. This property will be part of the definitions in the algebrization of §3 that follows.*

*In Lemmata 13 and 14 of §3 we shall also prove a few analogues of some extremely important topological facts regarding diagrams. Among them the analogue of the existence of oriented boundaries which assures that every (algebraic) main arc lies in exactly two cycles, defined as cyclically ordered sequences consisting of distinct main arcs where each arc is the left (or right) turn of the previous one.*

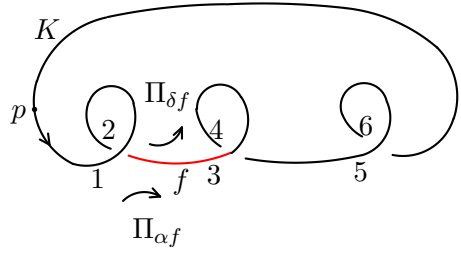

FIGURE 33. A diagram $K$ and a base point $p$ of it. The symbol $_pK$ is calculated according to the definition. Knowledge of $_pK$ alone, suffices to describe all oriented cycles (and oriented rooms) like for example the left and right cycles $\partial \Pi'_{\delta f}, \partial \Pi_{\alpha f}$ of the basic arc $f = (2)(3)$.

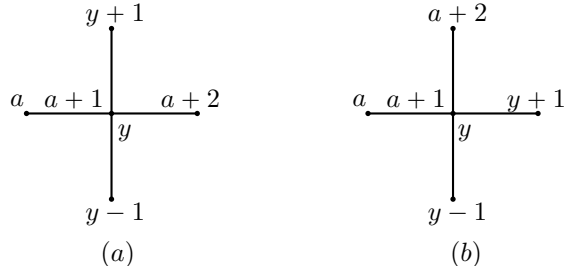

FIGURE 34. At any crossing of a diagram the intersection is transverse, thus (a) can happen but (b) cannot.

## 3. Algebrization

In previous sections we have seen how to assign to any given diagram $K$ a bunch of based symbols $_pK$ where $p$ are points of $K$ other than crossings. We have also seen that knowledge of $_pK$ actually determines the cycles of $K$, although it does not determine the exact topological placement of $K$ on the plane as well. And we have seen how $_pK$ changes when we change the base point $p$ or when we perform a topological move on $K$ (an Alexander-Briggs-Reidemeister move or an isotopy of the plane).

We start now the effort to work out the changes of the based symbols in a completely algebraic setting. We call $\Sigma_D$ the set of all based symbols of diagrams and we intend to describe algebraic moves $A \xrightarrow{\omega} B$ for $A, B \in \Sigma_D$ which correspond to the topological ones in the hope that they generate some equivalence relation in $\Sigma_D$. Succeeding in this would allow us (§4) to assign to any knot $\varkappa$, the equivalence class of any symbol of any of its representative diagrams.

The definition we are about to give is a long one and treats the based symbols as purely algebraic objects. It is based on the analysis we have done so far in the topological case of diagrams and their moves. We start the definition repeating the useful notion of $\mathbb{N}_{2n}$ and its cyclic ordering (Definition 14). We continue repeating the notion of our main object which is that of a based symbol (Definition 13). We then define analogues of basic arcs, cycles, $i$-gons etc. At the end we define some new symbols similar to the based ones which we call as numbered symbols. The set of numbered symbols includes the set of based symbols, and its role is auxiliary, namely, to ensure that the results of the algebraic moves have a temporary place to stay, since although they are proved eventually to be based symbols, this is not at all immediate from the very beginning. We define the algebraic moves in a bunch of separate definitions after we prove a few important results regarding the rest of the defined notions below.

**Definition 20.** We write $\mathbb{N}_0 = \emptyset$, and for $n > 0$ we call $\mathbb{N}_{2n} = \{1, 2, \ldots, 2n\}$ considered as an additive group with identity the element $2n$. For $n < m$, it is $\mathbb{N}_{2n} \subset \mathbb{N}_{2m}$. We consider $\mathbb{N}_{2n}$ cyclically ordered as $1 < 2$, $2 < 3, \ldots, 2n-1 < 2n$, $2n < 1$ and with no other pair of elements related in this order. We call any two comparable elements by $<$ as consecutive.

For any $m > 0$ and any $x \in \mathbb{N}_{2m}$ we write $\mathbb{N}_{2m,x} = \mathbb{N}_{2m} - \{x, x+1\}$.

For any $m > 0$ and any $x, y \in \mathbb{N}_{2m}$ with $|x-y| \neq 1$ we write $\mathbb{N}_{2m,x,y} = \mathbb{N}_{2m} - \{x, x+1, y, y+1\}\}$ for $|x-y| \geq 2$ and $\mathbb{N}_{2m,x,x} = \mathbb{N}_{2m} - \{x, x+1, x+2, x+3\}\}$ for $|x-y| = 0$.

For a diagram $K$ with $n > 0$ crossings we call as its based symbol at its base point $p$ the set $\{\Delta_1, \Delta_2, \ldots, \Delta_n\} = \{(i_1, i_2)^{\pi_1}, (i_3, i_4)^{\pi_1}, \ldots, (i_{2n-1}, i_{2n})^{\pi_n}\}$ of its crossings $\Delta_i$ equipped with all of their information. We conveniently denote the based symbol at $p$ as:

$$_pK = (i_1, i_2)^{\pi_1} (i_3, i_4)^{\pi_1} \ldots (i_{2n-1}, i_{2n})^{\pi_n},$$

where the order of the pairs is not important.

We denote the crossing $(a, b)^\pi$ without some of its information also as $a:b$, $b:a$, $(a, b)$, $(a:b)$, $(b:a)$, $(a:b)^\pi$ or $(b:a)^\pi$.

For a crossing $\Delta = (i, j)^\pi$ we define as height of $i$ and $j$ the numbers $v_i = 1, v_j = -1$. We recall that $\pi \in \{-1, +1\}$ is the sign of $\Delta$.

For $K$ the unknot, we define as its based symbol at $p$ the empty set $_pK = \emptyset$ for all $p \in K$. We call $\emptyset$ as the empty based symbol. This symbol contains no crossings by definition.

$\Sigma_\Delta$ denotes the set of all based symbols of all diagrams on the plane.

We define as order $\tau(A)$ of any $A \in \Sigma_D$, the number of crossings it contains.

We denote a crossing $(a, b)^\pi$ without some of its information also as $a:b, b:a, (a, b), (a:b), (b:a), (a:b)^\pi, (b:a)^\pi$.

Let us fix some $A \in \Sigma_D$ with $\tau(A) > 0$, say $A = (i_1, i_2)^{\pi_1}(i_3, i_4)^{\pi_2} \ldots (i_{2n-1}, i_{2n})^{\pi_n}$.

For disjoint subsets $A_1, A_2$ of $A$ we write $A = A_1 A_2$ allowing the subsets to be the empty set as well.

We call as basic arcs of $A$ the symbolic objects $e_1 = \overrightarrow{12}, e_2 = \overrightarrow{23}, \ldots, e_{n-1} = \overrightarrow{(2n-1)(2n)}, e_n = \overrightarrow{(2n)1}$. We denote an arbitrary basic arc as $e_i = \overrightarrow{i(i+1)}$ where $i \in \mathbb{N}$.

For $\emptyset \in \Sigma_D$ we do not define basic arcs.

For the basic arc $e_i = \overrightarrow{i(i+1)}$ the numbers $i$ and $i+1$ are called respectively as the first and last endpoints of $e$.

We call as negative basic arcs of $A$ the symbolic objects $-e_1 = -\overrightarrow{12}, -e_2 = -\overrightarrow{23}, \ldots, -e_{n-1} = -\overrightarrow{(2n-1)(2n)}, -e_n = -\overrightarrow{(2n)1}$. For $-e_i$, the numbers $i+1$ and $i$ are called respectively as the first and last endpoints.

We call the basic arcs and the negative basic arcs, with the common name of main arcs of $A$.

For a main arc $f$ with endpoints $a$ and $b$ it is certainly true that $|a-b| = 1 \bmod 2n$. We denote the arc as $f = ab = ba$ without specifying its first and last endpoints. If $a$ is the first endpoint and $b$ the last, then we define the sign of $f$ to be the number $f_e = b - a \in \{-1, +1\}$.

We call the letters $\alpha, \delta$ as negative and positive orientations regardless of the specific based symbol $A$. We call the integers $\pi_\alpha = +1, \pi_\delta = -1$ as the signs of the orientations.

For a basic arc $f = e = \overrightarrow{i(i+1)}$ or a negative basic arc $f = -e = -\overrightarrow{i(i+1)}$ of $A$, and for $\phi \in \{\alpha, \delta\}$, we define a main arc which we call the $\phi$-turn $\phi(f)$ of $f$ as follows:

If $f = e = \overrightarrow{i(i+1)}$, let $\Delta = (i+1:y) \in A$. If $f = -e = -\overrightarrow{i(i+1)}$, let $\Delta = (i:y) \in A$. Then $\phi(f)$ has first endpoint $y$ and last endpoint $y' = y - \pi_\phi \pi_f \pi_\Delta v_y$ (with this integer considered as an element of $\mathbb{N}_{2n}$). So if $y < y'$ it is $\phi(f) = \overrightarrow{yy'}$, whereas if $y' < y$ it is $\phi(f) = -\overrightarrow{yy'}$.

For a main arc $f$ of $A$, we write $C_{\alpha f} = (f_1 = f, f_2, f_3 \ldots)$ with $\alpha(f_i) = f_{i+1}, \forall i$ and $C_{\delta f} = (g_1 = f, g_2, g_3 \ldots)$ with $\delta(g_i) = g_{i+1}, \forall i$. We call $C_{\alpha f}, C_{\delta f}$ as the $\alpha$- and $\delta$-cycles of $f$ respectively. We also call them as (oriented) cycles of $A$. By Lemma 13 below, the cycles are sequences of finitely *many distinct arcs which are repeated indefinitely. So from now on we think of them as cyclically ordered sets $C_{\alpha f} = (f_1 = e, f_2, \ldots, f_m)_{cycl}, C_{\delta f} = (g_1 = e, g_2, \ldots, g_k)_{cycl}$ each with distinct elements, which we consider respectively as the $\alpha$-cycle of every $f_i$ and the $\delta$-cycle of every $g_i$. Thus we can write $C_{\alpha f_i} = C_{\alpha f_j}, \forall i, j$ and $C_{\delta g_i} = C_{\delta g_j}, \forall i, j$. We also write $f_i \in C_{\alpha f}, g_i \in C_{\alpha f}$.*

*Restricting attention to basic arcs, we say that two basic arcs $e, f$ have the same cycle sense (of orientation) in A whenever $e, f$ are elements of the same cycle of A; and we say that $e, f$ have opposite cycle senses (of orientation) in A whenever $e, -f$ are elements of the same cycle of A or equivalently (by Lemma 14(d) below) whenever $-e, f$ are elements of the same cycle of A. We write respectively $e \upuparrows f$ and $e \updownarrows f$. By Lemma 14(c) below, the same and opposite cycle senses are distinct notions. Then we can also define the number of similarity of cycle senses as $\rho_{e,f} = 1$ whenever $e, f$ have the same cycle sense, and $\rho_{e,f} = -1$ whenever $e, f$ have opposite cycle senses. Since by Lemma 13 cycles have a cyclic nature, it is true that both the definition of cycle sense as well as of the number of similarity of cycle senses are independent of the order of $f_1, f_2$; thus $\rho_{f_2, f_1} = \rho_{f_1, f_2}$.*

*Remaining to basic arcs, let $e_1$ and $e_2$ be two of them. For $e_1 \neq e_2$ we define the placement number of the ordered pairs $(e_1, e_2)$, $(e_2, e_1)$ to be two variables $\theta_{e_1, e_2}, \theta_{e_2, e_1}$ with values in $\{-1, +1\}$, so that $\theta_{e_2, e_1} = -\theta_{e_1, e_2}$; whenever $\theta_{e_i, e_j} = 1$ we say that $e_i$ is over $e_j$ and $e_j$ is under $e_i$ (since $\theta_{e_2, e_1} = -\theta_{e_1, e_2}$, the terminology is consistent). For $e_1 = e_2 = e$ we define the placement number of e with respect to itself to be a variable $\theta_{e,e} \in \{-1, 1\}$. For $\theta_{e,e+} = 1$ we say that e is placed over and then under itself, otherwise we say it is placed under and then over itself.*

*The empty based symbol $\varnothing$ does not contain basic arcs (since it contains no crossings), but we still define the placement number of the empty based symbol $\varnothing$ with respect to itself as a variable $\theta_{\varnothing, \varnothing} \in \{-1, 1\}$. For $\theta_{\varnothing, \varnothing} = +1$ we say that $\varnothing$ is placed over and then under itself, otherwise we say it is placed under and then over itself.*

*For a basic arc e of A, we define the twist-number of e as a variable $\theta_e$ with values in $\{-1, +1\}$. Whenever $\theta_e = +1$ we say that e is twisted over itself, otherwise we say it is twisted under itself.*

*For the empty based symbol $\varnothing$ we define its placement- or twist-number as a variable $\theta_\varnothing$ with values in $\{-1, +1\}$. Whenever $\theta_e = +1$ we say that $\varnothing$ is twisted over itself, otherwise we say it is twisted under itself.*

*Let $e = \overrightarrow{a(a+1)}$ be a basic arc of A. We say that e forms a 1-gon or a 1-loop of A whenever $\Delta = (a : a+1)$ is a crossing of A. We call $\Delta$ as the crossing of the 1-gon.*

*For an 1-gon e we define the placement- or twist-number of e as a variable $\theta_e^-$ with values in $\{-1, +1\}$.*

*Let $e_1 = \overrightarrow{a(a+1)}$, $e_2 = \overrightarrow{b(b+1)}$ be two distinct basic arcs of A non of which is a 1-gon. We say that $e_1, e_2$ form a 2-gon of A whenever $\Delta_1 = (a : b), \Delta_2 = (a+1 : b+1)$ or $\Delta_1 = (a : b+1), \Delta_2 = (a+1 : b)$ are crossings of A and the height of both $a, a+1$ is $+1$ or of both $-1$. We call the arc whose labels have heights $+1$ as the top arc and the other as the bottom one, and we say that the first is above the second (or that the second is below the first). We call $\Delta_1, \Delta_2$ as the crossing of the 2-gon.*

*For a 2-gon e we define the placement number of the ordered pairs $(e_1, e_2), (e_2, e_1)$ as two variables $\theta_{e_1, e_2}^-, \theta_{e_2, e_1}^-$ with values in $\{-1, +1\}$ so that $\theta_{e_2, e_1}^- = -\theta_{e_2, e_1}^-$.*

*Let $e_1 = \alpha\beta$, $e_2 = \gamma\delta$, $e_3 = \epsilon\zeta$ be three distinct basic arcs of A, non of them a 1-gon, and no two of them a 2-gon. We say that $e_1, e_2, e_3$ form a 3-gon or a triangle of A whenever $\Delta_1 = (\alpha : \gamma), \Delta_2 = (\beta : \epsilon), \Delta_3 = (\delta : \zeta)$ are crossings of A with $|\alpha - \beta| = |\gamma - \delta| = |\epsilon - \zeta| = 1$, and the labels of both endpoints of one of the arcs have height $+1$ (we call this arc as top arc and say it is above the other two); so then the endpoints of another arc have height $-1$ (we call it bottom arc and say it is below the other two), and one endpoint of the third arc has label $+1$ while the other $-1$ (and we call the arc as the middle one). We call $\Delta_1, \Delta_2, \Delta_3$ as the crossings of the 3-gon and we also* say that $\Delta_1, \Delta_2, \Delta_3$ form a 3-gon (or a triangle). We call the two labels of each of the crossings as companions; so $\alpha, \beta$ are companions, also $\gamma, \delta$ are companions and $\epsilon, \zeta$ are companions.

We call as a numbered symbol the empty set $\varnothing$ and any expression

$$(i_1, i_2)^{\pi_1}(i_3, i_4)^{\pi_2}\ldots(i_{2n-1}, i_{2n})^{\pi_n},$$

for some $n \in \mathbb{N}$, $n \geqslant 1$ when the set of elements $i_j$ is exactly $\mathbb{N}_{2n} = \{1, 2, \ldots, 2n\}$ and each pair is considered ordered. The order of the pairs in the notation is not important.

We call the pairs (if any) in a numbered symbol $A$ as crossings of $A$. So $\varnothing$ does not have any crossings.

We denote the set of all numbered symbols by $\Sigma$.

For $A \in \Sigma$ we define its order $\tau(A)$ to be 0 if $A = \varnothing$ and $n$ if $A = (i_1, i_2)^{\pi_1}(i_3, i_4)^{\pi_1}\ldots(i_{2n-1}, i_{2n})^{\pi_n}$.

For a crossing $\Delta = (a, b)^\pi$ of some numbered symbol $A$, we call $\pi$ as sign of $\Delta$, and the numbers $v_a = -1, v_b = +1$ as the height of $a, b$ respectively. We denote a crossing $(a, b)^\pi$ without some of its information also as $a : b, b : a, (a, b), (a : b), (b : a), (a : b)^\pi$ or $(b : a)^\pi$.

**Remark 3.** (a) The algebraic 1-gons, 2-gons and 3-gons defined here, do not exactly incarnate the corresponding topological ones because they do not contain any provision related to some kind of algebraic interior of the $i$-gons.

We will overpass this nuisance later on (§4.5, 4.6, 4.7) with the help of the isotopy move $\Omega_{iso\sigma}$ in the 2-sphere $S^2 = \rho \cup \{\infty\}$ (§1.5), and with the creation of some new topological moves $\Omega_{1\gamma}^-, \Omega_{2\gamma}^-, \Omega_{3\gamma}$ (§4.4).

(b) Of course $\Sigma_\Delta \subseteq \Sigma$ but it is not hard to show that $\Sigma_\Delta$ is actually a proper subset of $\Sigma$.

For the algebraic moves $A \xrightarrow{\omega} A'$ defined below in §3.1-3.5 it will be relatively easy to prove that $A' \in \Sigma$ (Lemmata 16, 17, 18, 19, 20). It will be considerably more demanding to prove that actually $A; \in \Sigma_D$ (Lemmata 22, 23, 25, 26, 27). This effort will also reveal that the algebraic moves are mirrored by the topological moves. This would imply eventually the desired correspondence between the equivalence of based symbols and the equivalence of knot diagrams.

In the following Lemmata we gather some important properties of cycles of based symbols. We expect then to hold since based symbols are derived by diagrams and since Definition 20 is nothing more than an algebraic synopsis to the various notions relevant to the diagrams introduced so far. We prove these properties here quite algebraically, assuming knowledge only of the symbol itself as an algebraic object without any mention to the topological nature of the original diagram. This reinforces our conviction that a number of important topological properties of the diagrams are determined by the purely algebraic properties of their based symbols.

**Lemma 12.** Let $A \in \Sigma_D$ be a given based symbol and $f, g$ main arcs of $A$. Then for any orientation $\phi \in \{\alpha, \delta\}$ it holds $\phi(f) = g \Rightarrow -f = (-\phi)(-g) = (-\phi)(-\phi(f))$.

*Proof.* We repeat the proof of Lemma 11 which as already noted was given so as to hold with no dependence on diagrams but rather only on based symbols. Only now we replace the topological notions of crossing point, sign of a crossing point, height of the labels of a crossing point, main arc, first and last point of a main arc, sing of a main arc, sign of an orientation, left- or right-turn of a main arc, which are used implicitly or explicitly in the proof, with the corresponding algebraic notions defined for based symbols as algebraic objects in Definition 20. □

**Lemma 13.** Let $A \in \Sigma_D$ be a based symbol and $f$ a main arc of $A$. Then each one of the two cycles $C_{\alpha f} = (f_1 = f, f_2, \ldots)$, $C_{\delta f} = (g_1 = f, g_2, \ldots)$ of $f$ is actually a finite sequence of distinct arcs which is repeated indefinitely. Thus we can write them as $C_{\alpha f} = (f_1 = e, f_2, \ldots, f_m)_{cycl}, C_{\delta f} = (g_1 = e, g_2, \ldots, g_k)_{cycl}$ for some $m, k > 1$ where $f_1, \ldots, f_m$ distinct and $\alpha(f_i) = f_{i+}, \forall i$, and similarly $g_1, \ldots, g_m$ distinct and $\delta(g_j) = g_{j+1}, \forall j$.

*Proof.* We repeat the proof of Corollary 1. Only now we replace all topological notions used in the proof of the Corollary with the corresponding in Definition 20. And we replace the use of Lemma 10 with the algebraic definition of the left- or right-turn of a main arc in Definition 20, and finally we replace Lemma 11 with its purely algebraic version, i.e. Lemma 12. □

**Lemma 14.** *Let $A \in \Sigma_D$ be a based symbol and $f$ a main arc of $A$. Then:*

*(a) The cycles $C_{\alpha f}, C_{\delta f}$ are distinct (they do not coincide as cyclically ordered sets), thus $f$ has exactly two cycles.*

*(b) The cycles $C_{\alpha f}, C_{\delta f}$ cannot contain $-f$.*

*(c) If $e$ is a main arc and $\phi \in \{\alpha, \delta\}$ then it is not possible that $e, -e \in C_{\phi f}$.*

*(d) If $C_{\phi f} = (f_1, f_2, \ldots, f_k)_{cyl}$ then $C_{-\phi(-e)} = (-f_k, -f_{k-1}, \ldots, -f_2, -f_1)_{cycl}$.*

*Proof.* (a) We fix an orientation $\phi \in \{\alpha, \delta\}$ and we work in $C_{\phi f}$.

Assume on the contrary that $-f \in C_{\phi f}$. Let $a$ be the first endpoint and $b$ the last endpoint of $e$. And let $\Delta = (b : y)$ be the crossing of $A$ with $b$ as one of its labels. Call $e$ the main arc that follows $f$ and $g$ the main arc which precedes $-f$ in $C_{\phi f}$. So $\phi(f) = e$ and $\phi(g) = -e$.

The first and last endpoints of $-f$ are $b$ and $a$ respectively.

Since $\phi(g) = -f$ and since $-f$ has $b$ as first endpoint, we have by the definition of $\phi(g)$ that the last endpoint of $\phi(g)$ should be $b - \pi_\phi \pi_g \pi_\Delta v_b$. Since $\phi(g) = -f$ we must have:
$$b - \pi_\phi \pi_g \pi_\Delta v_b = a. \tag{1}$$

By Lemma 12 we have $\phi(f) = e \Rightarrow (-\phi)(-e) = -f$. Since $\phi(f) = e$, the first endpoint of $e$ is $y$, so the last endpoint of $-e$ is $y$. The crossing containing $y$ as a label is $\Delta$ and the other label of $\Delta$ is $b$. Hence by the definition of $(-\phi)(-e)$, the last endpoint of $(-\phi)(-e)$ is $b - \pi_{-\phi}\pi_{-e}\pi_\Delta v_b$. Since $(-\phi)(-e) = -f$ we have:
$$b - \pi_{-\phi}\pi_{-e}\pi_\Delta v_b = a. \tag{2}$$

By (1), (2) and since $\pi_\phi = -\pi_{-\phi}$ we get:
$$\pi_g = -\pi_e. \tag{3}$$

But as $y$ is the last endpoint of $g$, its first endpoint can be only one of $y+1, y-1$. Similarly, since the first endpoint of $e$ is $y$, its last endpoint can be only one of $y+1, y-1$. Then the only possibilities for (3) to hold are $g = \overrightarrow{(y-1)y}, e = -\overrightarrow{(y-1)y}$ and $g = -\overrightarrow{y(y+1)}, e = \overrightarrow{y(y+1)}$. So it must be:
$$g = -e. \tag{4}$$

In the cyclic order of $C_{\phi f}$, let $f = f_1, f_2, \ldots, f_{k-1}, f_k = -f$ be the arcs between $f$ and $-f$. $f_2$ is the arc $e$ and $g$ is the arc $f_k$. (4) says that $f_2 = -f_{k-1}$. Repeating the same reasoning gives that moving in these arcs $t$ places forward from $f$ and $t$ places backwards from $-f$ we find opposite main arcs. In the middle places one of the following two can happen: i) there exists an arc $h$ for which $h = -h$, or ii) there exist two arcs $h_1, h_2$ so that $h_2 = -h_1$ and $\phi(h_1) = h_2$.

In i), the two endpoints of $h$ should coincide and this can never happen: even if the two endpoints form a crossing ($h$ would be a loop), the endpoints are always distinct labels (recall that each main arc has endpoints $a$ and $a+1$ modulo $2n$; here $n > 0$, so $2n > 2$, thus $a$ and $a+1$ cannot be the same modulo $2n$).

In ii), we get $\phi(h_1) = -h_1$. By the definition of $\phi(h_1)$, it must be that $\phi(h_1)$ has beginning point the label which together with the last point of $h_1$ consist the two labels of a crossing. Since $\phi(h_1)$ is $-h_1$, this means that the last point of $h_1$ should provide both labels of a crossing which of course cannot happen.

The contradiction proves our result.

(b) In $C_{\alpha f}$ the arc after $f$ is $\alpha(f)$ and in $C_{\delta f}$ the arc after $f$ is $\delta(f)$. But $\alpha(f)$ and $\delta(f)$ have distinct endpoints, thus they are distinct themselves: their last endpoints differ by $\pm 2$ modulo $2n$ and $n = \tau(A) > 0$ since by assumption $A$ has arcs; thus $2n > 2$ implying $\pm 2 \neq 0$. Hence the cyclically ordered sets $C_{\alpha f}, C_{\delta f}$ are distinct.

(c) If $e \in C_{\phi f}$ it will be $C_{\phi f} = C_{\phi e}$ and the result comes from part (a) of the Lemma.

(d) By Lemma 12, since $\phi(f_i) = f_{i+1}$, it is $(-\phi)(-f_{i+1}) = -f_i$ and the result is immediate. $\square$

**Remark 4.** *Part (b) of the last Lemma allows us to define for a based symbol A an analogue of the rooms for the diagrams: a room is any 2-member set with elements the two cycles $C_{\phi e}, C_{(-\phi)(-e)}$ for some basic arc e of A. But there will be no need to pursue this analogy any further in the sequel. Our main tool will be the oriented cycles of the symbols (Definition 20) which are the analogues of the oriented cycles (oriented boundaries of the rooms) of the diagrams (Definition 8). And clearly, the oriented boundaries of the rooms topologically determine the oriented rooms themselves.*

The Lemma that follows takes us outside the algebraic realm, and discusses the obvious relation of the algebraic and the topological cycles. It will be useful later on in the proof of Lemma 26 in §4.6.

**Lemma 15.** *For a diagram K with $\tau(K) > 0$, the $\phi$-cycles of a based symbol $_pK$ are the $\phi$-cycles of K. That is, they are the oriented boundaries of the rooms of K considered as cyclically ordered arrays of main arcs.*

Here we follow the convention of Definition 20 to consider a common notation $f = \overrightarrow{a(a+1)}$ for basic arcs of $K$ and their symbolic counterparts in $_pK$, and similarly we consider a common notation $-f = -\overrightarrow{} = \overleftarrow{a(a+1)}$ for their opposites.

*Proof.* Let $f = \overrightarrow{a(a+1)}$ be a main arc of $_pK$, and $C_{\phi f} = (f_1 = f, f_2, f_3 \ldots)_{cycl}$ a $\phi$-cycle of $f$. Then $\phi(f_i) = f_{i+1}, \forall i$.

$f$ is also a main arc of $K$. By Lemmata 3, 4 there exists a room $\Pi$ of $K$ containing $f$ in its boundary so that $\Pi_{\phi f}$. The $\phi$-cycle of $f$ in $K$ is $\partial \Pi_{\phi f}$, that is, the oriented boundary of $\Pi$ when the room is considered with its $\phi$ orientation. When we consider $\partial \Pi_{\phi f}$ as cyclically ordered arrays of main arcs $\partial \Pi_{\phi f} = (f_1 = f, g_2, \ldots, g_k)_{cycl}$ where $f_{i+1} = \phi(f_i), \forall i \in \mathbb{N}$ and $\phi(f_i)$ as described in Lemma 10. The description for the first and last point of $g_2 = \phi(f_1) = \phi(a(a+1))$ is the same as the description in Definition 20 for the first and last point of $f_2 = \phi(f_1) = \phi(a(a+1))$ when we consider $f_1$ algebraically as a based symbol. So $f_2 = g_2$.

Similarly $f_3 = g_3$ and so on until one of the two cycles $(f_1 = f, f_2, f_3 \ldots)_{cycl}, (f_1 = f, g_2, \ldots, g_k)_{cycl}$ is run out of elements. Then the other one has to run out as well since both cycles consist of distinct elements. Thus the two cycles coincide as wanted. $\square$

We now define the algebraic versions of the topological moves.

The first move $\Omega_\beta$ which we define does not have a topological analogue, but it is dictated by the fact that each diagram $K$ of order $n$, produces $2n$ based symbols; one for each basic arc chosen to hold the base point for $K$. The next move we define is the algebraic isotopy move $^{al}\Omega_{is}$ which is actually a special case of the $\Omega_\beta$ move and looks quite trivial, but which will nevertheless be very useful to us. This move stands at the core of mirroring all $\Omega_\beta$ moves by topological moves. So we for the most part we consider it as a separate move and keep it as such until the end, although from a point on we can just treat as a special case of the $\Omega_\beta$ move. The analogues of $\Omega_1, \Omega_2, \Omega_3$ moves are defined afterwards, and as expected the first two come in two versions.

After the definitions we provide some basic properties of the new moves. The main one is that if $K \xrightarrow{\Omega} K'$ is a usual topological move of diagrams then the corresponding algebraic move $^{al}\Omega$ is well defined and for some choices of based points $p \in K$, $q \in K'$ it is $_pK \xrightarrow{^{al}\Omega} {_qK'}$; thus topological moves are mirrored by algebraic ones.

3.1. $^{al}\Omega_\beta$ **moves.**

3.1.1. *Relabelings.* For $n \in \mathbb{N}, n > 0$ and $a \in \mathbb{N}_{2n} = \{1, 2, \ldots, 2n\}$, let $\mu_{2n,a} : \mathbb{N}_{2n} \to \mathbb{N}_{2n}$ be the 1-1 function of Figure 35 which we call as *relabeling function* for the choices $n$ and $a$. We call $a$ as a *new choice of base*.

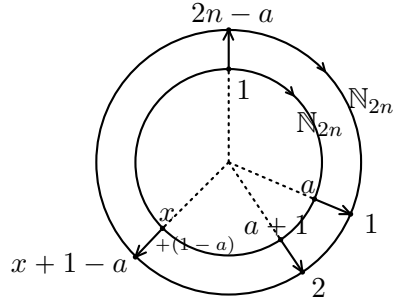

FIGURE 35. The relabeling function $\mu_{2n,a}$ that takes the "base" label 1 to some new label $a \in \mathbb{N}_{2n}$.

For a based symbol $A = (i_1, i_2)^{\pi_1}(i_3, i_4)^{\pi_2} \ldots (i_{2n-1}, i_{2n})^{\pi_n} \in \Sigma_D$ with $\tau(A) = n > 0$ and a relabeling function $\mu = \mu_{2n,a}$ we write $\mu(A) = (\mu(i_1), \mu(i_2))^{\pi_1}(\mu(i_3), \mu(i_4))^{\pi_2} \ldots (\mu(i_{2n-1}), \mu(i_{2n}))^{\pi_n}$.

For $n = 0$, we consider a unique relabeling function $\mu_{0,0} : \mathbb{N}_0 = \emptyset \to \emptyset = \mathbb{N}_0$ and we write $\mu_{0,0}(\emptyset) = \emptyset$.

3.1.2. *Definition of $^{al}\Omega_\beta$ moves.* For $A_1 \in \Sigma_D, A_2 \in \Sigma$ we say that $A_2$ is derived from $A_1$ by an algebraic move $^{al}\Omega_\beta$ whenever $A_2 = \mu_{2n,a}(A_1)$ for a relabeling function $\mu_{2n,a}$. We denote $A_1 \xrightarrow{^{al}\Omega_\beta} A_2$ without any mention to $\mu_{2n,a}$ which we call as the relabeling of the move. We call $^{al}\Omega_\beta$ as a *change of base move*.

Whenever one of $A_1, A_2$ is the empty set, the other is the empty set too and the move is $A_1 = \emptyset \xrightarrow{^{al}\Omega_\beta} \emptyset = A_2$ which is the only move among empty sets.

3.1.3. *Properties.* We gather in a single Lemma some important properties of $^{al}\Omega_\beta$ moves.

**Lemma 16.** *(a) $^{al}\Omega_\beta$ moves are well defined: for a based symbol $A_1 \in \Sigma_D$ and a relabeling function $\mu_{2n,a}$ it is indeed true that $A_2 = \mu_{2n,a}(A_1)$ is an element of $\Sigma$.*

*(b) If $A_1 \xrightarrow{^{al}\Omega_\beta} A_2$, then $\tau(A_1) = \tau(A_2)$.*

*(c) For any diagram $K$ and any pair of base points $p, q$ of $K$ it is $_pK \xrightarrow{^{al}\Omega_\beta} {_qK}$.*

*(d) For any diagram $K$ and any base point $p$ of $K$ it is $_pK \xrightarrow{^{al}\Omega_\beta} A_2 \Rightarrow A_2 \in \Sigma_D$. Namely, $A_2 = {_qK}$ for some base point $q$ of $K$.*

*(e) For a diagram $K$ with $\tau(K) = 0$ (that is: $K \approx S^1$), calling $-K$ a diagram with the same set of points as $K$ but with opposite orientation, any point of $K$ is a base point of both $K, -K$. Also, $_pK \xrightarrow{^{al}\Omega_\beta} {_q(-K)}$ for all pairs of base points $p, q$ of $K$ and $-K$.*

*Proof.* (a) For $A_1 = \emptyset$ the only relabeling function $\mu_{2n,(a)}$ defined on $A_1$ is $\mu_{0,0} : \mathbb{N}_0 = \emptyset \to \emptyset = \mathbb{N}_0$. Thus $A_2 = \emptyset$ and then $A_2 \in \Sigma$ as wanted.

For $A_1 = (i_1, i_2)^{\pi_1}(i_3, i_4)^{\pi_1} \ldots (i_{2n-1}, i_{2n})^{\pi_n}$, it is $A_2 = \mu_{2n,(a)}(A_1) = (\mu_{2n,(a)}(i_1), \mu_{2n,(a)}(i_2))^{\pi_1} (\mu_{2n,(a)}(i_3), \mu_{2n,(a)}(i_4))^{\pi_1} \ldots (\mu_{2n,(a)}(i_{2n-1}), \mu_{2n,(a)}(i_{2n}))^{\pi_n}$. Since the set of the $i_j$'s is $\mathbb{N}_{2n}$ and since $\mu_{2n,a}$ is a 1-1 function $\mathbb{N}_{2n} \to \mathbb{N}_{2n}$, the set of the $\mu_{2n,a}(i_j)$'s is again $\mathbb{N}_{2n}$. So by the definition of the numbered symbols we have $A_2 \in \Sigma$ as wanted.

(b) Immediate.

(c) The names of the crossings in $_pK$ are related to those in $_qK$ as in table 2.1 in §2.3.1. This relation is actually a relabeling function $\mu_{2n,a}$ as in our definition above. Moreover, table 2.1 shows that according to our definition above, it is $_qK = \mu_{2n,a}(_pK)$. So according to our definition for the algebraic $^{al}\Omega_\beta$ move we have $_pK \xrightarrow{^{al}\Omega_\beta} {_qK}$ as wanted.

(d) For $\tau(K) = 0$, it is $_pK = \emptyset$ and as explained in part (a) it is $A_2 = \emptyset$ as well. So $A_2 = {_qK}$ for any base point $q$ of $K$ and $A_2 \in \Sigma_D$ as wanted.

For $\tau(K) > 0$, let $\mu_{2n,a}$ be the relabeling function of the given $^{al}\Omega_\beta$ move, thus $A_2 = \mu_{2n,a}(_pK)$. We consider a base point $q$ in the basic arc $\overrightarrow{a(a+1)}$ of $K$. Then $_qK$ is the based symbol of $K$ coming from $_pK$ after the changes in the crossing names shown in table 2.1 of §2.3.1. This name changing coincides with the effect on $_pK$ of the relabeling $\mu_{2n,a}$. So $_qK = \mu_{2n,a}(_pK)$ which gives $A_2 = {_qK} \in \Sigma_D$ completing the proof.

(e) Any point of a diagram homeomorphic to a circle is a base point of the diagram. Since $K, -K$ share their points this implies that the first part of the result is true. Now recall that from part (d) we know $_pK = \emptyset$ and $_q(-K) = \emptyset$. So the move $\emptyset \xrightarrow{^{al}\Omega_\beta} \emptyset$ can be written $_pK \xrightarrow{^{al}\Omega_\beta} {_q(-K)}$ as wanted. □

**3.2. $^{al}\Omega_{iso}$ moves, an important special case of the $^{al}\Omega_\beta$ moves.** Let $A \in \Sigma_D$ be a based symbol. We say that $A$ is derived from $A$ by an algebraic isotopy move $^{al}\Omega_{iso}$ and we write $A \xrightarrow{^{al}\Omega_{iso}} A$.

In particular, for $A = \emptyset$ it is $\emptyset \xrightarrow{^{al}\Omega_{iso}} \emptyset$.

For a based symbol $A$ of positive order $\tau(A) = 2n > 0$, the $^{al}\Omega_{iso}$ move has the same effect as an $^{al}\Omega_\beta$ move with relabeling function $\mu_{2n,1}$ which also sends $A$ to $A$. Similarly, for $A = \emptyset$ and the isotopy move on $A$ has the same effect as the unique relabeling move $A = \emptyset \xrightarrow{^{al}\Omega_\beta} \emptyset$ on $A$.

So in effect, as already mentioned the algebraic isotopy move is not a genuinely new move but rather a special case of the relabeling move.

**Lemma 17.** (a) If $A_1 \in \Sigma_D$ and $A_1 \xrightarrow{^{al}\Omega_{iso}} A_2$ then $A_2 \in \Sigma_D$.

(b) If for the diagrams $K, K'$ it is $K \xrightarrow{\Omega_{iso}} K'$ then $_pK \xrightarrow{^{al}\Omega_{iso}} {_qK'}$ for suitable choices of base points $p$ and $q$; thus $_pK = {_qK'}$.

*Proof.* (a) By definition $A_1 \xrightarrow{^{al}\Omega_{iso}} A_2$ implies $A_2 = A_1$. But $A_1 \in \Sigma_D$ by hypothesis and we are done.

(b) Let $K \xrightarrow{^{al}\Omega_{iso}} K'$. We choose arbitrary some base point $p$ of $K$.

Let $F$ be the isotopy of the move $\Omega_{iso}$ on the plane $\rho$, and let $f$ be the last moment of $F$ which we think conveniently as $f : \rho \to \rho$ instead of $f : \rho \times 1 \to \rho \times 1$. It is $K' = f(K)$.

As mentioned in §1.3, the orientation of $K'$ is the one of $K$ induced via $f$. And the crossings of $K'$ are the images $\Delta' = f(\Delta)$ of the crossings $\Delta$ of $K$. So $q = f(p)$ is not a crossing of $K'$ and $_qK'$ has a meaning as a based symbol.

For $\tau(K) = 0$, we have $K \approx S^1$, and since $F_1$ is a homeomorphism we have $K' = f(K) \approx S^1$ as well. So $\emptyset = {_pK} = {_q K'}$ and then $_pK \xrightarrow{^{al}\Omega_{iso}} {_qK'}$ as wanted.

For $\tau(K) = n > 0$: Since $f$ is a homeomorphism, if we make a journey on $K$ starting at $q$ along its orientation, we meet the crossings $\Delta'_1 = f(\Delta_1), \ldots \Delta'_n = f(\Delta_n)$ of $K'$ in the same order as we meet the crossings $\Delta_1, \ldots, \Delta_n$ of $K$ in a journey on it starting at $p$ along its orientation. So the labels of any $\Delta'_j$ in the journey on $K'$ are the same as the labels of $\Delta_j$ in the journey on $K$. Thus if $\Delta = a : b$ is a crossing of $K$ then $f(\Delta) = \Delta' = a : b$ is a crossing of $K'$, and vice versa.

Regarding the signs, let $e_1, e_2$ be the over- and under-arc respectively of $K$ in a canonical neighborhood $U$ of $\Delta$. Then $U' = f(U)$ is a canonical neighborhood of $\Delta'$ and $e'_1 = f(e_1), e'_2 = f(e_2)$ are the two arcs of $K'$ in $U'$. According to the definition of the isotopy move in §1.3, $e'_1$ is the over-arc and $e'_2$ is the under-arc. Let $i$ be the label of $\Delta$ when we meet it moving along $e_1$ and let $j$ be the label of $\Delta$ when we meet it moving along $e_2$. Then $i$ will be the label of $\Delta'$ when we meet it moving along $e'_1$ and $j$ the label of $\Delta'$ when we meet it moving along $e'_2$. The sign $\pi$ of $\Delta$ is $-1$ or $+1$ according to if moving along $e_1$ we need to turn at $\Delta$ right or left to find ourselves on $e_2$ along its orientation and the crossing in $_pK$ is then $\Delta = (i:j)^{-1}, (i:j)^{+1}$ respectively. Similarly as we walk along $e'_1 = f(e_1)$ and turn at $\Delta' = f(\Delta)$ in order to find ourselves on $e'_2 = f(e_2)$ along its orientation, we need to make a right or left turn to get accordingly $\Delta' = (i:j)^{-1}$ or $\Delta' = (i:j)^{+1}$ in $_qK'$. Since $F$ is an isotopy starting with the identity of the plane, it is an orientation preserving isotopy, hence its last moment $f$ is an orientation preserving homeomorphism of the plane. This makes the turn of $e'_1$ to be right or left exactly whenever the turn of $e_1$ is right or left respectively. So we always have $\Delta = \Delta'$ as elements of $\Sigma_\Delta$. Thus $_pK \xrightarrow{al\Omega_{iso}} _qK'$ as wanted. $\square$

### 3.3. $^{al}\Omega_1$ moves.

3.3.1. *Relabelings.* For the algebraic $\Omega_1$ moves which we are about to define, we are going to use the relabeling functions of the topological $\Omega_1$ moves introduced in §2.3.3. We defined them in Figures 15, 17, and merged in Figure 18.

Recall that the relabelings are functions of one of the two following kinds:

(i) $\mu^+_{2n,a,(x)} : \mathbb{N}_{2n} \to \mathbb{N}_{2(n+1),(a+1)+(x+1)}$ where $n > 0$, $a \in \mathbb{N}_{2n}$ and $x = -1$ if $a < 2n$ whereas $x \in \{-1, 0, 1\}$ if $a = 2n$. For $n = 0$ all relabelings become equal to $\mu : \mathbb{N}_0 = \emptyset \to \emptyset = \mathbb{N}_{2,2+x}$, and we denote this relabeling as $\mu^+_{0,0}$.

(ii) $\mu^-_{2(n+1),z} : \mathbb{N}_{2(n+1),z} \to \mathbb{N}_{2n}$ where $n \geqslant 0$, $z \in \mathbb{N}_{2(n+1)}$. They are the inverse functions of the relabelings of the first kind:
$(\mu^+_{2n,a,(x)})^{-1} = \mu^-_{2(n+1),(a+1)+(x+1)} : \mathbb{N}_{2(n+1),(a+1)+(x+1)} \to \mathbb{N}_{2n}$, where $a \in \mathbb{N}_{2n}$, $x = -1$ if $a < 2n$ and $x \in \{-1, 0, 1\}$ if $a = 2n$.

For $n > 0$, any based symbol $A = (i_1, i_2)^{\pi_1}(i_3, i_4)^{\pi_2} \ldots (i_{2n-1}, i_{2n})^{\pi_n} \in \Sigma$ with $\tau(A) = n$ and any relabeling function $\mu = \mu^+_{2n,a,x}$, we write $\mu(A) = (\mu(i_1), \mu(i_2))^{\pi_1} (\mu(i_3), \mu(i_4))^{\pi_2} \ldots (\mu(i_{2n-1}), \mu(i_{2n}))^{\pi_n}$. Accordingly, for $n = 0$, for the based symbol $\emptyset$ and the relabeling function $\mu = \mu^+_{0,0}$, we write $\mu(\emptyset) = \emptyset$.

Let $A = (i_1, i_2)^{\pi_1}(i_3, i_4)^{\pi_2} \ldots (i_{2n+1}, i_{2n+2})^{\pi_{n+1}} \in \Sigma$ be a symbol with $\tau(A) = n + 1$ which contains the crossing of a 1-gon, say the crossing $(z : z + 1)$, $z \in \mathbb{N}_{2(n+1)}$ where with no loss of generality it is $\{z, z+1\} = \{i_{2n+1}, i_{2n+2}\}$. Then for the relabeling function $\mu = \mu^-_{2(n+1),z}$ we write $\mu(A) = (\mu(i_1), \mu(i_2))^{\pi_1} (\mu(i_3), \mu(i_4))^{\pi_2} \ldots (\mu(i_{2n-1}), \mu(i_{2n}))^{\pi_n}$.

Especially for $n = 0$, $A$ can only be $A = (1:2)^\pi \in \mathbb{N}_2$ ($\pi \in \{-1, 1\}$) and as $\mathbb{N}_{2,1} = \mathbb{N}_{2,2} = \mathbb{N}_0 = \emptyset$, the relabeling functions can only be $\mu = \mu^-_{2,1} = \mu^-_{2,2} : \emptyset \to \emptyset$; we write $\mu(A) = \emptyset$.

If $\mu(A) = A'$ is defined for some $A \in \Sigma_D$ and some relabeling $\mu$ as those in the current paragraph, then $A' \in \Sigma$ (Lemma 18 below).

3.3.2. *Definition of $^{al}\Omega_1^+$ moves.* Let $A_1 \in \Sigma_D$ with $\tau(A_1) = n > 0$. Let $e = \overrightarrow{a(a+1)}$ be a basic arc of $A_1$, where $a \in \mathbb{N}_{2n}$.

We say that the numbered symbol $A_2 \in \Sigma$ is derived from $A_1$ by an algebraic move $^{al}\Omega_1^+$ on $e$ whenever it is:

$$A_2 = \mu(A_1)\Delta$$
(3.1) $$\mu = \mu^+_{2n,a,(x)} \quad \text{(for some } x \in \{-1, 0, 1\})$$
$$\Delta = \left(\tfrac{1-\theta}{2}(a+1) + \tfrac{1+\theta}{2}(a+2), \tfrac{1+\theta}{2}(a+1) + \tfrac{1-\theta}{2}(a+2)\right)^{\pi_\phi},$$

where $\theta = \theta_e$ is a placement number of $e$, $\phi$ is some orientation and $\pi_\phi$ is the sign of $\phi$.

We call the move as an *under-* or *over-placement* (or twist) of $e$ depending on if $\theta_e = -1$ or $\theta_e = +1$ respectively. Suppressing most information, we denote the move simply as $A_1 \xrightarrow{^{al}\Omega_1^+} A_2$. We say this is an $^{al}\Omega_1^+$ move on $A_1$ or on the basic arc $e$ of $A_1$. We call $\mu$ as the relabeling of the move.

An $\Omega_1^+$ move can happen on any basic arc $e$ of $A_1$.

For $A_1 \in \Sigma_D$ with $\tau(A_1) = 0$ (that is for $A_1 = \varnothing$), we say that the numbered symbol $A_2 \in \Sigma$ is derived from $A_1$ by an algebraic move $^{al}\Omega_1^+$ whenever:

(3.2)
$$A_2 = \Delta$$
$$\Delta = \left(\tfrac{1-\theta}{2}1 + \tfrac{1+\theta}{2}2, \tfrac{1+\theta}{2}1 + \tfrac{1-\theta}{2}2\right)^{\pi\theta},$$

where $\theta_\varnothing$ is a placement number of $\varnothing$, and $\pi = \pi_\phi$ is the sign of some orientation $\phi$.

We call the move as *an under- or over-twist* of $\varnothing$ depending on if $\theta_\varnothing = -1$ or $\theta_\varnothing = +1$, and suppressing information we denote it as $A_1 \xrightarrow{^{al}\Omega_1^+} A_2$.

We also say that $A_1 \xrightarrow{^{al}\Omega_1^+} A_2$ is an $\Omega_1^+$ move on $A_1$.

3.3.3. *Definition of $^{al}\Omega_1^-$ moves.* Let $A_1 \in \Sigma_D$ with $\tau(A_1) = n+1 \geqslant 1$. Let $e = \overrightarrow{z(z+1)}$, $z \in \mathbb{N}_{2(n+1)}$ be a loop (1-gon) of $A_1$ whose crossing is $\Delta$.

We say that the numbered symbol $A_2 \in \Sigma$ is derived from $A_1$ by an algebraic move $^{al}\Omega_1^-$ on $e$ whenever $A_2$ is what we get deleting from $A_1$ the crossing $\Delta$ and relabeling the other crossings appropriately. More rigorously, whenever it is:

(3.3)
$$A_1 = A_1'\Delta = A_1'(z : z+1)$$
$$\mu = \mu_{2(n+1),z}^-$$
$$A_2 = \mu(A_1).$$

Then it is:

(3.4)
$$\mu = \mu_{2(n+1),z}^- = (\mu_{2n,a,(x)}^+)^{-1} \text{ for appropriate } a \in \mathbb{N}_{2n}, x \in \{-1, 0, 1\}$$
$$\Delta = \left(\tfrac{1-\theta}{2}(a+1) + \tfrac{1+\theta}{2}(a+2), \tfrac{1+\theta}{2}(a+1) + \tfrac{1-\theta}{2}(a+2)\right)^{\pi}$$
$$\theta = \theta_e \text{ a placement number of } e, \phi \text{ some orientation}, \pi = \pi_\phi \text{ the sign of } \phi.$$

We call $\mu$ as the relabeling of the move. Suppressing most information, we denote the move simply as $A_1 \xrightarrow{^{al}\Omega_1^-} A_2$. We say this is an $\Omega_1^-$ move on $A_1$ or on the loop $e$ of $A_1$.

3.3.4. *Properties.* In the Lemma below we have gathered a bunch of important properties of $^{al}\Omega_1$ moves.

**Lemma 18.** (a) *Let $A_1 \in \Sigma_D$ and $\mu$ a relabeling $\mu_{2n,a,(x)}^+$ or an inverse relabeling $\mu_{2(n+1),z}^-$. If $\mu(A_1) = A_2$ is defined then $A_2 \in \Sigma$.*

*So the $^{al}\Omega_1^+$ and $^{al}\Omega_1^-$ moves are well defined. Also, if $A_1 \xrightarrow{\Omega} A_2$ is one of them, then $A_1 \in \Sigma_D$, $A_2 \in \Sigma$.*

(b) *If $A_1 \xrightarrow{^{al}\Omega_1^+} A_2$, then $\tau(A_2) = \tau(A_1) + 1$, and if $A_1 \xrightarrow{\Omega_1^-} A_2$, then $\tau(A_2) = \tau(A_1) - 1$.*

(c) *If $K, K'$ are diagrams, then:*

- $K \xrightarrow{\Omega_1^+} K' \Rightarrow {}_pK \xrightarrow{^{al}\Omega_1^+} {}_qK'$ *for suitable choices of base points $p \in K, q \in K'$.*
- $K' \xrightarrow{\Omega_1^-} K \Rightarrow {}_qK' \xrightarrow{^{al}\Omega_1^-} {}_pK$ *for suitable choices of base points $p \in K, q \in K'$.*

*Proof.* (a) We only have to follow the definitions. In detail:

Consider first a relabeling $\mu = \mu^+_{2n,a,(x)}$

Since $\mu(A_1)$ is defined, it is $\tau(A_1) = n \geqslant 0$.

For $n = 0$: $A_1 = \emptyset$. The relabeling function is unique: $\mu = \mu^+_{0,0}$ and $\mu(A_1) = \emptyset$. Then $A_2 = (1:2) \in \Sigma_D$, so $A_2 \in \Sigma$ as wanted.

For $n > 0$: Say $A_1 = (i_1, i_2)^{\pi_1}(i_3, i_4)^{\pi_2} \ldots (i_{2n-1}, i_{2n})^{\pi_n}$. Since $\mu = \mu^+_{2n,a,(x)}$ it is $A_2 = \mu(A_1)(a + 1 : a + 2) = (\mu(i_1), \mu(i_2))^{\pi_1}(\mu(i_3), \mu(i_4))^{\pi_2} \ldots (\mu(i_{2n-1}), \mu(i_{2n}))^{\pi_n}(a + 1 : a + 2)$. By definition $\mu : \mathbb{N}_{2n} \to \mathbb{N}_{2(n+1),a} = \mathbb{N}_{2(n+1)} - \{a+1, a+2\}$ is 1-1 and onto. Hence the integers in the pairs of $A_2$ form the set $\mathbb{N}_{2(n+1)}$. Thus by the definition of the numbered symbols we have $A_2 \in \Sigma$ as wanted.

Now consider the relabeling $\mu = \mu^-_{2(n+1),z} = (\mu^+_{2n,a,(x)})^{-1}$ for the appropriate $a, x$.

Since $\mu(A_1)$ is defined, it is $\tau(A_1) = n + 1, n \geqslant 0, a \in \mathbb{N}_{2n}$ and $\overrightarrow{z(z+1)}$ must be a 1-gon in $A_1$. Say $A_1 = (i_1, i_2)^{\pi_1}(i_3, i_4)^{\pi_1} \ldots (i_{2n-1}, i_{2n})^{\pi_n}(z : z + 1)$.

Then $A_2 = \mu(A_1) = (\mu(i_1), \mu(i_2))^{\pi_1}(\mu(i_3), \mu(i_4))^{\pi_1} \ldots (\mu(i_{2n-1}), \mu(i_{2n}))^{\pi_n}$. By definition $\mu : \mathbb{N}_{2(n+1),a} \to \mathbb{N}_{2n}$ is 1-1 and onto. Hence the integers in the pairs of $A_2$ form the set $\mathbb{N}_{2n}$. Thus by the definition of the numbered symbols we have $A_2 \in \Sigma$ as wanted.

The claims for the moves are now immediate.

(b) Immediate.

(c) We only have to follow the definitions.

First we deal with the $K \xrightarrow{\Omega^+_1} K'$ case:

Let $\tau(K) = n \geqslant 0$.

For $n > 0$: Figure 12 presents the topological setting near the basic arc $e$ of the move as well as near the result $e'$ of the move along with choices of base points $p \in K$ and $q \in K'$. All cases (A)-(D) present the same topological setting differing only in the choices of $p, q$. So we choose any one of them to work with.

It is always $e = \overrightarrow{a(a+1)}$ for the integer $a \in \mathbb{N}_{2n}$ shown in the Figure, and the new crossing $\Delta$ formed is described in all detail in Relations (2.6), (2.7), (2.8), (2.9) for $\phi, \theta$ depending on the choices made in the move. Also, the labels of the crossings of $K$ and their new labels when considered as crossings of $K'$ are given in Figure 18 for $\mu^+_{2n,a,(x)}$ for the corresponding $x$ depending on the case we are working on.

Then the based symbols ${}_pK, {}_qK'$ are related exactly as in (3.1) for the same relabeling $\mu^+_{2n,a,(x)}$ in Figure 18 for the relabeling of based symbols. This is exactly what we ask for an algebraic move ${}_pK \xrightarrow{al\,\Omega^+_1} A$ on the 1-gon $e = \overrightarrow{a(a+1)}$ of the based symbol ${}_pK$ with the same $\phi, \theta$ as above; and then the result of this move, is a symbol $A$ with the same description as ${}_qK'$ for the diagram $K'$ in the case we are working and the choice of $q$ we made above. So $A = {}_qK'$. And $K \xrightarrow{\Omega^+_1} K' \Rightarrow {}_pK \xrightarrow{al\,\Omega^+_1} {}_pK'$ as wanted.

For $n = 0$, it is $K \approx S^1$ and ${}_pK = \emptyset$ for any choice of the base point $p$. The topological move $K \xrightarrow{\Omega^+_1} K'$ produces a crossing $\Delta_0$ for $K'$. No matter what the orientation of $K$, the possible results of the move produce all possible ordered pairs with all possible signs as crossing of $K'$ (Figure 11). The same happens with the crossings $\Delta$ in relation (3.2) whenever we consider all choices of $\phi, \theta = \theta_\emptyset$. So for any choice of the base point $q \in K'$ as in Figure 11, it is ${}_qK' = \Delta_0 = \Delta$ for some choices of $\phi, \theta_\emptyset$. Hence $K \xrightarrow{\Omega^+_1} K' \Rightarrow {}_pK \xrightarrow{al\,\Omega^+_1} {}_pK'$ in his case too, as wanted.

Next we deal with the $K' \xrightarrow{\Omega^-_1} K$ case:

Let $\tau(K') = n + 1,\ n \geqslant 0$.

The topological setting for the move is given in Figure 12 in cases (A)-(D) which differ only in the choices of the base points $q \in K'$ and $p \in K$. We choose any one we consider convenient to work with.

The full description of the deleted crossing $\Delta$ is given in relation (2.11) where $a \in \mathbb{N}_{2n}$, $x \in \{-1, 0, 1\}$, and $\phi, \theta$ depend on the choices made for the move. Also, the labels of the other crossings of $K'$ and their new labels when considered as crossings of $K$ are given in Figure 18 for $\mu^-_{2(n+1),z} = (\mu^+_{2n,a,(x)})^{-1}$, for the corresponding $x$ depending on the case we are working on.

Then the based symbols ${}_qK', {}_pK$ are related exactly as in Relations (??), (3.3) for the same relabeling $(\mu^+_{2n,a,(x)})^{-1}$, the same $\phi$ and the same $\theta$.

Let us also notice that the basic arc of the move, say $e'$, is a 1-gon of $K'$. Then $e'$ as a basic arc of ${}_qK'$ is an 1-gon too and we can perform on it an ${}^{al}\Omega^-_1$ move.

By (??) this move ${}_qK' \xrightarrow{{}^{al}\Omega^-_1} A$ produces a symbol $A$ for which Relation (??) holds, thus with the same description as ${}_qK'$ for the diagram $K'$ in the case we are working and the choice of $q$ we made above. So $A = {}_qK'$. And then $K' \xrightarrow{\Omega^-_1} K \Rightarrow {}_qK' \xrightarrow{{}^{al}\Omega^-_1} {}_pK$ as wanted. $\square$

3.4. ${}^{al}\Omega_2$ **moves.**

3.4.1. *Relabelings.* Recall that the relabelings are functions of one of the two following kinds:

(i) $\mu^+_{2n,a,b,(x),(y)} : \mathbb{N}_{2n} \to \mathbb{N}_{2(n+2),(a+1)+(y+1)(x+1),(b+3)+(y+1)(x+1)}$ where $n > 0$, $a, b \in \mathbb{N}_{2n}$, $a \leq b$ and $x = -1$ if $b < 2n$ whereas $x \in \{-1, 0, 1\}$ if $b = 2n$ and also $y = +1$ if $a = b$ whereas $y = 0$ if $a \neq b$. For $n = 0$ all relabelings become equal to $\mu : \mathbb{N}_0 = \emptyset \to \emptyset$ and we denote this relabeling as $\mu^+_{0,0,0}$.

(ii) $\mu^-_{2(n+2),z,w} : \mathbb{N}_{2(n+2),z,w} \to \mathbb{N}_{2n}$ where $n \geq 0$, $z, w \in \mathbb{N}_{2(n+2)}$ with $|w - z| \neq 1$. They are the inverse functions of the relabelings of the first kind:
$(\mu^+_{2n,a,b,(x),(y)})^{-1} = \mu^-_{2(n+2),(a+1)+(y+1)(x+1),b+3+(y+1)(x+1)} : \mathbb{N}_{2(n+2),(a+1)+(y+1)(x+1),b+3+(y+1)(x+1)}$
$\to \mathbb{N}_{2n}$, where $a, b \in \mathbb{N}_{2n}$, $a \leq b$. Here also $x = -1$ if $b < 2n$, whereas $x \in \{-1, 0, 1\}$ if $b = 2n$. And $y = +1$ if $a = b$, whereas $y = 0$ if $a \neq b$. For $n = 0$ the two relabeling functions coincide: $\mu^-_{0,1,3} = \mu^-_{0,2,4} : \emptyset \to \emptyset$.

Let $A = (i_1, i_2)^{\pi_1}(i_3, i_4)^{\pi_2} \ldots (i_{2n+1}, i_{2n+2})^{\pi_{n+1}}(i_{2n+3}, i_{2n+4})^{\pi_{n+2}} \in \Sigma_D$ be a based symbol with $\tau(A) = n + 2$ which contains the crossings of a 2-gon, say $(z : w), (z+1 : w+1)$ or $(z : w+1), (z+1 : w)$ (where $z, w \in \mathbb{N}_{2(n+1)}$, $z \neq w$) with heights of $z, z+1$ both $-1$ or both $+1$ in their corresponding pairs. Say $(i_{2n+1}, i_{2n+2})$ and $(i_{2n+3}, i_{2n+4})$ be the two pairs. Then for the relabeling function $\mu = \mu^-_{2(n+2),z+1,w+3}$ we write $\mu(A) = (\mu(i_1), \mu(i_2))^{\pi_1}(\mu(i_3), \mu(i_4))^{\pi_2} \ldots (\mu(i_{2n-1}), \mu(i_{2n}))^{\pi_n}$. So for example for $n = 0$ the possible based symbols are $A_1 = (1 : 2)^{\pi_1}(3 : 4)^{\pi_2}$ and $A_2 = (1 : 3)^{\pi_1}(2 : 4)^{\pi_2}$, and for $\mu = \mu^-_{0,1,3} = \mu^-_{0,2,4}$ it is $\mu(A_i) = \emptyset$, $i = 1, 2$.

If $\mu(A) = A'$ is defined for some $A \in \Sigma_D$ and some relabeling or inverse relabeling $\mu$ as those in the current paragraph, then $A' \in \Sigma$ (Lemma 19 below).

3.4.2. *Definition of ${}^{al}\Omega^+_2$ moves.* Let $A_1 \in \Sigma_D$ with $\tau(A_1) = n > 0$ and let $e_1 = \overrightarrow{a(a+1)}, e_2 = \overrightarrow{b(b+1)}$ for some $a, b \in \mathbb{N}_{2n}$, $a \leq b$ be two basic arcs of $A_1$ whose cycle senses can be compared, that is for which $e_1 \upuparrows e_2$ or $e_1 \updownarrows e_2$. We allow all possible choices of $e_1, e_2$, so we allow them to coincide as well in which case it is $a = b$. The asymmetry on $a, b$ dictates that we consider $e_1$ as the first arc and $e_2$ as the second among the two.

We say that the numbered symbol $A_2 \in \Sigma$ is derived from $A_1$ by an algebraic move $^{al}\Omega_2^+$ on $e_1, e_2$ whenever it is:

$$A_2 = \mu(A_1)\Delta_1\Delta_2$$
$$\mu = \mu^+_{2n,a,b,(x),(y)} \quad \text{(for some } x, y \in \{-1, 0, 1\}\text{)}$$

(3.5)
$$\Delta_1 = \left(\tfrac{1-\theta}{2}(a+1) + \tfrac{1+\theta}{2}\left(b + \tfrac{7+\rho}{2}\right), \tfrac{1+\theta}{2}(a+1) + \tfrac{1-\theta}{2}\left(b + \tfrac{7+\rho}{2}\right)\right)^{\pi\theta\rho}$$
$$\Delta_2 = \left(\tfrac{1-\theta}{2}(a+2) + \tfrac{1+\theta}{2}\left(b + \tfrac{7-\rho}{2}\right), \tfrac{1+\theta}{2}(a+2) + \tfrac{1-\theta}{2}\left(b + \tfrac{7-\rho}{2}\right)\right)^{-\pi\theta\rho},$$

where $\rho = \rho_{e_1,e_2}$ is the number of similarity of cycle sense of $e_1$ and $e_2$, $\theta_{e_1,e_2}$ is a placement number of the ordered pair $(e_1, e_2)$, $\phi$ is an orientation and $\pi_\phi$ is the sign of $\phi$.

Suppressing information we denote the move as $A_1 \xrightarrow{^{al}\Omega_2^+} A_2$. We say this is an $^{al}\Omega_2^+$ move on $A_1$ or on the basic arcs $e_1, e_2$ of $A_1$. We call $\mu$ as the relabeling of the move.

The only usefulness of the $a \leq b$ inequality is to make sense of the relabeling $\mu^+_{2n,a,b,(x),(y)}$ which was defined for $a, b$ in $\mathbb{N}_{2n}$ satisfying this inequality.

In the special case $A_1 \in \Sigma_D$ with $\tau(A_1) = 0$ (that is for $A_1 = \emptyset$), the only relabeling is $\mu_{0,0,0}$ and then the numbered symbol $A_2 \in \Sigma$ is derived from $A_1$ by an algebraic move $^{al}\Omega_2^+$ whenever:

$$A_2 = \Delta_1\Delta_2$$

(3.6)
$$\Delta_1 = \left(\tfrac{1-\theta}{2}(0+1) + \tfrac{1+\theta}{2}\left(0 + \tfrac{7+1}{2}\right), \tfrac{1+\theta}{2}(0+1) + \tfrac{1-\theta}{2}\left(0 + \tfrac{7+1}{2}\right)\right)^{\pi\theta\cdot 1}$$
$$\Delta_2 = \left(\tfrac{1-\theta}{2}(0+2) + \tfrac{1+\theta}{2}\left(0 + \tfrac{7-1}{2}\right), \tfrac{1+\theta}{2}(0+2) + \tfrac{1-\theta}{2}\left(0 + \tfrac{7-1}{2}\right)\right)^{-\pi\theta\cdot 1},$$

where $\theta = \theta_{\emptyset,\emptyset}$ is a placement number of the empty based symbol $A_1 = \emptyset$, $\phi$ is an orientation and $\pi_\phi$ is the sign of $\phi$.

We denote the move as $A_1 = \emptyset \xrightarrow{^{al}\Omega_2^+} A_2$.

3.4.3. *Definition of $^{al}\Omega_2^-$ moves.* Let $A_1 \in \Sigma_D$ with $\tau(A_1) = n + 2 \geq 2$. Let $e_1 = \overrightarrow{z(z+1)}$, $e_2 = \overrightarrow{w(w+1)}$ ($z, w \in \mathbb{N}_{n+2}$, $z \neq w$) be two basic arcs of $A_1$ which form a 2-gon of $A_1$ with crossings $\Delta_1, \Delta_2$.

We say that the numbered symbol $A_2 \in \Sigma$ is derived from $A_1$ by an algebraic move $^{al}\Omega_2^-$ on $e_1, e_2$ whenever $A_2$ is what we get deleting from $A_1$ the crossings $\Delta_1, \Delta_2$ of the 2-gon and relabeling the other crossings appropriately. More rigorously, whenever it is:

$$A_1 = A_1'\Delta_1\Delta_2 = A'(z:w)^{\pi_0}(z+1:w+1)^{-\pi_0} \text{ or } A'(z:w+1)^{\pi_0}(z+1:w)^{-\pi_0}$$

(3.7)
$$\mu = \mu^-_{2(n+1),z,w}$$
$$A_2 = \mu(A_1)$$

Then it is:

$$\mu = \mu^-_{2(n+1),z,w} = (\mu^+_{2n,a,b,(x),(y)})^{-1} \quad \text{(for appropriate } a, b \in \mathbb{N}_{2n}, x, y \in \{-1, 0, +1\}\text{)}$$

(3.8)
$$\Delta_1 = \left(\tfrac{1-\theta}{2}(a+1) + \tfrac{1+\theta}{2}\left(b + \tfrac{7+\rho}{2}\right), \tfrac{1+\theta}{2}(a+1) + \tfrac{1-\theta}{2}\left(b + \tfrac{7+\rho}{2}\right)\right)^{\pi\theta\rho}$$
$$\Delta_2 = \left(\tfrac{1-\theta}{2}(a+2) + \tfrac{1+\theta}{2}\left(b + \tfrac{7-\rho}{2}\right), \tfrac{1+\theta}{2}(a+2) + \tfrac{1-\theta}{2}\left(b + \tfrac{7-\rho}{2}\right)\right)^{-\pi\theta\rho}$$
$$\rho = \rho_{e_1,e_2} \text{ is the number of similarity of cycle sense of } e_1, e_2$$
$$\phi \text{ an orientation}, \pi = \pi_\phi \text{ the sign of } \phi.$$

Recall that $\rho_{e_1,e_2} \in \{1,1\}$ is uniquely determined for the basic arcs $e_1, e_2$.

Suppressing most information, we denote the move simply as $A_1 \xrightarrow{^{al}\Omega_2^-} A_2$. We say this is an $\Omega_2^-$ move on $A_1$ or on the 2-gon $e_1, e_2$ of $A_1$ and we call $\mu$ as its relabeling.

### 3.4.4. Properties.
We gather in a single Lemma some important properties of $^{al}\Omega_2$ moves.

**Lemma 19.** (a) Let $A_1 \in \Sigma_D$ and $\mu$ a relabeling $\mu_{2n,a,b,(x),(y)}^+$ or $\mu_{2(n+2),z,w}^-$. If $\mu(A_1) = A_2$ is defined, then $A_2 \in \Sigma$.

So the $\Omega_2^+$ and $\Omega_2^-$ moves are well defined. Also, if $A_1 \xrightarrow{\Omega} A_2$ is one of them, then $A_1 \in \Sigma_D$, $A_2 \in \Sigma$.

(b) If $A_1 \xrightarrow{^{al}\Omega_2^+} A_2$, then $\tau(A_2) = \tau(A_1) + 2$, and if $A_1 \xrightarrow{\Omega_2^-} A_2$, then $\tau(A_2) = \tau(A_1) - 2$.

(c) If $K, K'$ are diagrams, then:

- $K \xrightarrow{\Omega_2^+} K' \Rightarrow {}_pK \xrightarrow{^{al}\Omega_2^+} {}_pK'$ for suitable choices of base points $p \in K, q \in K'$.
- $K' \xrightarrow{\Omega_2^-} K \Rightarrow {}_qK' \xrightarrow{^{al}\Omega_2^-} {}_pK$ for suitable choices of base points $p \in K, q \in K'$.

*Proof.* (a) We follow the definitions:

Consider first a relabeling $\mu = \mu_{2n,a,b,(x),(y)}^+$.

Since $\mu(A_1)$ is defined, it is $\tau(A_1) = n \geq 0$.

For $n = 0$: $A_1 = \emptyset$. The relabeling function is unique: $\mu = \mu_{0,0,0}^+$ and $\mu(A_1) = \emptyset$. Then $A_2 = (1:4)(2:3) \in \Sigma_D$, so $A_2 \in \Sigma$ as wanted.

For $n > 0$: Say $A_1 = (i_1, i_2)^{\pi_1}(i_3, i_4)^{\pi_2} \ldots (i_{2n-1}, i_{2n})^{\pi_n}$. Since $\mu = \mu_{2n,a,b,(x),(y)}^+$ it is $A_2 = \mu(A_1)(a+1:b+3)(a+2:b+4) = (\mu(i_1), \mu(i_2))^{\pi_1}(\mu(i_3), \mu(i_4))^{\pi_2} \ldots (\mu(i_{2n-1}), \mu(i_{2n}))^{\pi_n}(a+1:b+3)(a+2:b+4)$ (or $(a+1:b+4)(a+2:b+3)$ in place of $(a+1:b+3)(a+2:b+4)$). By definition $\mu: \mathbb{N}_{2n} \to \mathbb{N}_{2(n+2),a,b} = \mathbb{N}_{2(n+2)} - \{a+1, a+2, b+3, b+4\}$ is 1-1 and onto. Hence the integers in the pairs of $A_2$ form the set $\mathbb{N}_{2(n+2)}$. Thus by the definition of the numbered symbols we have $A_2 \in \Sigma$ as wanted.

Now consider a relabeling $\mu = \mu_{2(n+2),a,b}^-$.

Since $\mu(A_1)$ is defined, it is $\tau(A_1) = n+2, n \geq 0, a, b \in \mathbb{N}_{2(n+1)}$ and $e_1 = \overrightarrow{z(z+1)}$, $e_2 = \overrightarrow{w(w+1)}$ ($z \neq w$) must form a 2-gon in $A_1$. Say without loss of generality $A_1 = (i_1, i_2)^{\pi_1}(i_3, i_4)^{\pi_1} \ldots (i_{2n-1}, i_{2n})^{\pi_n}(z:w)(z+1:w+1)$, and similarly then if $(z:w+1)(z+1:w)$ are the crossings in place of $(z:w)(z+1:w+1)$.

It is $A_2 = \mu(A_1) = (\mu(i_1), \mu(i_2))^{\pi_1}(\mu(i_3), \mu(i_4))^{\pi_1} \ldots (\mu(i_{2n-1}), \mu(i_{2n}))^{\pi_n}$. By definition $\mu: \mathbb{N}_{2(n+2),a,b} \to \mathbb{N}_{2n}$ is 1-1 and onto. Hence the integers in the pairs of $A_2$ form the set $\mathbb{N}_{2n}$. Thus by the definition of the numbered symbols we have $A_2 \in \Sigma$ as wanted.

The claims for the moves are now immediate.

(b) Immediate.

(c) We follow the definitions again:

First we deal with the $K \xrightarrow{\Omega_2^+} K'$ case:

Let $\tau(K) = n \geq 0$.

For $n > 0$: Figure 20 presents the topological setting near the basic arcs $e_1, e_2$ of $K$ for the move as well as the result of the move (the arcs $e_1', e_2'$) along with choices of base points $p \in K$ and $q \in K'$. Cases (A),(C),(E),(G) present the same topological setting, (B),(D),(F),(H) present another and (I)-(L) yet another. The cases in each setting differ only in the choices of $p, q$. For any one of the settings it suffices to work in one of its cases. Let's say we have chosen one.

It is always $e_1 = \overrightarrow{a(a+1)}$, $e_2 = \overrightarrow{b(b+1)}$ for $a, b \in \mathbb{N}_{2n}$ shown in Figure 20 and the new crossings $\Delta_1, \Delta_2$ are described in all detail in Relations (2.14) for these $a, b \in \mathbb{N}_{2n}$ and the corresponding $\phi, \theta, \rho$ depending on the choices made in the move. Also, the labels of the crossings of $K$ and their new labels when considered as crossings of $K'$ are given in Figure 25 for $\mu^+_{2n,a,b,(x),(y)}$ and the corresponding $x, y$ depending on the case we are working on.

Then the based symbols $_pK, {}_qK'$ are related exactly as in Relation (3.5) for the same relabeling $\mu^+_{2n,a,b,(x),(y)}$. This is exactly what we ask asks for an algebraic move $_pK \xrightarrow{al\Omega^+_2} {}_qK'$ on the 2-gon $e_1 = \overrightarrow{a(a+1)}$, $e_2 = \overrightarrow{b(b+1)}$ of the based symbol $_pK$ with the same $\phi, \theta, \rho$ as above; and then the result of this move, is a symbol $A$ with the same description as $_qK'$ for the diagram $K'$ in the case we are working and the choice of $q$ we made above. So $A = {}_qK'$. And $K \xrightarrow{\Omega^+_2} K' \Rightarrow {}_pK \xrightarrow{al\Omega^+_2} {}_pK'$ as wanted.

For $n = 0$, it is $K \approx S^1$ and $_pK = \varnothing$ for any choice of the base point $p$. We choose the base point $q$ of $K'$ as in Figure 19. The topological move $K \xrightarrow{\Omega^+_2} K'$ produces two crossings $\Delta_{01}, \Delta_{02}$ for $K'$. No matter what the orientation of $K$, the possible results of the move produce all possible ordered pairs $(1:4), (2:3)$ with all possible pairs of opposite signs as the crossing of $K'$. The same happens with the crossings $\Delta_1, \Delta_2$ in Relation (3.6) whenever we consider all choices of $\phi, \theta = \theta_{\varnothing, \varnothing}$. So it is $_qK' = \Delta_{01}\Delta_{02} = \Delta_1\Delta_2$ for some choices of $\phi, \theta_{\varnothing, \varnothing}$. Hence $K \xrightarrow{\Omega^+_1} K' \Rightarrow {}_pK \xrightarrow{al\Omega^+_2} {}_pK'$ in his case too, as wanted.

Next we deal with the $K' \xrightarrow{\Omega^-_2} K$ case:

Let $\tau(K') = n+2$, $n \geq 0$.

The topological setting for the move is given in Figure 20 along with choices of base points $p \in K$ and $q \in K'$. Cases (A),(C),(E),(G) present the same topological setting, (B),(D),(F),(H) present another and (I)-(L) yet another. The cases in each setting differ only in the choices of $p, q$. For any one of the settings it suffices to work in one of its cases. Let's say we have chosen one.

The full description of the deleted crossings $\Delta_1, \Delta_2$ is given in Relation (2.18) where $a, b \in \mathbb{N}_{2n}$, $x, y \in \{-1, 0, 1\}$, $\rho = \rho_{e'_1, e'_2}$ is the boundary sense of the arcs $e'_1, e'_2$ of $K'$ on which we perform the move, and $\phi, \theta$ depend on the choices made in the move. Also, the labels of the other crossings of $K'$ and their new labels when considered as crossings of $K$ are given in Figure 25 for $\mu^-_{2(n+2),z,w} = (\mu^+_{2n,a,b,(x),(y)})^{-1}$, for the corresponding $x, y$ depending on the case we are working on.

Then the based symbols $_qK', {}_pK$ are related exactly as in Relations (3.7), (3.8) for the same relabeling $(\mu_{2n,a,b,(x),(y)})^{-1}$, the same $\rho_{e'_1, e'_2}$, the same $\phi$ and the same $\theta$.

And let us note that since the basic arcs $e'_1, e'_2$ of the move form a 2-gon of $K'$, then as basic arcs of $_qK'$ they also form a 2-gon and we can perform on them an $^{al}\Omega^-_2$ move. By (3.7) this move $_qK' \xrightarrow{al\Omega^-_2} A$ produces a symbol $A$ for which Relation (3.8) holds, thus with the same description as $_qK'$ for the diagram $K'$ in the case we are working and the choice of $q$ we made above. So $A = {}_qK'$. And then $K' \xrightarrow{\Omega^-_2} K \Rightarrow {}_qK' \xrightarrow{al\Omega^-_2} {}_pK$ as wanted. □

3.5. $^{al}\Omega_3$ **moves.** Let $A$ be a based symbol with $\tau(A) \geq 3$ and $e_1, e_2, e_3$ a 3-gon of $A$. If $\Delta_1, \Delta_2, \Delta_3$ are the three crossings of the 3-gon, then:

$\Delta_1 = (\alpha, \gamma)^{\pi_1}, \Delta_2 = (\beta, \epsilon)^{\pi_2}, \Delta_3 = (\delta, \zeta)^{\pi_3}$, $\alpha, \beta, \gamma, \delta, \epsilon, \zeta \in \mathbb{N}_{2n}$, $|\alpha - \beta| = |\gamma - \delta| = |\epsilon - \zeta| = 1$.

We say that the numbered symbol $A' \in \Sigma$ is derived from $A$ by an algebraic move $^{al}\Omega_3$ on the 3-gon $e_1, e_2, e_3$ whenever:

$$\begin{aligned} A &= \cdots (x:y)^\pi \cdots (\alpha, \gamma)^{\pi_1} (\beta, \epsilon)^{\pi_2} (\delta, \zeta)^{\pi_3} \\ A' &= \cdots (x:y)^\pi \cdots (\beta, \delta)^{\pi_1} (\alpha, \zeta)^{\pi_2} (\gamma, \epsilon)^{\pi_3}. \end{aligned}$$
(3.9)

If for each arc $e_i$ we call as companion the labels of its endpoints, then $A'$ is described as the numbered symbol we get from $A$ by interchanging in the crossings $\Delta_1, \Delta_2, \Delta_3$ the places of the companion labels.

We denote the move as $A \xrightarrow{al\Omega_3} A'$.

Some important properties of $^{al}\Omega_3$ moves are gathered below:

**Lemma 20.** *(a) Let $A \xrightarrow{al\Omega_3} A'$. Then $A' \in \Sigma$, so the $^{al}\Omega_3$ moves are well defined.*

*(b) If $A_1 \xrightarrow{al\Omega_3} A_2$, then $\tau(A_2) = \tau(A_1)$.*

*(c) If $K, K'$ are diagrams, then: $K \xrightarrow{\Omega_3} K' \Rightarrow {}_pK \xrightarrow{al\Omega_3} {}_pK'$ for suitable choices of base points $p \in K, q \in K'$.*

*Proof.* (a) The integers in the pairs of $A'$ form the set $\mathbb{N}_{2n}$. Thus by the definition of the numbered symbols we have $A_2 \in \Sigma$ as wanted.

(b) Immediate.

(c) Figure 27 presents the topological setting near the basic arcs $e_1, e_2, e_3 \in K$ and $e'_1, e'_2, e'_3 \in K'$ for the given $K \xrightarrow{\Omega_3} K'$ move. The base points $p, q$ are chosen either outside the 2-disk of the move or on some $e_i$ and $e'_i$ respectively. We fix a choice to work with.

For the chosen base points, relation (2.19) explains how the crossings inside the 2-disk of the move change and relation (2.20) explains how ${}_pK$ and ${}_qK'$ are related.

Let us note that since the three basic arcs $e_1, e_2, e_3$ of $K$ form a 3-gon of $K$, then as basic arcs of ${}_pK$ they form a 3-gon of it too. So we can perform an $^{al}\Omega_3$ move on it ${}_pK \xrightarrow{al\Omega_3} A$ which satisfies Relation (3.9). This and Relations (2.19) and (2.20) for ${}_qK'$ make $A$ and ${}_qK'$ have the same description. Thus $A = {}_qK'$ and ${}_pK \xrightarrow{al\Omega_3} {}_q K'$ as wanted.

Also, (2.19) and (2.20) imply relation (3.9). It follows that there exists an algebraic $\Omega_3$ move ${}_pK \xrightarrow{al\Omega_3} A$. The numbered symbol □

## 4. Algebraic moves are mirrored by topological moves

In section §3 we gave an algebrization of some topological notions involved in the description of the well-known equivalence of diagrams defined via the topological moves of Alexander, Briggs and Reidemeister. Notably we defined based symbols of diagrams and algebraic moves that send them to similar symbols called numbered symbols. Each topological move (including the isotopies of the plane) has a corresponding algebraic one which we proved mirrors it in the sense that if $K \to [\Omega]K'$ then ${}_pK \to [{}^{al}\Omega]_q K'$ for some suitable choices of the base points $p, q$.

In order to completely replace topological moves by algebraic ones in our study of knots, we need to prove the converse mirroring as well. And this is our next goal. At the same time we'll show that each algebraic move results to an element of $\Sigma_\Delta$ rather than of the larger set $\Sigma$, making the moves an internal operation on $\Sigma_\Delta$.

### 4.1. $\Omega_d$ uncrossing moves, crossing and uncrossing disks, decompositions of diagrams.
This section is auxiliary to the next one where we prove that two based symbols $pK, {}_qK'$ which coincide imply that $K, K'$ are equivalent diagrams.

We start by describing an alteration on a diagram which we call a move and which removes a crossing from the diagram. And then we describe another kind of move that adds a crossing to a diagram. For an exact description of these moves we need to consider base points on the diagrams.

**Definition 21.** *(Nice disks, uncrossing disks and moves, crossing disks, survival and heredity) Let $K$ be a diagram and $p$ a base of it. We denote this based diagram as $K, p$.*

(**Nice disks**) *We call as nice disk of $K$ any 2-disk $U$ on the plane which intersects $K$ transversely and which contains exactly two arcs $\ell_1, \ell_2$ of the parametrized curve $K$ (Figure 36). We call the endpoints of $\ell_1, \ell_2$ as nodes of $K$ in $U$. We give an ordering to these endpoints: we call them as first, second, third and fourth node of $U$ in $K, p$, according to if they are the first entrance, first exit, second entrance and second exit point of $U$ in the order we meet them in a trip on $K$ from its base point $p$ along its orientation. The first and second nodes are respectively the first and second endpoints of the arc $\ell_1$. The third and fourth nodes are respectively the first and second endpoints of the arcs $\ell_2$. We call $\ell_1, \ell_2$ as the first and second leg of $K, p$ in $U$ respectively. It is convenient to call them respectively as the first and second leg of $U$ in $K, p$ as well.*

*$K$ as a parametrized curve splits into three parts with respect to $U$: one is its inner part $U \cap K = \ell_1 \cup \ell_2$, one is its arc $\gamma_1(K)$ from the first exit to the second entrance (second and third nodes) and one is its arc $\gamma_2(K)$ from the second exit to the first entrance (fourth and first nodes). We call $\gamma_1(K), \gamma_2(K)$ as the first and second arm of $K, p$ outside $U$. We also call them as the first and second arm of $U$ in $K, p$.*

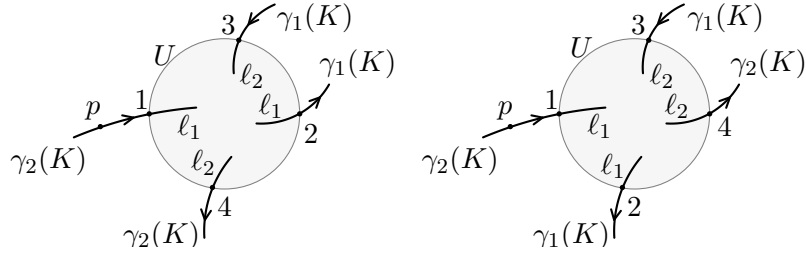

FIGURE 36. Nice disks $U$ for based diagrams $K, p$, their ordered nodes $1, 2, 3, 4$, their first and second legs $\ell_1, \ell_2$, and their first and second arms $\gamma_1(K), \gamma_2(K)$.

(**Uncrossing disks**) *The canonical neighborhoods of the crossings of $K$ (ch. §1.2) are nice 2-disks which hereafter we shall also call as uncrossing disks of $K$ (Figure 37).*

*The two legs of $K$ inside*

*If in an uncrossing disk $U$ of a crossing point the leg $\ell_1$ is over $\ell_2$ at the crossing, then we assign to $U$ the number $\epsilon_U = +1$, otherwise the number $\epsilon_U = -1$. We call $\epsilon_U$ as the kind-number of $U$ for (the based diagram) $K, p$ and we write $U^\epsilon$ for the disk accompanied by its kind-number $\epsilon$.*

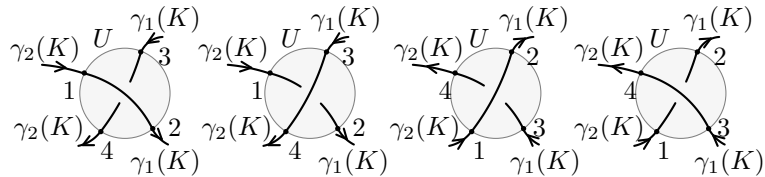

FIGURE 37. Uncrossing disks. The first, second, third and fourth node (entrance and exit points) are named by the numbers 1,2,3,4 respectively. All possible topological placements on the plane are shown here. In all cases the arm $\gamma_1(K)$ is placed to the right of the disk. The base point $p$ is not shown, but it has to lie in $\gamma_2(K)$.

(**Uncrossing moves**) *From a diagram $K$ with base point $p$, and for any uncrossing disk $U$ of a crossing point, we get a new diagram $K_d$ in which the crossing point is removed as follows: we replace the legs of $K$ in $U$ as in Figure 38, and we change the orientation of the first arm $\gamma_1(K)$ to the opposite. Similarly to the usual moves on diagrams, we allow such an alteration up to an* *isotopy of the plane that fixes throughout everything outside the disk $U$. Necessarily, the isotopy fixes the points on the boundary of the disk as well.*

*We consider the result $K_d$ as an oriented closed curve with the orientation on its arms and legs given to them in Definition 21. Then $K_d$ is a diagram and by its definition it is uniquely defined up to an isotopy of the plane. We call $K_d$ as an uncrossing diagram of $K$ in $U$ and the alteration as an uncrossing of $K$ in $U$. We write $K, p \xrightarrow[U^\epsilon]{\Omega_d} K_d, p$, or $K, p \xrightarrow[U]{\Omega_d} K_d, p$, or just $K, p \xrightarrow{\Omega_d} K_d, p$.*

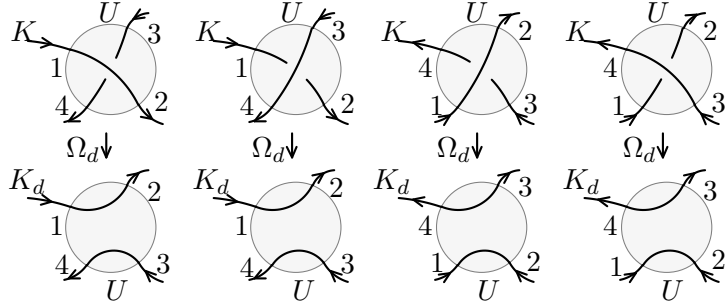

FIGURE 38. The uncrossing move $K, p \xrightarrow[U]{\Omega_d} K_d, p$. The base point $p$ is not shown here.

**(Crossing disks)** *We call a nice disk $U$ of $K$ as a crossing disk of $K$ whenever the legs $\ell_1, \ell_2$ of $K$ in $U$ are disjoint and oppositely orientated. This means that for $\ell_1, \ell_2$ with the induced orientations of $K$, the pair $(U, U \cap K) = (U, \ell_1 \cup \ell_2)$ is homeomorphic to the one in Figure 39 (a). We call the parts $\Gamma_1, \Gamma_3$ in the Figure as the ears of $U$ and the part $\Gamma_2$ as the body of $U$.*

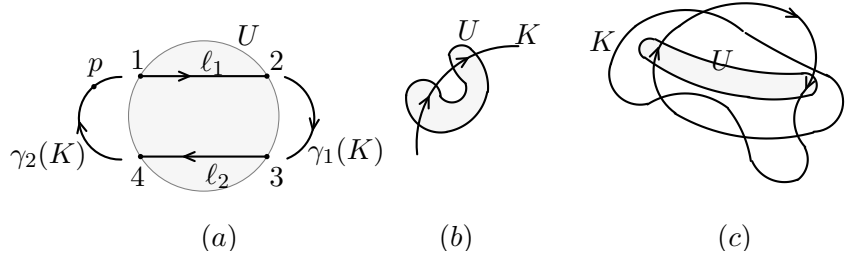

FIGURE 39. (a) The topological model for a crossing disk $U$ of a diagram $K$. The first, second, third and fourth node of $U$ is named here as $1, 2, 3, 4$, respectively. (b), (c) Examples of nice disks in diagrams (the base points are not shown).

*For use in the crucial Lemma 22 later on, we consider here isotopies of our diagrams (they lie on $\rho$) happening not only on the plane $\rho$, but also on the sphere $S^2 = \rho \cup \{\infty\}$. This sphere lies in $S^3 = \mathbb{R}^3 \cup \{\infty\}$ (ch. §1.5).* **(Survival and heredity)** *Let $K$ be a diagram, $p$ a base point of it and $K, p \xrightarrow[U_0]{\Omega_d} K', p$ an uncrossing move. If $U_0$ is the disk of the move and $U = U_0$ or $U$ is a nice disk of $K$ disjoint from $U_0$, then $U$ is a nice disk of $K'$ as well, and we say that $U$ survives in $K'$. We also say that $p$ survives in $K'$. In the case $U = U_0$, the disk becomes a crossing disk of $K'$ (there exist exactly two disjoint legs of $K'$ in it), and if its kind-number in $K, p$ is $\epsilon$, then we say that $\epsilon$ survives as the (inherited) ghost kind-number of $U$ in $K', p$.*

*Let $K$ be a diagram, $p$ a base point of it, and $K \xrightarrow{\Omega_{iso} \text{ or } \Omega_{iso\sigma}} K'$ an isotopy move on the plane or on the sphere with isotopy map $F$, and let $F_1$ be the last moment of $F$. We say that $p$ survives in $K'$ as $p' = F_1(p)$. If $U$ is a nice disk of $K$, the position $U' = F_1(U)$ of $U$ at the end of the move* is a nice disk of $K'$ and we say that $U$ survives in $K'$ as $U'$. If the disk $U$ is an uncrossing disk with kind-number $\epsilon$, then $U'$ is an uncrossing disk of $K'$ with kind-number $\epsilon$ in $K', p'$ as well, so we say that $\epsilon$ survives in $K'$ as kind-number of $U'$. If the disk $U$ is a crossing disk equipped with a ghost kind-number $\epsilon$ in $K, p$, then $U'$ is a crossing disk of $K'$ and we say that $\epsilon$ survives as ghost kind-number of $U'$ in $K', p'$.

**Remark 5.** *Clearly the position of the base point p affects the order of the nodes of a nice disk, as well as the order of its legs and arms.*

*The net result of an uncrossing move is to produce a new diagram $K_d$ which as a point set coincides with $K$ except in $U$ where now it fails to cross itself. $U$ is transformed to a crossing disk of the new diagram with the same points as nodes but with their order changed, namely the old $2, 3$ nodes interchange names becoming the $3, 2$ nodes respectively $(1 \to 1, 2 \to 3, 3 \to 2, 4 \to 4)$.*

*For a given diagram $K$ and a given base point $p$ of it, the kind-number of any uncrossing disk is uniquely determined.*

*The role of the kind numbers around a collection of disjoint uncrossing disks for the crossing points of a diagram is destined below to be that of the recovery instructions of the diagram after it is completely uncrossed in all the disks. And the exact purpose of the last part in Definition 23, is to help us talk about decompositions and reconstructions of diagrams. The definition involves the sphere $S^2 = \rho \cup \{\infty\}$ because this will be a much easier surface than the plane $\rho$ to work with in Lemma 22. Nevertheless we can stick with the plane throughout if we wish.*

The following Definition describes a way to change a given diagram in steps so as to get rid off all its crossings.

**Definition 22.** *(Uncrossing decompositions of diagrams) Let $K$ be a diagram of order $\tau(K) = n \geqslant 1$ and $p$ a base point of it. Let $(\Delta_1, \ldots, \Delta_n)$ be a fixed order on the crossings of $K$.*

*We call as uncrossing sequence or as a decomposition of $K, p$, any sequence of alterations of $K$ as follows:*

*Choose an ordered collection $S = (U_1, U_2, \ldots, U_n)$ of uncrossing disks of $K$ which is nice for $K, p$ and $(\Delta_1, \ldots, \Delta_n)$, meaning that $\Delta_i \in U_i \forall i$, the disks are disjoint and none contains $p$.*

*We call $K = K_1$ and we perform an uncrossing move $K_1, p \xrightarrow[U_1^{\epsilon_1}]{\omega_{1d}} K_2, p$. Then perform an uncrossing move $K_2, p \xrightarrow[U_2^{\epsilon_2}]{\omega_{1d}} K_3, p$. And so on until a last uncrossing move $K_n, p \xrightarrow[U_n^{\epsilon_n}]{\omega_{nd}} K_{n+1}, p$. $\epsilon_i \in \{-1, 1\}$ are the uniquely defined kind-numbers of the disks $U_i$ for the based diagrams $K_i, p$.*

*We remind (ch. Definition 21) that after the $\omega_{id}$ move, the uncrossing disk $U_i$ of $K_i$ becomes a crossing disk of $K_{i+1}$ and its kind-number $\epsilon_i$ for $K_i, p$ becomes its ghost kind-number in $K_{i+1}, p$. After each one of the subsequent moves $\omega_{(i+1)d}, \ldots, \omega_{nd}$, the disk $U_i$ remains a crossing disk of $K_j, p$ ($j > i$), possibly with its nodes changing order. By its definition, the ghost kind-number of $U_i$ in $U_j, p$ does not change during these moves. Let $1_i, 2_i, 3_i, 4_i$ be the first, second, third and fourth node of $U_i$ as a crossing disk of the resulting diagram $K_{n+1}, p$ (we keep the base point $p$ fixed throughout during the moves).*

*We call the ordered collection of topological spaces:*
$D_{K,p,S} = (K_{n+1}, p, U_1^{\epsilon_1}, 1_1, 2_1, 3_1, 4_1, U_2^{\epsilon_2}, 1_2, 2_2, 3_2, 4_2, \ldots, U_n^{\epsilon_n}, 1_n, 2_n, 3_n, 4_n)$
*as an uncrossing decomposition of the based diagram $K, p$. We also write for it:*
$D_{K,p,S} = (K_{n+1}, p, U_1, 1_1, 2_1, 3_1, 4_1, U_2, 1_2, 2_2, 3_2, 4_2, \ldots, U_n, 1_n, 2_n, 3_n, 4_n)^{(\epsilon_1, \epsilon_2, \ldots, \epsilon_n)}$.

*We call the $(5n+2)$-ad of topological spaces in the last notation as the base of the decomposition and the $n$-ad $(\epsilon_1, \epsilon_2, \ldots, \epsilon_n)$ as its exponent. We call the ordered collection of topological spaces $(K = K_1, K_2, \ldots, K_{n+1})$ as an uncrossing sequence of diagrams of $K, p$. We call $K_{n+1}, p$ as an uncrossing of $K, p$ or as a final (based) diagram in the decomposition of $K, p$.*

*According to Lemma 21 that follows, the diagrams $K_2, \ldots, K_{n+1}$ in the above process are indeed defined. Also, $D_{K,p,S}$ is well defined up to isotopy on the plane for the spaces of its base, regardless of the chosen nice n-ad $S$ of disks (for $K, p$ for the ordered crossings $(\Delta_1, \ldots, \Delta_n)$) or of the exact location of the legs in the $U_i$'s after the $\omega_{id}$ moves. So we write $D_{K,p}$ for any member in this isotopy class and call it as the uncrossing decomposition of $K, p$ (instead of just "an" uncrossing decomposition). We also call $K_{n+1}, p$ as the uncrossing of $K, p$ (instead of just "an" uncrossing) or as the final (based) diagram in the decomposition of $K, p$.*

*Notice that the isotopy is an operation that by the definition of survival and heredity, fixes the exponent in the sense that for all $i$ and for the diagram at hand at any moment of the isotopy, the ghost kind-number $\epsilon_i$ remains always the same.*

*Whenever for all $i$ the crossing $\Delta_i$ is the one in $K_i$ which has 1 as one of its labels, then we call the uncrossing decomposition $D_{K,p,S}$ of $K, p$ as direct and the corresponding sequence as a direct sequence. We denote it also as $\mathbb{D}_{K,p,S}$. As a special case to the general one, we have that all direct decompositions are isotopic on the plane and we write $\mathbb{D}_{K,p}$ for any member in this isotopy class and call it as the unique direct uncrossing decomposition of $K, p$. We call $K_{n+1}, p$ as the direct uncrossing of $K, p$ or as the final (based) diagram in the direct uncrossing of $K, p$.*

*For a diagram $K$ with $\tau(K) = 0$ it is $K \approx S^1$ and there exist no crossings on it. So there exist no non-empty ordered collections $S = (U_1, U_2, \ldots, U_n)$ of uncrossing disks of $K$ which are nice for $K, p$, where $p$ is an arbitrary base point of $K$. We say $S = \varnothing$ is the only nice ordered collection of $K, p$ and we define the only uncrossing decomposition of $K, p$ to be the ordered pair of topological spaces $D_{K,p,\varnothing} = (K, p)^\varnothing$. We also call this as the unique direct uncrossing decomposition $\mathbb{D}_{K,p,\varnothing}$ of $K, p$.*

*Finally, let $D_{K,p,S} = B^E$ and $D_{K',p',S'} = (B')^{E'}$ be two uncrossing decompositions. We call them isotopic by an isotopy $F$ on the plane or the sphere which sends $D_{K,p,S}$ onto $D_{K',p',S'}$, whenever $E = E'$ and $F_1(B) = B'$ ($F_1$ is the last moment of $F$). Then $F^{-1}$ is an isotopy that sends $D_{K',p',S'}$ onto $D_{K,p,S}$, thus we just say that the two decompositions are isotopic.*

**Remark 6.** *Although in an uncrossing decomposition the disk $U_i$ of the i-th move survives as a nice disk of $K_j$ during the previous moves $\omega_{jd}$, $(j < i)$, its kind-number first comes to life as a disk of $K_i$ and immediately after $\omega_{id}$ it dies surviving as a ghost kind-number for the rest of the diagrams until the last one $K_{n+1}$.*

*It will not harm to stress the following detail in the definition of an uncrossing decomposition: the order of the nodes $1_i, 2_i, 3_i, 4_i$ of a crossing disk $U_i$ of $K$ as first, second, third, fourth respectively, is the one they have when $U_i$ is considered as a crossing disk of the end result $K_{n+1} \approx S^1$. This order should not be confused with the order of the same points when $U_i$ is considered as an uncrossing disk of the given diagram $K$. Neither with the order they have for $U$ as a crossing or an uncrossing disk for any other diagram $K_j$ in the uncrossing sequence, except of course for the last one, $K_{n+1}$.*

*The exponent part of an uncrossing decomposition is just a set of over/under instructions which does not makes sense all at once for the $(5n + 2)$-ad of the topological spaces in the base. But it instructs us to perform inversely one after the other the "crossing moves" $K_{n+1}, p \xrightarrow[U_n^{\epsilon_n}]{\omega_n^{\epsilon_n}} K_n, p \cdots K_3, p \xrightarrow[U_2^{\epsilon_n}]{\omega_2^{\epsilon_2}} K_2, p \xrightarrow[U_1^{\epsilon_1}]{\omega_1^{\epsilon_1}} K_1, p$ recovering an isotopic (on the plane) version of $K, p$. The crossing moves can be defined as an inverse alteration of that in the uncrossing moves, but we do not define them officially here, because we shall not be really interested in them in what follows. The exact way of reconstruction of the diagram $K$ is contained in the proof of part (d) of Lemma 22 below.*

*As we can easily check through specific examples, the order of the crossings is important in forming an uncrossing decomposition of a diagram. The direct uncrossings are in a way the easiest*

to handle in respect to a given position of the base point p. They are at least powerful enough to help us prove Lemma 22 in the next section.

The last part of the Definition involves the 2-sphere $S^2$ for the same reason (ch. Remark 5) this happens in the previous Definition.

**Lemma 21.** *Let $K$ be a diagram of order $\tau(K) = n$, let $p$ be a base point of it, and $(\Delta_1, \ldots, \Delta_n)$ be a fixed order on the crossings of $K$. Let also $S = (U_1, U_2, \ldots, U_n)$ be a nice n-ad for $K, p$ for the ordered crossings $(\Delta_1, \ldots, \Delta_n)$; recall this means the $U_i$'s are disjoint uncrossing neighborhoods of $K$ which do not contain $p$ and such that $\Delta_i \in U_i$ for all $i$.*

*(a) The uncrossing sequence of diagrams $(K = K_1, K_2, \ldots, K_{n+1})$ of $K, p$, corresponding to the n-ad $S = (U_1, U_2, \ldots, U_n)$ is indeed defined. That is, the uncrossing moves $K_i, p \xrightarrow[U_i]{\omega_{1d}} K_{i+1}, p$ are defined.*

*(b) $K_{n+1}$ is a circle.*

*(c) An uncrossing decomposition $D_{K,p,S} = (K_{n+1}, p, U_1, 1_1, 2_1, 3_1, 4_1, U_2, 1_2, 2_2, 3_2, 4_2, \ldots, U_n, n_1, n_2, n_3, n_4)^{(\epsilon_1, \epsilon_2, \ldots, \epsilon_n)}$ of $K, p$ corresponding to the n-ad $S = (U_1, U_2, \ldots, U_n)$ remains the same up to isotopy of the plane regardless of the exact location of the legs in the $U_i$'s after the $\omega_{id}$'s.*

*(d) The unique (by (c)) up to isotopy on the plane uncrossing decompositions $D_{K,p,S_U}$ and $D_{K,p,S_V}$ of $K, p$ corresponding to two nice n-ads $S_U = (U_1, U_2, \ldots, U_n)$ and $S_V = (V_1, V_2, \ldots, V_n)$ (for $K, p$ for the ordered crossings $(\Delta_1, \ldots, \Delta_n)$), are isotopic on the plane.*

*(e) If $(K = K_1, K_2, \ldots, K_{n+1})$, $n \geq 2$ is an uncrossing sequence for a direct uncrossing decomposition $\mathbb{D}_{K,p,S}$ of $K, p$ where $S = (U_1, \ldots, U_n)$, then $(K_2, \ldots, K_{n+1})$ is an uncrossing sequence for a direct uncrossing decomposition $\mathbb{D}_{K_2, p, S'}$ ($S' = (U_2, \ldots, U_n)$) of $K, p$. And if*

$$\mathbb{D}_{K,p,S} = (K_{n+1}, p, U_1, 1_1, 2_1, 3_1, 4_1, U_2, 1_2, 2_2, 3_2, 4_2, \ldots, U_n, 1_n, 2_n, 3_n, 4_n)^{(\epsilon_1, \epsilon_2, \ldots, \epsilon_n)},$$

*then:*

$$\mathbb{D}_{K_2, p, S'} = (K_{n+1}, p, U_2, 1_2, 2_2, 3_2, 4_2, \ldots, U_n, 1_n, 2_n, 3_n, 4_n)^{(\epsilon_2, \ldots, \epsilon_n)},$$

*where for $n = 1$, $S' = \varnothing$ and $D_{K_{n+1}, p, S'} = D_{K_2, p, \varnothing} = (K_2, p)^\varnothing$.*

*(f) If $S = (U_1, \ldots, U_n)$ is an ordered nice sequence for a direct decomposition of $K, p$, then for each crossing point $\Delta = (a : b)^\pi$ the index $m$ for which $\Delta \in U_m$ is determined by any of $a, b$.*

In other words, all based diagrams with the same based symbol, have a common 1-1 correspondence between crossing points as we meet them in a trip on the diagrams (from the base points and along the orientations) and uncrossing disks in the ordered nice n-ads of the direct decompositions.

*(g) Let $S = (U_1, \ldots, U_n)$ be an ordered nice sequence for a direct decomposition of $K, p$, and let $K_{n+1}$ be the final diagram in the decomposition. For each $i = 1, 2, \ldots, n$, the k-th node ($k = 1, 2, 3, 4$) of $U_i$ as an uncrossing disk of $K, p$ becomes the $k_j$-th node for some $k_j \in \{1, 2, 3, 4\}$ of $U_i$ as a crossing disk of $K_{n+1}, p$. Then $k_j$ is uniquely determined by $i$ and $k$.*

In other words, if the k-th node of the i-th uncrossing disk $U_i$ of a based diagram $K, p$ becomes the $k_j$ node of $U_i$ in the direct decomposition $\mathbb{D}_{K,p}$, then the k-th node of the i-th uncrossing disk $U'_i$ of a based diagram $K', q$ for which $_pK = {_q}K'$, also becomes the $k_j$-th node of $U'_i$ in the direct decomposition $\mathbb{D}_{K',q}$.

*(h) Let $S = (U_1, \ldots, U_n)$ be an ordered nice sequence for a direct decomposition of $K, p$, and let $(K_1 = K, K_2, \ldots, K_n)$ the sequence of uncrossing diagrams in the direct decomposition $\mathbb{D}_{K,p,S}$. Then the first, second, third and fourth node of $U_1$ in $K_2, p$ remain respectively as the first, second, third and fourth node of $U_1$ in each one of the subsequent based diagrams $K_3, p$ and $K_4, p$, and $\ldots, K_{n+1}, p$.*

More generally, for each $i = 1, 2, \ldots, n$, the first, second, third and fourth node of $U_i$ in $K_{i+1}, p$ remain as nodes with this order for $U_i$ in all subsequents based diagrams, up to the last one $K_{n+1}, p$.

The typical topological picture of a direct uncrossing decomposition of some diagram is shown in Figure 40.

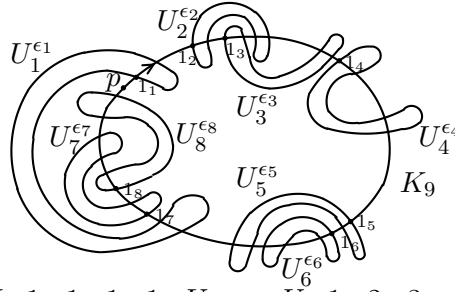

$$\mathbb{D}_{K,p} = (K_9, U_1, 1_1, 1_2, 1_3, 1_4, U_2, \ldots, U_8, 1_8, 2_8, 3_8, 4_8)^{(\epsilon_1,\ldots,\epsilon_8)}$$

FIGURE 40. A typical topological picture of a direct uncrossing decomposition of a based diagram $K, p$. The diagram $K = K_1$ was given here with 8 crossing points, thus the final diagram in its decomposition is named $K_9$. By the definition of the decomposition, $K_9$ is a circle and contains the base point $p$. For each uncrossing disk $U_i$ in the decomposition considered as a nice disk of $K_9, p$, its first, second, third and fourth node are denoted as $1_i, 2_i, 3_i, 4_i$; in the figure only the first node $1_i$ is depicted. The kind-numbers $\epsilon_i$ belong to $\{-1, +1\}$ and they are determined by $K, p$. In the depicted example the disks $U_1, U_2, U_7$ have their bodies in the interior of $K_9$ and their ears in the exterior, whereas the converse holds for the rest of the disks. This remark has a meaning whenever $K_9$ is considered on the plane $\rho$. Whenever $K_9$ is viewed on the sphere $S^2$, a way to replace the remark is to say that the bodies of $U_1, U_2, U_7$ are in the same region of the sphere with respect to its circle $K_9$, which is different from the region in which lie the bodies of the rest of the disks.

*Proof.* (a) For $n = 0$ the sequence contains only $K_1 = K$ and the result holds.

Let $n > 0$. Since $U_1$ is an uncrossing disk of $K_1$, the move $K_1, p \xrightarrow[U_1]{\omega_{1d}} K_2, p$ is defined, so $K_2$ is defined. By the definition of the move $\omega_{1d}$, everything outside $U_1$ remains unaltered, and since $U_1, U_2$ are disjoint, $U_2$ is a crossing disk of $K_2$ as well. So $K_2, p \xrightarrow[U_2]{\omega_{2d}} K_3, p$ is defined, thus $K_3$ is defined, and so on, implying all $K_i$ are defined inductively.

(b) $K_{n+1}$ is a diagram with no crossings, thus it is a circle.

(c) For $\tau(K) = 0$, there exist no moves in a sequence so no choices for them in forming the base of the decomposition. So $\mathbb{D}_{K,p} = (K, p)^\emptyset$ and there exists no pair of topological spaces resulting as decomposition of $K, p$ other than $(K, p)$. Of course $(K, p)$ is isotopic on the plane to itself, for example by the constant isotopy, and the desired result holds.

Let now $\tau(K) > 0$.

First, let us prove that the first, second, third and fourth node $1_i, 2_i, 3_i, 4_i$ of $U_i$ are well defined points regardless of the choices:

By definition, there exist a lot of positions of $D_{K,p,S}$ on the plane depending on the choices we make during each $\omega_{id}$ move for the location of the legs in $U_i$ after the move is over. No matter what the choices are, $U_i$ remains an uncrossing disk after the moves $\omega_{jd}$ with $j < i$ and it turns to a crossing disk after the moves $\omega_{jd}$ with $j \geqslant i$.

Consider some $i$. After a move $\omega_{jd}$ with $j < i$, the disk $U_i$ remains pointwise fixed. So its two legs remain pointwise fixed, still crossing each other but they may have exchanged orders, and also, the endpoints of each may have exchanged orders.

Such changes are due to orientation alterations of the arms of the resulting diagram $K_{j+1}$ after the move. Nevertheless, the orientation of $K_{j+1}$ does not depend on the choices we make in $U_j$ regarding the exact location of the legs after $\omega_{jd}$. So the order of the nodes of $U_i$ as a crossing disk of $K_{j+1}$ is fixed. Also, for the move $\omega_{id}$ the disk changes as the move dictates making the legs lose their intersection. The nodes change their order, but this order will be the same regardless of the exact location of the legs in $U_i$ after this move. So the order of the nodes of $U_i$ as a crossing disk of $K_{i+1}$ is fixed. Finally, after a move $\omega_{jd}$ with $j > i$ no changes happen to $U_i$ except possibly for the order of its nodes and as in the case for moves with $j < i$ this order does not depend on the choice for the location of the legs in the disk $U_j$ of the move. So the order of the nodes of $U_i$ as a nice disk of the final diagram $K_{n+1}$ is well defined independently of the choices in forming this diagram. In other words $i_1, i_2, i_3, i_4$ are well defined points as wanted.

Next let us prove that the ghost kind-numbers $\epsilon_i$ of the $U_i$'s are well defined regardless of our choices:

Now for any $i$, let $\ell_{i1}, \ell_{i2}$ be the choice for the legs of $U_i$ in an uncrossing decomposition $D_{K,p,S}$ of $K, p$ ($S = (U_1, \ldots, U_n)$), and say $K_{n+1}$ is the final diagram. And let $\ell'_{i1}, \ell'_{i2}$ be the choice for the legs of $U_i$ in another uncrossing decomposition $D'_{K,p,S}$ of $K, p$ for the same $S = (U_1, \ldots, U_n)$. And let this time $K'_{n+1}$ be the final diagram.

Consider an $i$. $\epsilon_i$ first comes to existence as a kind-number for $K_i$. It becomes a ghost kind-number of $K_{i+1}$ after the $\omega_{id}$ move and continues to exist as a ghost kind number retaining its value (by its definition) after every subsequent move. Now according to what we proved just above, the choices we make for the moves $\omega_{jd}$ with $j < d$ do not alter the order of the nodes of $U_i$ as a crossing disk of $K_i$. So the first and second leg of $K_i$ in $U$ are uniquely determined. Then the over/under information of $K = K_1$ inside $U_i$ which did not change up until $K_i$ came up in the uncrossing sequence of diagrams, gives us the kind-number $\epsilon_i$ of $K_i$ in $U_i$ in a unique way. So $\epsilon_i$ are well defined regardless of our choices as we wanted to prove.

We are now in position to prove the desired result:

We consider the above setting where for $S = (U_1, \ldots, U_n)$, we get the uncrossing decomposition $D_{K,p,S}$ and $D'_{K,p,S}$ of $K, p$, by choosing in $U_i$ for all $i$ the legs $\ell_{i1}, \ell_{i2}$ for the first decomposition and the legs $\ell'_{i1}, \ell'_{i2}$ for the second. Say $K_{n+1}$ and $K'_{n+1}$ be the final diagrams respectively.

According to what we proved above, it is:
$D_{K,p,S} = (K_{n+1}, p, U_1, 1_1, 2_1, 3_1, 4_1, U_2, 1_2, 2_2, 3_2, 4_2, \ldots, U_n, 1_n, 2_n, 3_n, 4_n)^{(\epsilon_1, \epsilon_2, \ldots, \epsilon_n)}$
and
$D'_{K,p,S} = (K'_{n+1}, p, U_1, 1_1, 2_1, 3_1, 4_1, U_2, 1_2, 2_2, 3_2, 4_2, \ldots, U_n, 1_n, 2_n, 3_n, 4_n)^{(\epsilon_1, \epsilon_2, \ldots, \epsilon_n)}$.

Let $F$ be the constant isotopy of the the complement of the interiors of the $U_i$'s. And let $F_i$ an isotopy of each $U_i$ into itself which fixes throughout the boundary and brings the pair $(\ell_{i1}, \ell_{i2})$ onto the pair $(\ell'_{i1}, \ell'_{i2})$; such an isotopy inside each $U_i$ exists by the definition of the $\omega_{id}$ move. Then $F$ agrees with the $F_i$'s on the boundaries and we can assemble all these isotopies to an isotopy $G$ of the whole plane. It is $G(U_i) = U_i$ (although not point-wise) and $G(K_{n+1}) = K'_{n+1}$. Since $G$ keeps all points of $\partial U_i$ fixed, it keeps fixed their nodes $i_1, i_2, i_3, i_4$ with respect to either $K_{n+1}, K'_{n+1}$, i.e. $G(i_j) = i_j$ for all $i = 1, 2, \ldots, n$ and $j = 1, 2, 3, 4$. Since $p$ lies outside the disks, we also have $G(p) = p$. Thus we have $G(D_{K,p,S}) = D'_{K,p,S}$ as wanted.

(d) Let
$D_{K,p,S_U} = (K_{(n+1)U}, p, U_1, 1_{1U}, 2_{1U}, 3_{1U}, 4_{1U}, U_2, 1_{2U}, 2_{1U}, 3_{2U}, 4_{2U}, \ldots, U_n, 1_{nU}, 2_{nU}, 3_{nU}, 4_{nU})$
and
$D_{K,p,S_V} = (K_{(n+1)V}, p, V_1, 1_{1V}, 2_{1V}, 3_{1V}, 4_{1V}, V_2, 1_{2V}, 2_{1V}, 3_{2V}, 4_{2V}, \ldots, V_n, 1_{nV}, 2_{nV}, 3_{nV}, 4_{nV})$.

For all $i$ we consider disks $W_i$ small enough to be contained in both $U_i, V_i$ and we choose them to be uncrossing disks for the corresponding crossing points of $K$. Then $(W_1, \ldots, W_n)$ is a nice sequence of uncrossing disks for $K, p$ and for the given order of the crossings. Let this gives the following uncrossing decomposition of $K, p$:
$D_{K,p,S_W} = (K_{(n+1)W}, p, W_1, 1_{1W}, 2_{1W}, 3_{1W}, 4_{1W}, W_2, 1_{2W}, 2_{1W}, 3_{2W}, 4_{2W}, \ldots, W_n, 1_{nW}, 2_{nW}, 3_{nW}, 4_{nW})$.

Since $W_i \subset U_i$ for all $i$, the diagrams $K_{(n+1)U}, K_{(n+1)W}$ coincide outside the union of the $U_i$'s. In each $U_i$ they differ as in Figure 41.

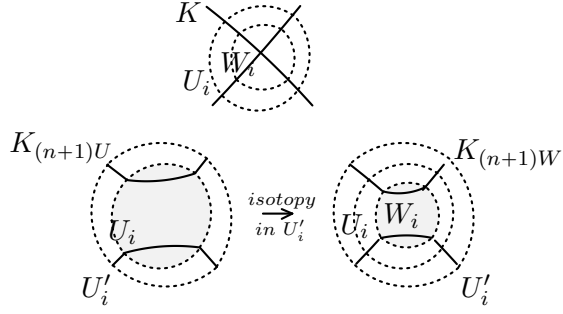

FIGURE 41. Two uncrossing disks $U_i, W_i$ of the same crossing of $K$, one inside the other. A disk $U'_i$ slightly bigger than $U_i$. An isotopy in $U'_i$ fixing $\partial U'_i$ brings $U_i$ onto $W_i$ and the legs of $U_i$ in $K_{(n+1)U}$ onto the legs of $W_i$ in $K_{(n+1)W}$.

We consider a slightly bigger uncrossing disk $U'_i$ with $U_i$ in its interior (Figure 41). Then there exists an isotopy in $U'_i$ that fixes the boundary, brings $U_i$ to $W_i$ and the legs of $U_i$ (in $K_{(n+1)U}, p$) to the legs of $W_i$ (in $K_{(n+1)W}, p$). Assembling all these isotopies we get an isotopy $F$ of the whole plane that fixes everything outside the $U'_i$'s and which brings the $(n+1)$-ad $(K_{(n+1)U}, U_1, \ldots, U_n)$ to the $(n+1)$-ad $(K_{(n+1)W}, W_1, \ldots W_n)$. For each $i$ this isotopy actually brings the first, second, third and fourth node of $U_i$ (in $K_{(n+1)U}, p$) to the corresponding node of $W_i$ (in $K_{(n+1)W}, p$) for all $i$. This is readily checked inductively for each pair of corresponding moves in the two uncrossing procedures. Similarly there exists an isotopy $G$ of the plane which brings the $(n+1)$-ad $(K_{(n+1)V}, V_1, \ldots, V_n)$ to the $(n+1)$-ad $(K_{(n+1)W}, W_1, \ldots W_n)$ and the first, second, third and fourth node of $V_i$ to the corresponding node of $W_i$ for all $i$. Since $p$ lies outside the disks, we also have $F(p) = p$. Then the isotopy $F^{-1} \circ G$ brings $D_{K,p,S_V}$ onto $D_{K,p,S_U}$ as wanted.

(e) For the direct uncrossing decomposition $\mathbb{D}_{K,p,S}$ of $K, p$ we perform in a row the uncrossing moves $K = K_1, p \xrightarrow[U_1]{\omega_{1d}} K_2, p \xrightarrow[U_2]{\omega_{2d}} K_3, p \cdots \xrightarrow[U_n]{\omega_{nd}} K_{n+1}, p$. Since the sequence is the direct one, the uncrossing disk $U_i$ contains the crossing for which 1 is a label. Then if $n > 1$, the sequence $K_2, p \xrightarrow[U_2]{\omega_{2d}} K_3, p \cdots \xrightarrow[U_n]{\omega_{nd}} K_{n+1}, p$ is an uncrossing sequence of $K_2$ where in each $K_i$ when considered with base the point $p$, the uncrossing disk $U_i$ contains the crossing for which 1 is a label. So by definition this is the direct uncrossing sequence of $K_2$. Since this sequence is part of the one for $K = K_1$, the order of the nodes in each $U_i, i \geq 2$ does not change from the order they had in the bigger sequence. If on the other hand $n = 1$, then $K_2$ is a circle and it is also the last diagram $K_{n+1}$ in the uncrossing sequence. Then $S' = \emptyset$ and $\mathbb{D}_{K_{n+1},p,S'} = \mathbb{D}_{K_2,p,\emptyset} = (K_2,p)^\emptyset$ is the unique uncrossing sequence of $K_2$. So we have the result in all cases.

(f) We prove it inductively on $n = \tau(K)$.

For $n = 0$ it holds vacuously, whereas for $n = 1$ it holds immediately.

Let it be true for some $n \geq 1$. And let $K, p$ be a based diagram with $\tau(K) = n + 1$. Let also, $S = (U_1, U_2, \ldots, U_{n+1})$ be a nice $(n+1)$-ad of uncrossing disks for a direct uncrossing decomposition of $K$, and let $(K_1 = K, K_2, \ldots K_{n+1})$ be the corresponding uncrossing sequence of diagrams.

According to part (e) of the Lemma, $S_2 = (U_2, \ldots, U_{n+1})$ is a nice $n$-ad of uncrossing disks for a direct uncrossing decomposition of $K_2$, with $(K_2, \ldots K_{n+1})$ as its corresponding uncrossing sequence of diagrams.

$K_2$ contains all crossings of $K_1$ except the first which disappears because of the first uncrossing move of $K$. Point $p$ lies in $K_2$. The labels of the crossings in $K_2, p$ come by the relabeling of the same crossings in $K, p$ given in the following table:

$$\begin{array}{cccccccccc}
K = K_1, p & : & 1 & 2 & \ldots & a-1 & a & a+1 & \ldots & 2n+2 \\
K_2, p & : & & a-2 & \underset{\leftarrow}{\ldots} & 1 & & a-1 & \underset{\rightarrow}{\ldots} & 2n
\end{array}$$

Let $\Delta = (x : y)$ be a crossing point of $K, q$ other than the first one. As a crossing point of $K_2, p$, it is $\Delta = (x' : y')$ with the labels $x', y'$ determined uniquely be $x, y$ according to the table above. Since $\tau(K_2) = n$, the induction hypothesis holds for $K_2, p$. So the index $m \in \{2, \ldots, n+1\}$ for which $\Delta \in U_m$ in the nice $n$-ad $(U_2, \ldots, U_{n+1})$ of $K_2, p$ is determined by either $x', y'$, thus by either $x, y$. Hence the index $m \in \{2, \ldots, n+1\}$ for which $\Delta \in U_m$ in the nice $n$-ad $(U_1, U_2, \ldots, U_{n+1})$ of $K, p$ is determined by either $x, y$.

Finally, for the first crossing $\Delta = (1 : a)$ of $K_p$ it is always $\Delta \in U_1$. For all based diagrams with the same based symbol, their first crossing has the same description $\Delta = (1 : a)$ for the same $a$. So the index 1 of $U_1$ in which the first crossing point $\Delta$ lies, is determined uniquely by either label 1 or $a$. And the desired result have been proved.

(g) We prove it inductively on $n = \tau(K)$.

For $\tau(K) = 0$ the result holds vacuously.

For $\tau(K) = 1$: passing from $K = K_1, p$ to $K_2, p$ in the direct decomposition, we know the exact way the node $k$ of $U_1$ in $K, p$ becomes the node $k'$ of $U_1$ in $K_2, p$: $1 \to 1, 2 \to 3, 3 \to 2, 4 \to 4$. So $k$ and $i = 1$ determine uniquely $k'$.

Let the result be true for some $n \geq 1$. And let $K, p$ be a based diagram with $\tau(K) = n + 1$. Let also, $S = (U_1, U_2, \ldots, U_{n+1})$ be a nice $(n+1)$-ad of uncrossing disks for a direct uncrossing decomposition of $K$, and let $(K_1 = K, K_2, \ldots K_{n+1})$ be the corresponding uncrossing sequence of diagrams.

$K_2$ contains $p$ and all crossings of $K_1$ except the first one. The labels of the crossings in $K_2, p$ come by the relabeling of the same crossings in $K, p$ given in the table of part (e) of this Lemma.

The table says that the crossing point of $U_i$, $i = 2, 3, \ldots, n+1$ has two descriptions $(x, y)^\pi$ and $(x', y')^\pi$ when considered in $K, p$ and $K_2, p$ respectively, so that $x$ determines $x'$ and $y$ determines $y'$ and vice versa. In $K, p$, the unique relation $x < y$ or $y < x$ that holds for $x, y$ determines which one of the legs of $U_i$ gives entrance and exit points in $U_i$ as first and second node, and which as third and fourth node. So knowledge of $i$, means knowledge of $x, y$ and this means knowledge of the order of the nodes of $U_i$ in $K, p$. Similarly in $K_2, p$, the unique relation $x' < y'$ or $y' < x'$ that holds for $x', y'$ determines which one of the legs of $U_i$ gives entrance and exit points in $U_i$ as first and second node, and which as third and fourth node. So knowledge of $i$, means knowledge of $x', y'$ and this means knowledge of the order of the nodes of $U_i$ in $K_2, p$. Since $x, y$ determine $x', y'$ we have that $i$ and each $k$-node of $U_i$ in $K, p$ determine the name $k'$ of the same point as a node of $U_i$ in $K_2, p$. Then since $\tau(K) = n$, the induction hypothesis holds for $K_2$ and we have that $i$ and $k'$ determine the order of this node in $U_i$ as a disk in the resulting direct uncrossing based diagram $K_{n+1}, p$. Since this is also the resulting based diagram in the direct uncrossing decomposition of $K, p$ as well (ch. part (e) of the current Lemma), we have the desired result.

It remains to show the result for $U_1$ as well: The order $1, 2, 3, 4$ of the nodes of $U_1$ in $K, p$, changes when $U_1$ is considered a disk in $K_2, p$ as follows: $1 \to 1, 2 \to 3, 3 \to 2, 4 \to 4$. This is the order these points have after each one of the subsequent uncrossing moves, so this is the order they have in $K_{n+1}, p$ as well: Indeed, the arc $\ell$ of $K_2$ between $p$ and node 2 of $U_1$ (Figure 42), also contains node 1 of $U_1$ as well as the leg $\ell_1$ of $U_1$ in $K_2, p$ between 1 and 2. Say $\ell_2$ is the second leg of $U_1$ in $K_2, p$. The arcs $\ell, \ell_1, \ell_2$ are not disturbed point-wise by any subsequent move in the uncrossing sequence of $K$. As $p$ always lies in the second arm of each $U_i$ in all diagrams produced during the uncrossings, the orientation of $\ell$ is not changed. Thus if we move in each of $K_2, \ldots, K_{n+1}$ along $\ell$ (starting at $p$ and keeping in the orientation of the diagram), we meet $U_1$ and then we exit this disk always at the same points, thus these are always the first and second node of $U_1$. Hence $\ell_1$ remains unchanged also orientation-wise through the moves. Since $\ell_1, \ell_2$ are always oppositely oriented in $U_1$, it follows that $\ell_2$ remains always unchanged orientation-wise, thus its first and second endpoints remain alway the nodes 3 and 4 respectively of $U_1$ in all $K_2, p$ and $K_3.p$ up to $K_{n+1}, p$. Now as explained above, the order of the nodes in $K_2, p$ is determined by the one they had for $U$ in $K, p$ (recall the change $1 \to 1, 2 \to 3, 3 \to 2, 4 \to 4$ as we pass from $K = K_1, p$ to $K_2, p$). So the numbers $i = 1$ and $k = 1, 2, 3, 4$ determine $k' \in \{1, 2, 3, 4\}$, in this case as well. We are done!

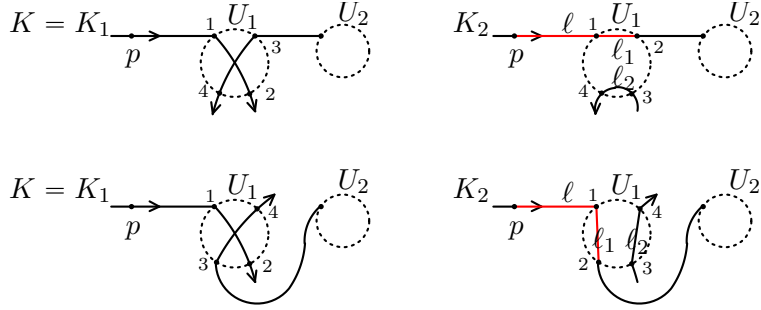

FIGURE 42. Depicted are the first uncrossing disk $U_1$ in a direct decomposition of a diagram $K, p$ and the position of the nodes $1, 2, 3, 4$ of $U_1$ when considered as a disk of either $K_1, p$ or $K_2, p$ ($K_2$ is the result of uncrossing $K = K_1$ in $U_1$). $\ell$ is the oriented arc in $K_2$ between $p$ and node 2, with first endpoint $p$ and the induced orientation of $K_2$. $\ell, \ell_1, \ell_2$ are not affected by any subsequent uncrossing move in the decomposition sequence of $K = K_1$, both point-wise and orientation-wise. Thus the order of the nodes of $U_1$ in $K_2, p$ is the same with their order for $U_1$ in $K_{n+1}, p$.

(h) A proof for $i = 1$ is included in the last paragraph above, in the proof of part (g) of the Lemma. This immediately implies the truth of the result for $i = 2$, and then this the truth for $i = 3$ and so on, as wanted. □

4.2. $^{al}\Omega_{iso}$ **mirroring.** As a first concrete step in showing that algebraic moves are mirrored by topological moves, we prove in this section that the algebraic isotopy move implies that both symbols involved in the move are realized as based symbols of diagrams instead of merely being numbered symbols (we already know this holds for the first symbol in the move by the very definition of the move). And we also prove that although we may need more than one topological move to get equivalence of the two diagrams in the equivalence relation $\underset{top}{\sim}$ (in the set $\mathbb{D}$ of diagrams on the plane), the two diagrams are indeed equivalent.

Most auxiliary work for the proof was done in the previous section. It is completed formally in part (d) of the Lemma that follows, after some extra needed work done in parts (a)-(c). The proof of parts (a), (b) is a bit tedious, and it is exactly for their proof that we developed the material in §4.1. In their proof we shall also need to work out some purely topological arguments. The presented arguments do not claim to be the most elegant ones, but at least they fulfill their purpose.

**Lemma 22.** *Let $K, K'$ be diagrams and $p, q$ be base points of them. The following hold:*

*(a) If $_pK = {}_qK'$ then the direct uncrossing decompositions $\mathbb{D}_{K,p}$ and $\mathbb{D}_{K',q}$ are isotopic on the sphere $S^2 = \rho \cup \{\infty\}$ (a subset of $S^3 = \mathbb{R}^3 \cup \{\infty\}$).*

*(b) If the direct uncrossing decompositions $\mathbb{D}_{K,p}$ and $\mathbb{D}_{K',q}$ are isotopic on the sphere $S^2$, then $K$ as a diagram is isotopic on $S^2$ to $K'$ as a diagram (isotopic as sets, and at the corresponding crossings they have the same over/under information).*

*(c) If $_pK = {}_qK'$ then $K \underset{top}{\sim} K'$.*

(d) If $_pK \xrightarrow{^{al}\Omega_{iso}} A \in \Sigma$, then $A = {_q}K'$ for some diagram $K'$ and some base point $q$ of it, so $A \in \Sigma_D$. Moreover, for any diagram $K'$ and base point $q$ of it for which $A = {_q}K'$, it is $K \underset{top}{\sim} K'$.

*Proof.* (a) $_pK = {_q}K'$ implies $\tau(K) = \tau(K')$.

We prove the result by induction on $\tau(K) = \tau(K') = n$.

For $n = 0$, the diagrams $K, K'$ are circles and $_pK = {_q}K'$ are indeed equal to each other for all choices of base points $p, q$ (both diagrams equal to $\varnothing$). The direct decompositions of $K, p$ and $K', q$ are $\mathbb{D}_K = (K, p)^\varnothing$ and $\mathbb{D}_{K'} = (K', q)^\varnothing$. But the bases $(K, p)$ and $(K', q)$ are isotopic on $S^2$ no matter what the orientations of $K, K'$ might be. So $\mathbb{D}_K = (K, p)^\varnothing$ and $\mathbb{D}_{K'} = (K', q)^\varnothing$ are isotopic on $S^2$.

Let the result be true for some $n \geq 0$.

For $n + 1$: Let $K, K'$ be diagrams with $\tau(K) = \tau(K') = n + 1$ and $p, q$ bases on them so that $_pK = {_q}K'$.

Let
$$\mathbb{D}_{K,p} = (K_{n+2}, p, U_1, 1_1, 2_1, 3_1, 4_1, U_2, 1_2, 2_2, 3_2, 4_2, \ldots, U_{n+1}, 1_{n+1}, 2_{n+1}, 3_{n+1}, 4_{n+1})^{(\epsilon_1, \epsilon_2, \ldots, \epsilon_{n+1})}$$
be a direct uncrossing decomposition of $K, p$, and
$$\mathbb{D}_{K',q} = (K'_{n+2}, q, U'_1, 1'_1, 2'_1, 3'_1, 4'_1, U'_2, 1'_2, 2'_2, 3'_2, 4'_2, \ldots, U'_{n+1}, 1'_{n+1}, 2'_{n+1}, 3'_{n+1}, 4'_{n+1})^{(\epsilon'_1, \epsilon'_2, \ldots, \epsilon'_{n+1})}.$$
be a direct uncrossing decomposition of $K', p'$.

We wish to show that $\mathbb{D}_{K,p}$ and $\mathbb{D}_{K',q}$ are isotopic on the sphere $S^2 = \rho \cup \{\infty\}$. So we seek for an isotopy which sends the base of one of them to the base of the other, and we also need to have that the two exponents coincide.

To start observe that the first moves in the sequence of the uncrossings in these two decompositions are $K = K_1, p \xrightarrow[U_1^{\epsilon_1}]{\omega_{1d}} K_2, p$ and $K' = K'_1, q \xrightarrow[U_1^{\prime\epsilon'_1}]{\omega'_{1d}} K'_2, q$. Since $_pK = {_q}K'$ and the crossings $\Delta_1 \in U_1$ and $\Delta'_1 \in U'_1$ each have 1 as one of their labels, it is $\Delta_1 = (1, a)^{\pi_1} = \Delta'_1$ or else $\Delta_1 = (a, 1)^{\pi_1} = \Delta'_1$ for some $a \in \{2, 3, \ldots, 2(n+1)\}$.

Claim 1: $\epsilon_1 = \epsilon'_1$.

Indeed:

Recall that $\epsilon_1$ is the kind-number of $U_1$ in $K, p$. To get its value, we make a trip on $K$ along its orientation starting at $p$ until we arrive at $U_1$ and we continue going through $U_1$ along the first leg of $U_1$ observing if this leg is the over- or under-crossing of the crossing point $\Delta_1$ to get $\epsilon = +1$ or $-1$ respectively. The possibilities are shown in Figure 43. Similarly for $\epsilon'_1$.

If $\Delta' = \Delta = (1, a)^\pi$, then the first leg of $K, p$ and $K', q$ respectively in $U_1, U'_1$ is over the second, so $\epsilon_1 = \epsilon'_1 = +1$. Similarly, if $\Delta' = \Delta = (a, 1)^\pi$, then the first leg of $K, p$ and $K', q$ respectively in $U_1, U'_1$ is under the second, so $\epsilon_1 = \epsilon'_1 = -1$. Hence $\epsilon_1 = \epsilon'_1$ in all cases, and Claim 1 is proved.

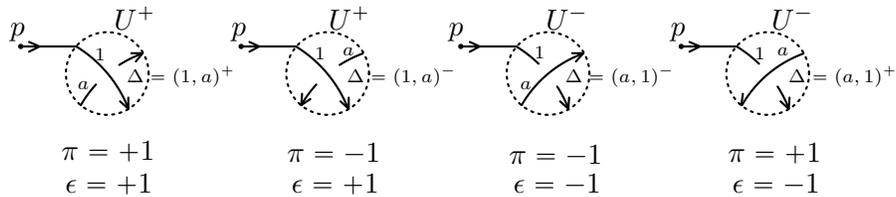

FIGURE 43. Let $U$ be an uncrossing disk of the first crossing $\Delta$ which we meet after the base point $p$ when we move along the orientation of a diagram $K$. The sign of $\Delta = (1 : a)^\pi$ and the kind-number $\epsilon$ of $U^\epsilon$ in $K, p$ are as shown.

Claim 2: It holds $_pK_2 = {_q}K'_2$.

Indeed:

$K_2$ contains all crossings of $K_1$ except $\Delta_1$ which disappears because of the uncrossing move $\omega_{1d}$. Similarly for $K_2'$. The points $p$ and $q$ lie in $K_2, K_2'$ respectively. The labels of the crossings in $K_2, p$ and $K_2', q$ come by some relabeling of the same crossings in $K, p$ and $K', q$. In both cases, the relabeling is done according to the following table:

$$
\begin{array}{lcccccccc}
(K = K_1, p), (K' = K_1', q) & : & 1 & 2 & \ldots & a-1 & a & a+1 & \ldots & 2n+2 \\
(K_2, p), (K_2', q) & : & & a-2 & \underset{\leftarrow}{\ldots} & 1 & & a-1 & \underset{\rightarrow}{\ldots} & 2n
\end{array}
$$

Thus if $\mu$ is the relabeling as a function from the top to the bottom row (1 and $a$ are not in the domain of this function), then any old crossing $\Delta = (x, y)^\pi \in {}_pK = {}_qK'$ other than $\Delta_1$ becomes $\Delta_d = (\mu(x), \mu(y))^\epsilon \in {}_pK_2$ and $\Delta_d' = (\mu(x), \mu(y))^{\epsilon'} \in {}_pK_2'$ for some signs $\epsilon, \epsilon'$. Actually it is $\epsilon = \epsilon'$:

By construction (the position of $p$, the definition of the arms and the orientations of $K_1, K_2$) it is $\gamma_1(K) = -\gamma_1(K_2)$, $\gamma_2(K) = \gamma_2(K_2)$ and similarly $\gamma_1(K') = -\gamma_1(K_2')$, $\gamma_2(K') = \gamma_2(K_2')$. Also, because of the equality ${}_pK = {}_qK'$, the arms $\gamma_1(K), \gamma_1(K')$ contain crossings with the same description (signs, labels and height of labels): the arms contain the crossings with a label $2, 3, \ldots, a-1$ in this order as we traverse them along the orientation on each diagram, and then for each crossing the second label is determined and the heights of the two labels and the sign of the crossing are determined as well, one after the other in both diagrams, and coincide in both diagrams since ${}_pK = {}_qK'$. Similarly for the arms $\gamma_2(K), \gamma_2(K')$.

So if $\Delta$ is a crossing of $\gamma_1(K)$ with itself, then $\Delta'$ is a crossing of $\gamma_1(K')$ with itself, thus both $\Delta_d, \Delta_d'$ are crossings of $-\gamma_1(K)$ with itself. Similarly if $\Delta$ is a crossing of $\gamma_2(K)$ with itself then $\Delta_d$ is a crossings of $\gamma_2(K)$ with itself and $\Delta_d'$ a crossings of $\gamma_2(K)$ with itself. In the second case it is immediate that the signs of $\Delta_d, \Delta_d'$ is the same as the sign of $\Delta$ so they are equal to each other. In the first case the orientations of the over and under arcs in the uncrossing disks of $\Delta, \Delta'$ change simultaneously when we consider them as disks for the same crossing points for $K_d, K_d'$. This makes the signs of $\Delta_d, \Delta_d'$ to be again the same as the sign of $\Delta$, thus again equal to each other.

If finally $\Delta$ is a crossing of $\gamma_1(K)$ with $\gamma_2(K)$ then (as above) $\Delta_d$ is a crossings of $-\gamma_1(K)$ with $\gamma_2(K)$ and $\Delta_d'$ a crossings of $-\gamma_1(K')$ with $\gamma_2(K')$. Then in the uncrossing disks of $\Delta, \Delta'$ the sign of one of the over or under arc changes when the disk is considered for $K_d, K_d'$ around $\Delta_d$ and $\Delta_d'$. This makes the signs of $\Delta_d, \Delta_d'$ to be opposite to the ones of $\Delta, \Delta'$, thus again equal to each other.

So it is always true that $\epsilon = \epsilon'$ as claimed. This implies $\Delta_d = \Delta_d'$. Then ${}_pK_2 = {}_qK_2'$, hence Claim 2 is proved.

Now observe that $\tau({}_pK_2) = \tau({}_pK_2') = n$ and ${}_pK_2 = {}_qK_2'$. According to the induction assumption, the direct uncrossing decompositions of $K_2, p$ and $K_2', q$ are isotopic in $S^2$. According to Lemma 21 these decompositions are respectively:

$\mathbb{D}_{K_2, p} = (K_{n+2}, p, U_2, 1_2, 2_2, 3_2, 4_2, \ldots, U_{n+1}, 1_{n+1}, 2_{n+1}, 3_{n+1}, 4_{n+1})^{(\epsilon_2, \ldots, \epsilon_{n+1})}$

and

$\mathbb{D}_{K_2', q} = (K_{n+2}', q, U_2', 1_2', 2_2', 3_2', 4_2', \ldots, U_{n+1}', 1_{n+1}', 2_{n+1}', 3_{n+1}', 4_{n+1}')^{(\epsilon_2', \ldots, \epsilon_{n+1}')}$.

Since $\mathbb{D}_{K_2', q}, \mathbb{D}_{K_2', q}$ are isotopic in $S^2$, say by some isotopy $F$ we have by definition $(\epsilon_1, \epsilon_2, \ldots, \epsilon_{n+1}) = (\epsilon_1', \epsilon_2', \ldots, \epsilon_{n+1}')$. Combined with $\epsilon_1 = \epsilon_1'$ proved above, we have that the exponents of $D_{K, p}, D_{K', q}$ coincide.

Thus we shall finish with the desired proof of the induction hypothesis for $n + 1$ and the whole induction argument, after we prove Claim 3 that follows. In the proof we'll utilize isotopy $F$ to produce a desired isotopy $G$ of $S^2$ that sends the base of $D_{K, p}$ to that of $D_{K', q}$.

Claim 3: There exists an isotopy $G$ of $S^2$ that sends the base of a direct decomposition $\mathbb{D}_{K, p}$ to that of $\mathbb{D}_{K', q}$.

Indeed:

For $n = 0$ it is $\tau(K) = \tau(K') = n + 1 = 1$, and Figure 44 shows the truth of the claim whenever $_pK = {_qK'} = (1,2)^{-1}$ by exhibiting all possible topological placements of a based diagram on the plane with the single crossing $(1,2)^{-1}$, along with its decompositions.

Now let us observe that the 7-ads $(K_2, p, U_1, 1_1, 2, 3_1, 4_1)$, $(K_2', p', U_1', 1_1', 2', 3_1', 4_1')$, $(K_2'', p'', U_1'', 1_1'', 2'', 3_1'', 4_1'')$ of all diagrams in the Figure are clearly isotopic to each other on $S^2$. These sequences are the bases of the direct decompositions of the diagrams, thus we have the result for diagrams with the single crossing $(1,2)^{-1}$. Similarly the claim holds for the diagrams with a single crossing $(1,2)^{+1}$ or $(2,1)^{-1}$ or $(2,1)^{+1}$.

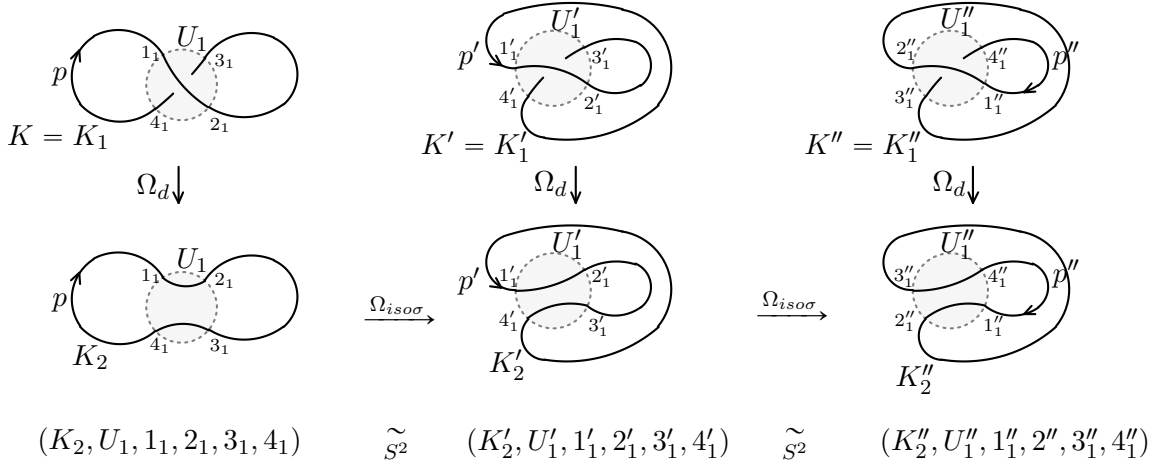

FIGURE 44. All possible topological placements on the plane of a based diagram with the single crossing $(1,2)^{-1}$, along with the corresponding direct uncrossing decompositions. The bases of all these decompositions are clearly isotopic in $S^2$.

For $n > 0$ it is $\tau(K) = \tau(K') = n + 1 > 1$ and there exists at least two crossing points and corresponding uncrossing disks $U_1, U_2, \ldots$ and $U_1', U_2', \ldots$ for the two diagrams.

All topological placements on the plane for the first two crossings in a direct uncrossing decomposition of a diagram $L, p$ with at least two crossings, are given in Figure 45. According to the given calculations, no two of the placements share all of $\pi_1, \epsilon_1, \pi_2, \epsilon_2$.

Since $_pK_1 = {_qK_1'}$, the first crossing points of the two based diagrams have the same description, thus the same sign $\pi_1 = \pi_1'$ and then the kind-numbers $\epsilon, \epsilon'$ of $U_1, U_1'$ in $K_1, p$ and $K_1', q$ are the same (argue as in Claim 1).

Also, the uncrossing disks $U_2$ and $U_2'$ are the ones which contain the first crossing points of $K_2, p$ and $K_2', q$. These points have to have the same description in $_pK_2$ and $_qK_2'$ (recall $_pK_2 = {_qK_2'}$, and the points both contain 1 as a label), say $\Delta = (1, b)^\pi, \Delta' = (1, b)^\pi$ or $\Delta = (b, 1)^\pi, \Delta' = (b, 1)^\pi$. Then the kind-numbers $\epsilon, \epsilon'$ of $U_2, U_2'$ in $K_2, p$ and $K_2', q$ are the same (argue as in Claim 1).

So for both $K, p$ and $K', q$ the same case among I-XVI of Figure 45 holds. The argument that follows shows the desired isotopy of the bases of $\mathbb{D}_{K,p}, \mathbb{D}_{K',q}$ when case I holds, but it can be repeated for all other cases $II - XVI$. So the completion of the argument will at the same time complete the inductive proof of part (b).

Let us assume then that case I holds for $K, p$ and $K', q$.

We pay attention to the first two uncrossing disks $U_1, U_2$ and $U_1', U_2'$ in the direct uncrossing decompositions of $K, p$ and $K, q'$. Let $\Delta_1 = (1 : a) \in {_pK}$ and $\Delta_1' = (1 : a) \in {_qK'}$ be the first intersection points of $K, p$ and $K', q$. We also pay attention to the disks $U_m, U_m'$ of Figure 46 which contain the crossing points in the two based diagrams whose one label is $a + 1$. $U_m$ and $U_m'$ are the first uncrossing disks that we meet after the fourth node of $U_1$ and $U_1'$ as we move along the orientations of $K, K'$. The use of a common index $m$ for the two disk is justified by the part (f) of Lemma 21. Let us notice that it may happen that these disks coincide with $U_2, U_2'$, but this does not affect the arguments that follow.

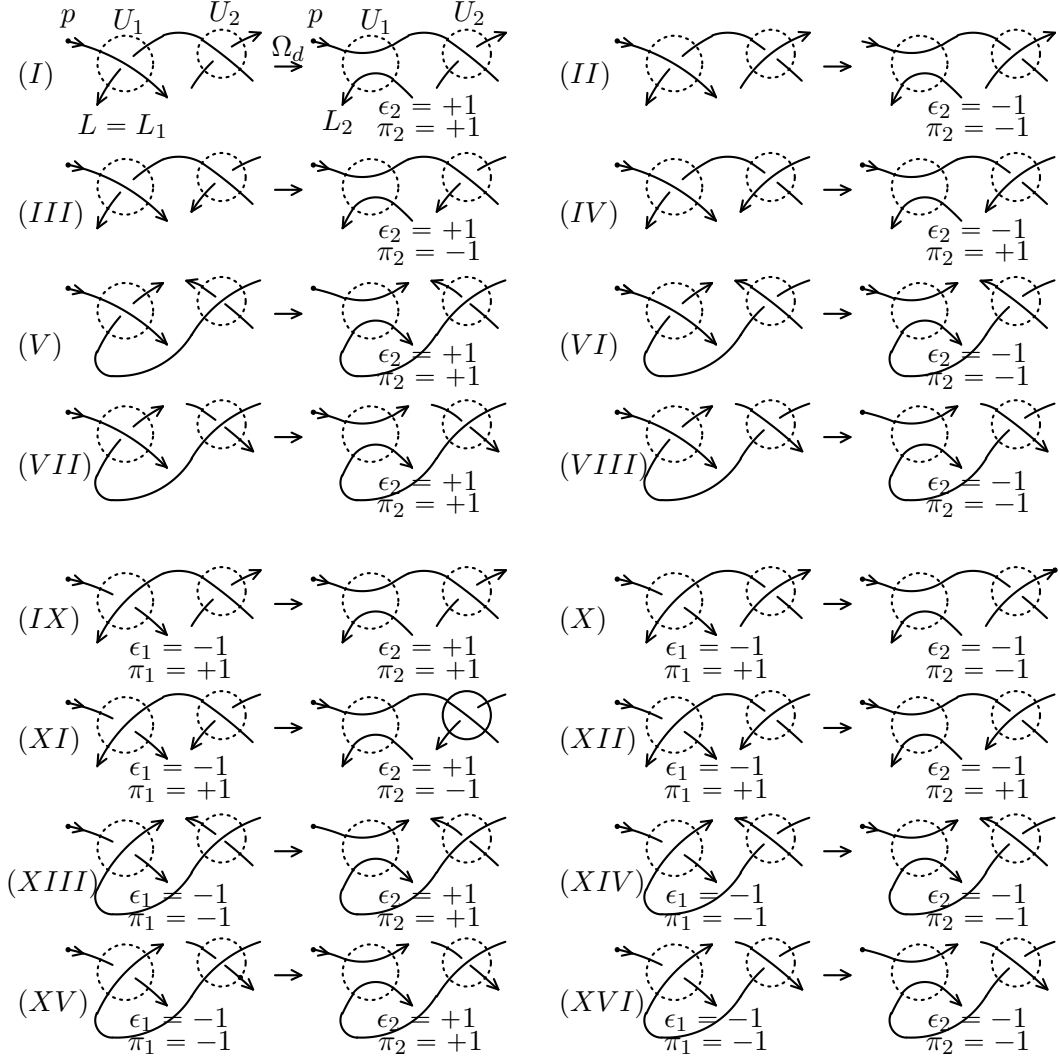

FIGURE 45. Depicted are all possible topological placements on the plane of the crossings appearing as first and second in the direct uncrossing decomposition of a based diagram $L, p$ with at least two crossings. For each case, the sign $\pi_1$ of the first crossing point and also the kind-number $\epsilon_1$ of the first uncrossing disk $U_1$ in $L = L_1, p$ are calculated. Also shown here is the first uncrossing move $L = L_1, p \xrightarrow[U_1]{\Omega_d} L_2, p$ in the decomposition. For each case the sign $\pi_2$ of the second crossing and also the kind-number $\epsilon_2$ of the second uncrossing disk $U_2$ in $L_2, p$ are calculated.

In Figure 46 we denote $t_1, t_2, t_3, t_4$ respectively the first, second, third and fourth node of $U_1$ in $K_2$. We denote $A$ the body of $U_1$ as a nice disk of $K_2, p$ and $A_1, A_2$ its two ears, and $a, b, c, d$ the first,second,third and fourth node of $U_2, p$ as a crossing disk of $K_2, p$. We denote $x$ the first node of $U_m$ as a crossing disk of $K_2, p$. We call $B$ a 2-disk containing $A$ with boundary the arcs $\sigma_1 = pa$ and $\sigma_3 = dx$ of $K_2$, the arc $\sigma_2 = ad$ of $\partial U_2$ and an arc $\sigma_4$ joining $p$ and $x$ and avoiding the rest of $K$ and the $U_i$'s.

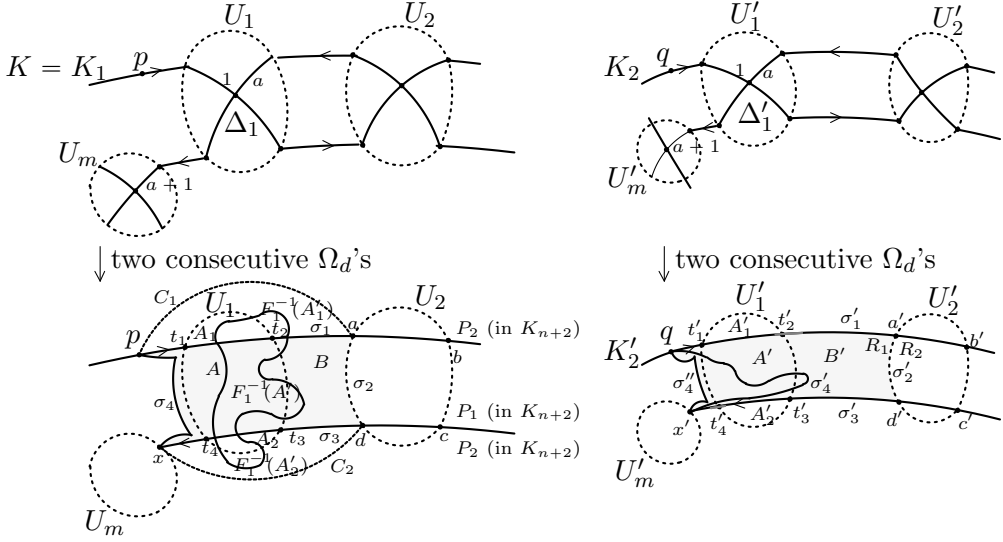

FIGURE 46. The schema shows two based diagrams $K, p$ and $K', q$ for which $_pK = {_q}K$, along with the resulting diagrams $K_2, p$ and $K'_2, q$ after we perform the first two uncrossing moves in the direct uncrossing decompositions of $K, p$ and $K', q$. Regarding the relative position of $U_1, U_2$ and $U'_1, U'_2$ we assume here that case I of Figure 45 holds. $\Delta_1, \Delta'_1$ are the first crossing points in the two decompositions. $U_m, U'_m$ are the first uncrossing disks we meet moving along the orientations of the diagrams, starting at the fourth node of $U_1, U'_1$ respectively. The bodies $A, A'$ of $U_1, U'_1$ in $K_2, p$ and $K'_2, q$ respectively, are parts of the 2-disks $B, B'$ depicted in gray color. These discs remain unaffected by the whole sequence of uncrossing moves for $K$ and $K'$, thus each one of them remains in one region with respect to the final uncrossing diagram $K_{n+1}$ and $K'_{n+1}$ respectively. The disk $C = B \cup C_1 \cup C_2$ is constructed so that an isotopy fixed point-wise in the exterior of $C$, brings $(A, A_1, A_2)$ onto $(A', A'_1, A'_2)$ and $(\sigma_1, \sigma_2, \sigma_3)$ onto $(\sigma'_1, \sigma'_2, \sigma'_3)$. $R_1, R_2$ indicate angles near $a'$: $R_2$ lies in $U'_2$, whereas $R_1$ lies outside $U'_2$. $R_2$ is connected to $A'$ with arcs not intersecting $K'_2$. After the decomposition is completed, the resulting 2-disk $K_{n+2}$ will contain the top and bottom of the obvious horizontal zone (at the middle of the Figure) in one of its interior or exterior $P_1, P_2$, say in $P_2$, and it will contain the horizontal zone in $P_1$.

Similarly for the regions, arcs and points $t'_1, t'_2, t'_3, t'_4, A', A'_1, A'_2, a', b', c', d', \sigma'_1, \sigma'_2, \sigma'_3$ in Figure 46. $\sigma'_4, B'$ and some more regions and arcs will be defined in a moment.

The nodes $a, d$ are the first and fourth nodes of $K_2, p$ and also $a', d'$ are the first and fourth nodes of $K'_2, q$. By part (g) of Lemma 21 we know that $a, a'$ get the same name $\bar{a} = \bar{a}' \in \{1, 2, 3, 4\}$ as nodes in $K_{n+2}, p$ and $K'_{n+2}, q$, and also $d, d'$ get the same name $\bar{d} = \bar{d}' \in \{1, 2, 3, 4\}$ as nodes in $K_{n+2}, p$ and $K'_{n+2}, q$. The same also holds for the order of $b$ and $c$ in $K_{n+2}, p$ compared with the orders of $b'$ and $c'$ in $K'_{n+2}, q$. Similarly, the points $x$ and $x'$ are the first nodes of $U_m$ and $U'_m$ in $K, p$ and $K', q$, thus they get the same name $\bar{x} = \bar{x}' \in \{1, 2, 3, 4\}$ as nodes in $K_{n+2}, p$ and $K'_{n+2}, q$.

Now recall the isotopy $F$ on $S^2$ between $\mathbb{D}_{K_2,p}$ and $\mathbb{D}_{K'_2,q}$, and let $F_1$ the last moment of $F$. By definition $F_1$ maps $p$ to $q$, the nodes $\bar{a}, \bar{d}$ of $U_2$ to $\bar{a}, \bar{d}$ of $U'_2$ and the node $\bar{x}$ of $U_m$ to the node $\bar{x}'$ of $U'_m$. So it maps the points $p, a, d, c$ onto the points $q, a', d', x'$. Since $F_1$ maps the oriented diagram $K_{n+1}$ onto the oriented diagram $K'_{n+1}$, it follows the it maps the arc $\sigma_1$ onto $\sigma'_1$ and the arc $\sigma_3$ onto $\sigma'_3$. Since $F_1$ maps $U_2$ onto $U'_2$, thus $\partial U_2$ onto $\partial U'_2$, it follows that $F_1$ sends the arc $\sigma$ of $\partial U_2$ onto $\sigma_1'$ of $\partial U_2'$ (because of the orders of the nodes on the two boundaries). Let $\sigma_3' = F_1(\sigma_3)$ and $B' = F_1(B)$.

Since $F_1$ is a homeomorphism, $B'$ is a 2-disk and $\partial B' = \sigma_1' \cup \sigma_2' \cup \sigma_3' \cup \sigma_4'$. $B'$ contains points in the region $R_1$ in Figure 46 near $a'$ or else points in the region $R_2$ near $a'$ in the interior of $U_2'$. The second case cannot hold: if it did, by Jordan's Theorem for curves on a plane (or a 2-sphere) the interior of $B'$ should contain some point $y'$ on the boundary of $U_2'$. Then there would be some point $y$ in the interior of $B$ mapped by $F_1$ onto $y$. By definition of $F$, the points of $S^2$ that map onto $\partial U_2'$ are those on $\partial U_2$. Thus $y$ should lie on $\partial U_2$. Then the interior of $B$ would intersect $\partial U_2$, a contradiction since by the construction of $B$ its interior does not have common points with $\partial U_2$.

If $\sigma_2'$ intersects $A'$, we modify $F$ locally in a tubular neighborhood of $\sigma_4$ so as its image $s_4''$ avoids $A'$. We continue calling $F$ the whole isotopy of $S^2$. Calling once again as $F_1$ the last moment of $F$, it is $F_1(B) = B'$, a 2-disk which contains $A'$, whose boundary is $\partial B' = \sigma_1' \cup \sigma_2' \cup \sigma_3' \cup \sigma_4''$ and whose interior does not intersect $K_2'$ or the boundaries of the $U_i$'s with $i \geqslant 2$. Also, this $F$ is an isotopy which coincides with the old one outside the interior of $B$ so, it is an isotopy on $S^2$ which continues to send the base of the direct uncrossing decomposition of $K_2, p$ to the base of the direct uncrossing decomposition of $K_2', q$.

Let $P_1, P_2$ be the two 2-disks of $S^2$ with boundaries the circle $K_{n+1}$. And let $B$ lies in $P_1$.

By construction $F_1^{-1}(A') \subseteq B$, whereas $F_1^{-1}(A_1'), F_1^{-1}(A_2')$ are small 2-disks with an arc of their boundaries on the arcs $\sigma_1$ and $\sigma_3$ of $\partial B$ respectively. The two disks $F_1^{-1}(A_1'), F_1^{-1}(A_2')$ lie on $P_2$. Let $C_1$ be a 2-disk in $P_2$, intersecting $B$ in $\sigma_1$ and containing $A_1, F_1^{-1}(A_1')$, small enough so that it does not have common points outside $\sigma_1$ with $K_2$ and the $U_i$'s. Similarly, let $C_2$ be a 2-disk in $P_2$, intersecting $B$ in $\sigma_3$ and containing $A_2, F_1^{-1}(A_2')$, small enough so that it does not have common points outside $\sigma_3$ with $K_2$ and the $U_i$'s, and also small enough so that it has no common points with $C_1$. Then $C = B \cup C_1 \cup C_2$ is a 2-disk in $S^2$.

By Schonflies Theorem for disks on a plane (or 2-sphere) we can trivially construct an isotopy $G$ of $S^2$ fixed point-wise outside $C$ which sends $B, C_1, C_2$ onto $B, C_1, C_2$ so that $\sigma_1$ goes onto $\sigma_1$, also $\sigma_3$ onto $\sigma_3$, and $A, A_1, A_2$ onto $F^{-1}(A), F^{-1}(A_1), F^{-1}(A_2)$ respectively. Then $H = G \cup F$ is an isotopy of $S^2$ which similarly to $F$, sends the base of the direct uncrossing decomposition $\mathbb{D}_{K_2, p}$ to the base of the direct uncrossing decomposition $\mathbb{D}_{K_2', q}$. Moreover, $p, t_1, t_2, a$ in $\sigma_1$ have the order $(p, t_1, t_2, a)$, and $q, t_1', t_2', a'$ in $\sigma_1'$ have the order $(q, t_1', t_2', a')$ and by construction $H$ sends $t_1, t_2$ onto $t_1', t_2'$ respectively. Similarly, $H$ sends $t_3, t_4$ onto $t_3', t_4'$ respectively. By part (h) of Lemma 21 we know that the first, second, third and fourth node of $U_1$ in $K_{n+2}, p$ are the first, second, third and fourth node of $U_1$ in $K_2, p$, that is, the points $t_1, t_2, t_3, t_4$ respectively. Similarly, we know that the first, second, third and fourth node of $U_1$ in $K_{n+2}, p$ are the points $t_1', t_2', t_3', t_4'$ respectively.

So except possibly for $U_1$, the isotopy $H$ sends all spaces in the base of the direct uncrossing decomposition $\mathbb{D}_{K,p}$ to the base of the direct uncrossing decomposition $\mathbb{D}_{K_2', q}$. But $H$ sends the four points $t_1, t_2, t_3, t_4$ on $\partial U_1$ to the correct four points $t_1', t_2', t_3', t_4'$ on $\partial U_1'$. Perturbing $H$ in a small 2-disk which contains $U_1$ and is contained in $C$ (so we fix it point-wise outside this 2-disk) we easily manage to have the resulting isotopy $h$ of $S^2$ send $\partial U_1$ to $\partial U_1'$, also $t_1, t_2, t_3, t_4$ onto $t_1', t_2', t_3', t_4'$ respectively and $K_{n+1}$ onto $K_{n+1}'$. Then $h \circ H$ is an isotopy of $S^2$ that sends the base of the direct uncrossing decomposition $\mathbb{D}_{K,p}$ to the base of the direct uncrossing decomposition $\mathbb{D}_{K_2', q}$. And we are done!

(b) The isotopy between the decompositions implies there exists the same number of uncrossing neighborhoods for the two diagrams, thus the order of $K, K'$ is the same, say $n$.

For $n = 0$ the result is immediate.

For $n > 1$, let:
$$D_{K,p} = (K_{n+1}, p, U_1, 1_1, 2_1, 3_1, 4_1, U_2, 1_2, 2_2, 3_2, 4_2, \ldots, U_n, 1_n, 2_n, 3_n, 4_n)^{(\epsilon_1, \epsilon_2, \ldots, \epsilon_n)}$$
be a direct uncrossing decomposition of $K, p$, and:

$$D_{K',q} = (K'_{n+1}, q, U'_1, 1'_1, 2'_1, 3'_1, 4'_1, U'_2, 1'_2, 2'_2, 3'_2, 4'_2, \ldots, U'_n, 1'_n, 2'_n, 3'_n, 4'_n)^{(\epsilon'_1, \epsilon'_2, \ldots, \epsilon'_n)}.$$

be a direct uncrossing decomposition of $K', q$, and let $F$ be an isotopy of $S^2$ which sends $D_{K,p}$ onto $D_{K',q}$. In particular, by the definition of $F$ it is $(\epsilon_1, \epsilon_2, \ldots, \epsilon_n) = (\epsilon'_1, \epsilon'_2, \ldots, \epsilon'_n)$.

Also, let $(K_1 = K, K_2, \ldots, K_{n+1})$ and $(K'_1 = K', K'_2, \ldots, K'_{n+1})$ be the sequence of the diagrams in the two direct uncrossing decompositions.

By Lemma 21, we know that for each $i = 1, 2, \ldots$ the points $1_i, 2_i, 3_i, 4_i$ are not only the first, second, third and fourth node respectively of $U_i$ in $K_{n+1}$, but they are also the first, second, third and fourth node respectively of $U_i$ in $K_{i+1}$. Then by the definition of an uncrossing move, we have that $1_i, 3_i$ are respectively the first and second endpoint of the first leg of $U_i$ in $K_i, p$, and $2_i, 4_i$ are respectively the first and second endpoint of the second leg of $U_i$ in $K_i, p$ (here we consider $\ell_{i1}, \ell_{i2}$ as directed arcs of $K_i$), thus $1_i, 3_i, 2_i, 4_i$ are the first, second, third and fourth node of $U_i$ in $K_i, p$. We also know that the arcs $\ell_{i1}, \ell_{i2}$ remain unaltered in $K = K_1, K_2, \ldots K_{i-1}$ since up to that point of the uncrossing sequence only the disks $U_1, U_2, \ldots, U_{i-1}$ were used. And we know that the sign $\epsilon_i = +1$ or $-1$ informs us about if respectively $\ell_{i1}$ is over $\ell_{i2}$ or $\ell_{i1}$ is under $\ell_{i2}$ whenever these arcs are considered in $K_i$. Then the same happens in $K = K_1$ as well.

So let us construct a diagram $\overline{K}$ from the circle $K_{n+1}$ by keeping everything fixed outside the $U_i$'s (orientations are fixed as well), and altering $K_{n+1}$ inside the $U_i$'s so that $1_i, 3_i$ are joined by an arc $\overline{\ell}_{i1}$, and $2_i, 4_i$ are joined by an arc $\overline{\ell}_{i2}$ so that they intersect transversely in a single point. We give to $\overline{\ell}_{i1}, \overline{\ell}_{i2}$ orientations so that their first endpoints are $1_i$ and $2_i$ respectively. We also give to the crossing point the over/under information $\epsilon_i$ ($\ell_{i1}$ is over $\ell_{i2}$ exactly when $\epsilon_i = 1$).

Since $K_{n+1}$ is the final diagram in the uncrossing decomposition $D_{K,p}$ of $K$ and because of the remarks made above, $K$ and $\overline{K}$ differ only on their points inside the $U_i$'s, but in any case they agree on the over/under information for their intersection points in each $U_i$. But their first legs $\ell_{i1}, \overline{\ell}_{i1}$ have the same first and second endpoints, and similarly for their second legs $\ell_{i2}, \overline{\ell}_{i2}$. Since the intersection of the legs in the two cases is transverse in a single point, it follows that we can send the pair $(\ell_{i1}, \ell_{i2})$ onto the pair $(\overline{\ell}_{i1}, \overline{\ell}_{i2})$ by an isotopy $G_i$ of $S^2$ keeping everything fixed point-wise outside $U_i$ (also the boundary of $U_i$ as well). Then $G = G_n \circ \cdots \circ G_1$ is an isotopy of $S^2$ that brings $K$ point-wise onto $\overline{K}$. Since the over/under information of $K$ and $\overline{K}$ at their corresponding points are the same, and as isotopies retain this information, $G$ brings $K$ onto $\overline{K}$ as diagrams.

In exactly the same manner we construct a diagram $\overline{K'}$ and an isotopy $H$ of $S^2$ that brings $K'$ onto $\overline{K'}$ as diagrams. Let $(\ell'_{i1}, \ell'_{i2})$ be the first and second leg which we consider for the disk $U'_i$ in $K'_{i+1}$ in the construction of this diagram $\overline{K'}$.

Recall that $F$ sends $(K_{n+1}, p, U_1, 1_1, 2_1, 3_1, 4_1, U_2, 1_2, 2_2, 3_2, 4_2, \ldots, U_n, 1_n, 2_n, 3_n, 4_n)$ onto $(K'_{n+1}, q, U'_1, 1'_1, 2'_1, 3'_1, 4'_1, U'_2, 1'_2, 2'_2, 3'_2, 4'_2, \ldots, U'_n, 1'_n, 2'_n, 3'_n, 4'_n)$.

So $F$ sends the first and second endpoints of $(\ell_{i1}, \ell_{i2})$ onto the first and second endpoints of $(\ell'_{i1}, \ell'_{i2})$. Thus we can modify $F$ inside the $U_i$'s so that it sends $(\ell_{i1}, \ell_{i2})$ onto $(\ell'_{i1}, \ell'_{i2})$. Call $f$ the modified isotopy. Then $f$ brings $\overline{K}$ point-wise onto $\overline{K'}$. Since $\epsilon_i = \epsilon'_i$, $f$ brings $\overline{K}$ onto $\overline{K'}$ as diagrams.

Hence $H^{-1} \circ f \circ G$ is an isotopy of $S^2$ which brings $K$ onto $K'$ as diagrams, as wanted.

(c) According to part (a) of the Lemma, the uncrossing decompositions $\mathbb{D}_{K,p}$ and $\mathbb{D}_{K',q}$ are isotopic in the sphere $S^2$. Then according to part (b), $K'$ is isotopic to $K$ on $S^2$. By Theorem 3, we have then that $K \underset{top}{\sim} K'$ as wanted.

(d) If $_pK \xrightarrow{^{al}\Omega_{iso}} A$, then by the definition of the $^{al}\Omega_{iso}$ move we have $A = {_pK}$ and the first part of the Lemma holds, setting for example $K' = K$ and $q = p$.

Let now $A = {_qK'}$ for a diagram $K'$ and a base point $q$ of it. Then $_pK = A = {_qK'}$ and we get the desired result by part (c) of the Lemma. $\square$

**Remark 7.** *As already mentioned in Remark 6, the proof of part (d) of the last Lemma contains as a byproduct the way to reconstruct a diagram $K$ from a direct uncrossing decomposition of it.*

4.3. $^{al}\Omega_\beta$ **mirroring.** As in the special case of the $^{al}\Omega_{iso}$ moves, the $^{al}\Omega_\beta$ moves actually happen between elements of $\Sigma_D$, and moreover any such move is mirrored by topological moves between diagrams. The following Lemma 22 asserts these facts:

**Lemma 23.** *If $K$ is a diagram, $p$ a base point of $K$ and $_pK \xrightarrow{^{al}\Omega_\beta} A \in \Sigma$, then $A = {_qK'}$ for some diagram $K'$ and base point $q$ of it, thus $A \in \Sigma_D$. Moreover $K \underset{top}{\sim} K'$.*

*Proof.* If $_pK \xrightarrow{^{al}\Omega_\beta} A \in \Sigma$ then by part (d) of Lemma 16 we have $A = {_{p'}K}$ for some choice $p'$ of base point for $K$, thus $A \in \Sigma$. Now if $A = {_qK'}$ for some diagram $K'$ and base point $q$ of it, then $_{p'}K = {_qK'}$ and by part (c) of Lemma 22 we have $K \underset{top}{\sim} K'$. □

4.4. $\Omega_{1\sigma}^-, \Omega_{2\sigma}^-, \Omega_{2\mu}^-, \Omega_{3\mu}$ **long-range moves, and** $\Omega_{1\gamma}^-, \Omega_{2\gamma}^-, \Omega_{3\gamma}$ **generic moves.** This section is auxiliary to the next three sections where we prove that $^{al}\Omega_1, ^{al}\Omega_2$ and $^{al}\Omega_3$ moves are mirrored by topological ones. To facilitate the proof, we define here a bunch of long-range topological moves generalizing the usual Alexander-Briggs-Reidemeister moves. We do not define any similarly exotic version of $\Omega_1^+, \Omega_2^+$ moves since as we shall see, these are adequate to express the mirroring of their algebraic counterparts (§4.5,4.6).

First let us notice that if in a diagram $K$ the union of a sequence of some basic arcs $e_1, e_2, \ldots, e_k$ forms a simple closed curve $\gamma$ on the plane (each arc shares its endpoints with exactly two of the other arcs in the sequence) then by Jordan's Theorem there exists an interior $\Pi_1$ and an exterior $\Pi_2$ of $\gamma$ on the plane. And $K - \gamma$ lies in $\Pi_1 \cup \Pi_2$. So $K - \gamma$ is either contained entirely in one of $\Pi_1, \Pi_2$ or else it has non-empty parts in both of $\Pi_1, \Pi_2$. Whenever $\gamma$ is an 1-gon, 2-gon or 3-gon on which an $\Omega_1^-, \Omega_2^-$ or $\Omega_3$ move is performed, then $K - \gamma$ lies always in the exterior $\Pi_2$. Below we define moves in the rest of the cases in regards to the location of $K - \gamma$.

**Definition 23.** *Let $K$ be a diagram.*

*We call as a generic 1-gon of $K$ any basic arc $e$ of $K$ for which the two endpoints coincide. We call as a generic 2-gon $e, f$ of $K$ with $e$ over $f$ (and $f$ under $e$), any two basic arcs $e, f$ of $K$, none of which is a 1-gon and which share their two distinct endpoints, so that $e$ is over $f$ in both of them. And we call as a generic 3-gon $e, f, g$ of $K$ with $e$ as the over-arc, $f$ as the middle arc and $g$ as the under-arc, any three distinct basic arcs $e, f, g$ of $K$, none of which is a 1-gon and no two of which form a 2-gon, so that each shares an endpoint with another arc and the second endpoint with the third arc, and so that $e$ is over $f, g$ at its endpoints, whereas $g$ is under $e, f$ at its endpoints.*

- *Let $\Pi_1$ be the interior and $\Pi_2$ the exterior on the plane, of a generic 1-gon $e$. All possible topological settings (up to reflections) are shown in Figure 47. Then:*

    *If $K - e$ lies in $\Pi_1$ we call as an $\Omega_1^-$ move on the sphere $S^2 = \rho \cup \{\infty\}$ the alteration of $K$ to $K'$ as shown in Figure 47. The name comes from the fact that we can think of this move as a sequence $K \xrightarrow{\Omega_{iso\sigma}} K_0 \xrightarrow{\Omega_1^-} K'$ on the sphere $S^2 = \rho \cup \{\infty\}$. We denote the move as $\Omega_{1\sigma}^-$, we write $K \xrightarrow{\Omega_{1\sigma}^-} K'$, and we say the move is performed on the generic 1-gon $e$.*

    *A generic $\Omega_1^-$ move is any one of $\Omega_1^-, \Omega_{1\sigma}^-$ and we denote it as $\Omega_{1\gamma}^-$.*

- *Let $\Pi_1$ be the interior and $\Pi_2$ the exterior on the plane, of a generic 2-gon $e, f$. All possible topological settings (up to reflection) are shown in Figure 48. Then:*

    *If $K - (e \cup f)$ lies in $\Pi_1$ we call as an $\Omega_2^-$ move on the sphere the alteration of $K$ to $K'$ as shown in Figure 48 (a). The name comes from the fact that we can think of the move as*

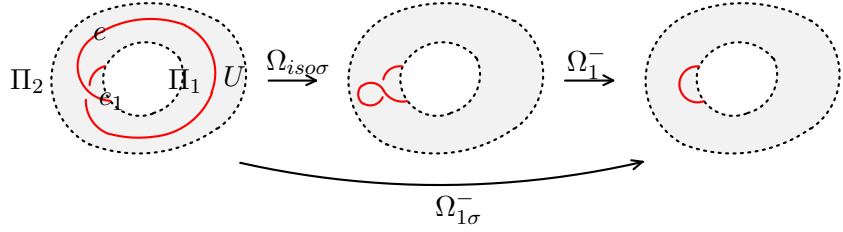

FIGURE 47. The $\Omega_{1\sigma}$ move. $e$ is a loop in the interior of $U$. $e_1$ is $e$ together with the two legs touching the boundary of $U$.

a sequence $K \xrightarrow{\Omega_{iso\sigma}} K_0 \xrightarrow{\Omega_2^-} K'$ of moves on the sphere. We denote the move as $\Omega_{2\sigma}^-$, we write $K \xrightarrow{\Omega_{2\sigma}^-} K'$, and we say the move is performed on the generic 2-gon $e, f$.

If $K - (e \cup f)$ lies in both $\Pi_1, \Pi_2$ then we call as an $\Omega_{2\mu}^-$ move the alteration of $K$ to $K'$ as shown in Figure 48 (b). We call $\Omega_{2\mu}^-$ as a big $\Omega_2^-$ move. We write $K \xrightarrow{\Omega_{2\mu}^-} K'$ and we say it is performed on the generic 2-gon $e, f$.

A generic $\Omega_2^-$ move is any one of $\Omega_2^-, \Omega_{2\sigma}^-, \Omega_{2\mu}^-$ and we denote it as $\Omega_{2\gamma}^-$.

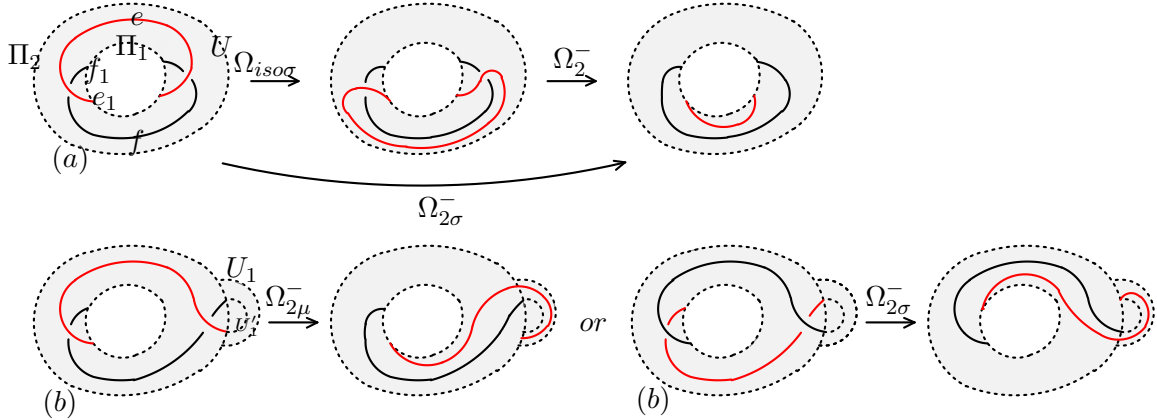

FIGURE 48. The $\Omega_{2\sigma}, \Omega_{2\mu}$ moves.

- Let $\Pi_1$ be the interior and $\Pi_2$ the exterior on the plane, of a generic 3-gon $e, f, g$. If $K - (e \cup f \cup g)$ lies in $\Pi_1$ or in both $\Pi_1, \Pi_2$ then all possible topological placements (up to reflections) are shown in Figure 49. Then:

  We call as an $\Omega_{3\mu}$ move the alteration of $K$ to $K'$ as shown in each one of the cases. We call $\Omega_{3\mu}$ as a big $\Omega_3$ move. We write $K \xrightarrow{\Omega_{3\mu}} K'$.

  A generic $\Omega_3$ move is any one of $\Omega_3, \Omega_{3\mu}$ and we denote it as $\Omega_{3\gamma}$.

In all cases the annulus $U = [-1, 1] \times \gamma$ is a tubular (bicolar) neighborhood on the plane of the simple closed curve $\gamma$ of the $i$-gon. $\Pi_1, \Pi_2$ are the interior and exterior of $\gamma$ on the plane. $U$ is small enough to intersect $K$ in arcs $e_1$ or $e_1, f_1$ or $e_1, f_1, g_1$ depending on the case, slightly larger than $e, f, g$ containing no more intersection points on them, as is actually expected from the usual $\Omega_1^-, \Omega_2^-, \Omega_3$ cases.

In case (b) of the $\Omega_{2\mu}^-$ move the 2-disk $U_1'$ is outside $U$ and contains the part of $K - \gamma$ in $\Pi_2$ (such a disk exists since $K - \gamma$ is bounded). $U_1$ is a slightly bigger 2-disk outside $U$ so that $\partial U_1'$ touches $\partial U_1'$ only in a small arc on $\partial U$. We can similarly choose the disks $U_1', U_1$ inside $P_1$ (we do not depict this in the Figure).

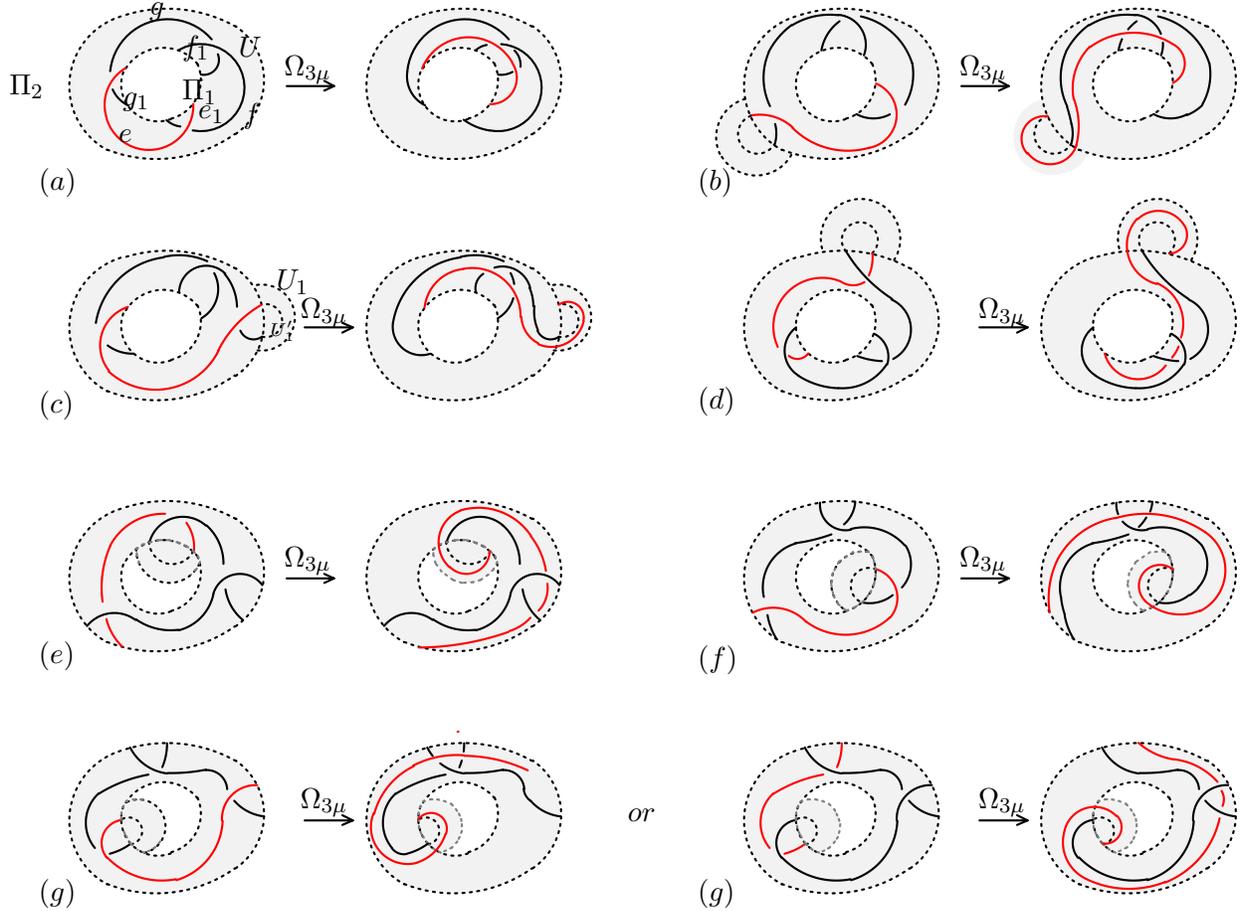

FIGURE 49. The $\Omega_{3\mu}$ move.

In a similar way we consider the 2-disks $U_1', U_1$ in cases (b)-(g) of the $\Omega_{3\mu}$ move. We choose them so that they contain an endpoint, say $x$, of the top or of the bottom arc, whenever the boundary component of $U$ on which $x$ lies is touched upon only by one of the other two arcs. In case (g) each one of the top and bottom arcs satisfies the conditions, thus there exist two versions of the move (both of them depicted).

The moves are defined by altering the red arc of $K$ as shown. In more detail:

For an 1-gon, we delete the part $e$ of $e_1$, and if we are in the smooth case we also smooth the remaining part of $e_1$ at its old crossing point.

For a 2-gon the arc that changes is the red arc. In case (a) of Figure 48 the alteration comes also as the result of an isotopy move on the 2-sphere $S^2 = \rho \cup \{\infty\}$ which moves $e_1$ inside the 2-disk $\Pi_2 \cup \{\infty\}$ to a position on which a usual $\Omega_2^-$ move can be performed. The position of the red arc after the move is outside $U_1'$ and inside $U \cup U_1$. Similarly for the not depicted variation where $U_1', U_1$ are inside $\Pi_1$.

For a 3-gon the arc that changes is the red one. Its position after the move is inside $U \cup U_1$ and outside $U_1'$.

*If we insist to assign a 2-disk for each move as in the case of the usual moves, then we can say that for the* $\Omega_{1\sigma}^-, \Omega_{2\sigma}^-$ *moves the 2-disk is* $U \cup \Pi_1$, *and for the* $\Omega_{2\mu}^-, \Omega_{3\mu}$ *moves the 2-disk is* $\Pi \cup U \cup U_1$.

*In exact analogy to what happens in the usual* $\Omega_1^-, \Omega_2^-, \Omega_3$ *moves, the final position of the arcs* $e_1, f_1, g_1$ *inside the 2-disk of the move for all moves here is not exact but considered up to isotopy of the disk.*

The location of the disks $U_1', U_1$ is dictated by our wish to use them delete crossings (for $\Omega_{2\mu}$ moves) or relocate crossings (for $\Omega_{3\mu}$ moves), and at the same time produce a diagram that is equivalent to the given one (Lemma 24 below). If we put them on the side of $U$ that contains the end points of all three arcs of a (generic) 3-gon, then arcs of $K - \gamma$ joining the endpoints of the other two arcs may go into the way, prohibiting the movement of our chosen top or bottom arc.

An obvious remark is that a generic i-gon ($i = 1, 2, 3$) is not necessarily an $i$-gon (Definitions 16, 18, 19) since the simple closed curve formed by the arcs of the $i$-gon can contain points of the diagram $K$ in its interior on the plane. Of course the converse is true: the i-gons are generic igons as well. Also, it is clear that we can avoid mentioning the 2-sphere altogether in the last Definition.

The $\Omega_{1\sigma}^-, \Omega_{2\sigma}^-$ moves clearly do not change the equivalence class of the given diagram $K$ (by Theorem 3). Also clearly, the choice of the tubular neighborhood $U$ or of the 2-disks $U_1, U_2$ does not alter the equivalence class of the result of an $\Omega_{2\mu}^-$ or an $\Omega_{3\mu}$ move (the results for any two choices are isotopic on the plane), thus regarding the topological equivalence of diagrams these new moves are well-defined. As we shall see in a moment (Lemma 24), these moves also do not change the equivalence class of the given diagram. So all the exotic moves we define here are legitimate, that is, the result of each one in any given diagram is expressed by a sequence of the usual Alexander-Briggs-Reidemeister moves and of isotopies. Since the $\Omega_1^-, \Omega_2^-, \Omega_3$ moves are special cases of the $\Omega_{1\gamma}^-, \Omega_{2\gamma}^-, \Omega_{3\gamma}$ moves we get:

**Corollary 2.** *The moves* $\Omega_{iso}, \Omega_1^+, \Omega_{1\gamma}^-, \Omega_2^+, \Omega_{2\gamma}^-, \Omega_{3\gamma}$ *generate the same equivalence classes of diagrams as the usual moves* $\Omega_{iso}, \Omega_1^+, \Omega_1^-, \Omega_2^+, \Omega_2^-, \Omega_3$.

**Lemma 24.** *Let $K$ be a diagram with $\tau(K) = n + i$, $n \geqslant 0$, $i = 1, 2, 3$. And let $\gamma$ be the simple closed curve formed by the union of arcs of a generic i-gon of $K$. Say $\gamma = e$ for $i = 1$, $\gamma = e \cup f$ for $i = 2$, and $\gamma = e \cup f \cup g$ for $i = 3$. Let $\Pi_1$ be the interior and $\Pi_2$ the exterior of the simple closed curve $\gamma$ on the plane.*

*(a) For $i = 1, 2$ and according to the definitions, we can destroy the $i$ in number crossings lying on $\gamma$ (keeping the remaining crossings of $K$ intact) by an $\Omega_i^-$ or an $\Omega_{i\sigma}^-$ or an $\Omega_{i\mu}^-$ move depending on if $K - \gamma$ respectively (I) lies in $\Pi_2$, (II) lies in $\Pi_1$, (III) contains non-empty parts in both $\Pi_1, \Pi_2$. Then in all cases the resulting diagram $K'$ is equivalent to $K$, that is, $K \underset{top}{\sim} K'$.*

*Summing up: there always exists a generic move $\Omega_{i\gamma}^-$ so that $K \xrightarrow{\Omega_{i\gamma}^-} K'$ and $K \xrightarrow{\Omega_{i\gamma}^-} K' \Rightarrow K \underset{top}{\sim} K'$.*

*(b) For $i = 2$, let $p$ be a base point of $K$ not on $e, f$. Then $e, f$ as basic arcs of $_pK$ are expressed uniquely as $e = \overrightarrow{z(z+1)}, f = \overrightarrow{w(w+1)}$ ($z, w \in \mathbb{N}_{2(n+2)}, z \neq w$). And there exist unique $a, b \in \mathbb{N}_{2n}$ and $x, y \in \{-1, 0, 1\}$ such that $\mu_{2(n+2),w,z}^- = (\mu_{2n,a,b,(x),(y)}^+)^{-1}$. Let us choose a base point $q$ of $K'$ in each case (I), (II), (III) of part (a) as follows: in Figure 20 rename $K', K$ to $K, K'$ respectively and rename $p, q$ to $q, p$ respectively. Depending on $a, b$, choose the base point $q$ as the position it has in the figure after these renamings. Then:*

*The for the moves $\Omega_{i\gamma}^-$ of part (a), the renaming of the labels of the crossings of $_pK$ when considered as crossings of $_qK'$ is given by the relabeling function $\mu_{2(n+2),w,z}^- = (\mu_{2n,a,b,(x),(y)}^+)^{-1}$.*

(c) For $i = 3$, an $\Omega_3$ or $\Omega_{3\mu}$ move on a generic 3-gon $e, f, g$ exchanges it with a new generic 3-gon $e', f', g'$, namely, if $e_1, f_1, g_1$ are the slightly bigger arcs than $e, f, g$ taking part in the move, and if $e'_1, f'_1, g'_1$ are their final positions, then $e', f', g'$ are the basic arcs of the resulting diagram which lie on $e'_1, f'_1, g'_1$ respectively. In all cases, the resulting diagram $K'$ is equivalent to $K$, that is, $K \underset{top}{\sim} K'$.

Summing up: there always exists a generic move $\Omega_{3\gamma}$ so that $K \xrightarrow{\Omega_{3\gamma}} K'$ and $K \xrightarrow{\Omega_{3\gamma}} K' \Rightarrow K \underset{top}{\sim} K'$.

(d) For $i = 3$, let $\Omega_3$ or $\Omega_{3\mu}$ be a move on a 3-gon $e, f, g$. Let $p$ be a base point of $K$ either not inside the 2-disk of the move or else on one of the arcs $e, f, g$. The three crossings appearing as endpoints of $e, f, g$ are expressed in ${}_pK$ as:

$$\Delta_1 = (\alpha, \gamma)^{\pi_1}, \Delta_2 = (\beta, \epsilon)^{\pi_2}, \Delta_3 = (\delta, \zeta)^{\pi_3},$$

for $\alpha, \gamma, \epsilon, \beta, \delta, \epsilon$ distinct integers in $\mathbb{N}_{2(n+3)}$ and $|\alpha - \beta| = |\gamma - \delta| = |\epsilon - \zeta| = 1$ modulo $2(n+3)$.

Let us choose a base point $q$ of $K'$ as $p$ itself whenever $p$ does not lie in the 2-disk of the move, or else we choose it as a point on the final position $e', f', g'$ of $e, f, g$ depending on if $p$ lies on $e, f, g$ respectively.

Then the three crossings appearing as endpoints of $e', f', g'$ are expressed in ${}_qK'$ as:

$$\Delta'_1 = (\beta, \delta)^{\pi_1}, \ \Delta'_2 = (\alpha, \zeta)^{\pi_2}, \ \Delta'_3 = (\gamma, \epsilon)^{\pi_3}$$

and ${}_pK$ becomes ${}_qK'$ as:

$${}_pK = \cdots (x, y)^{\pi} \cdots \Delta_1, \Delta_2, \Delta_3 \xrightarrow{\Omega_3} {}_qK' = \cdots (x, y)^{\pi} \cdots \Delta'_1, \Delta'_2, \Delta'_3.$$

*Proof.* (a) We prove it for $i = 2$ and quite similarly it is proved for $i = 1$ as well.

For an $\Omega_2^-$ move this is already known but the following argument covers all cases: Let $\kappa$ be a knotted loop in $\mathbb{R}^3$ with $K$ as its diagram on the plane $\rho$. And let $\overline{e}, \overline{f}$ the two arcs of $\kappa$ that project onto $e, f$ respectively. The over/under information at the endpoints of $e, f$ imply that the endpoints of $\overline{e}$ are higher than those of $\overline{f}$, that is, in a higher distance from the plane $\rho$ (with respect to our point of view, which is in one half-space of $\rho$). So for some arcs $\overline{e}', \overline{f}'$ of $\kappa$ slightly larger than $\overline{e}, \overline{f}$, we can perform a vertical movement in space high enough for $\overline{e}'$ and low enough for $\overline{f}'$ (fixing throughout their endpoints). Enough mans that the new position of all points of $\overline{e}$ are higher than all other points of $K$ (the endpoints of $\overline{e}$ are at the same distance), and similarly, the new position of all points of $\overline{f}$ are lower than all other points of $K$ (the endpoints of $\overline{f}$ are at the same distance). Subsequently, we can move parallel to $\rho$ the points of the new positions of $\overline{e}', \overline{f}'$ (fixing the endpoints throughout). We take care this parallel move to make the final positions of $\overline{e}, \overline{f}$ project on $e', f'$.

The final position $\kappa'$ of $\kappa$ has the correct projection $K'$. Since the two moves described above can be implemented by isotopies of the space, we have that $\kappa, \kappa'$ are isotopic in space. So $K \underset{top}{\sim} K'$ as wanted.

(b) The given description of $e, f$ concerns basic arcs, thus the numbers $z, w$ are unique. The definition of the relabelings $\mu^-_{2(n+2),w,z}$ and of the relabelings $(\mu^+_{2n,a,b,(x),(y)})^{-1}$ are given in common in Figure 25. By these definitions it follows immediately that for the numbers $x, y$ at hand, there exist unique $a, b \in \mathbb{N}_{2n}$ and $x, y \in \{-1, 0, 1\}$ such that $\mu^-_{2(n+2),w,z} = (\mu^+_{2n,a,b,(x),(y)})^{-1}$, as stated here.

If our move on the 2-gon $e, f$ is an $\Omega_2^-$ move, then for these $a, b, x, y$ and depending on which of $e, f$ is over the other, there exists a unique case for the topological picture near the 2-gon $e, f$, among the ones (cases (A)-(L)) depicted in Figure 20. In it, the base point $q$ of $K'$ and $p$ of $K$ are as described here in the statement of the Lemma (with the reversing in the roles of $K, K'$ and of $p, q$). And of course then the renaming of the labels of the crossings of $_pK$ when considered as crossings of $_qK'$ is given by the relabeling function $\mu^-_{2(n+2),w,z} = (\mu^+_{2n,a,b,(x),(y)})^{-1}$ as stated.

Finally, if it happens our move on the 2-gon $e, f$ to be an $\Omega_{2\sigma}^-$ or $\Omega_{2\mu}$ move, let us notice that the labeling of the crossings in the 12 cases (A)-(L) of Figure 20 is not affected if we modify the pictures near the two crossings to indicate that any one of the parts of $K - (e \cup f)$ lies inside $\Pi_1$ or $\Pi_2$ at will, but keeping the over/under information intact. So then, for the values of $a, b, x, y$ at hand, and depending on which of $e, f$ is over the other, there exists a unique case for the topological picture near the 2-gon $e, f$ among the 12 cases (A)-(L) in the (modified) Figure 20. As above, in this unique case the base point $q$ of $K'$ and $p$ of $K$ are as described here in the statement of the Lemma (with the reversing in the roles of $K, K'$ and of $p, q$). And of course then again, the renaming of the labels of the crossings of $_pK$ when considered as crossings of $_qK'$ is given by the relabeling function $\mu^-_{2(n+2),w,z} = (\mu^+_{2n,a,b,(x),(y)})^{-1}$ as stated here, proving the statement.

(c) The claim for the exchange of a 3-gon by another holds immediately by the definition of the $\Omega_3, \Omega_{3\mu}$ moves.

The claim $K \underset{top}{\sim} K'$, holds automatically for $\Omega_3$ moves. For an $\Omega_{3\mu}$ move the desired result is proved quite similarly to part (b): we move in space an arc of a knot projecting to $K$. Namely we move the arc $s$ that projects to the red arc in Figure 49. $e_1$ is the top or the bottom arc among the three in the 3-gon. So $s$ can be moved freely in space in the desired way.

(d) The arcs $e_1, f_1, g_1$ of the move that contain $e, f, g$ respectively, intersect at the endpoints of $e, f, g$ and have the same over/under relation at these points as $e, f, g$. By their definition, the moves $\Omega_3$ or $\Omega_{3\mu}$ retain this information at the final position $e'_1, f'_1, g'_1$ of $e_1, f_1, g_1$, thus this information is inherited by the basic arcs $e', f', g'$ on $e'_1, f'_1, g'_1$, and by definition this makes the triad $e', f', g'$ to a generic 3-gon.

For an $\Omega_3$ move, the noted description (in the statement of the current Lemma) for the crossings at the endpoints of $e, f, g$ is correct, as was seen in §2.3.5 (ch. Relation 2.19). This description is not affected by the region $\Pi_1$ or $\Pi_2$ in which each part of $K - (e \cup f \cup g)$ belongs, so it is a correct description for $\Omega_{3\mu}$ moves as well.

Also, for the chosen $p$ and $q$ and for an $\Omega_3$ move, the noted description of the crossings at the endpoints of $e', f', g'$ is a correct one as was seen in §2.3.5 (ch. Relation 2.20). This description is not affected by the region $\Pi_1$ or $\Pi_2$ in which each part of $K - (e \cup f \cup g)$ belongs, so it is the correct description for $\Omega_{3\mu}$ moves as well, which finishes the proof. □

4.5. $^{al}\Omega_1$ **mirroring.** In this section we deal with the mirroring of the algebraic $\Omega_1$ moves via the topological moves of diagrams. Namely, we prove:

**Lemma 25.** (a) If $_pK \xrightarrow{^{al}\Omega_1^+} A \in \Sigma$, then $A = {_qK'}$ for some diagram $K'$ and base point $q$ of it, thus $A \in \Sigma_D$. Moreover $K \xrightarrow{\Omega_1^+} K'$, thus $K \underset{top}{\sim} K'$.

(b) If $_pK \xrightarrow{^{al}\Omega_1^-} A \in \Sigma$, then $A = {_qK'}$ for some diagram $K'$ and base point $q$ of it, thus $A \in \Sigma_D$. Moreover $K \xrightarrow{\Omega_1^-} K'$ or $K \xrightarrow{\Omega_{iso\sigma}} K_0 \xrightarrow{\Omega_1^-} K'$ for some diagram $K_0$. Equivalently, $K \xrightarrow{\Omega_1^-} K'$ or $K \xrightarrow{\Omega_{1\sigma}^-} K'$, that is, always $K \xrightarrow{\Omega_{1\gamma}^-} K'$. Hence it is always $K \underset{top}{\sim} K'$.

*Proof.* (a) Let $\mu(_pK) = \mu(K) = n \geqslant 0$.

For $n > 0$, let the given algebraic move ${}^{al}\Omega_1^+$ be performed on $e = \overrightarrow{a(a+1)}$ (ch. §3.3.2) of ${}_pK = \Delta_1 \cdots \Delta_n$. Then (ch. relation (3.1)):

(4.1)
$$A = \mu({}_pK)\Delta$$
$$\mu = \mu_{2n,a,(x)}^+ \quad \text{(for some } x \in \{-1, 0, 1\})$$
$$\Delta = \left(\tfrac{1-\theta}{2}(a+1) + \tfrac{1+\theta}{2}(a+2), \tfrac{1+\theta}{2}(a+1) + \tfrac{1-\theta}{2}(a+2)\right)^{\pi\theta},$$

where $\theta = \theta_e$ is a placement number of $e$, $\phi$ is some orientation and $\pi_\phi$ is the sign of $\phi$.

We consider the basic arc $\overline{e} = \overrightarrow{a(a+1)}$ of $K$ corresponding to the algebraic arc $e$ of ${}_pK$.

By Lemma 4 there exists a room $\Pi$ of $K$ containing $\overline{e}$ in its boundary so that $\Pi_{\phi\overline{e}}$ (that is: $\Pi$ equipped with the orientation $\phi$ induces on the part $\overline{e}$ of its boundary, the orientation which $\overline{e}$ has as a basic arc of $K$).

Inside $\Pi$ we perform a topological move $K \xrightarrow{\Omega_1^+} K'$ of diagrams on a subarc $\ell$ of $\overline{e}$. We position a base point $q$ of $K'$ as in Figure 12, depending on $a$ and $x$. Then the relabeling function of our $\Omega_1^+$ topological move is $\mu = \mu_{2n,a,(x)}^+$. Thus the crossing $\Delta_i = (x_i, y_i)^{\pi_i}, i = 1, \ldots, n$ of $K$ expressed in the base point $p$, becomes the crossing $(\mu(x_i), (\mu(y_i)))^{\pi_i}$ of $K'$ expressed in the base $q$.

For our topological move, we moreover choose a twist (placement) of $\overline{e}$ (ch. Definition 15) so that its placement number is $\theta_{\overline{e}} = \theta_e$. This means that the final position $\ell'$ of $\ell$ after the move lies in $\Pi$ as in Figure 13 for the value $\theta_{\overline{e}} = \theta_e$. Let $\overline{\Delta}$ be the new crossing created by the move.

It is trivially checked that the relations $\Pi_{\phi\overline{e}}, \theta_{\overline{e}} = \theta_e$ and the fact that we perform the topological move inside $\Pi$ always imply that the sign of $\overline{\Delta}$ is $\pi_\phi\theta$ and its two labels are $\tfrac{1-\theta}{2}(a+1) + \tfrac{1+\theta}{2}(a+2), \tfrac{1+\theta}{2}(a+1) + \tfrac{1-\theta}{2}(a+2)$ with heights $-1$ and $1$ respectively. So:
$$\overline{\Delta} = \left(\tfrac{1-\theta}{2}(a+1) + \tfrac{1+\theta}{2}(a+2), \tfrac{1+\theta}{2}(a+1) + \tfrac{1-\theta}{2}(a+2)\right)^{\pi_\phi\theta} = \Delta.$$
Hence:
${}_pK' = (\mu(x_1), (\mu(y_1)))^{\pi_1} \cdots \mu(\Delta_0) \cdots (\mu(x_n), (\mu(y_n)))^{\pi_n} \overline{\Delta} = \mu({}_pK)\Delta = A$ which is the first required result, in particular $A \in \Sigma_D$. Moreover, by construction $K \xrightarrow{\Omega_1^+} K'$, thus $K \underset{top}{\sim} K'$ as wanted.

For $n = 0$, a similar argument shows the result. In this case it is $K \approx S^1$ and the given algebraic move ${}^{al}\Omega_1^+$ is not performed to some basic arc, since such arcs are not defined for circles. Nevertheless it holds (ch. §3.3.2 and relation 3.2):
$$A_2 = \Delta$$
$$\Delta = \left(\tfrac{1-\theta}{2}1 + \tfrac{1+\theta}{2}2, \tfrac{1+\theta}{2}1 + \tfrac{1-\theta}{2}2\right)^{\pi\theta},$$
where $\theta_\varnothing$ is a placement number of $\varnothing$, and $\pi = \pi_\phi$ is the sign of some orientation $\phi$.

This time we perform a topological move $K \xrightarrow{\Omega_1^+} K'$ on a subarc $e$ of the circle $K$ which does not contain the base point $p$. The circle $K$ has exactly two rooms $\Pi, \Pi'$ (its interior and exterior). We perform the move in the room which equipped with the orientation $\phi$ induces on its boundary $K$ the orientation it has as a diagram. Moreover, we choose a twist of $K$ for this move so that its placement number is $\theta_e = \theta_\varnothing$.

If we call $\overline{\Delta}$ the crossing created by the move, it is trivially checked that its sign is always $\pi_\phi\theta$. So:
$$\overline{\Delta} = \left(\tfrac{1-\theta}{2}(a+1) + \tfrac{1+\theta}{2}(a+2), \tfrac{1+\theta}{2}(a+1) + \tfrac{1-\theta}{2}(a+2)\right)^{\pi\theta} = \Delta,$$
hence:
${}_pK' = \overline{\Delta} = \Delta = A$ which is the first required result(the claimed $q$ is $p$), thus in particular $A \in \Sigma_D$. Moreover, by construction $K \xrightarrow{\Omega_1^+} K'$, thus $K \underset{top}{\sim} K'$ as wanted.

(b) The argument is similar to that for part (a) with some extra complication regarding the exact topological moves that mirror the given algebraic one.

So let $\mu(_pK) = \mu(K) = n + 1 \geqslant 1$. And let the given algebraic move be performed on the loop (1-gon) $e = \overrightarrow{z(z+1)}$ ($z \in \mathbb{N}_{2(n+1)}$) of $_pK$.

Then $A$ is what we get deleting from $_pK$ the crossing of the loop and relabeling the other crossings appropriately (ch. Relations 3.3, 3.4):

$_pK = (x_1, y_1)^{\pi_1} \cdots (x_n, y_n)^{\pi_n} \Delta$

$\Delta = (z : z+1) = \left(\frac{1-\theta}{2}(a+1) + \frac{1+\theta}{2}(a+2), \frac{1+\theta}{2}(a+1) + \frac{1-\theta}{2}(a+2)\right)^\pi$ (for some $a \in \mathbb{N}_{2n}$)

$A = \mu(_pK)$

$\mu = \mu^-_{2(n+1),z} = (\mu^+_{2n,a,(x)})^{-1}$ (for the appropriate $x \in \{-1, 0, 1\}$),

where $\theta = \theta_e$ is a placement number of $e$, $\phi$ is some orientation and $\pi_\phi$ is the sign of $\phi$.

We consider the basic arc $\bar{e} = \overrightarrow{z(z+1)}$ of $K$. Since $e$ is a loop of $_pK$, the labels $z, z+1$ concern the same crossing of $_pK$. So the two endpoints of the arc $\bar{e}$ of the diagram $K$ coincide say to a point $\overline{\Delta}$, making it a closed curve. Since $\bar{e}$ is a basic arc of $K$, it contains no self-intersections in its interior, thus it is a simple closed curve.

At this point we can continue with the generic 1-gon terminology developed in the previous section, or we can conveniently consider once again the placement of the diagrams on the 2-sphere $S^2 = \rho \cup \{\infty\}$ instead on the plane $\rho$ (ch. §1.5). We go with the second option: By Schonflies Theorem the complement of $\bar{e}$ on the sphere consists of two open disks $\Pi_1$ and $\Pi_2$.

Let us observe that $K - \bar{e}$ is a connected set, and it has to lie entirely in $\Pi_1$ or entirely in $\Pi_2$. Call $\Pi$ this $\Pi_i$. We consider a small collar neighborhood $U \approx \bar{e} \times [0, 1]$ of $\bar{e}$ in $\Pi$ (Figure 47), which intersects $K - \bar{e}$ in two sub arcs stemming from the crossing point $\Delta$ of $K$ and with no other common points. Let $\bar{e}_1$ be the part of the given diagram $K$ lying on the 2-disk $U \cup \Pi$.

In case $\Pi = \Pi_2$ we perform a topological move $K \xrightarrow{\Omega^-_1} K'$ of diagrams on $\bar{e}'$ inside the 2-disk $U \cup \Pi$ of the plane making $\Delta$ disappear. In case $\Pi = \Pi_1$ we perform an isotopy move $\Omega_{iso\sigma}$ on the sphere $S^2$, moving $\bar{e}_1$ inside the 2-disk $\Pi$ bringing $\bar{e}_1$ to the position of the middle diagram in Figure 47. Let $K \xrightarrow{\Omega^s_{iso}} K_0$ be this move. We consider as base point of $K_0$ the final position $p_0$ of $p$ after the isotopy. Then $_{p_0}K_0 = {_pK} =$. Consequently we perform on $K_0$ a move $K_0 \xrightarrow{\Omega^-_1} K'$ inside the 2-disk $U \cup \Pi$ making $\Delta$ disappear.

In both cases, we choose a base point $q$ of $K'$ with the help of Figure 12 as follows: In the figure we consider $K$ or $K_0$ as the diagram $K'$, also $p_0$ as the point $q$, and $K'$ as the diagram $K$. Then depending on $a$, we chose for $q$ the place where point $p$ is located.

Then the relabeling function of our $\Omega^-_1$ topological move is $(\mu^+_{2n,a,(x)})^{-1} = \mu^-_{2(n+1),z} = \mu$. Thus the crossing $\Delta_i = (x_i, y_i)^{\pi_i}, i = 1, \ldots, n$ of $K$ expressed in the base point $p$, becomes the crossing $(\mu(x_i), (\mu(y_i)))^{\pi_i}$ of $K'$ expressed in the base $q$. So:

$_qK' = (\mu(x_1), (\mu(y_1)))^{\pi_1} \cdots (\mu(x_n), (\mu(y_n)))^{\pi_n} = \mu(_pK) = A$ which is the first required result, thus in particular $A \in \Sigma_D$. Moreover, by construction $K \xrightarrow{\Omega^-_1} K'$ (whenever $\Pi = \Pi_2$) or $K \xrightarrow{\Omega_{iso\sigma}} K_0 \xrightarrow{\Omega^-_1} K'$ (whenever $\Pi = \Pi_1$). Thus $K \underset{top}{\sim} K'$ (ch. Theorem 3) as wanted. □

4.6. $^{al}\Omega_2$ **mirroring.** Now we prove the mirroring of the algebraic $\Omega_2$ moves by topological moves:

**Lemma 26.** (a) If $_pK \xrightarrow{^{al}\Omega^+_2} A \in \Sigma$, then $A = {_qK'}$ for some diagram $K'$ and base point $q$ of it, thus $A \in \Sigma_D$. Moreover $K \xrightarrow{\Omega^+_2} K'$, so $K \underset{top}{\sim} K'$.

(b) If $_pK \xrightarrow{^{al}\Omega_2^-} A \in \Sigma$, then $A = {_qK'}$ for some diagram $K'$ and base point $q$ of it, thus $A \in \Sigma_D$. Moreover $K \xrightarrow{\Omega_2^-} K'$ or $K \xrightarrow{\Omega_{2\sigma}^-} K'$ or $K \xrightarrow{\Omega_{2\mu}^-} K'$, that is, always $K \xrightarrow{\Omega_{2\gamma}^-} K'$. Hence by Lemma 24 it is always $K \underset{top}{\sim} K'$.

*Proof.* (a) We shall work as in part (a) of Lemma 25. Some extra care is needed due to the fact that the two arcs of a topological $\Omega_2^+$ move must lie on the boundary of a common room of $K$ in order for the move to be performed.

So, let $\tau(_pK) = \tau(K) = n \geqslant 0$.

For $n > 0$, let the given algebraic move $^{al}\Omega_2^+$ be performed on two basic arcs $e_1 = \overrightarrow{a(a+1)}$, $e_2 = \overrightarrow{b(b+1)}$ of $_pK = \Delta_1 \cdots \Delta_n = (x_1, y_1)^{\pi_1} \cdots (x_n, y_n)^{\pi_n}$, for some $a, b \in \mathbb{N}_{2n}$, $a \leqslant b$ (§3.4.2). By definition of the $^{al}\Omega_2^+$ move, the cycle senses of $e_1, e_2$ can be compared, and this by definition (Definition 20) implies in particular that $e_1, e_2$ or $e_1, -e_2$ are parts of the same oriented $\phi$-cycle of the based symbol $_pK$. Also, by the definition of $^{al}\Omega_2^+$ moves it holds (ch. Relation (3.5)):

$$A = \mu(_pK)\Delta_{n+1}\Delta_{n+2}$$
$$\mu = \mu^+_{2n,a,b,(x),(y)} \quad \text{(for some } x, y \in \{-1, 0, 1\})$$
$$\Delta_1 = \left(\tfrac{1-\theta}{2}(a+1) + \tfrac{1+\theta}{2}\left(b + \tfrac{7+\rho}{2}\right), \tfrac{1+\theta}{2}(a+1) + \tfrac{1-\theta}{2}\left(b + \tfrac{7+\rho}{2}\right)\right)^{\pi\theta\rho}$$
$$\Delta_2 = \left(\tfrac{1-\theta}{2}(a+2) + \tfrac{1+\theta}{2}\left(b + \tfrac{7-\rho}{2}\right), \tfrac{1+\theta}{2}(a+2) + \tfrac{1-\theta}{2}\left(b + \tfrac{7-\rho}{2}\right)\right)^{-\pi\theta\rho},$$

where $\rho = \rho_{e_1,e_2}$ is the number of similarity of cycle sense of $e_1, e_2$, $\theta = \theta_{e_1,e_2}$ is a placement number of the ordered pair $(e_1, e_2)$, $\phi$ is an orientation and $\pi_\phi$ is the sign of $\phi$.

We consider the basic arcs $\bar{e}_1 = \overrightarrow{a(a+1)}$, $\bar{e}_2 = \overrightarrow{b(b+1)}$ of $K$. We noted above that the basic arcs $e_1 = \overrightarrow{a(a+1)}$, $e_2 = \overrightarrow{b(b+1)}$ of the based symbol $_pK$, or the main arcs $e_1, -e_2$ are parts of the same $\phi$-cycle $C$ of $_pK$. We know by Lemma 15 that the $\phi$-cycles of $_pK$ are the $\phi$-cycles of $K$, that is, they are the oriented boundaries of the rooms of $K$ considered with their $\phi$ orientations. Thus accordingly, the basic arcs $\bar{e}_1, \bar{e}_2$ of $K$ or the main arcs $\bar{e}_1, -\bar{e}_2$ are parts of the oriented boundary $\partial\Pi$ of the same $\phi$-room $\Pi_\phi$ of $K$. In particular this implies that $\rho_{e_1,e_2} = \rho_{\bar{e}_1,\bar{e}_2}$, that is, the similarity or not of senses of $e_1, e_2$ on their $\phi$-cycle is the same as the similarity or not of the senses of $\bar{e}_1, \bar{e}_2$ on theirs. Since $e_1$ is an arc of the considered $\phi$-cycle $C$ of $_pK$, it also follows that $\bar{e}_1$ is an arc of $\partial\Pi$, in other words that $\Pi_{\phi\bar{e}_1}$.

Inside $\Pi$ we perform a topological move $K \xrightarrow{\Omega_2^+} K'$ of diagrams on subarcs $\ell_1, \ell_2$ of $\bar{e}_1, \bar{e}_2$. Since $a \leqslant b$, arc $\bar{e}_1$ is the first among the two in this move. We chose as placement number for this move the number $\theta$ (that is: we place $\bar{e}_1$ over $\bar{e}_2$ whenever $\theta = 1$ and under $\bar{e}_2$ whenever $\theta = -1$). We position a base point $q$ of $K'$ as in Figure 20, depending on $a, b, x, y$. Then the relabeling function of our $\Omega_2^+$ topological move is $\mu = \mu^+_{2n,a,b,(x),(y)}$. Thus the crossing $\Delta_i = (x_i, y_i)^{\pi_i}, i = 1, \ldots, n$ of $K$ expressed in the base point $p$, becomes the crossing $(\mu(x_i), (\mu(y_i)))^{\pi_i}$ of $K'$ expressed in the base $q$.

Let $\overline{\Delta}_{n+1}, \overline{\Delta}_{n+2}$ be the two new crossings created by the topological move. For this move, we moreover choose a placement of $\bar{e}_1, \bar{e}_2$ with placement number $\theta_{\bar{e}_1,\bar{e}_2} = \theta_{e_1,e_2}$ (ch. Definition 17).

As is trivially checked, in all possible topological placements (Figure 24 of Lemma 8) the relations $\Pi_{\phi\bar{e}}$, $\rho_{e_1,e_2} = \rho_{\bar{e}_1,\bar{e}_2}$, $\theta_{\bar{e}_1,\bar{e}_2} = \theta_{e_1,e_2}$ and the fact that we perform the topological move inside $\Pi$ always imply (and no matter if $a < b$ or $a = b$) that one of $\overline{\Delta}_{n+1}, \overline{\Delta}_{n+2}$ has the description of $\Delta_1$ given above, and the other has the description of $\Delta_2$.

Hence:

$_qK' = (\mu(x_1),(\mu(y_1)))^{\pi_1} \cdots (\mu(x_n),(\mu(y_n)))^{\pi_n} \overline{\Delta}_{n+1} \overline{\Delta}_{n+2} = \mu(_pK)\Delta_{n+1}\Delta_{n+2} = A$ which is the first required result, and in particular $A \in \Sigma_D$. Moreover, by construction $K \xrightarrow{\Omega_2^+} K'$, thus $K \underset{top}{\sim} K'$ as wanted.

For $n = 0$ we argue similarly. In this case it is $K \approx S^1$. The given algebraic move $^{al}\Omega_2^+$ is performed on $_pK = \varnothing$ and it holds (ch. §3.4.2 and Relation (3.6)):

$$A = \Delta_1 \Delta_2$$
$$\Delta_1 = \left(\tfrac{1-\theta}{2}(0+1) + \tfrac{1+\theta}{2}\left(0+\tfrac{7+1}{2}\right), \tfrac{1+\theta}{2}(0+1) + \tfrac{1-\theta}{2}\left(0+\tfrac{7+1}{2}\right)\right)^{\pi\theta \cdot 1}$$
$$\Delta_2 = \left(\tfrac{1-\theta}{2}(0+2) + \tfrac{1+\theta}{2}\left(0+\tfrac{7-1}{2}\right), \tfrac{1+\theta}{2}(0+2) + \tfrac{1-\theta}{2}\left(0+\tfrac{7-1}{2}\right)\right)^{-\pi\theta \cdot 1},$$

where $\theta = \theta_{\varnothing,\varnothing}$ is a placement number of the empty based symbol $A_1 = \varnothing$, $\phi$ is an orientation and $\pi_\phi$ is the sign of $\phi$.

This time we perform a topological move $K \xrightarrow{\Omega_2^+} K'$ on two disjoint subarcs $\ell_1, \ell_2$ of an arc $e$ of the circle $K$, none of which contains the base point $p$. Let in the arc $e$ oriented by the induced orientation of $K$, be $e_1$ the arc appearing before $e_2$ as we travel from the first point of $e$ along its orientation. We consider $e_1$ as the first among the two arcs for our move. And we take care to use as placement number of the move the number $\theta$. The circle $K$ has exactly two rooms $\Pi, \Pi'$ (its interior and exterior). We perform the move in the room which equipped with the orientation $\phi$ induces on its boundary $K$ the orientation it has as a diagram.

If we call $\overline{\Delta}_1, \overline{\Delta}_2$ the two crossings created by the move, it is trivially checked that their description is the one for $\Delta_1, \Delta_2$ given above. So:
$\{\overline{\Delta}_1, \overline{\Delta}_2\} = \{\Delta_1, \Delta_2\}$,
hence:
$_pK' = \overline{\Delta}_1, \overline{\Delta}_2 = \Delta_1, \Delta_2 = A$ which is the first required result (here actually it is $q = p$), thus in particular $A \in \Sigma_D$. Moreover, by construction $K \xrightarrow{\Omega_2^+} K'$, thus $K \underset{top}{\sim} K'$ as wanted.

(b) The argument is similar to that for part (a) with some extra complication regarding the exact topological moves that mirror the given algebraic one. This time we shall use the generic 2-gon terminology developed in §4.4.

So let $\tau(_pK) = \tau(K) = n+2 \geq 2$. And let the given algebraic move be performed on the 2-gon $e_1 = \overrightarrow{z(z+1)}$, $e_2 = \overrightarrow{w(w+1)}$ of $_pK = \Delta_1 \ldots \Delta_{n+2}$, where $1 \leq z, w \leq 2(n+2)$, $z \neq w$ (ch. §3.4.3). And let $a, b \in \mathbb{N}_{2n}$ and $x, y \in \{-1, 0, 1\}$ be such that $\mu^-_{2(n+2),w,z} = (\mu^+_{2n,a,b,(x),(y)})^{-1}$ (ch. §2.3.4 and Figure 25).

Then $A$ is what we get deleting from $_pK$ the two crossings of the 2-gon and relabeling the other crossings appropriately (ch. Relations 3.7, 3.8):

$$_pK = (x_1, y_1)^{\pi_1} \cdots (x_n, y_n)^{\pi_n} \Delta_{n+1} \Delta_{n+2}$$
$$\Delta_{n+1} = \left(\tfrac{1-\theta}{2}(a+1) + \tfrac{1+\theta}{2}\left(b+\tfrac{7+\rho}{2}\right), \tfrac{1+\theta}{2}(a+1) + \tfrac{1-\theta}{2}\left(b+\tfrac{7+\rho}{2}\right)\right)^{\pi\theta\rho} \text{ (for some } a,b \in \mathbb{N}_{2n})$$
$$\Delta_{n+2} = \left(\tfrac{1-\theta}{2}(a+2) + \tfrac{1+\theta}{2}\left(b+\tfrac{7-\rho}{2}\right), \tfrac{1+\theta}{2}(a+2) + \tfrac{1-\theta}{2}\left(b+\tfrac{7-\rho}{2}\right)\right)^{-\pi\theta\rho}$$
$$A = \mu(_pK)$$
$$\mu = \mu^-_{2(n+1),z,w} = (\mu^+_{2n,a,b,(x),(y)})^{-1} \quad \text{(for appropriate} \quad x, y \in \{-1, 0, +1\})$$

where $\rho = \rho_{e_1,e_2}$ is the (uniquely determined) number of similarity of cycle sense of $e_1, e_2$, $\theta = \theta_{e_1,e_2}$ is a placement number of the ordered pair $(e_1, e_2)$, $\phi$ is an orientation and $\pi_\phi$ is the sign of $\phi$.

We consider the basic arcs $\bar{e}_1 = \overrightarrow{z(z+1)}$, $\bar{e}_2 = \overrightarrow{w(w+)}$ of $K$ corresponding to the algebraic arcs $e_1, e_2$ of $_pK$. Since $e_1, e_2$ is a 2-gon of $_pK$, the arcs $\bar{e}_1, \bar{e}_2$ share their two endpoints, so the union $\bar{e}_1 \cup \bar{e}_2$ is a closed curve. Their endpoints are the points $\Delta_{n+1}, \Delta_{n+2}$ on the plane and they are distinct. Since $\bar{e}_1, \bar{e}_2$ are basic arcs of $K$, they do not contain any intersection of $K$ in their interior, and since their endpoints are distinct, the union $\bar{e}_1 \cup \bar{e}_2$ is a simple closed curve. Also, since one of $e_1, e_2$ is over the other at both of their two common endpoints, the same happens for $\bar{e}_1, \bar{e}_2$. So $\bar{e}_1, \bar{e}_2$ is a generic 2-gon of $K$.

By Lemma 24 we can perform a topological $\Omega_{2\gamma}^-$ move $K \xrightarrow{\Omega_{2\gamma}^-} K'$ on the generic 2-gon $\bar{e}_1, \bar{e}_2$ which destroys $\Delta_{n+1}, \Delta_{n+2}$ retaining the other crossings $\Delta_1, \ldots, \Delta_n$. By the same lemma we can also choose a base point $q$ of $K'$ so that the relabeling of the labels for the crossings of $K$ with base point $p$ which survive in $K'$ is given by the function $\mu = \mu_{2(n+1),z,w}^- = (\mu_{2n,a,b,(x),(y)}^+)^{-1}$. Thus the crossing $\Delta_i = (x_i, y_i)^{\pi_i}$, $i = 1, \ldots, n$ of $K$ expressed in the base point $p$, becomes the crossing $(\mu(x_i), (\mu(y_i)))^{\pi_i}$ of $K'$ expressed in the base point $q$.

So:
$_qK' = (\mu(x_1), (\mu(y_1)))^{\pi_1} \cdots (\mu(x_n), (\mu(y_n)))^{\pi_n} = \mu(_pK) = A$ which is the first required result, thus in particular $A \in \Sigma_D$. Moreover, by construction $K \xrightarrow{\Omega_{2\gamma}^-} K'$ and by Lemma 24 it is $K \underset{top}{\sim} K'$, as wanted. $\square$

### 4.7. $^{al}\Omega_3$ mirroring.
We now prove the mirroring of $^{al}\Omega_3$ moves:

**Lemma 27.** *If $_pK \xrightarrow{^{al}\Omega_3} A \in \Sigma$, then $A = {}_qK'$ for some diagram $K'$ and base point $q$ of it, thus $A \in \Sigma_D$. Moreover $K \xrightarrow{\Omega_3} K'$ or $K \xrightarrow{\Omega_{3\sigma}} K'$ or $K \xrightarrow{\Omega_{3\mu}} K'$, that is, always $K \xrightarrow{\Omega_{3\gamma}} K'$. Hence by Lemma 24 it is always $K \underset{top}{\sim} K'$.*

*Proof.* Let $\mu(_pK) = \mu(K) = n + 3 \geqslant 3$ and let the given algebraic move $^{al}\Omega^3$ be performed on a 3-gon $e, f, g$ of $_pK = \Delta_1 \ldots \Delta_{n+2}$.

Then it holds (ch. Relation 3.9):

(4.2)
$$_pK = (x_1, y_1)^{\pi_1} \cdots (x_n, y_n)^{\pi_n} (\alpha, \gamma)^{\pi_1} (\beta, \epsilon)^{\pi_2} (\delta, \zeta)^{\pi_3}$$
$$A \cdots (x_1, y_1)^{\pi_1} \cdots (x_n, y_n)^{\pi_n} (\beta, \delta)^{\pi_1} (\alpha, \zeta)^{\pi_2} (\gamma, \epsilon)^{\pi_3}.$$

We consider the basic arcs $\bar{e}_1 = \alpha\beta$, $\bar{e}_2 = \gamma\delta, \bar{e}_3 = \epsilon\zeta$ of $K$. In this notation $\alpha, \beta$ are labels of the crossings at the endpoints of $\bar{e}_1$ with respect to the base point $p$, but they do not necessarily denote respectively the first and second endpoint of $\bar{e}_1$ (it could be the other way around). Since $e_1, e_2, e_3$ is a 3-gon of $_pK$, each one of the arcs $\bar{e}_1, \bar{e}_2, \bar{e}_2$ shares an endpoint with one of the other two arcs and its other endpoint with the third arc. So the union $\bar{e}_1 \cup \bar{e}_2 \cup \bar{e}_3$ is a closed curve. Their endpoints are $\Delta_{n+1} = \alpha : \gamma$, $\Delta_{n+2} = \beta : \epsilon$, $\Delta_{n+3} = \delta : \zeta$ and they are distinct. Since $\bar{e}_1, \bar{e}_2, \bar{e}_3$ are basic arcs of $K$, they do not contain any intersection of $K$ in their interior, and since the point $\Delta_{n+1}, \Delta_{n+2}, \Delta_{n+3}$ of the plane are distinct, the union $\bar{e}_1 \cup \bar{e}_2$ is a simple closed curve. Also, since one of $e_1, e_2, e_3$ is over the rest at both of its endpoints and another one of $e_1, e_2, e_3$ is under the rest at both of its endpoints, the same happens for $\bar{e}_1, \bar{e}_2, \bar{e}_3$. So $\bar{e}_1, \bar{e}_2, \bar{e}_1$ is a generic 3-gon of $K$.

By Lemma 24 we can perform a topological $\Omega_{3\gamma}$ move $K \xrightarrow{\Omega_{3\gamma}} K'$ on the generic 3-gon $\bar{e}_1, \bar{e}_2, \bar{e}_3$ which replaces the generic 3-gon with another one. The three crossings appearing as endpoints of the new 3-gon are expressed in $_qK'$ as:

$$\Delta'_{n+1} = (\beta, \delta)^{\pi_1}, \ \Delta'_{n+2} = (\alpha, \zeta)^{\pi_2}, \ \Delta'_{n+3} = (\gamma, \epsilon)^{\pi_3}$$

and $_pK$ becomes $_qK'$ as:

$$_pK = \cdots (x_i, y_i)^{\pi_i} \cdots \Delta_{n+1}, \Delta_{n+2}, \Delta_{n+3} \xrightarrow{\Omega_{e\gamma}} {}_qK' = \cdots (x,y)^\pi \cdots \Delta'_{n+1}, \Delta'_{n+2}, \Delta'_{n+3}.$$

So:
$_qK' = \cdots (x_i, y_i)^{\pi_i} \cdots \Delta'_{n+1}\Delta'_{n+2}\Delta'_{n+3} = \mu(_pK) = A$ which is the first required result, thus in particular $A \in \Sigma_D$. Moreover, by construction $K \xrightarrow{\Omega_{3\gamma}} K'$ and by Lemma 24 it is $K \underset{top}{\sim} K'$, as wanted. □

## 5. The symbol invariant $S$

We now assemble all work done so far and we finally define our first invariant.

By Lemmata 22, 23, 25, 26, 27 we immediately have:

**Corollary 3.** *The algebraic moves* $^{al}\Omega_\beta, {}^{al}\Omega_{iso}, {}^{al}\Omega_1^+, {}^{al}\Omega_1^-, {}^{al}\Omega_2^+, {}^{al}\Omega_2^-, {}^{al}\Omega_3$ *happen within the set $\Sigma_D$ of diagrams on the plane. That is: if $A \xrightarrow{\omega} A'$ is any one of the algebraic moves, then both $A, A' \in \Sigma_D$.*

Although $^{al}\Omega_{iso}$ moves are special cases of the $^{al}\Omega_\beta$ moves, we continue bellow to consider them as a separate entity. It will be no harm to forget altogether about them.

The 7 algebraic moves actually define an equivalence relation in the set $\Sigma_D$:

**Definition 24.** *For $A, A' \in \Sigma_D$ we say that $A'$ is derived from $A$ via algebraic moves, whenever there exists a finite sequence of such moves $A = A_1 \xrightarrow{\omega_1} A_2 \xrightarrow{\omega_2} A_3 \cdots \xrightarrow{\omega_{n-1}} A_n = A'$ starting with $A$ and ending with $A'$. We write $A \underset{al}{\sim} A'$.*

**Lemma 28.** $\underset{al}{\sim}$ *is an equivalence relation in $\Sigma_D$.*

*Proof.* Reflexivity holds by Lemma 16. Transitivity holds immediately from the definition of $\underset{al}{\sim}$ (the fact that the sequences relating two based symbols are finite is important here).

For the commutativity:

Let $A = {}_pK \underset{al}{\sim} {}_qK' = A'$. Then there exists a sequence $_pK = {}_{p_1}K_1 \xrightarrow{\omega_1} {}_{p_2}K_2 \xrightarrow{\omega_2} {}_{p_3}K_3 \cdots \xrightarrow{\omega_{n-1}} {}_{p_n}K_n = {}_qK'$ ($\omega_i$'s are algebraic moves) connecting the two symbols, where we called $K$ as $K_1$ and $K'$ as $K_n$. By Lemmata 22, 23, 25, 26, 27, we have $K_1 \underset{top}{\sim} K_2 \underset{top}{\sim} K_3 \cdots \underset{top}{\sim} K_n$, so $K_1 \underset{top}{\sim} K_n$. Thus $K_n \underset{top}{\sim} K_1$ as well (since $\underset{top}{\sim}$ is an equivalence relation). So there exists a finite sequence $K_n = L_1 \xrightarrow{f_1} L_2 \xrightarrow{f_2} L_3 \cdots \xrightarrow{f_{m-1}} L_m = K_1$ ($f_j$'s are topological moves). Then by Lemmata 16, 17, 18, 19, 20, there exist suitable base points of the $L_j$'s and algebraic moves between the corresponding based symbols: $_{q_1}L_1 \xrightarrow{f_1} {}_{q_2}L_2 \xrightarrow{f_2} {}_{q_3}L_3 \cdots \xrightarrow{f_{m-1}} {}_{q_m}L_m$. So $_{q_1}L_1 \underset{al}{\sim} {}_{q_m}L_m$, that is, $_{q_1}K_n \underset{al}{\sim} {}_{q_m}K_1$. For the base points $q \in K_n = K'$ and $p \in K_1 = K$, Lemma 16 implies $_qK_n \xrightarrow{{}^{al}\Omega_\beta} {}_{q_1}K_n$ and $_{q_m}K_1 \xrightarrow{{}^{al}\Omega_\beta} {}_pK_1$ for suitable $^{al}\Omega_\beta$ moves. So $_qK_n \xrightarrow{{}^{al}\Omega_\beta} {}_{q_1}K_n \underset{al}{\sim} {}_{q_m}K_1 \xrightarrow{{}^{al}\Omega_\beta} {}_pK_1$ which by the definition of the relation $\underset{al}{\sim}$ gives $_qK_n \underset{al}{\sim} {}_pK_1$. Since $_qK_n = A', {}_pK_1 = A$ we have proved $A' \underset{al}{\sim} A$, so the required reflexivity property completing the proof. □

**Definition 25.** *For a based symbol $A \in \Sigma_D$ we denote by $[\![A]\!]_{al}$ its equivalence class in the relation $\underset{al}{\sim}$ defined in the set $\sigma_D$ of diagrams. We call the classes simply as symbols. We denote as $[\![\Sigma_D]\!]_{al}$ the set of all symbols (the classes of $\underset{al}{\sim}$).*

*We recall that for a diagram $K \in D$ we have denoted earlier by $[\![K]\!]_{top}$ its equivalence class in the relation $\underset{top}{\sim}$.*

*For a knotted loop $\kappa$ we denote by $[\![\kappa]\!]$ its equivalence class in space under the isotopies on the plane, that is, the knot $\kappa$ it represents.*

The last proof distilled also gives:

**Lemma 29.** *(a) Let $A = {}_pK$, $A' = {}_qK'$ be two based symbols. If $A \underset{al}{\sim} A'$ then $K \underset{top}{\sim} K'$.*

*(b) Let $K, K'$ be two equivalent diagrams: $K \underset{top}{\sim} K'$. Then for arbitrary base points $p \in K, q \in K'$, it is ${}_pK \underset{al}{\sim} {}_qK'$, so $[\![{}_pK]\!]_{al} = [\![{}_qK']\!]_{al}$.*

*(c) Let $K$ be any diagram. Then for arbitrary base points $p, q$ of $K$ it is ${}_pK \underset{al}{\sim} {}_qK$. So $[\![{}_pK]\!]_{al} = [\![{}_qK]\!]_{al}$ as well.*

*Proof.* (a) Since $A = {}_pK \underset{al}{\sim} {}_qK' = A'$, there exists a sequence ${}_pK \xrightarrow{\omega_1} {}_{p_2}K_2 \xrightarrow{\omega_2} {}_{p_3}K_3 \cdots \xrightarrow{\omega_{n-1}} {}_qK'$ ($\omega_i$'s are algebraic moves). By Lemmata 22, 23, 25, 26, 27, we have $K \underset{top}{\sim} K_2 \underset{top}{\sim} K_3 \cdots \underset{top}{\sim} K_{n-1} \underset{top}{\sim} K'$, so $K \underset{top}{\sim} K'$.

(b) $K \underset{top}{\sim} K'$ implies that there exists a finite sequence $K = K_1 \xrightarrow{f_1} K_2 \xrightarrow{f_2} K_3 \cdots \xrightarrow{f_{n-1}} K_n = K'$ ($f_i$'s are topological moves). Then by Lemmata 16, 17, 18, 19, 20, there exist suitable base points of the $K_i$'s and algebraic moves between the corresponding based symbols: ${}_{p_1}K_1 \xrightarrow{f_1} {}_{p_2}K_2 \xrightarrow{f_2} {}_{p_3}K_3 \cdots \xrightarrow{f_{n-1}} {}_{p_n}K_n$. So ${}_{p_1}K_1 \underset{al}{\sim} {}_{p_n}K_n$, that is, ${}_{p_1}K_1 \underset{al}{\sim} {}_{p_n}K_n$. For the base points $p \in K_1 = K$ and $q \in K_n = K'$, Lemma 16 implies ${}_pK_1 \xrightarrow{{}^{al}\Omega_\beta} {}_{p_1}K_1$ and ${}_{q_n}K_n \xrightarrow{{}^{al}\Omega_\beta} {}_qK_n$ for suitable ${}^{al}\Omega_\beta$ moves. So ${}_pK_1 \xrightarrow{{}^{al}\Omega_\beta} {}_{p_1}K_1 \underset{al}{\sim} {}_{p_n}K_n \xrightarrow{{}^{al}\Omega_\beta} {}_qK_n$ which by the definition of the relation $\underset{al}{\sim}$ gives ${}_pK_1 \underset{al}{\sim} {}_qK_n$, that is ${}_pK \underset{al}{\sim} {}_qK'$ as wanted. Then the two equivalent classes in $\underset{al}{\sim}$ are the same, thus $[\![{}_pK]\!]_{al} = [\![{}_qK']\!]_{al}$ holds.

(c) Since $K \underset{top}{\sim} K$, the result comes immediately from part (a) of the current Lemma. $\square$

**Definition 26.** *According to Lemma 29 (c), a given diagram $K$ defines a unique symbol $[\![{}_pK]\!]_{al}$, no matter what its base point $p$ might be. We call this symbol as the symbol of $K$ and we denote it as $S(K)$.*

*According to Lemma 29 (b), the symbols of any two equivalent diagrams $K, K'$ are the same: $S(K) = S(K')$. This way, the class $[\![K]\!]_{top}$ of any diagram corresponds to a unique symbol which we call as the symbol of the class. We denote this symbol as $S([\![K]\!]_{top})$. So we define it as: $S([\![K]\!]_{top}) = S(K)$.*

Essentially the Last Lemma says that equivalence classes of diagrams correspond to symbols (equivalent classes of based symbols). So then knots in space correspond to symbols. This correspondence is 1-1 and onto, that is:

**Theorem 4.** *Let $\mathbb{K}$ be the set of (smooth or pl) knots in the space ($\mathbb{R}^3$ or $S^3$). Let the map*

(5.1) $\quad\quad\quad\quad \begin{array}{rl} Symb : \mathbb{K} & \to [\![\Sigma_D]\!]_{al} \\ \varkappa = [\![\kappa]\!] & \mapsto Symb(\varkappa) = S([\![K]\!]_{top})(= S(K) = [\![{}_pK]\!]_{al}) \end{array}$

*where $\kappa$ is a (smooth or pl) knotted loop in space representing the knot $\varkappa$, and $K$ is the projection of $\kappa$ on a plane ($p$ is some base point of $K$). Then Symb is 1-1 and onto.*

*Proof.* For the injectivity: Let $\varkappa, \varkappa'$ be two knots with $Symb(\varkappa) = Symb(\varkappa')$. Let $\kappa, \kappa'$ be knotted loops representing $\varkappa, \varkappa'$ respectively, and let $K, K'$ be the respective diagrams of $\kappa, \kappa'$ on the plane. By definition, the equality $Symb(\varkappa) = Symb(\varkappa')$ means that $S(K) = S(K')$. If $p, q$ are some base points of $K, K'$ respectively, then by definition $S(K) = [\![_p K]\!]_{al}$ and $S(K') = [\![q'_K]\!]_{al}$. So it is $[\![p_K]\!]_{al} = [\![q'_K]\!]_{al}$. Then by Lemma 29 (a) it is $K \underset{top}{\sim} K'$, so the equivalent classes of the two diagrams in relation $\underset{top}{\sim}$ coincide: $[\![K]\!]_{top} = [\![K']\!]_{top}$. And due to the work of Alexander, Briggs and Reidemeister this implies $\kappa, \kappa'$ are isotopic in space, thus they represent the same knot. Hence $\varkappa = \varkappa'$ which proves $Symb$=1-1.

For the surjectivity: If $X \in [\![\Sigma_D]\!]_{al}$, then $X = [\![A]\!]_{al}$ for some $A \in \Sigma_D$. And $A = {}_p K$ for a diagram $K$ and a base point $p$ of it. Now, there exists a knotted loop $\kappa$ in space with $K$ as its diagram (this is clear and we commented on this very early in the paper at the end of §1.2). Let $\varkappa$ be the knot class of $\kappa$. Then by definition, $Symb(\varkappa) = [\![_p K]\!]_{al}$, and since $A = {}_p K$ we have $Symb(\varkappa) = [\![A]\!]_{al} = X$, proving the surjectivity. $\square$

**Definition 27.** *(The symbol invariant of knots) For the (smooth or pl) knots $\varkappa$ in space ($\mathbb{R}^3$ or $S^3$) we call as the symbol of $\varkappa$ and we write $S(\varkappa)$, the value $Symb(\kappa)$ of the function in the last Theorem.*

*$Symb(\kappa)$ is an equivalence class of based symbols, in the equivalence relation $\underset{al}{\sim}$ defined in the set $\Sigma_D$ of all based symbols. More specifically it is:*

$$Symb(\varkappa) = S(\varkappa) = S(K) = S([\![K]\!]_{top}) = [\![_p K]\!]_{al},$$

*where $\kappa$ is a (smooth or pl) knotted loop in space respresenting the knot $\varkappa$, $K$ is the projection of $\kappa$ on a plane, and $p$ is some base point of $K$.*

*We also define the symbol of the knotted loop $\kappa$ to be the symbol of its diagram $K$ and we write $S(\kappa)$ for it. So we have a homogeneous notation for all 1-dimensional objects involved. In order to have as a homogeneous notation as possible, we also denote $[\![A]\!]_{al} = S(A)$ for the class of any based symbol $A$ in the relation $\underset{al}{\sim}$. With this notation it holds:*

(5.2) $$Symb(\varkappa) = S(\varkappa) = S(\kappa) = S(K) = S([\![K]\!]_{top}) = S({}_p K) = [\![_p K]\!]_{al}.$$

*According to Theorem 4, the symbol of knots is a complete knot invariant in the sense that it distinguishes any two oriented knots. We call it as the symbol invariant of knots. For each knot it is a set of based symbols and it is calculated from any knot diagram of a given knot.*

The completeness of this invariant has the expected consequences. A noted example is the following:

**Corollary 4.** *The symbol invariant $S$ of knots detects chirality and amphichirality.*

**Example 3.** *As an example, let us compute the symbol invariant for the two trefoil knots $\varkappa, \varkappa'$.*

*Let $K, K'$ be the diagrams of $\varkappa, \varkappa'$ respectively given in Figure 50. We have ${}_p K = (1,4)^{+1} (2,5)^{+1} (3,6)^{+1}$ and ${}_q K' = (4,1)^{-1} (5,2)^{-1} (6,3)^{-1}$.*

*So:*

$$Symb(\varkappa) = S(\varkappa) = S(\kappa) = S(K) = S([\![K]\!]_{top}) = S({}_p K) = [\![_p K]\!]_{al} = [\![(1,4)^{+1} (5,2)^{+1} (3,6)^{+1}]\!]_{al}$$
$$Symb(\varkappa') = S(\varkappa') = S(\kappa') = S(K') = S([\![K']\!]_{top}) = S({}_q K') = [\![_q K']\!]_{al} = [\![(1,4)^{-1} (5,2)^{-1} (3,6)^{-1}]\!]_{al}.$$

## 6. The reduced symbol invariant $R$

The symbol invariant $S$ constructed in the previous section assigns to each knot $\varkappa$ an appropriate equivalence class $S(\varkappa)$ of based symbols $\Sigma_D$. The class concerns the equivalence relation $\underset{al}{\sim}$ defined in $\Sigma_D$ which is derived by the algerbraic moves $^{al}\Omega_\beta, ^{al}\Omega_{iso}, {}^{al}\Omega_1^+, {}^{al}\Omega_1^-, \Omega_2^+, {}^{al}\Omega_2^-, {}^{al}\Omega_3$ on elements of $\Sigma_D$. $S(\varkappa)$ contains infinitely many representatives, because we can always perform $^{al}\Omega_1^+$ and $^{al}\Omega_2^+$ moves on any given representative, increasing the order of the representatives indefinitely. In this respect the $S$ invariant does not immediately inform us about if two given knots $\varkappa_1, \varkappa_2$ are the same or not, unless of course we can somehow decide about the coincidence or non-coincidence of the two sets of representatives of $S(\varkappa_1), S(\varkappa_2)$.

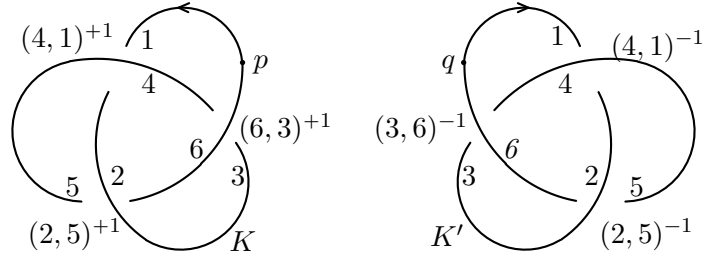

FIGURE 50. The knot diagrams $K, K'$ of some representatives $\kappa, \kappa'$ of the two trefoil knots $\varkappa, \varkappa'$, and some base points $p, q$ on them.

To overcome this difficulty, we now construct a new invariant $R$ based on the symbol invariant $S$, which can immediately decide about if two knots are the same or knot; so it can also be used effectively to classify all knots. $R(\varkappa)$ will be defined as a finite subset of the representatives of $S(\varkappa)$. This set will be computed algorithmically by any given representative of $S(\varkappa)$. We define $R(\varkappa)$ as the set of the representatives with minimum order. Each member in $R(\varkappa)$ will be called a reduced based symbol and $R$ as the reduced symbol invariant.

We provide the preliminary work in the section that follows leaving the main Theorem for the subsequent section. The main Theorem assures that $R$ is a knot invariant and that it can be computed effectively algorithmically. For the line of reasoning presented here, it will be necessary to define for each based symbol $A$ not only the set $R(A)$ of equivalent and reduced symbols, but also a similar set $\overline{R}(A)$ of equivalent symbols to which we can arrive from $A$ via moves that decrease the order of the symbols up to a point that further decrease is not possible. Eventually the two sets will be proved to be the same, allowing as a corollary the main Theorem of §6.2.

6.1. **Reduced based symbols, heredity of crossings and i-gons.**

**Definition 28.** *We call the $^{al}\Omega_1^+$, $^{al}\Omega_2^+$ as positive moves, the $^{al}\Omega_1^-$, $^{al}\Omega_2^-$ as negative moves and the $^{al}\Omega_\beta$, $^{al}\Omega_{iso}$, $^{al}\Omega_3$ as neutral moves.*

*We call any finite sequence of algebraic moves as negative or neutral depending on if all the moves are negative or neutral respectively. If the last symbol in a negative sequence does not accept any negative move, then we call the sequence as a full negative sequence.*

*We call the symbols which do not accept any negative move as momentary reduced symbols. So full negative sequences have to end with momentary reduced symbols.*

*For an $A \in \Sigma_D$ we denote $R(A)$ the set of all equivalent to $A$ (in the equivalence $\underset{al}{\sim}$) based symbols which have the minimum possible order among the equivalent symbols to $A$. In other words, $R(A)$ is the set of all those representatives of $[\![A]\!]_{al}$ of minimum possible order. We call $R(A)$ as the reduced set of $A$ and we call each member of $R(A)$ as a reduced based symbol.*

*For an $A \in \Sigma_D$ we denote $\overline{R}(A)$ the set of all based symbols $B$ to which we can arrive starting from $A$ after a full negative sequence of moves, followed by a neutral sequence of moves. Thus $B \in \overline{R}(A)$ exactly when there exists a sequence $A = A_0 \xrightarrow{\omega_1} A_1 \xrightarrow{\omega_2} \cdots \xrightarrow{\omega_n} A_n \xrightarrow{\omega_{n+1}} A_{n+1} \cdots \xrightarrow{\omega_m} A_m = B$ with $\omega_1, \cdots, \omega_n$ negative, $\omega_{n+1}, \omega_m$ neutral and $A_n$ accepting no negative moves. We allow full* *negative sequences to be empty as well, that is, we allow beginning symbols $A$ which have no negative moves. We call each member of $\overline{R}(A)$ as locally reduced based symbol.*

In the definition we allowed empty full negative sequences so that $\overline{R}(A)$ can be defined for reduced and momentary reduced symbols $A$ as well.

If $B$ is a based symbol to which we arrive from $A$ after a full negative sequence of moves then $B$ is momentary reduced and since we can extend the sequence for example by the neutral move $B \xrightarrow{\omega_1} B$, we have $B \in \overline{R}(A)$. This is the reason that we did not include in the last definition the possibility of empty neutral sequences in the second part of the demanded moves for the elements of $\overline{R}(A)$.

The locally reduced and the reduced based symbols resemble in a sense respectively the usual local and total minima of real graphs. But unlike the case of graphs, we intend to show that here the two notions of minima coincide!

Our next definition concerns the fate of the crossings of some based symbol after we perform an algebraic move on it.

The exact description of the crossings after a move is given in the definitions of the moves (§3.1, 3.2, 3.3, 3.4, 3.5 respectively). Here it is convenient to describe the resulting based symbol more loosely as follows:

After a positive move one or two new crossings appear, whereas each old crossing $\Delta = (a, b)^\pi$ is conserved suitably relabeling its labels $a, b$. After a negative move, one or two crossings are eliminated, whereas each one of the remaining crossings is conserved suitably relabeling its labels. After an $^{al}\Omega_{iso}$ move all crossings remain as they are. After an $^{al}\Omega_\beta$ move, each crossing is conserved suitably relabeling its labels. After an $^{al}\Omega_3$ move three crossings interchange their labels suitably, whereas each one of the remaining crossings remains as it is. More rigorously, but still relatively relaxed:

**Definition 29.** *For each algebraic move $A \xrightarrow{\omega} A'$ and each crossing $\Delta$ of $A$ we say that: (a) $\Delta$ is eliminated whenever the move is performed on an $i$-gon ($i = 1, 2$) with $\Delta$ as a crossing, (b) $\Delta$ is exchanged whenever the move is performed on a 3-gon with $\Delta$ as a crossing, (c) $\Delta$ is conserved otherwise. For $\omega$ positive, we say that the one or two crossings described in Relations 3.1, 3.5 respectively, are new crossings of $A'$ created by the move.*

*For the exchanged and the conserved crossings of $A$ we say that they are crossings of $A'$ as well, and that they are inherited to $A'$. Whenever $\Delta$ is conserved we allow $\Delta$ to be also the name of the crossing of $A'$ that takes the place of $\Delta$ (that is, the one derived from $\Delta$ by relabeling its labels as described in the definition of the move in §3.1, 3.2, 3.5). In the case of a positive or negative move, we have defined a relabeling function $\mu$ respectively (§3.3.2, 3.3.3, 3.4.2, 3.4.3).*

*For the other kinds of moves we have not defined some kind of relabeling function $\mu$, but we can do so in a fast track fashion beneath. Then for an exchanged or preserved crossing $\Delta \in A$ in a move with relabeling function $\mu$ we shall write either $\Delta$ or $\mu(\Delta)$ for the corresponding crossing of $B$ after the move.*

*So let $A \xrightarrow{\omega} A'$ be an algebraic move and $\tau(A) = 2n$. Let $c$ be a circle with some $2n$ places $M_1, \ldots, M_{2n}$ on it. Although not necessary, we can imagine the points chosen so that as the indices increase, the points move in the negative orientation on the circle; then we say the places for the points are negatively ordered. And again although not necessary, we can put the numbers $1, 2, \ldots, 2n$ on $M_1, \ldots, M_{2n}$ respectively. Whenever $n = 0$, no places are put on $c$. And let $c'$ be a copy of $c$. $c$ corresponds to $A$ and $c'$ to $A'$.*

*- For an $\omega = {}^{al}\Omega_{iso}$ move we do nothing on $c'$ and we say that $c$ and $c'$ with the numbering of all the places $M_i$ represents the identity relabeling $\mu$ of the move that sends each integer $k$ to $k$.*

*- For an $\omega = {}^{al}\Omega_\beta$ move we put the number $1$ on some other place $M_i$ (currently assigned to some integer possibly other than $1$), as demanded in the definition of the ${}^{al}\Omega_\beta$ move. Subsequently,* *if the points of the places are imagined negatively ordered, we assign in the negative orientation all the other integers $2, \ldots, 2n$ to the other places $M_{i+1}, M_{i+2} \ldots, M_{i-1}$ of $c'$ (Figure 51). We say that $c$ and $c'$ with the labels for their places represent the relabeling $\mu$ of $\omega$, which sends the label of any $M_j$ in $c$ to the label of $M_j$ in $c'$. This relabeling $\mu$ is the one we met in §3.1.2.*

*- For an $\omega = {}^{al}\Omega_3$ move we indicate (say by green color) on $c, c'$ the arcs $e_1 = M_{i_1}M_{i_1+1}$, $e_2 = M_{i_2}M_{i_2+1}$, $e_3 = M_{i_3}M_{i_3+1}$ (indices mod $2n$) whose endpoints label the three crossings of the 3-gon on which the move is performed. We then change on $c'$ the labels of the six endpoints of the $e_i$'s as demanded in the definition of the ${}^{al}\Omega_3$ move. There is no need to be exact in our schematic symbolism in Figure 51; indicating the arcs $e_1, e_2, e_3$ on the two circles is enough. We say that $c$ and $c'$ with the new labeling of all places $M_i$ represents the relabeling $\mu$ of $\omega$, which sends the old label of any $M_i$ to the new one.*

*We conveniently extend to negative moves as well:*

*- For an $\omega = {}^{al}\Omega_1^-$ move we indicate on $c, c'$ the arc $e_1 = M_{i_1}M_{i_1+1}$ (indices mod $2n$) whose endpoints label the crossing of the 1-gon on which the move is performed. On $c$ we denote the arc green to indicate its availability for a move on it. We then eliminate on $c'$ the points $M_{i_1}, M_{i_1+1}$, we put the integer 1 at one of the remaining places $M_i$ as demanded in the definition of the ${}^{al}\Omega_1^-$ move (there is no need to be exact in our schematic symbolism in Figure 51), and if the points for the places are imagined negatively ordered, we assign in the negative orientation the integers $2, 3, \ldots, 2(n-1)$ in turn to the other places remaining on $c$. We denote the elimination of the arc by red color, indicating that it was already used in a negative move. We say that $c$ and $c'$ (with arc $e_1$ deleted) with the labels of their places represent the relabeling $\mu$ of $\omega$. This relabeling $\mu$ is the one we met in §3.3.3.*

*- For an $\omega = {}^{al}\Omega_2^-$ move we indicate on $c, c'$ the arcs $e_1 = M_{i_1}M_{i_1+1}, e_2 = M_{i_2}M_{i_2+1}$ (indices mod $2n$) whose endpoints label the crossings of the 2-gon on which the move is performed. On $c$ we denote the arcs green to indicate their availability for a move on them. We then eliminate on $c'$ the points $M_{i_1}, M_{i_1+1}, M_{i_2}, M_{i_2+1}$ (Figure 51), we put the integer 1 at one of the remaining places $M_i$ as demanded in the definition of the ${}^{al}\Omega_2^-$ move (there is no need to be exact in our schematic symbolism in Figure 51), and if the points for the places are imagined negatively ordered, we assign in the negative orientation the integers $2, 3, \ldots, 2(n-2)$ in turn to the other places remaining on $c'$. We denote the elimination of the arcs by red color, indicating that they were already used in a negative move. We say that $c$ and $c'$ (with the arcs $e_1, e_2$ deleted) with the labels of all their places represents the relabeling $\mu$ of $\omega$, which sends the old label of any such $M_i$ to its new label. This relabeling $\mu$ is the one we met in §3.4.3.*

*For each kind of algebraic move $A \xrightarrow{\omega} A'$ from the ones above, we say that the cyclic representation of the relabeling $\mu$ with the help of the circles $c, c'$ also represents the move itself. For the negative moves and the ${}^{al}\Omega_3$ moves, we say that the arcs $e_1$ or $e_1, e_2$ or $e_1, e_2, e_3$ respectively of the $k$-gon on which $\omega$ is performed, represent the $k$-gon itself. If the crossings of a $k$-gon survive in a move, we say that the $k$-gon survives as well, and that it is a $k$-gon of the resulting based symbol.*

*In accordance to the above terminology, we also say that a negative move is performed by deleting the appropriate arcs $e_i$ (or deleting just the places at the endpoints of $e_i$) on the representing circles, and by moving the base integer 1 to an appropriate new place among the remaining ones after the deletions.*

There is no reason to work out a similar representation for the positive moves, since the Lemma 30 below is all we need for our purposes regarding these moves.

After a move $A \xrightarrow{\omega} A'$ the placement of the base integer 1 on the remaining places $M_i$ on $c'$ is important for determining the correct based symbol $A'$. Nevertheless, all based symbols that can be derived from $c'$ by changing the placement of 1, are connected to each other by a change of base move $^{al}\Omega_\beta$, and this is all that is important to us. So in what follows we shall not try to accurately indicate the placements of 1.

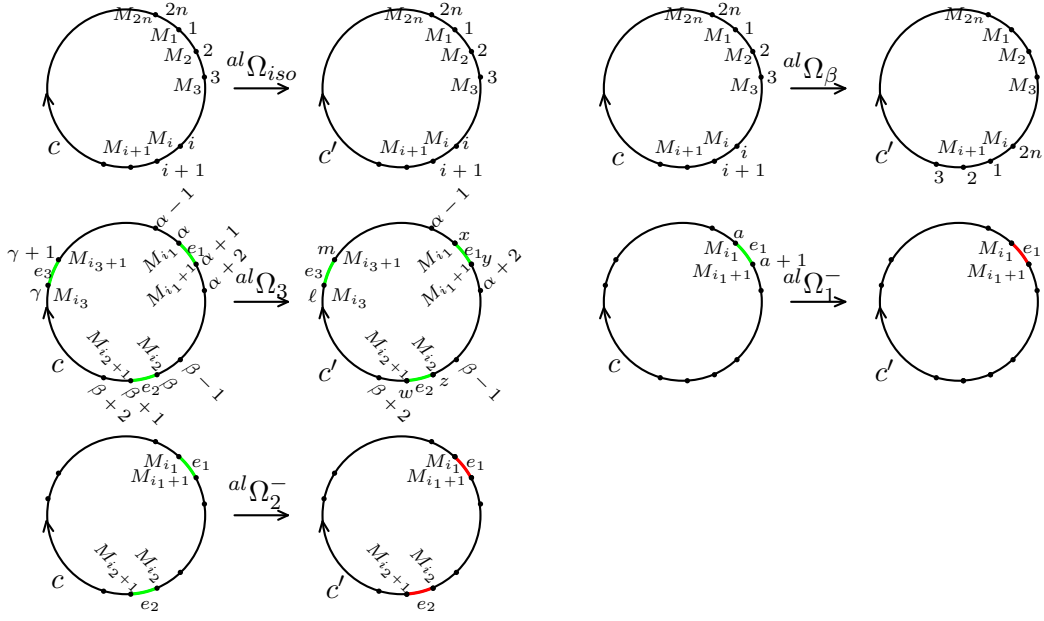

FIGURE 51. Cyclic representation of the relabeling functions of the algebraic moves.

The terminology of the last Definition aims mainly to let us talk freely about $k$-gons that do or do not disappear after a pair of consecutive moves one of which is negative, or after a bunch of consecutive negative moves. Also, it allows us to have a schematic picture of the $k$-gon on each one of the circles corresponding to the based symbols involved in these moves. The picture is a set of $k$ in number small arcs on each circle, always at the same place when we superimpose them. According to the durability of the $k$-gon, its color might change from green to red, never becoming green again.

The Lemma that follows concerns some rather obvious connections between algebraic moves. They hold immediately because of the definitions of the moves. The proof included here makes use of the terminology of the last Definition.

**Lemma 30.** *(a) Let $A \xrightarrow{\omega} A'$ be an algebraic move on $i$-gon $x$ and $c, c'$ the circles of $A, A'$ for this move in the cyclic representation of the move. Iy $y$ is a $k$-gon of $A$ other than $x$, then $y$ is inherited to $A'$ as a $k$-gon on $c'$.*

*(b) If $A \xrightarrow{\omega} A'$ is a positive $^{al}\Omega_i^+$ move, then there exists a negative $^{al}\Omega_i^-$ move $A' \xrightarrow{\omega'} A$ as well.*

*(c) If $A \xrightarrow{\omega} A'$ is a neutral move ($^{al}\Omega_\beta$ or $^{al}\Omega_3$), then there exists a neutral move $A' \xrightarrow{\omega'} A$ of the same type as well.*

*(d) If $A \xrightarrow{\omega} A'$ is a neutral move ($^{al}\Omega_\beta$ or $^{al}\Omega_3$), and $A$ is momentarily reduced (it does not accept a negative move), then $A'$ is momentarily reduced as well.*

*Proof.* (a) By the definition of the $1, 2, 3$-gons, any two of them as sets of basic arcs (or as sets of crossings) are either disjoint or coincide. Since $y$ is given distinct from $x$, it is disjoint from it as a set of crossings. So $y$ retains all its crossings after $\omega$, of course relabeled by the relabeling $\mu$ of the move.

If $x$ is a 3-gon, then $\mu$ does not change the labels of the crossings of $y$, thus they continue to satisfy the conditions in the definition of a $k$-gon. So the image $y'$ of $y$ in $c'$ is a $k$-gon as well.

If $x$ is an $i$-gon for $i = 1, 2$, $\mu$ transforms all remaining labels cyclically after the appropriate determination of the place of $c'$ on which 1 will be put. This simultaneous cyclic relabeling of the labels of the crossings of $y$ trivially retains the relations between them demanded in the definition of the $k$-gon, for any value $k = 1, 2, 3$. So the image $y'$ of $y$ in $c'$ is a $k$-gon as well.

(b) Let $\tau(A) = n$. Then $\tau(A') = n + i$. Let $\mu$ be the relabeling function of $\omega$. By the definition of the relabeling functions, there exists a relabeling function $\lambda = \mu^{-1}$.

The move creates i in number new crossings for $A'$ and also assigns to it each one of the crossings of $A$ with new labels: each label $x$ of any crossing of $A$ now becomes $\mu(x)$.

For $i = 1$ let $\Delta = (a : a+1)$ be the new crossing created by the move. According to the definition, $\Delta$ is the crossing of an 1-gon $T$ of $A'$.

For $i = 2$ let $\Delta_1 = (a : b), \Delta_2 = (a + 1 : c = b \pm 1)$ be the two new crossings created by the move. None of $(a : a+1), (b, c)$ is a crossing of $A'$, thus none of them is 1-gon and according to the description of the crossing, the heights of $a, a+1$ are the same, thus $\Delta_1, \Delta_2$ are the crossings of a 2-gon $T$ of $A'$.

In both cases $i = 1, 2$, we can perform an $\omega' = {}^{al}\Omega_i^-$ move on $T$. The relabeling function of $\omega'$ is $\mu^{-1} = \lambda$. Let $A' \xrightarrow{\omega'} A''$. The result of $\omega'$ on the crossings of $A'$ is to eliminate $\Delta$ or $\Delta_1, \Delta_2$ respectively, and relabel the labels $\mu(x)$ of each one of the remaining crossings to $\lambda(\mu(x)) = \mu^{-1}(\mu(x)) = x$, i.e. to the values they had in the crossings of $A$. So $A'', A$ have the same crossings, thus they are the same based symbol and the move $\omega'$ becomes $A' \xrightarrow{\omega'} A$ which proves the result.

(c) With the cyclic representation of moves developed above the result is immediate. Or we can say that for a relabeling move that sends 1 onto $a$, the relabeling move that sends $a$ to 1 (thus it sends 1 to $2n - a + 2$), is the required move. And if $\omega$ iss an ${}^{al}\Omega_3$ move replaces crossings $(\alpha, \gamma)^{\pi_1} (\beta, \epsilon)^{\pi_2} (\delta, \zeta)^{\pi_3}$ with crossings $(\beta, \delta)^{\pi_1} (\alpha, \zeta)^{\pi_2} (\gamma, \epsilon)^{\pi_3}$, then we define $\omega'$ as the ${}^{al}\Omega_3$ move that replaces crossings $(\beta, \delta)^{\pi_1} (\alpha, \zeta)^{\pi_2} (\gamma, \epsilon)^{\pi_3}$ with crossings $(\alpha, \gamma)^{\pi_1} (\beta, \epsilon)^{\pi_2} (\delta, \zeta)^{\pi_3}$.

(d) By part (a) of the current Lemma, there exists a move $A' \xrightarrow{\omega'} A$ of the same type as $\omega$. Let in the cyclic representation of the moves be $c', c$ respectively the circles of $A', A$ for the move $\omega'$. If there exists a negative move on $A'$, let $x'$ be its $i$-gon.

If $\omega' = {}^{al}\Omega_3$ and $y'$ the 3-gon of the move, then $x'$ and $y'$ are disjoint since different kinds of move are performed on them. So $\omega'$ keeps $x'$ as a copy, say $x$, on $c$ without changing the labels of the places of $x'$. If on the other hand $\omega' = {}^{al}\Omega_\beta$ then at once we have that $\omega'$ keeps $x'$ as a copy, say $x$, on $c$ with an appropriate relabeling of the places of $x'$.

In both cases the labels of the places of $x$ are those of an $i$-gon. Thus a negative move can be performed on $x$ for $A$, a contradiction which proves the result. $\square$

Now we prove a Lemma which we use to expedite subsequent arguments.

**Lemma 31.** *In each one of the following situations, the horizontal move on $A$ happens on an $i$-gon $x$ and the horizontal move on $A'$ happens on the inherited (after the vertical move) $i$-gon $x'$ of $x$ in $A'$; so in particular in cases (b), (c) we assume that the the horizontal and vertical moves on $A$ are performed on disjoint $j$-gons $x$ and $y$ respectively.*

$$
(a) \quad \begin{array}{c} A \xrightarrow{{}^{al}\Omega_i^-} A_1 \\ {}^{al}\Omega_\beta \downarrow \\ A' \xrightarrow{{}^{al}\Omega_i^-} A'_1 \end{array} \qquad (b) \quad \begin{array}{c} A \xrightarrow{{}^{al}\Omega_i^-} A_1 \\ {}^{al}\Omega_3 \downarrow \\ A' \xrightarrow{{}^{al}\Omega_i^-} A'_1 \end{array} \qquad (c) \quad \begin{array}{c} A \xrightarrow{{}^{al}\Omega_i^-} A_1 \\ {}^{al}\Omega_k^- \downarrow \\ A' \xrightarrow{{}^{al}\Omega_i^-} A'_1 \end{array}
$$

*Then there exist vertical moves which "complete the square" as shown:*

$$
(a)\ \begin{array}{ccc} A & \xrightarrow{^{al}\Omega_i^-} & A_1 \\ ^{al}\Omega_\beta \downarrow & & \downarrow ^{al}\Omega_\beta \\ A' & \xrightarrow[^{al}\Omega_i^-]{} & A'_1 \end{array} \quad (b)\ \begin{array}{ccc} A & \xrightarrow{^{al}\Omega_i^-} & A_1 \\ ^{al}\Omega_3 \downarrow & & \downarrow ^{al}\Omega_3 \\ A' & \xrightarrow[^{al}\Omega_i^-]{} & A'_1 \end{array} \quad (c)\ \begin{array}{ccc} A & \xrightarrow{^{al}\Omega_i^-} & A_1 \\ ^{al}\Omega_k^- \downarrow & & \downarrow ^{al}\Omega_k^- \\ A' & \xrightarrow[^{al}\Omega_i^-]{} & A'_1 \end{array}
$$

*In cases (b) and (c), the right vertical move is performed on the inherited $k$-gon $y_1$ of $y$ in $A_1$ (after the top horizontal move).*

*Proof.* Let us call $\omega_1$ the top horizontal move, $\omega'_1$ the bottom horizontal move and $\omega_0$ the vertical move. Let in the cyclic representation of the moves be $c, c_1$ the circles of $A, A_1$ for $\omega_1$, $c', c'_1$ the circles of $A', A'_1$ for $\omega'_1$, and $c, c'$ the circles of $A, A'$ for $\omega_0$. $c', c_1$ are copies of $c$ with some of the places of $c_1$ deleted, and $c'_1$ is a copy of $c'$ with some of the places deleted:

In all cases the $i$-gon $x$ is represented on $c$ by one or two arcs (depending on if $i = 1$ or 2). Similarly, in cases (b) and (c) the $k$-gon $y$ is represented by one or two arcs (depending on if $k = 1$ or 2), and by hypothesis $x, y$ are disjoint. $c_1$ is $c$ with the places on arc $x$ deleted. $c'$ is $c$ with the places on arc $y$ deleted. $c'_1$ is $c'$ with the places on arc $x$ deleted. So $c'_1$ is $c$ with the places on the disjoint arcs $x, y$ deleted.

In case (a) each one of the circles $c_1, c'_1$ contains all places of $c$. After each one the moves $\omega_1, \omega'_1, \omega_0$ the place of the integer 1 may change from its original position 1 on $c$ due to the fact that it may be contained in one of $x, y$. So the label $a$ of a place in $c$ that remains a place of $c_1, c'_1$, say becomes $a', a''$ respectively, with $a', a''$ possibly differing. The exact description of them is as $a' = (\mu(\omega_1))(a)$, $a'' = (\mu(\omega'_1) \circ \mu(\omega_0))(a)$, where $\mu(\omega)$ denotes the relabeling of any move $\omega$ as in Definition 29. So $a'' = (\mu(\omega'_1) \circ \mu(\omega_0)) \circ \mu(\omega_1)^{-1}(a')$. It is a triviality to check that this is the relabeling function of a relabeling move on $A'_1$ (the one that puts the integer 1 of the circle $c_1$ to the place of the integer 1 on the circle $c'_1$). This move sends $A_1$ to a based symbol with the same crossings as $A'_1$, so it is the symbol $A'_1$ and we are done.

In case (b) we perform an $\omega'_0 = s^{al}\Omega_3$ move on arc $y$ of $A_1$, let $A_1 \xrightarrow{\omega'_0} A''_1$. The circle $c''_1$ of $A''_1$ in this move is $c_1$ with the places on $y$ relabeled as demanded by the move. So $c''_1$ is $c$ with the places on $x$ deleted. Also $c'_1$ is $c'$ with the places on $y$ relabeled as demanded by the $^{al}\Omega_3$ move on its 3-gon $y$. So $c'_1$ is also $c$ with the places on $x$ deleted. In other words, $c'_1, c''_1$ contain the same places of $c$. After each one the moves $\omega_1, \omega'_1, \omega_0$ the place of the integer 1 may change from its original position 1 on $c$. So the label $a$ of a place in $c$ that remains a place of $c'_1, c''_1$, say becomes $a', a''$ respectively, with $a', a''$ possibly differing. The exact description of them is as $a'' = (\mu(\omega'_0) \circ \mu(\omega_1))(a)$, $a' = (\mu(\omega'_1) \circ \mu(\omega_0))(a)$, where $\mu(\omega)$ denotes the relabeling of any move $\omega$ as in Definition 29. So $a'' = (\mu(\omega_0) \circ \mu(\omega'_1)) \circ ((\mu(\omega'_1) \circ \mu(\omega_0)))^{-1}(a')$. It is a triviality to check that this is the identity function. So $A''_1$ has the same crossings as $A'_1$, so it is the symbol $A'_1$ and we are done.

In case (c) we perform an $\omega'_0 = ^{al}\Omega_k^-$ move on arc $y$ of $A_1$, let $A_1 \xrightarrow{\omega'_0} A''_1$. The circle $c''_1$ of $A''_1$ in this move is $c_1$ with the places on $y$ deleted. So $c''_1$ is $c$ with the places on $x$ and $y$ deleted. This makes $c'_1, c''_1$ contain the same places of $c$. Also, after each one the moves $\omega_1, \omega'_1, \omega_0$ the place of the integer 1 may change from its original position 1 on $c$. So the label $a$ of a place in $c$ that remains a place of $c'_1, c''_1$, say becomes $a', a''$ respectively with $a', a''$ possibly differing. The exact description of them is as $a'' = (\mu(\omega'_0) \circ \mu(\omega_1))(a)$, $a' = (\mu(\omega'_1) \circ \mu(\omega_0))(a)$, where $\mu(\omega)$ denotes the relabeling of the move $\omega$ as in Definition 29. So $a'' = (\mu(\omega'_0) \circ \mu(\omega_)) \circ ((\mu(\omega'_1) \circ \mu(\omega_0)))^{-1}(a')$. It is a triviality to check that this is the identity function. So $A''_1$ has the same crossings as $A'_1$, so it is the symbol $A'_1$ and we are done. $\square$

We need one more Lemma similar to the last one:

**Lemma 32.** *In each one of the following situations, the vertical move on $A$ happens on a $k$-gon $x$ and the vertical move on $A'$ happens on the inherited (after the horizontal move) $k$-gon $x_1$ of $x$ in $A_1$; so in particular in cases (b), (c) we assume that the the horizontal and vertical moves on $A$ are performed on disjoint $j$-gons $x$ and $y$ respectively.*

$$
\text{(a)} \quad \begin{array}{ccc} A & \xrightarrow{^{al}\Omega_\beta} & A_1 \\ {}^{al}\Omega_k^-\downarrow & & \downarrow {}^{al}\Omega_k^- \\ A' & & A'_1 \end{array} \quad\quad \text{(b)} \quad \begin{array}{ccc} A & \xrightarrow{^{al}\Omega_3} & A_1 \\ {}^{al}\Omega_k^-\downarrow & & \downarrow {}^{al}\Omega_k^- \\ A' & & A'_1 \end{array} \quad\quad \text{(c)} \quad \begin{array}{ccc} A & \xrightarrow{^{al}\Omega_i^-} & A_1 \\ {}^{al}\Omega_k^-\downarrow & & \downarrow {}^{al}\Omega_k^- \\ A' & & A'_1 \end{array}
$$

*Then there exist horizontal moves which "complete the square" as shown:*

$$
\text{(a)} \quad \begin{array}{ccc} A & \xrightarrow{^{al}\Omega_\beta} & A_1 \\ {}^{al}\Omega_k^-\downarrow & & \downarrow {}^{al}\Omega_k^- \\ A' & \xrightarrow[^{al}\Omega_\beta]{} & A'_1 \end{array} \quad\quad \text{(b)} \quad \begin{array}{ccc} A & \xrightarrow{^{al}\Omega_3} & A_1 \\ {}^{al}\Omega_k^-\downarrow & & \downarrow {}^{al}\Omega_k^- \\ A' & \xrightarrow[^{al}\Omega_3]{} & A'_1 \end{array} \quad\quad \text{(c)} \quad \begin{array}{ccc} A & \xrightarrow{^{al}\Omega_i^-} & A_1 \\ {}^{al}\Omega_k^-\downarrow & & \downarrow {}^{al}\Omega_k^- \\ A' & \xrightarrow[^{al}\Omega_i^-]{} & A'_1 \end{array}
$$

*In cases (b) and (c), the bottom horizontal move is performed on the inherited $i$-gon $y'$ of $y$ in $A'$ (after the left vertical move).*

*Proof.* Similar to the proof of Lemma 32. □

**Lemma 33.** *(a) If $A \xrightarrow{\omega} A'$ is an algebraic move then $\overline{R}(A) = \overline{R}(A')$.*
*(b) If $A \underset{al}{\sim} A'$ then $\overline{R}(A) = \overline{R}(A')$.*

*Proof.* (a) We prove it case by case for each one of the five kinds of algebraic moves.

(I) For $\omega = {}^{al}\Omega_{iso}$:

It is $A = A'$ and the truth of the result is immediate. Since ${}^{al}\Omega_{iso}$ moves is a special case of ${}^{al}\Omega_\beta$ moves, the results can also come from the $\omega = {}^{al}\Omega_\beta$ case proved immediately below.

(II) For $\omega = {}^{al}\Omega_\beta$:

If no negative move can be made on $A$, then for the arbitrary $B \in \overline{R}(A)$ there exists by the definition a neutral sequence from A to B. Let us denote it $A \xrightarrow{\omega_1,\ldots,\omega_n} B$. Also, by Lemma 30 the existence of $A \xrightarrow{\omega} A'$ implies the existence of another change of base move $A' \xrightarrow{\omega'} A$. Thus $A' \xrightarrow{\omega'} A \xrightarrow{\omega_1,\ldots,\omega_n} B$ is a neutral sequence showing $B \in \overline{R}(A')$. Hence $\overline{R}(A) \subseteq \overline{R}(A')$. Similarly $\overline{R}(A') \subseteq \overline{R}(A)$, so $\overline{R}(A) = \overline{R}(A')$.

If there exist negative moves on $A$, then for the arbitrary $B \in \overline{R}(A)$ let $A = A_0 \xrightarrow{\omega_1} A_1 \xrightarrow{\omega_2} A_2 \cdots \xrightarrow{\omega_n} A_n \xrightarrow{\overline{\omega}_1,\ldots,\overline{\omega}_m} B$ be a full negative sequence of moves $\omega_i$ followed by a sequence of neutral moves $\overline{\omega}_j$ as asked by the definition of $\overline{R}(A)$. We shall show the existence of a negative sequence, say $f$, from $A'$ to some $A'_n$, so that $A'_n$ is connected to $A_n$ with an $\omega' = {}^{al}\Omega_\beta$ move $A'_n \xrightarrow{\omega'} A_n$. Then $A' = A'_0 \xrightarrow{f} A'_n \xrightarrow{\omega'} A_n \xrightarrow{\overline{\omega}_1,\ldots,\overline{\omega}_m} B$ would be a negative sequence followed by a neutral sequence from $A'$ to $B$. The $A'_n \xrightarrow{\omega'={}^{al}\Omega_\beta} A_n$ move and the momentarily reduced character of $A_n$ imply by Lemma 30 that $A'$ is also momentarily reduced. So the above negative sequence leading to $B$ is full, showing that $B \in \overline{R}(A')$. Then as above, we would have the required equality $\overline{R}(A) = \overline{R}(A')$. For the claim about the existence of $f$ and $\omega'$:

In the cyclic representation of the moves $\omega_i$ according to Definition 29, let $c_0, c_1, \ldots, c_n$ be the circles corresponding to $A_0, A_1, \ldots, A_n$ respectively. Let $x_0, x_1, \ldots, x_{n-1}$ be the k-gons of $c_0, c_1, \ldots, c_{n-1}$ on which the $\omega_1, \ldots, \omega_n$ move is performed. And let for definiteness that they are $k_0, k_1, \ldots, k_{n-1}$-gons respectively.

The given move $\omega = {}^{al}\Omega_\beta$ is represented by the circle $c_0$ of $A = A_0$ and a copy $c'_0$ of $c$ for $A' = A'_0$. $c'_0$ hosts a copy $x'_0$ of $x_0$ of relabeled places which is a $k_0$-gon of $A'_0$. We perform on $x'_0$ an $\omega'_1 = {}^{al}\Omega^-_{k_0}$ move $A'_0 \xrightarrow{\omega'_1} A'_1$. By Lemma 31 there exists an $\Omega_\beta$ move $A_1 \xrightarrow{\omega_{(1)}} A'_1$.

We repeat then this reasoning and construct moves $\omega'_2, \omega_{(2)}, \omega'_3, \omega_{(3)}, \ldots$ in turn so that $\omega'_i = {}^{al}\Omega^-_{k_i}, \omega_{(i)} = {}^{al}\Omega_\beta$:

$$\begin{array}{ccccccccc}
A = A_0 & \xrightarrow{\omega_1} & A_1 & \xrightarrow{\omega_2} & A_2 & \cdots & \xrightarrow{\omega_n} & A_n & \xrightarrow{\overline{\omega}_1,\ldots,\overline{\omega}'_m} B \\
\downarrow \omega & & \downarrow \omega_{(1)} & & \downarrow \omega_{(2)} & & & \downarrow \omega_{(n)} & \\
A' = A'_0 & \xrightarrow{\omega'_1} & A'_1 & \xrightarrow{\omega'_2} & A'_2 & \cdots & \xrightarrow{\omega'_n} & A'_n &
\end{array}$$

By Lemma 30 and the fact that $A_n$ is momentarily reduced, we have that $A'_n$ is momentarily reduced as well. The same lemma provides the existence of the inverse $\omega'_{(n)}$ of the last move $\omega_{(n)}$ assuring it is again an ${}^{al}\Omega_\beta$ move. Then the sequence of moves $\omega'_1, \ldots, \omega'_n, \omega'_{(n)}, \overline{\omega}_1, \ldots, \overline{\omega}'_m$ shows the result.

(III) For $\omega = {}^{al}\Omega_3$:

We work as in case (II):

If no negative move can be made on $A$, then for the arbitrary $B \in \overline{R}(A)$ there exists by the definition a neutral sequence from A to B. Let us denote it $A \xrightarrow{\omega_1,\ldots,\omega_n} B$. Also, the existence of $A \xrightarrow{\omega} A'$ implies by Lemma 30 the existence of another ${}^{al}\Omega_3$ move $A' \xrightarrow{\omega'} A$. Thus $A' \xrightarrow{\omega'} A \xrightarrow{\omega_1,\ldots,\omega_n} B$ is a neutral sequence showing $B \in \overline{R}(A')$. Hence $\overline{R}(A) \subseteq \overline{R}(A')$. Similarly $\overline{R}(A') \subseteq \overline{R}(A)$, so $\overline{R}(A) = \overline{R}(A')$.

If there exist negative moves on $A$, then for the arbitrary $B \in \overline{R}(A)$ let $A = A_0 \xrightarrow{\omega_1} A_1 \xrightarrow{\omega_2} A_2 \cdots \xrightarrow{\omega_n} A_n \xrightarrow{f} B$ be a full sequence of negative moves $\omega_i$ followed by a sequence $\overline{\omega}_1, \ldots, \overline{\omega}_m$ of neural moves as asked by the definition of $\overline{R}(A)$.

In the cyclic representation of the moves $\omega_i$, let $c_0, c_1, \ldots, c_n$ be the circles corresponding to $A_0, A_1, \ldots, A_n$ respectively. Let $x_0, x_1, \ldots, x_{n-1}$ be the k-gons of $c_0, c_1, \ldots, c_{n-1}$ on which the $\omega_1, \ldots, \omega_n$ move is performed. And let for definiteness that they are $k_0, k_1, \ldots, k_{n-1}$-gons respectively.

Let the given move $\omega = {}^{al}\Omega_3$ be performed on the 3-gon $y = y$. $\omega$ is represented by the circle $c_0$ of $A = A_0$ and a copy $c'_0$ of $c$ for $A' = A'_0$ so that the places on the copy $y'$ of $y$ are relabeled as the move demands. Since the moves $\omega_1, \omega$ which can be performed on $A$ are of different type, their corresponding $i$-gons are disjoint (by the definition of the $i$-gons). That is: $x_0, y$ are disjoint. Then the copy $x'_0$ of $x_0$ in $c'_0$ is disjoint from $y'$ as well. The labels of the endpoints of $x'_0$ are those of the corresponding endpoints of $x_0$, thus $x_0$ is a $k_0$-gon of $A'_0$. We perform on $x'_0$ an $\omega'_1 = {}^{al}\Omega^-_{k_0}$ move $A'_0 \xrightarrow{\omega'_1} A'_1$. By Lemma 31 there exists an ${}^{al}\Omega_3$ move $A_1 \xrightarrow{\omega_{(1)}} A'_1$.

We repeat then this reasoning and construct moves $\omega'_2, \omega_{(2)}, \omega'_3, \omega_{(3)}, \ldots$ in turn so that $\omega'_i = {}^{al}\Omega^-_{k_i}, \omega_{(i)} = {}^{al}\Omega_3$ and:

$$\begin{array}{ccccccccc}
A = A_0 & \xrightarrow{\omega_1} & A_1 & \xrightarrow{\omega_2} & A_2 & \cdots & \xrightarrow{\omega_n} & A_n & \xrightarrow{\omega'_1,\ldots,\omega'_m} B \\
\downarrow \omega & & \downarrow \omega_{(1)} & & \downarrow \omega_{(2)} & & & \downarrow \omega_{(n)} & \\
A' = A'_0 & \xrightarrow{\omega'_1} & A'_1 & \xrightarrow{\omega'_2} & A'_2 & \cdots & \xrightarrow{\omega'_n} & A'_n &
\end{array}$$

The ${}^{al}\Omega_3$ move $A_n \xrightarrow{\omega_{(n)}} A'_n$ implies the existence of an ${}^{al}\Omega_3$ move $A'_n \xrightarrow{\omega'_{(n)}} A_n$, and then the sequence $A' = A'_0 \xrightarrow{\omega'_1} A'_1 \xrightarrow{\omega'_2} A'_2 \cdots \xrightarrow{\omega'_n} A'_n \xrightarrow{\omega'_{(n)}} A_n \xrightarrow{\overline{\omega}_1,\ldots,\overline{\omega}_m} B$ is a negative sequence followed by a neutral sequence from $A'$ to $B$. Moreover since $A_n$ is momentarily reduced, Lemma 30 says $An'$ is also momentarily reduced, thus the sequence of negative moves from $A'$ to $A'_n$ is full. Hence $\overline{R}(A) \subseteq \overline{R}(A')$. Similarly we get $\overline{R}(A') \subseteq \overline{R}(A)$, thus $\overline{R}(A) = \overline{R}(A')$ as wanted.

(IV) For $\omega = {}^{al}\Omega_1^-, {}^{al}\Omega_2^-$:

The proof is quite similar to that of cases (II), (III). In detail:

There exist negative moves on $A$ ($\omega$ is such a one), so for the arbitrary $B \in \overline{R}(A)$ let $A = A_0 \xrightarrow{\omega_1} A_1 \xrightarrow{\omega_2} A_2 \cdots \xrightarrow{\omega_n} A_n \xrightarrow{\overline{\omega}_1,\ldots,\overline{\omega}_m} B$ be a full sequence of negative moves $\omega_i$ followed by a sequence $\overline{\omega}_1,\ldots,\overline{\omega}_m$ of neural moves as asked by the definition of $\overline{R}(A)$.

In the cyclic representation of the moves $\omega_i$, let $c_0, c_1, \ldots, c_n$ be the circles corresponding to $A_0, A_1, \ldots, A_n$ respectively. Let $x_0, x_1, \ldots, x_{n-1}$ be the k-gons of $c_0, c_1, \ldots, c_{n-1}$ on which the $\omega_1, \ldots, \omega_n$ move is performed. And let for definiteness that they are $k_0, k_1, \ldots, k_{n-1}$-gons respectively.

Let the given move $\omega = {}^{al}\Omega_k$ be performed on the k-gon $y = y_0$. If $y$ is not $x_0$ then it is disjoint from it (by the definition of the $i$-gons) and it is inherited in $c_1$ as a copy $y_1$. If $y$ is not $x_1$ then it is inherited in $c_2$ as a copy $y_2$ and so on. If it happens that $y$ survives to $A_n$, then we can perform a negative move on this k-gon of $A_n$, a contradiction to our assumptions. So $y$ does not survive to $A_n$ which means that for some $t$, $0 \leq t \leq n-1$ the copy $y_t$ of $y$ in $c_t$ is actually $x_t$: $y_t = x_t$.

If it happens that $t = 0$, then $y_0 = x_0$, so $A' = A_1$ and then $A' = A_1 \xrightarrow{\omega_2} A_2 \cdots \xrightarrow{\omega_n} A_n \xrightarrow{f} B$ is a full negative sequence followed by a neutral sequence from $A'$ to $B$ proving $\overline{R}(A) \subseteq \overline{R}(A')$. Similarly we get $\overline{R}(A') \subseteq \overline{R}(A)$, thus $\overline{R}(A) = \overline{R}(A')$ as wanted.

If $t \geq 1$, the moves $\omega_1, \omega$ which can be performed on $A$ are different, hence their corresponding $i$-gons are disjoint (by the definition of the $i$-gons). That is: $x_0, y$ are disjoint. Then the copies $x'_0, y'_0$ respectively of $x_0, y_0$ in $c'_0$ are disjoint as well. $x'_0$ is the only k-gon deleted from $c'_0$, thus $y'_0$ remains on $c'_0$. So we can perform on $y'_0$ an $\omega'_1 = {}^{al}\Omega_{k_0}^-$ move $A'_0 \xrightarrow{\omega'_1} A'_1$. By Lemma 31 there exists an $\Omega_k^-$ move $A_1 \xrightarrow{\omega_{(1)}} A'_1$ performed on $y_1$.

We repeat then this reasoning and construct moves $\omega'_2, \omega_{(2)}, \omega'_3, \omega_{(3)}, \ldots, \omega_{(t)}$ in turn so that $\omega'_i = {}^{al}\Omega_{k_i}^-, \omega_{(i)} = {}^{al}\Omega_k^-$:

$$
\begin{array}{ccccccccccccc}
A = A_0 & \xrightarrow{\omega_1} & A_1 & \xrightarrow{\omega_2} & A_2 & \cdots & \xrightarrow{\omega_t} & A_t & \xrightarrow{\omega_{t+1}} & A_{t+1} & \xrightarrow{\omega_{t+2}} & \cdots A_n & \xrightarrow{\overline{\omega}_1,\ldots,\overline{\omega}_m} & B \\
\downarrow \omega & & \downarrow \omega_{(1)} & & \downarrow \omega_{(2)} & & & \downarrow \omega_{(t)} & & & & & & \\
A' = A'_0 & \xrightarrow{\omega'_1} & A'_1 & \xrightarrow{\omega'_2} & A'_2 & \cdots & \xrightarrow{\omega'_t} & A'_t = A_{t+1} & & & & & &
\end{array}
$$

with the vertical moves performed on the $x_i$'s and the horizontal moves performed on the $y_i$'s. Since $x_t = y_t$, it is $A'_t = A_{t+1}$. Then the sequence $A' = A'_0 \xrightarrow{\omega'_1} A'_1 \xrightarrow{\omega'_2} A'_2 \cdots \xrightarrow{\omega'_t} A'_t = A_{t+1} \xrightarrow{\omega_{t+2}} A_{t+2} \cdots \xrightarrow{\omega'_{(n)}} A_n \xrightarrow{\overline{\omega}_1,\ldots,\overline{\omega}_m} B$ is a full negative sequence ($A_n$ is momentarily reduced) followed by a neutral sequence from $A'$ to $B$ which implies as above that $\overline{R}(A) = \overline{R}(A')$ as wanted.

(V) By Lemma 30, the ${}^{al}\Omega_i^+$ move $A \xrightarrow{\omega} A'$ implies the existence of an ${}^{al}\Omega_i^-$ move $A' \xrightarrow{\omega} A$. Then from the (IV) case just above we have $\overline{A'} = \overline{A}$ which is the desired result.

(b) Since for any two equivalent symbols $A, A'$ there exists a finite sequence of moves starting with one and ending with the other, the result comes immediately from part (a) of the Lemma. □

According to Lemma 33, the locally reduced set $\overline{R}(A)$ is the same for all based symbols equivalent to $A$, so it is a notion of the whole class of $A$:

**Definition 30.** *If $C$ is any class in the equivalence $\underset{al}{\sim}$, then we define the locally reduced set $\overline{R}(C)$ of $C$ as $\overline{R}(C) = \overline{R}(A)$ for any based symbol $A \in C$ (that is: for any based symbol $A$ for which $C = [\![A]\!]_{al}$).*

**Lemma 34.** *(a) For any based symbol $A$ and any $B \in \overline{R}(A)$, there exists a finite sequence of non-positive algebraic moves from $A$ to $B$.*
*(b) Any momentarily reduced based symbol is a reduced symbol.*

*Proof.* (a) For any $B \in \overline{R}(A)$ it is $A \underset{al}{\sim} B$, hence there exists a finite sequence of algebraic moves $A = A_0 \xrightarrow{\omega_1} A_1 \xrightarrow{\omega_2} A_2 \cdots \xrightarrow{\omega_n} A_n = B$.

If no positive move exists among the $\omega_i$'s then this above sequence shows the result.

If there exist positive moves, then there is a last one among them, say $A_\ell \xrightarrow{\omega_{\ell+1}} A_{\ell+1}$. This move creates a $k$-gon $y$ on $B_{\ell+1}$ for some $k = 1$ or $2$. $y$ is inherited to each subsequent $A_i$ up to the point when an $^{al}\Omega_k^-$ move is performed on it to make it disappear: if no such move exists among the subsequent $\omega_i$'s, then $y$ is inherited to $A_n = B$ and $B$ can reduce its order further by an $^{al}\Omega_k^-$ move on $y$, a contradiction to the reduced character of $B$.

So let $A_m \xrightarrow{\omega_{m+1} = {}^{al}\Omega_k^-} A_{m+1}$ be a move ($m > \ell$) performed in the inherited $k$-gon of $y$ in $A_m$; say $y_{l+2}, \ldots, y_m$ be the $k$-gon $y$ as inherited in $A_{\ell+1}, A_{\ell+2}, \ldots, A_m$ respectively. The $\omega_i$ moves from $A_\ell$ to $A_{m+1}$ start with a positive one, they end with a negative one and they are neutral in between. We shall show that we can eliminate the first and the last one and replace each neutral by a new one of the same type, still starting from $A_\ell$ and ending to $A_{m+1}$. Then inductively we can eliminate all positive moves and we are done. Indeed:

On the $k$-gon $y = y_{\ell+1}, y_{l+2}, \ldots, y_m$ we perform an $^{al}\Omega_k^-$ move and get:

$$
A_\ell \xrightarrow{\omega_{\ell+1} = {}^{al}\Omega_i^+} \begin{array}{c} A_{\ell+1} \\ \downarrow \omega_{(\ell+1)} \\ A'_\ell = A_\ell \end{array} \xrightarrow{\omega_{\ell+2}} \begin{array}{c} A_{\ell+2} \\ \downarrow \omega_{(\ell+2)} \\ A'_{\ell+3} \end{array} \xrightarrow{\omega_{\ell+3}} \begin{array}{c} A_{\ell+3} \\ \downarrow \omega_{(\ell+3)} \\ A'_{\ell+3} \end{array} \cdots \xrightarrow{\omega_m} \begin{array}{c} A_m \\ \downarrow \omega_{(m)} \\ A'_m \end{array} \xrightarrow{\omega_{m+1} = {}^{al}\Omega_k^-} A_{m+1}
$$

The negative $\omega_{(\ell+1)}$ move on $A_{\ell+1}$ is performed on the $k$-gon $y$ created by the positive move $A_\ell \xrightarrow{\omega_{\ell+1}} A_{\ell+1}$, thus returning us to the symbol $A_\ell$ as if the two moves were never performed. So $A'_\ell = A_\ell$ as shown in the diagram.

Now Lemma 32 implies the existence of the lower horizontal moves in the following diagram. They are of the same type as the corresponding top horizontal moves (which in particular says that the new moves are non-positive):

$$
A_\ell \xrightarrow{\omega_{\ell+1} = {}^{al}\Omega_i^+} \begin{array}{c} A_{\ell+1} \\ \downarrow \omega_{(\ell+1)} \\ A'_\ell = A_\ell \end{array} \xrightarrow{\omega_{\ell+2}} \begin{array}{c} A_{\ell+2} \\ \downarrow \omega_{(\ell+2)} \\ A'_{\ell+3} \end{array} \xrightarrow{\omega_{\ell+3}} \begin{array}{c} A_{\ell+3} \\ \downarrow \omega_{(\ell+3)} \\ A'_{\ell+3} \end{array} \cdots \xrightarrow{\omega_m} \begin{array}{c} A_m \\ \downarrow \omega_{(m)} \\ A'_m = A_{m+1} \end{array} \xrightarrow{\omega_{m+1} = {}^{al}\Omega_k^-} A_{m+1}
$$

The negative $\omega_{(m)}$ move on $A_m$ is performed on the inherited $k$-gon $y_m$ of $y$. The same happens for the negative $\omega_{m+1}$ move on $A_m$. So they are the same move, thus resulting to the same based symbol. So $A'_m = A_{m+1}$ as shown in the diagram.

Then the sequence $A_\ell \xrightarrow{\omega'_{\ell+2}} A'_{\ell+2} \xrightarrow{\omega'_{\ell+3}} A'_{\ell+3} \to \cdots \xrightarrow{\omega'_m} A'_m = A_{m+1}$ from $A_\ell$ to $A_{m+1}$ consists of non-positive moves, proving our claim and the proposition.

(b) Let $A$ be a momentarily reduced based symbol, that is one that accepts no negative moves. And let $B$ be some element in $\overline{R}(A)$. By part (a) there exists a sequence $f$ of non-positive moves from $A$ to $B$. If there exists some negative move among them, let the first one happen to the symbol $A'$ (a symbol among those produced by the sequence). And let this negative move happen on a $k$-gon $x$ of $A'$. Then the moves from $A$ to $A'$ are neutral. But each neutral move reverses to another neutral move by Lemma 30. So there exists a sequence of neutral moves from $A'$ to $A$. Each one of them inherits $x$ as a $k$-gon to the next symbol as a $k$-gon each time, so eventually $A$ contains $x$ as a $k$-gon and we can perform on it a negative move, a contradiction. So no negative move exists among those in $f$. This implies $\tau(A) = \tau(B)$ and since $B$ is reduced this implies $A$ is reduced too. □

We are now in position to prove our original goal that reduced and locally reduced symbols are the same.

**Lemma 35.** *For any based symbol $A$ it holds $R(A) = \overline{R}(A)$.*

*Proof.* Let $B$ be an arbitrary element of $R(A)$. $B$ as a reduced symbol accepts no negative moves, and since $B \xrightarrow{\Omega_{iso}} B$ is a neutral sequence starting from $B$ and ending to $B$. We have $B$ belongs to $\overline{R}(B)$. Now since $A \underset{al}{\sim} B$ we have $\overline{R}(A) = \overline{R}(B)$ by Lemma 33. So then $B \in \overline{R}(B) = \overline{R}(A)$ which implies $R(A) \subseteq \overline{R}(A)$.

Let now $C$ be an arbitrary element of $\overline{R}(A)$. Then there exists a full negative sequence $f$ from $A$ to some momentarily reduced symbol $A'$ followed by a neutral sequence $g$ from $A'$ to $C$. By Lemma 34 $A'$ is a reduced symbol. Since the sequence $g$ is neutral, it follows that $\tau(A') = \tau(C)$. Since $A'$ is reduced, this equality of orders implies that $C'$ is reduced as well. So $C' \in R(A)$. Hence $\overline{R}(A) \subseteq R(A)$. Pairing this with $R(A) \subseteq \overline{R}(A)$ we get $R(A) = \overline{R}(A)$ as wanted. □

Lemma 35 essentially says that among all equivalent symbols to a given one, say $A$, the set of the reduced symbols $R(A)$ is the same as the set of the locally reduced symbols $\overline{A}$ and also the same as the set of momentarily reduced symbols which we did not name so far, say $MR(A)$, that is: $R(A) = \overline{A} = MR(A)$. Also, that we can go from any element of this set to any other via a sequence of neutral moves (ch. Corollary 5 in a moment).

According to Lemmata 35 and 33, the reduced set $\overline{R}(A)$ is the same for all based symbols equivalent to $A$, so it is a notion of the whole class of $A$ and actually coincides with that of the locally reduced set:

**Definition 31.** *If $C$ is any class in the equivalence $\underset{al}{\sim}$, then we define the reduced set $R(C)$ of $C$ as $R(C) = R(A)$ for any based symbol $A \in C$ (that is: for any based symbol $A$ for which $C = [\![A]\!]_{al}$). More specifically it is:*

$$R(C) = R([\![A]\!]_{al}) = R(A) = \overline{R}(C) = \overline{R}([\![A]\!]_{al}) = \overline{R}(A).$$

As an immediate corollary of the above work we have:

**Corollary 5.** *(a) For any based symbol $A$ and any $B \in R(A) = \overline{R}(A)$, there exists a finite sequence of non-positive moves from $A$ to $B$. Such a sequence contains no negative moves exactly when $A$ is reduced. $A$ is reduced exactly when no negative move exists on it, i.e. whenever it is momentarily reduced.*

*(b) For any based symbol $A$, the set $R(A)$ is finite.*

*(c) For any based symbol $A$ which is not reduced, there exists a negative sequence from $A$ to a reduced symbol in $R(A)$.*

*(d) For any reduced based symbol $A$ and any element $B$ of $R(A)$, there exists a neutral sequence from $A$ to $B$.*

*(e) For any based symbol $A$, we can calculate $R(A)$ algorithmically in finitely many steps as follows:*

*step 1.* Check for k-gons $x$ of $A$ for $k = 1, 2$. If there exists some, perform a negative move $A \xrightarrow{^{al}\Omega_k^-} A'$ on it, replace $A$ with $A'$ and repeat. If there exists no $x$'s, go to step 2.

*step 2.* Set $i = 1$, $R_i = \{A'\}$ where $A'$ is the symbol to which we arrived in at the end of step 1. For the unique element $A'$ of $R_i$ perform all neutral moves $\omega = {}^{al}\Omega_\beta, {}^{al}\Omega_3$ $A \xrightarrow{\omega} A'$. Replace $i$ by $i + 1$, set $R_{i+1} = R_i \cup R_i'$ where $R_i'$ is the set containing all the results $B$ of the moves $\omega$, and repeat. Stop whenever for first time $R_{m+1} = R_m$.

*step 3.* Write $R(A) = R_1 \cup R_2 \cdots \cup R_m$.

*Proof.* (a) The existence of non-positive moves from $A$ to $B$ is part (a) of Lemma 34. Whenever $A$ is reduced, such a sequence of course cannot contain a negative move since such a move would reduce the order below that of $A$ which is minimum.

If $A$ is momentarily reduced and $B$ some element of $R(A)$, then there exists a full sequence of moves from $A$ to some $A'$ followed by a neutral sequence ending at $B$. Since no negative move exists for $A$, the full negative sequence is empty, so the order of $A$ is that of $B$, thus $A$ is reduced too.

(b) By definition, each member of $R(A)$ has the minimum order a based symbol equivalent to $A$ can have, so all elements of $R(A)$ have the same order, say $n$. If $n = 0$. then $R(A) = \{\emptyset\}$, a finite set, whereas if $n > 0$, let us just note that there exist finitely only based symbols of order $n$: each based symbol of order $n$ is a collection of $n$ crossings and each crossing is a pair of distinct integers in $\mathbb{Z}_{2n}$ with all integers in all crossings giving $\mathbb{Z}_{2n}$ as their union. Such partitions of $\mathbb{Z}_{2n}$ can be done in only finitely many ways. The choice of sign $+1$ or $-1$ for each crossing that completes the description of the crossings, keeps the possibilities still finite in number.

(c) If $A$ contained no negative move, then it would be reduced by part (b), a contradiction by our hypothesis for $A$. Performing a negative move on $A \xrightarrow{\omega_1} A'$ we get a new symbol $A'$ which by the same reasoning is either a reduced symbol or else it accepts some negative move $\omega_2$ and so on. Since each negative move decreaces the order of the symbol and since the orders are greater or equal to 0, we eventually arrivve by a finite sequence $\omega_1, \omega_2, \ldots$ to some symbol with minimum order which of course is reduced and lies in $R(A)$ as wanted.

(d) By part (b) there exists a sequence of non-positive moves from $A$ to $B$, so for all symbols $A = A_0, A_1, \ldots, A_n = B$ in this sequence, it is $\tau(A_i) \leq \tau(A)$. Since $\tau(A)$ is the smallest possible order for the equivalent symbols to $A$, it must be $\tau(A_i) = \tau(A)$ which means that all $A_i$ (including $A_n = B$) must be reduced.

(e) Step 1 will be over after a finite number of repetitions because of part (c). Let $A'$ the based symbol we have arrived to.

Let us now observe that each based symbol has only finitely many 3-gons, so only finitely many ${}^{al}\Omega_3$ moves can be performed on it. Also, it has only finitely many labels, so only finitely many change of base moves. These, along with the fact that $R_1$ is finite imply $R_2$ is finite. Similarly then $R_3$ is finite and inductively, each $R_i$ is finite. This says that each repetition of step 2 terminates after a finite number of performed moves $\omega$.

We claim that after a finite number of repetitions of step 2, we arrive at a set $R_{i+1}$ equal to the one at hand (that is to $R_i$): each element of $R_2$ is connected by a neutral sequence to the symbol $A'$ at the end of step 1, thus the order of each one of them is that of $A'$. Similarly each element of $R_2$ has the order of $A'$ and so on for each element of any $R_i$. But in part (b) of the current Lemma we observed that there exist only finitely many based symbols of a given order. We conclude that inevitably, for an index $i$ it will be $R_{i+1} = R_i$. So step 2 will indeed terminate after a finite number of repetitions.

The symbol $A'$ we have arrived the moment step 1 terminates, contains no negative moves (by the definition of the termination), thus by Lemma 35 it is a reduced symbol, hence $A' \in R(A)$.

It is $R_{i+1} = R_1 \cup \cdots \cup R_i \Rightarrow R_i \subseteq R_{i+1}$, $\forall i$. Let $m$ be the smallest value of $i$ for which $R_{i+1} = R_i$. Then $R_{m+2} = R_1 \cup \cdots \cup R_m \cup R_{m+1} = R_1 \cup \cdots \cup R_m \cup R_m = R_1 \cup \cdots \cup R_m = R_{m+1}$ and inductively $R_m = R_{m+1} = R_{m+2} = \cdots$. So $R_1 \cup R_2 \cdots \cup R_m$ is the set of all based symbols which we can get from $A'$ by neutral sequences (finite sequences of neutral moves) starting with $A'$. Since $A'$ is reduced, we know by part (b) that we can arrive at each element of $R(A')$ by a neutral sequence. Hence $R(A') \subseteq R_1 \cup R_2 \cdots \cup R_m$. But $A, A'$ are equivalent based symbols, so $R(A) = R(A')$ according to Lemmata 33 and 35. Since each element of $R_1 \cup R_2 \cdots \cup R_m$ has the order of $A'$, and since $A'$ is reduced, all elements of $R_1 \cup R_2 \cdots \cup R_m$ are reduced, thus $R_1 \cup R_2 \cdots \cup R_m \subseteq R(A') = R(A)$. Hence $R(A) = R_1 \cup R_2 \cdots \cup R_m$ as wanted. $\square$

6.2. **Definition of the invariant $R$ and its algorithmic computability.**

**Definition 32.** *(The reduced symbol invariant of knots) For the (smooth or pl) knots $\varkappa$ in space ($\mathbb{R}^3$ or $S^3$) we call as the reduced symbol of $\varkappa$ and we write $R(\varkappa)$, the set $R(A)$ for any representative based symbol $A$ for the knot $\varkappa$.*

*If $\kappa$ is a (smooth or pl) knotted loop in space respresenting the knot $\varkappa$, $K$ the projection of $\kappa$ on a plane, and $p$ some base point of $K$ then:*

$$R(\kappa) = R([\![_p K]\!]_{al}) = R(_p K).$$

*For any base point $p$ of $K$ we also define as $R(_p K)$ the reduced symbol $R(\kappa)$ of the knotted loop $\kappa$, the reduced symbol $R(K)$ of the diagram $K$, and the reduced symbol $R([\![K]\!]_{top})$ of the equivalence $[\![K]\!]_{top}$ in $\underset{top}{\sim}$. With this notation it holds:*

$$R(\varkappa) = R(\kappa) = R(K) = R([\![K]\!]_{top}) = R([\![_p K]\!]_{al}) = R(_p K).$$

*According to Lemma 36 that follows, the reduced symbol set is a complete knot invariant in the sense that it distinguishes any two oriented knots. We call it as the reduced symbol invariant of knots. For each knot it is a set of based symbols calculated from any knot diagram of a given knot algorithmically in finite time, so we call the invariant as effectively computable.*

As we commented earlier, the symbol invariant $S(\varkappa)$ of a knot $\varkappa$ is a really huge set containing infinitely many representatives (each a based symbol of a diagram representing the knot), due to the fact that we can always perform an endless sequence of positive moves, even if we start from the simplest of all based symbols, the empty set. The reduced symbol $R(\varkappa)$ on the other hand does not have this handicap:

**Lemma 36.** *The symbol of knots is a complete knot invariant in the sense that it distinguishes any two oriented knots. For distinct knots, their reduced symbols are disjoint sets of based symbols. It is also algorithmically computable according to the algorithm of Corollary 5, which ends in a finite number of steps. For convenience we repeat the algorithm here:*

*Consider a knot diagram $K$ of $\kappa$ and a base point $p$ on it.*

*step 1. Check for k-gons $x$ of $_p K$ for $k = 1, 2$. If there exists some, perform a negative move $_p K \xrightarrow{al \Omega_k^-} A'$ on it, replace $_p K$ with $A'$ and repeat. If there exists no $x$'s, go to step 2.*

*step 2. Set $i = 1$, $R_i = \{A'\}$ where $A'$ is the symbol to which we arrived in at the end of step 1. For the unique element $A'$ of $R_i$ perform all neutral moves $\omega = {}^{al}\Omega_\beta, {}^{al}\Omega_3$ $A \xrightarrow{\omega} A'$. Replace $i$ by $i + 1$, set $R_{i+1} = R_i \cup R'_i$ where $R'_i$ is the set containing all the results $B$ of the moves, and repeat. Stop whenever for first time $R_{m+1} = R_m$.*

*step 3. Write $R(\kappa) = R(_p K) = R_1 \cup R_2 \cdots \cup R_m$.*

*Proof.* By Theorem 4 the invariant $S$ distinguishes knots, that is, two knots $\varkappa_1, \varkappa_2$ are the same or distinct according to if their symbols $S(\varkappa_1), S(\varkappa_2)$ are the same or distinct. In case the symbols are not the same, then as sets they are disjoint, since a based symbol lying in both $S(\varkappa_1), S(\varkappa_2)$ would make all elements in both sets to be equivalent in $\underset{al}{\sim}$. But by definition $R(\varkappa) \subset S(\varkappa)$ for any knot $\varkappa$. Hence $R(\varkappa_1), R(\varkappa_2)$ are disjoint for distinct knots, which proves the completeness of $R$.

The truth of the claim for the Algorithm stems from the definition of $R(\varkappa)$ as $R(_pK)$ and Corollary 5. $\square$

Similarly to the symbol invariant $S$, the completeness of the reduced symbol invariant $R$ has the expected consequences. For example:

**Corollary 6.** *The reduced symbol invariant $R$ of knots detects chirality and amphichirality.*

As a byproduct of the above work:

**Corollary 7.** *The crossing number of any knot $\varkappa$ is computed as follows: consider a knot diagram $K$ of $\varkappa$ and a base point $p$ of it. For the based symbol $_pK$ consider consecutive negative moves up to the point no such move exists. The order of the based symbol we arrive at, is the crossing number of $\varkappa$.*

We close with an example.

**Example 4.** *Let us compute the reduced symbol invariant for the two trefoil knots $\varkappa, \varkappa'$.*

*Let $K, K'$ be the diagrams of $\varkappa, \varkappa'$ respectively given in Figure 50. We have $_pK = (1,4)^{+1} \ (5,2)^{+1} \ (3,6)^{+1}$ and $_qK' = (1,4)^{-1} \ (5,2)^{-1} \ (3,6)^{-1}$.*

*Neither one of $_pK, _qK'$ accepts any negative $k$-move for $k = 1, 2$. So each one of them is a reduced based symbol. Neither accepts any 3-move as well. Thus only $^{al}\Omega_\beta$ moves exist on them.*

*All renamings of $_pK$ are derived by cyclically changing the integers $1, 2, 3, 4, 5, 6$ (putting $1$ to any one of them and keeping the cyclic order), thus they are:*

$$(1,4)^{+1} \ (5,2)^{+1} \ (3,6)^{+1}, \quad (2,5)^{+1} \ (6,3)^{+1} \ (4,1)^{+1}, \quad (3,6)^{+1} \ (1,4)^{+1} \ (5,2)^{+1},$$
$$(4,1)^{+1} \ (2,5)^{+1} \ (6,3)^{+1}, \quad (5,2)^{+1} \ (3,6)^{+1} \ (1,4)^{+1}, \quad (6,3)^{+1} \ (4,1)^{+1} \ (2,5)^{+1}.$$

*Similalry all renamings of $_qK'$ are:*

$$\{(1,4)^{-1} \ (5,2)^{-1} \ (3,6)^{-1}, \quad (2,5)^{-1} \ (6,3)^{-1} \ (4,1)^{-1}, \quad (3,6)^{-1} \ (1,4)^{-1} \ (5,2)^{-1},$$
$$(4,1)^{-1} \ (2,5)^{-1} \ (6,3)^{-1}, \quad (5,2)^{-1} \ (3,6)^{-1} \ (1,4)^{-1}, \quad (6,3)^{-1} \ (4,1)^{-1} \ (2,5)^{-1}.$$

*So:*

$$R(\varkappa) = \{(1,4)^{+1} \ (5,2)^{+1} \ (3,6)^{+1}, \quad (4,1)^{+1} \ (2,5)^{+1} \ (6,3)^{+1}\}.$$
$$R(\varkappa') = \{(1,4)^{-1} \ (5,2)^{-1} \ (3,6)^{-1}, \quad (4,1)^{-1} \ (2,5)^{-1} \ (6,3)^{-1}\}.$$

*Since $R(\varkappa), R(\varkappa')$ are distinct, $\varkappa, \varkappa'$ are distinct.*

**Remark 8.** *As a final remark, let us observe that a great deal of the technicalities in this paper concerns the relabelings. We introduced them in order to have a description of the crossings in a based diagram after a change of the base point or after a topological $\Omega_1^\pm$ or $\Omega_2^\pm$ move. Subsequently we incorporated this language to the algebraic situation regarding the $^{al}\Omega_\beta, ^{al}\Omega_1^\pm$ or $^{al}\Omega_2^\pm$ moves.*

*Actually we did not work out a complete set of descriptions for all possible base points for each one of the above topological moves. Although one can provide such a full set, this will change nothing regarding the equivalence classes of diagrams and of based symbols, thus it would change nothing regarding our invariants. The change of base move $^{al}\Omega_\beta$ takes care of things even when in shortage of relabelings. Hence it is a matter of taste which set of relabelings one chooses to work with.*

Department of Mathematics, National Technical University of Athens, Zografou campus, GR-15780 Athens, Greece.

*Email address*: `dkodokostas@math.ntua.gr`